\documentclass[11pt]{amsart}
\usepackage[english]{babel}
\usepackage[utf8]{inputenc}
\usepackage{fullpage}
\usepackage{enumerate}
\usepackage{amsaddr}

\usepackage{accents}
\newcommand*{\dt}[1]{%
  \accentset{\mbox{\large\bfseries .}}{#1}}
\newcommand*{\ddt}[1]{%
  \accentset{\mbox{\large\bfseries .\hspace{-0.25ex}.}}{#1}}

\usepackage[toc]{appendix}

\usepackage{amsmath,amsthm,amssymb,amsfonts}
\usepackage{mathtools}
\usepackage{esint}
\usepackage{appendix}

\usepackage{tikz}

\usepackage{float}
\usepackage{graphics,graphicx}
\usepackage[backref=page]{hyperref}
\hypersetup{
	plainpages=false,
	unicode=false,          % non-Latin characters in Acrobat\u2019s bookmarks
	pdftoolbar=true,        % show Acrobat\u2019s toolbar?
	pdfmenubar=true,        % show Acrobat\u2019s menu?
	pdffitwindow=true,      % page fit to window when opened
	pdftitle={My title},    % title
	pdfauthor={Author},     % author
	pdfsubject={Subject},   % subject of the document
	pdfcreator={Creator},   % creator of the document
	pdfproducer={Producer}, % producer of the document
	pdfkeywords={keywords}, % list of keywords
	pdfnewwindow=true,      % links in new window
	colorlinks=true,       % false: boxed links; true: colored links   %%%%%%%%%%% MODIFIE
	linkcolor=blue,          % color of internal links
	citecolor=blue,        % color of links to bibliography
	filecolor=blue,      % color of file links
	urlcolor=magenta,          % color of external links  
	urlbordercolor=0 1 1
}

\usepackage[hang,small]{caption}
\usepackage[]{subcaption}

\graphicspath{{./flagfoldsFigures/}}

\usepackage{faktor} % quotient

\newcommand{\R}{\mathbb R}
\newcommand{\N}{\mathbb N}

\newcommand{\cH}{\mathcal{H}}

\newcommand{\cE}{\mathcal{E}}
\newcommand{\cF}{\mathcal{F}}

\newcommand{\cI}{\mathcal{I}}
\newcommand{\cJ}{\mathcal{J}}
\newcommand{\cU}{\mathcal{U}}
\newcommand{\cL}{\mathcal{L}}
\newcommand{\cW}{\mathcal{W}}

\newcommand{\cV}{\mathcal{V}}

\def\xC{{\rm C}}
\def\xLip{ {\rm Lip} }

\def\xL{{\rm L}}

\def\xdiv{{\rm div}}
\def\xdiam{{\rm diam}}

\def\xM{{\rm M}}
\def\xO{{\rm O}}
\def\xspan{{\rm span \:}}
\def\xSym{{\rm Sym}}
\def\xSkew{{\rm Skew}}
\def\xdiag{{\rm diag \:}}

\def\xLM{L_{\overline{M}}}
\def\xdM{d_{\overline{M}}}
\def\xLWF{L_{\cW\cF}}
\def\xdWF{d_{\cW\cF}}

\newcommand{\V}{\|V\|}

\newtheorem{theorem}{Theorem}[section]
\newtheorem{proposition}[theorem]{Proposition}

\newtheorem*{theorem*}{Theorem}
\newtheorem{definition}[theorem]{Definition}

\theoremstyle{remark}
\newtheorem{remk}[theorem]{Remark}

\newtheorem{xmpl}[theorem]{Example}

\DeclareMathOperator*{\diam}{diam}
\DeclareMathOperator*{\tr}{tr}

\newcommand{\one}{\mathbf{1}}
\renewcommand{\phi}{\varphi}
\renewcommand{\epsilon}{\varepsilon}

\newcommand{\G}{G_{d,n}}

\title{Flagfolds: an approach to multi-dimensional varifolds}

\author{Blanche \textsc{Buet}}
\address{Universit{\'e} Paris-Saclay, Inria, Cnrs, Laboratoire de math{\'e}matiques d'Orsay (Orsay, France)}

\author{Xavier \textsc{Pennec}}
\address{Université Côte d'Azur and Inria, Epione team (Sophia-Antipolis, France)
}
\date{\today}
\thanks{B. Buet acknowledges support from the French
National Research Agency (ANR) under grant ANR-21-CE40-0013-01 (project GeMfaceT) and grant ANR-24-CE40-2216 (project STOIQUES).\\
X. Pennec was funded by the European Research Council (ERC) under the European Union’s Horizon 2020 research and innovation program (grant agreement Nr. 786854 G-Statistics). He was also supported by the French government through the 3IA Côte d’Azur Investments ANR-19-P3IA-0002 \& 23-IACL-0001 managed by the French National Research Agency(ANR).}


\begin{document}

\definecolor{gray}{rgb}{0.5, 0.5, 0.5}

\begin{abstract}
By interpreting the product of the Principal Component Analysis, that is the covariance matrix, as a sequence of nested subspaces naturally coming with weights according to the level of approximation they provide, we are able to embed all $d$--dimensional Grassmannians into a stratified space of covariance matrices. We observe that Grassmannians constitute the lowest dimensional skeleton of the stratification while it is possible to define a Riemannian metric on the highest dimensional and dense stratum, such a metric being compatible with the global stratification. With such a Riemannian metric at hand, it is possible to look for geodesics between two linear subspaces of different dimensions that do not go through higher dimensional linear subspaces as would euclidean geodesics. Building upon the proposed embedding of Grassmannians into the stratified space of covariance matrices, we generalize the concept of varifolds to what we call flagfolds in order to model multi-dimensional shapes.
\end{abstract}

     \maketitle

\section*{Introduction}

When analyzing data in large dimensions, one often looks for the optimal lower dimension to project the data without losing important information. However, the multiscale nature of the data
often prevents the identification of a single, well-defined intrinsic dimension: we rather obtain a set of dimensions depending on the quality of the approximation that we are allowing on the data.  
Moreover,  assuming that data live on a submanifold of fixed dimension (``the manifold hypothesis'') is often false: tree-like structures, for instance, live on stratified spaces \cite{billera_geometry_2001}. Likewise, quotient spaces are often stratified and their dimension vary at the points where the isotropy group changes. In practice, the local dimension of the data may vary with the location but also the scale at which we look at them. In astrophysics, for instance, the large scale structure of the universe has a web appearance with structures aggregated in dense compact clusters connected by elongated filaments and  sheetlike walls \cite{Yepes} while the smaller scale structures display galaxies and individual stars. Similarly, white matter tracts regrouping axons in the brain can locally have a tubular shape that evolves into a  thin sheet-like structures at some places when measured using diffusion MRI \cite{yushkevich_structure-specific_2008}. 
Such complex models can be seen at different scales and there is not a single notion of dimension, even locally. 

\emph{Flags to encode multi-dimensionality.} An interesting idea to tackle this variability is the notion of flags (Definition~\ref{dfn:flag}), which are series of properly embedded subspaces. Flags naturally arise implicitly in many statistical analyses such as principal component analysis (see Section~\ref{subsec:FlagPCA}). For instance, finding the local dimension of a dataset is often performed by selecting a  neighborhood and performing a local Principal Component Analysis (PCA) to determine the number of components and to find a basis of the local tangent space. This is for instance the essence of the Local Tangent Space Approximation (LTSA) manifold learning algorithms \cite{zhang_principal_2002}.  Although one often selects a fixed dimension to approximate data,  PCA  constructs in fact a series of nested subspaces that iteratively improves the approximation of data. 
Such an increasing sequence of properly embedded vector subspaces is a  filtration of subspaces called a flag: starting from a zero-dimensional subspace (the mean), the flag  is generated by adding dimensions along the successive eigenspaces of the covariance matrix with increasing eigenvalues. If all the eigenvalues have multiplicity one, one can parameterize the flag by the ordered set of eigenvectors. With a multiplicity larger than none, we should consider adding the subspace generated by all the eigenvectors of the eigenvalue. The importance of each subspace may be measured by the additional variance which is explained by each subspace, that is, the difference $\mu_i = i (\lambda_i - \lambda_{i+1})$ of the successive eigenvalues $\lambda_1 \geq \lambda_2 \geq \ldots \geq \lambda_n \geq \lambda_{n+1} = 0$. This gives birth to the notion of \emph{weighted flags} that couples a flag with a sequence $(\mu_1, \ldots, \mu_n)$ of associated weights (see Section~\ref{subsec:weigthedFlagTopo}). The completeness or sparsity of a weighted flag then reads on the zero weights. This is called the \emph{type} of the flag (Definition~\ref{dfn:typeDelta}).
In non-linear spaces, \cite{damon_backwards_2013} have argued that the nestedness of approximation spaces is one of the most important characteristics of PCA. This characteristics was one of the main features of the Barycentric Subspace Analysis method recently proposed to generalize PCA to Riemannian manifolds and geodesic spaces \cite{pennec:hal-01343881}. Thus, weighted flags of subspaces seem to be the natural mathematical objects to encode  hierarchically embedded approximation spaces with multiple dimensions. 

Flags generalize Grassmannians (linear subspaces of fixed dimensions) and moreover guaranty that higher dimensional subspaces include the lower dimensional ones, so that approximations done at different levels remain consistent. 
%Flags provide a convenient mathematical concept to encode the multiplicity and flexibility of dimensions which has not been thoroughly investigated so far in data science. 
Since the flag locally approximating the data may vary from one point of the embedding space to another, one is led to consider fields of flags. \emph{Flagfolds} are then distributions of fields of (weighted) flags, defined as linear forms integrating them, similarly to the way varifolds integrate fields of Grassmannians.

\emph{From varifolds to flagfolds.}
From the geometric measure theory point of view, $d$--varifolds provide a generalized notion of $d$--dimensional surface and can be seen as a non-oriented counterpart of currents (see \cite{Allard72,Almgren}). Technically, a $d$--varifold in $\mathbb{R}^n$ is a positive Radon measure on the product $\mathbb{R}^n$ times the $d$--Grassmannian $\G$. It can be also understood as a positive measure in $\mathbb{R}^n$, whose support is the geometric shape of interest (the measure being the multiplicity of the object or more generally a weight on the object)
coupled with a probability measure on the $d$--Grassmannian given at almost any point $x$ of the object and varying with the position $x$, typically, the probability for given a $d$--subspace to be the tangent space at $x$.
Although classical examples of varifolds are associated with surfaces or rectifiable sets, the flexible structure is also well-suited for ``discrete-like'' objects as it has been evidenced in \cite{trouve_charon, BuetLeonardiMasnou}.

As it is clear from the name itself, a $d$--varifold comes with a given dimension $d$ that is fixed from its definition as a measure in $\R^n \times \G$. In order to consider objects with a varying dimension or with several dimensions depending on the scale, one can substitute the $d$--Grassmannian $\G$ with a larger subset of symmetric matrices $\xSym(n)$. In this paper, we propose to replace $\G$ with the compact set of \emph{weighted flags} $\cW\cF(n)$, or equivalently in terms of topological space, with the set of positive semi definite matrices of trace $1$, $\xSym_+^1(n)$, that contains the Grassmannians $\G$ for $d = 1$ to $d = n$ (identifying $\G$ with rank $d$ orthogonal projectors). We call \emph{flagfolds} the Radon measures on $\R^n \times \cW\cF(n)$. 
In the same spirit, it has been proposed in \cite{ambrosioSoner} to define generalized varifolds as measures in $\R^n \times A_{d,n}$ where $A_{d,n} \subset \xSym(n)$  is compact and $A \in A_{d,n}$ if and only if $-n I_n \leq A \leq I_n$ and $\tr A = d$, and in particular $\G \subset A_{d,n}$ though only for the given $d$ defining $A_{d,n}$. Such generalized varifolds allow to exhibit a limit Brakke flow of a sequence of Ginzburg-Landau energies. In \cite{tinarrage} the author directly replaces $\G$ with the non compact whole space $\xSym(n)$, working with measures in $\R^n \times \xSym(n)$. Last but not least, in \cite{Fomenko}, the authors develop a theory of multivarifolds to solve a multi-dimensional Plateau problem: they glue the Grassmannians to form $\cup_d \G$ and consider measures on the product $\R^n \times \cup_d \G$ so that multivarifolds identify with sums of varifolds and such a construction does not seem to model transitions between objects of different dimensions.

%an increasing sequence of subspaces at each point (a flag) and we need to define a measure giving the probability of a flag to be the the increasing sequence of subspaces approximating the data at each point.
%

\subsection*{Main contributions of the paper and perspectives}
The aim of this paper is twofold: in a first part, we focus on the structure that can be given to the set of \emph{weighted flags} $\cW\cF(n)$ (see Definition~\ref{dfn:weightedFlags}) while in a second part we define and investigate properties of \emph{flagfolds} (see Definition~\ref{dfn:flagfolds}).
More precisely, Section~\ref{section:WeightedFlags} to \ref{section:RiemannianStructureFlags} delve into finer and finer structural properties of $\cW\cF(n)$ and of its \emph{strata}: global topological quotient, stratification with respect to the type of the flag, manifold and Riemannian structure of each stratum and eventually a global length structure on $\cW\cF(n)$. The main contribution of this first part is to define a Riemannian metric $g$ (see Proposition~\ref{prop:WFriemannianMetrics}) on the dense stratum $M(n)$ of complete weighted flags (i.e. of type $(1, \ldots, 1)$) that is compatible with the global quotient structure in the following sense: we prove in Theorem~\ref{thm:completion} (and Proposition~\ref{prop:xdMxdWF}) that the completion of $(M(n) , g)$ is a metric space (and even a length space) $(\cW\cF(n) , \xdWF)$ and $\xdWF$ is compatible with the quotient topology in $\cW\cF(n)$ (see Proposition~\ref{prop:dMtopoEquivalent}). This compatibility seems to translate into concrete properties of the geodesics: Section~\ref{section:geodesic} proposes a basic numerical method to compute numerical approximations of geodesics starting from a given point with given initial speed in the stratum $(M(n) , g)$.  Loosely speaking, experiments show that along a numerical geodesic, the complexity of the weighted flag does not increase: more precisely, if a weight $\mu_i$ is close to $0$ at both ends of the geodesic then it remains small along the geodesic itself, which is not true along euclidean geodesics (see Figures~\ref{figNumGeod1} and \ref{figEuclGeod1} for a comparison on a concrete example). Such property is exactly the point of the definition of the Riemannian metric $g$: ensuring that a geodesic from two flags that both are nearly a line does not pass through more ``diffuse'' flags but rather rotates the direction of the line, and more generally rotates the nested subspaces. In other words, even though it is only numerical evidence performed in a simple low-dimensional setting, we hope that $g$ consistently extends the Riemannian metrics on each $\G$, $d = 1, \ldots, n$ in the sense that geodesics connecting endpoints that are both close to some $\G$ or more generally to some elementary cell (of flags of a given type) remain close to it in between. Of course, such property should depend on the function $f$ which is involved in Proposition~\ref{prop:WFriemannianMetrics} to pinch the canonical Riemannian metric $g^{\R^n} + g^{(1,\ldots,1)}$ in $M(n)$ and define $g$ in the spirit of a warped metric.

Sections~\ref{section:flagfoldsDef} and \ref{section:flagfoldFirstVariation} are devoted to the investigation of flagfolds. As already mentioned, flagfolds are Radon measures in $\R^n \times \cW\cF(n)$ and generalize $d$--varifolds: more precisely, it is possible to associate a flagfold $\widehat{V}$ with a $d$--varifold $V$ (for any $d = 1 \ldots n$) and the resulting embedding of $d$--varifolds into flagfolds preserves the mass measure $\| \widehat{V} \| = \| V \|$, while conversely, one can associate with a given flagfold $W$ a $d$--varifold $V_d$ for each $d = 1 \ldots n$ satisfying $\| W \| = \| V_1 \| + \ldots \| V_n \|$ though $W = V_1 + \ldots + V_n$ does not hold in general since weighted flags contains more than $\cup_{d=1}^{n} \G$. In other words, there is a one-to-one correspondence between family of varifolds $(V_1, \ldots, V_n)$ and flagfolds with support in $\R^n \times\cup_{d=1}^{n} \G$ (see Section~\ref{subsec:flagfoldsEmbeddingVarifolds}). 
Flagfolds then allow to model diffused approximations of lower dimensional structures such as a transition from a $3d$--thin tube to a $1d$--line as detailed in Example~\ref{xmpl:diffusedFlagfold} which was not possible in the varifold framework as these two objects were not even living in the same space of varifolds because of their respective dimensions, such an example concretely arises from MRI diffusion data of the brain. We also emphasize that from a numerical perspective, the covariance matrix is very often the natural information arising from data, working with the flagfold structure avoids the additional truncation step performed to obtain a $d$--plane from the covariance matrix, that hence requires the a priori dimension $d$: the flagfold structure is closer to the information extracted from usual data than the varifold structure.
Eventually, Section~\ref{section:flagfoldFirstVariation} shows that it possible to extend the notion of first variation from varifolds to flagfolds (Definition~\ref{dfn:flagfoldFirstVariation}) consistently with the aforementioned embedding: for a varifold $V$, the first variation of $V$ and $\widehat{V}$ coincide. In a very similar way as it is done for varifolds, it is possible to infer a monotonicity formula (see Proposition~\ref{prop:monotonicity}) and a structure theorem (see Theorem~\ref{thm:rectifiability}) from the control of the first variation of a flagfold: under natural assumptions on $d$--dimensional lower densities of each layer $V_d$ of $W$ and supposing that the first variation of $W$ is locally bounded, we establish that $W$ decomposes as the sum of its $d$--dimensional varifold layers $W = V_1 + \ldots + V_d$. Yet, we are not able to prove their rectifiability under such assumptions.
The well-definition of the first variation of a flagfold leads the path to the definition of associated generalized mean curvature, second fundamental form, mean curvature flow and in particular multi-dimensional mean curvature flow as well as their approximate versions (see for instance \cite{Allard72,Hutchinson,Brakke,BuetLeonardiMasnou,BuetLeonardiMasnou2,BuetRumpf} when dealing with varifolds) as we intend to explore in the near future.

%Recall that in order to define $d$--varifolds, one only needs to provide $\G$ with a topology, similarly, the definition of flagfolds only requires to provide $\cW\cF(n)$ with a topology, nevertheless, the understanding of the structure ``stratification'' with respect to the type of the flag is crucial to extend classical concepts from varifolds to flagfolds.

\section*{Notations}
We give hereafter general notations and additionally, for more specific notations, we give the first occurrence where they are defined in the paper.

We fix integers $d$, $n$ satisfying $1 \leq d \leq n$.
\begin{itemize}
 \item $(e_1, \ldots, e_n)$ denotes the canonical basis of $\R^n$.

 \item $\xSym(n) = \{ A \in \xM_n(\R) \: : \: A^T = A \}$ is the subspace of $n \times n$ symmetric matrices.
 \item $\xSkew(n) = \{ A \in \xM_n(\R) \: : \: A^T + A = 0 \}$ is the subspace of $n \times n$ skew symmetric matrices.
 \item $\xSym_+(n) \subset \xSym(n)$ is the set of symmetric positive semi-definite matrices.
 \item $\xSym_+^1(n) = \{ A \in \xSym_+(n) \: : \: \tr(A) = 1 \}$.
 \item $\xO(n) = \{ P \in \xM_n(R) \: : \: P^T P = I_n \}$ is the orthogonal group.
 \end{itemize}
For $I = (p_1, p_2, \ldots, p_r)$ with positive integers $r, p_1, \ldots, p_r$,
 \begin{itemize}
 \item $\xO(I) = \left\lbrace {\rm diag}(R_1, \ldots, R_r) \: : \: \forall i = 1, \ldots, r, \, R_i \in \xO(p_i) \right\rbrace \simeq \xO(p_1) \times \xO(p_2) \times \ldots \times \xO(p_r)$.
 \item $\xSkew(I)  = \left\lbrace {\rm diag}(A_1, \ldots, A_r) \: : \: \forall i = 1, \ldots, r, \, A_i \in \xSkew(p_i) \right\rbrace \simeq \xSkew(p_1) \times \ldots \times \xSkew(p_r)$.
 \item $\cF_I = \faktor{\xO(n)}{\xO(I)}$ with quotient map $\pi_I : \xO(n) \rightarrow \cF_I$.
 \item Given $U \in \xM_{l,n}(\R)$, we write $U = (u_1, \ldots, u_n)$ meaning that $u_1, \ldots, u_n \in \R^l$ are the columns of $U$.
 \item Given a $d$--subspace $E$ of $\R^n$, $\Pi_E \in \xSym_+(n)$ is the (rank $d$) orthogonal projector onto $E$.
 \item $\begin{array}{lcl}
 \G & = & \{ E \subset \R^n \: : \: E \text{ is a vector subspace of } \R^n, \: \dim E = d \} \\
 & \simeq & \{ P \in \xSym(n) \: : \: P^2 = P \text{ and } \tr P = d \} \end{array}$.
\end{itemize}
Section~\ref{section:WeightedFlags}:
\begin{itemize}
\item $p_{J \to I} : \cF_J \rightarrow \cF_I$ for $J \preccurlyeq I$ in Definition~\ref{dfn:typeCoarser}.
\item $\Delta(n) = \{ x \in [0,1]^n \: : \: \sum_{i=1}^n x_i = 1 \}$ is the unit simplex of $\R^n$ and\\
$\cW(n) = \{ x \in [0,1]^n \: : \: x_1 \geq x_2 \geq \ldots \geq x_n, \: \sum_{i=1}^n x_i = 1 \}$.
\item Given $A \in \xSym_1^+(n)$, 
\begin{itemize}
    \item $\lambda_1(A) \geq \lambda_2(A) \geq \ldots \geq \lambda_n(A)$ are the ordered eigenvalues of $A$,
    \item $E_\lambda(A)$ is the eigenspace associated with the eigenvalue $\lambda \in [0,1]$, $E_{\lambda_k(A)}(A)$ is often referred as $E_k(A)$.
    \item $\begin{array}{lcl}
     \mu(A) = (\mu_1(A), \ldots, \mu_n(A)) & =&  (\lambda_1(A) - \lambda_2(A), 2(\lambda_2(A) - \lambda_3(A)),
     \\ & &\ldots , k (\lambda_k(A) - \lambda_{k+1}(A)), \ldots , n \lambda_n(A)) \in \Delta(n),  \mbox{\eqref{eq:lambdaToMu}.}
    \end{array}$
\end{itemize}
\item Weighted flags $\cW\cF(n) = \faktor{\Delta(n) \times \xO(n)}{ \sim }$ homeomorphic to $\xSym_+^1(n)$, see Definition~\ref{dfn:weightedFlags} and Proposition~\ref{prop:homeo}.
\item For $\alpha \in \Delta(n)$, the type $\tau(\alpha)$ is introduced in Definition~\ref{dfn:typeDelta}, then for $A \in \xSym_+^1(n)$, $\tau(A) = \tau ((\mu_1(A), \ldots, \mu_n(A))$.
\end{itemize}
Section~\ref{section:distancesWF}: given $r \in \{1, \ldots, n\}$ and $K = \{k_1, \ldots, k_r\} \subset \{1, \ldots, n\}$, 
\begin{itemize}
\item $\displaystyle \mathring{\Delta} (n \: ; K) = \left\lbrace \mu \in \Delta(n) \: \left| \begin{array}{l}
\mu_j > 0 \text{ for } j \in K   \\
\mu_j = 0 \text{ for } j \notin K  
\end{array} \right. \right\rbrace$, see \eqref{eq:stratifCell}.
\item $\displaystyle 
M(r \: ; K) = \left\lbrace ( \mu, W) \in \cW\cF(n) \: \left| \begin{array}{l}
\mu_j > 0 \text{ for } j \in K   \\
\mu_j = 0 \text{ for } j \notin K  
\end{array} \right. \right\rbrace 
$ and $M(n) \simeq \mathring{\Delta}(n) \times \cF_{(1, \ldots, 1)}$.
\end{itemize}
Section~\ref{section:SmoothAndRiemannianF_I}: given $I = (p_1, \ldots, p_r)$ and $U \in \xO(n)$,
\begin{itemize}
\item $\mathfrak{m}_{I}$ is the orthogonal complement of $\xSkew(I)$: 
$\xSkew(n) = \mathfrak{m}_{I} \: \oplus \:  \xSkew(I)$, see \eqref{eq:mI}.
\item $H_U^I = U \mathfrak{m}_{I}$, $V_U^I = \ker T_U \pi_I$ are the horizontal and vertical spaces at $U$: $T_U \xO(n) = H_U^I \oplus V_U^I$.
\item $d_I$ and $L_I$ are the distance and length associated with a Riemannian metric $g^I$ in $\cF_I$, see Proposition~\ref{prop:flagRiemannianStruct}, \eqref{eq:flagMetric} and \eqref{eq:lengthFI}.
\end{itemize}
Section~\ref{section:RiemannianStructureFlags}:
\begin{itemize}
\item Given $I = (p_1, \ldots, p_r)$, $\displaystyle
X_I = \{ (i,j) \: : \: p_1 + \ldots + p_{k} + 1 \leq i < j \leq p_1 + \ldots + p_{k+1}  \text{ for some } k \in \{ 1, \ldots, r-1 \} \}$ is the set of block diagonal indices, see \eqref{eq:diagonalIndicesI}.
\item For $\mu \in \Delta(n)$, $\mu_{i \to j} = (\underbrace{0, \ldots, 0}_{\in \R^{i-1}}, \underbrace{ \mu_i, \mu_{i+1}, \ldots, \mu_{j-1} }_{\in \R^{j-i}} , \underbrace{0, \ldots, 0}_{\in \R^{n-j+1}})$.
\item $g$ is defined in Proposition~\ref{prop:WFriemannianMetrics} and provides each $M (r \: ; K)$ with a Riemannian metric, length $L_g$ and distance $d_g$.
\item $\xLM$ is a length structure inducing the distance $\xdM$ in $\cW\cF(n)$, see \eqref{eq:lengthStructure} and \eqref{eq:RiemannianDistance}.
\item $\xLWF$ is another length structure in $\cW\cF(n)$, see Definition~\ref{dfn:piecewiseC1pathII}, and it induces the distance $\xdWF = \xdM$ in $\cW\cF(n)$, see Proposition~\ref{prop:xdMxdWF}.
\end{itemize}
Section~\ref{section:flagfoldsDef}:
\begin{itemize}
\item For $\mu \in \Delta(n)$, $\overline{d}(\mu) = \sum_{k=1}^n k \mu_k$ and for $A \in \xSym_+^n$, $\overline{d}(A) = \overline{d}(\mu(A))$, see \eqref{eq:dimensionMu}.
\item $i : E \in \G \rightarrow \frac{1}{d} \Pi_E \in \xSym_+^1(n)$ and for $S \in \xSym_+^1(n)$, $\overline{S} = \sum_{k=1}^n \mu_k(S) \Pi_{E_k(S)}$; for $E \in \G$, $\overline{i(E)} = \Pi_E$, see \eqref{eq:GtoWF} and \eqref{eq:WFtoG}.
\item From a $d$--varifold $V$ to the flagfold $\widehat{V} = ({\rm Id}, i)_\# V$, see \eqref{eq:varifoldToFlagfold}.
\item From the flagfold $W$ to the varifolds $(V_1, \ldots, V_d)$, $V_d = \mu_d \: ({\rm Id}, E_d)_\# W$, see \eqref{eq:flagfoldsToVarifolds}.
\end{itemize}

\section{Weighted flags: a topological quotient space}
\label{section:WeightedFlags}

The purpose of this section is to introduce and investigate the topological structure of what we call \emph{weighted flags} that couple a sequence of nested vector subspaces, i.e. a flag, with a sequence of respective weights. To this end, in Section~\ref{secFlag}, we begin with some known facts concerning flags: we introduce the quotient set $\cF_I$ (see \eqref{eq:flagHomogeneous}) of all flags of a fixed type $I$ and we characterize its topology in Proposition~\ref{prop:flagCV}. Section~\ref{subsec:FlagPCA} explains how the geometric information resulting from a Principal Component Analysis (PCA) can be naturally represented with a flag whose successive subspaces have different relative importance which naturally leads to Definition~\ref{dfn:weightedFlags} of \emph{weighted flags} as a quotient space denoted by $\cW\cF(n)$. Section~\ref{subsec:weigthedFlagTopo} furthermore characterizes the topology of $\cW\cF(n)$ in Proposition~\ref{prop:weightedFlagsCV} and checks that the eigen decomposition induces an homeomorphism between $\xSym_1^+(n)$ and $\cW\cF(n)$ in Proposition~\ref{prop:homeo}.

\subsection{Flags of a fixed type} \label{secFlag}

First of all, let us recall what a flag of $\R^n$ is:

\begin{definition}[Flag] \label{dfn:flag}
A \emph{flag} of $\R^n$ is an increasing sequence of vector subspaces of $\R^n$
\[
\{ 0 \} = E_0 \subset E_1 \subset E_2 \subset \ldots \subset E_r = \R^n \: .
\]
Increasing implicitly meaning distinct so that denoting $d_i = \dim E_i \in \{ 0, \ldots, n \}$, we have $0 = d_0 < d_1 < \ldots < d_r = n$ and $(d_0, d_1, \ldots, d_r)$ is called the \emph{signature} of the flag. Setting $p_i = d_i - d_{i-1} \in \{1, \ldots, n\}$ for $i \in \{1, \ldots, r\}$, we also introduce the \emph{type} $(p_1, \ldots , p_r)$ of the flag, satisfying $p_1 + \ldots + p_r = n$.
\end{definition}

Let us fix some type $I = (p_1, \ldots, p_r)$ and let $\cF_I$ denote the set of all flags of type $I$. We pass from the type $(p_1, \ldots, p_r)$ to the signature through $d_0 = 0$ and $d_i = p_1 + \ldots + p_i$. The orthogonal group $\xO(n)$ acts transitively on $\cF_I$ 
\[
\left\lbrace \begin{array}{lll}
\xO(n) \times \cF_I & \rightarrow & \cF_I \\
(q, (E_0, E_1, \ldots, E_r)) & \mapsto & (q(E_0), q(E_1), \ldots, q(E_r))
\end{array}
\right.
\]
Moreover, defining $F_0 = \{0\}$ and for $i \in \{0, \ldots, r \}$, $F_i = \xspan (e_1, \ldots, e_{d_i})$, where $\{e_1,\ldots e_n\}$ is the canonical basis of $\R^n$, one obtain a flag of type $I$, called \emph{standard flag} of type $I$. As the stabilizers of flags in the same orbit are conjugated to each other and the action is transitive, it is enough to compute the stabilizer of the standard flag above. It is not difficult to check that if $q$ stabilizes the standard flag of type $I$, it has a block--diagonal matrix with sizes $p_1$, $p_2$ to $p_r$ and then the stabilizer is isomorphic to $\xO(I) = \xO(p_1) \times \xO(p_2) \times \ldots \times \xO(p_r)$, leading to the identification
\begin{equation} \label{eq:flagHomogeneous}
\cF_I = \faktor{\xO(n)}{\xO(I)} \: .
%\cF_I = \frac{\xO(n)}{\xO(p_1) \times \ldots \times \xO(p_r)}
\end{equation}

\begin{xmpl}[Grassmannian]
It is possible to identify the $d$--dimensional Grassmannian of $\R^n$, 
\[
\G = \{ E \subset \R^n \: : \: E \text{ is a vector subspace of } \R^n, \: \dim E = d \}
\]
with the set of flags of type $I = (d, n-d)$.
\end{xmpl}
Identification \eqref{eq:flagHomogeneous} allows to give a rich structure to $\cF_I$, starting with a topology. Let us be more precise.
We fix hereafter a type $I = (p_1, \ldots, p_r)$. We recall that $\xO(n)$ is a Lie group and $\xO(I)$ is a closed (Lie) subgroup of $\xO(n)$. The right action of $\xO(I)$ on $\xO(n)$
\[
\begin{array}{ccc}
\xO(n)  \times \xO(I) & \rightarrow & \xO(n) \\
(U,  R) & \mapsto & U R
\end{array} \quad \text{with} \quad R = {\rm diag} (R_1, \ldots, R_r) \text{ and } (R_1, \ldots, R_r)  \in \xO(I)
\]
is continuous (and smooth) and $\cF_I = \faktor{\xO(n)}{\xO(I)}$ is the associated topological quotient. We recall the following topological properties of $\cF_I$ in the next proposition (see $7.12$ in \cite{Boothby} for details).

\begin{proposition}
We denote by $\pi_I : \xO(n) \rightarrow \faktor{\xO(n)}{\xO(I)}$ the canonical projection. Then $\pi_I$ is an open map and $\faktor{\xO(n)}{\xO(I)}$ is Hausdorff.
\end{proposition}

Given $U \in \xO(n)$, we denote by $\pi_I(U)$ (or simply $\pi(U)$ when there is no ambiguity) the class of $U$ in $\cF_I$. We then have
\[
\pi_I(U) = \{ U R \: : \: R = {\rm diag} (R_1, \ldots, R_r) \in \xO(I) \} \: .
\]
In order to explicit the convergence in the quotient $\cF_I$, we start with the case of $\G = \cF_{(d,n-d)}$. We do not prove the following assertions and we refer to the proof of Proposition~\ref{prop:flagCV} for details. Let $(V^{(m)})_{m \in \N}, V \in \cF_{(d,n-d)}$ and let $(U^{(m)})_{m \in \N}, U \in \xO(n)$ such that $\pi \left( U^{(m)} \right) = V^{(m)}$ and $\pi \left( U \right) = V$ (with $\pi = \pi_{(d, n-d)}$). Let us write $U = (u_1, \ldots, u_n)$ and for all $m \in \N$, $(U^{(m)}) = (u_1^{(m)}, \ldots, u_n^{(m)})$. Then $(V^{(m)})_{m \in \N}$ converges to $V$ in $\cF_{(d,n-d)}$ if and only if 
\begin{align*}
\Longleftrightarrow & \, \forall m, \exists R^{(m)} = {\rm diag}(R_1^{(m)}, R_2^{(m)}) \in \xO(I) \text{ such that } U^{(m)} R^{(m)} \xrightarrow[m \to \infty]{} U \text{ in } \xO(n)\\
\Longleftrightarrow & \, \forall m, \exists (R_1^{(m)}, R_2^{(m)}) \in \xO(I) \text{ such that } \left\lbrace
\begin{array}{ll}
 (u_1^{(m)}, \ldots, u_d^{(m)}) {R_1^{(m)}} & \xrightarrow[m \to \infty]{} (u_1, \ldots, u_d)\\
(u_{d+1}^{(m)}, \ldots, u_n^{(m)}) {R_2^{(m)}} & \xrightarrow[m \to \infty]{} (u_{d+1}, \ldots, u_n)
\end{array} \right. \\
\Longleftrightarrow & \, \text{the following orthogonal projectors converge} \\
& \quad\quad\quad
\left\lbrace
\begin{array}{ll}
 (u_1^{(m)}, \ldots, u_d^{(m)}) (u_1^{(m)}, \ldots, u_d^{(m)}) ^T & \xrightarrow[m \to \infty]{} (u_1, \ldots, u_d)(u_1, \ldots, u_d)^T\\
 (u_{d+1}^{(m)}, \ldots, u_n^{(m)})  (u_{d+1}^{(m)}, \ldots, u_n^{(m)})^T & \xrightarrow[m \to \infty]{} (u_{d+1}, \ldots, u_n)(u_{d+1}, \ldots, u_n)^T
\end{array} \right. \\
\Longleftrightarrow & \,
\Pi_{\xspan (u_1^{(m)}, \ldots, u_d^{(m)}) } 
\xrightarrow[m \to \infty]{} \Pi_{\xspan(u_1, \ldots, u_d)}
\end{align*}
This last characterization of the convergence in $\cF_{(d,n-d)}$ corresponds to a different identification of the Grassmannian with orthogonal projectors:
\[
\begin{array}{ccl}
\G & \simeq & \{ P \in \xSym(n) \: : \: P^2 = P \text{ and } \tr P = d \} \\
E & \mapsto & \Pi_E \: .
\end{array}
\]
We subsequently observe that both identifications of $\G$ with $\cF_{(d,n-d)}$ and orthogonal projectors of rank $d$, induce the same topology in $\G$. Given $d$--dimensional subspaces $\left( E^{(m)} \right)_{m \in \N}$, $E$ in $\G$ with respective orthonormal basis $(u_1^{(m)}, \ldots, u_n^{(m)})$, $(u_1, \ldots, u_n)$, we can define
\[
E^{(m)} \xrightarrow[m \to \infty]{\G} E \quad \Longleftrightarrow \quad \Pi_{E^{(m)}} \xrightarrow[m \to \infty]{} \Pi_E
\quad \Longleftrightarrow \quad \pi \left( u_1^{(m)}, \ldots, u_n^{(m)} \right) \xrightarrow[m \to \infty]{\cF_{(d,n-d)}} \pi \left( u_1, \ldots, u_n \right)  \: .
\]
The next proposition characterize the convergence in $\cF_I$, generalizing what we observed in the case of $\cF_{(d,n-d)}$.

\begin{proposition}[Convergence of fixed type flags] \label{prop:flagCV}
Let $I = (p_1, \ldots, p_r)$ and for all $k = 1 \ldots r$, $d_k = p_1 + \ldots + p_k$.
Let $(V^{(m)})_{m \in \N}, V \in \cF_I$ and let $(U^{(m)})_{m \in \N}, U \in \xO(n)$ such that $\pi \left( U^{(m)} \right) = V^{(m)}$ and $\pi \left(U \right) = V$. We introduce for $m \in \N$ and $k = 1, \ldots, r$:
\[
\begin{array}{ccccccccl}
U & = & (u_1, \ldots, u_n), &  F_k & = & \xspan \left(u_{d_{k-1} + 1}, \ldots , u_{d_k} \right) , &   E_k & = & \bigoplus_{i=1}^k F_i  \\
 (U^{(m)}) & = & (u_1^{(m)}, \ldots, u_n^{(m)}) , &  F_k^{(m)} & = & \xspan \left(u_{d_{k-1} + 1}^{(m)}, \ldots , u_{d_k}^{(m)} \right)  , & E_k^{(m)} & = & \bigoplus_{i=1}^k F_i^{(m)}  
\end{array}
\]
Then $(V^{(m)})_{m \in \N}$ converges to $V$ in $\cF_I$ and we write $V^{(m)} \xrightarrow[m \to \infty]{\cF_I} V$ if and only if 
\begin{equation*}
\forall k = 1, \ldots, r, \, F_k^{(m)} \xrightarrow[m \to +\infty]{{\rm G}_{k,n}} F_k \quad \Longleftrightarrow \quad \forall k = 1, \ldots, r, \, E_k^{(m)} \xrightarrow[m \to +\infty]{{\rm G}_{k,n}} E_k
\end{equation*}
\end{proposition}
\begin{proof}
The last equivalence is a direct consequence of the orthogonality of the subspaces.

\noindent Let $\| \cdot \|$ be a norm in $\xM_n(\R)$ and let $U \in \xO(n)$. Consider open balls $B(U,\epsilon) = \{ U' \in \xO(n) \: : \: \| U' - U \| < \epsilon \}$ of center $U \in \xO(n)$ and radius $\epsilon> 0$ for the induced distance in $\xO(n)$.
Then the sets 
$ \displaystyle
\mathcal{U}_\epsilon = \pi_I (B(U, \epsilon))$ are open, recalling that $\pi_I$ is an open map.
%= \bigcup_{U' \in B(U,\epsilon)} \pi_I(U') 
% = \bigcup_{U' \in B(U,r)} \{ U'R \: : \: R \in \xO(I) \}
%\text{ are open sets since } \pi_I^{-1} (\cU_\epsilon) = \bigcup_{R \in \xO(I)} B(U,\epsilon) R \text{  is open in } \xO(n) \: .
%\]

\noindent Assume that $V^{(m)} = \pi_I(U^{(m)})  \xrightarrow[m \to \infty]{\cF_I} V = \pi_I(U)$ and let $\epsilon> 0$, then $\cU_\epsilon$ is an open neighbourhood of $V$ and there exists $N = N_\epsilon \in \N$ such that for all $m \geq N$, $V^{(m)} \in \cU_\epsilon$. Consequently, for all $m \geq N$, $U^{(m)} \in \pi_I^{-1}(\cU_\epsilon) = \bigcup_{R \in \xO(I)} B(U,\epsilon) R$ and there exists $R^{(m),\epsilon} \in \xO(I)$ such that
$U^{(m)} \in B(U,\epsilon) {R^{(m),\epsilon}}^T$ i.e. $\| U^{(m)} {R^{(m), \epsilon}} - U \| < \epsilon$.
With $\epsilon = \frac{1}{l}$, for $l \in \N \setminus \{0\}$, we can for instance set $R^{(m)} = R^{(m),\epsilon_l}$ for $N_{\epsilon_l} \leq m < N_{\epsilon_{l+1}}$.
Letting $l \to \infty$, we infer that 
\[
U^{(m)} R^{(m)} \xrightarrow[m \to +\infty]{} U \quad \text{in } \xM_n(\R) \: .
\]
Let $k \in \{ 1, \ldots, r\}$. We recall that $R^{(m)} = \xdiag (R_1^{(m)}, \ldots, R_r^{(m)})$, which implies that the columns $d_{k-1} + 1$ to $d_k$ of $U^{(m)}$ and $U^{(m)} R^{(m)}$ span the same subspace $F_{k}^{(m)}$ and then $F_{k}^{(m)} \xrightarrow[m \to +\infty]{{\rm G}_{k,n}} F_k$.
%the orthogonal projection matrix onto $F_{k}^{(m)}$ is
%\begin{align*}
%\Pi_{F_{k}^{(m)}} & = U^{(m)} R^{(m)} \left( \begin{array}{c}
%                                    0 \\
%                                    \hline\\
%                                    I_{p_k}\\
%                                    \hline \\
%                                    0
%                                   \end{array} \right) \left[ U^{(m)} R^{(m)} \left( \begin{array}{c}
%                                    0 \\
%                                    \hline\\
%                                    I_{p_k}\\
%                                    \hline \\
%                                    0
%                                   \end{array} \right) \right]^T
%                                   = U^{(m)} \left( \begin{array}{ccc}
%                                    0 &  0 & 0 \\
%                                    0 & I_{p_k} & 0\\
%                                    0 & 0 & 0
%                                   \end{array} \right) {U^{(m)}}^T \\
%    & \xrightarrow[m \to +\infty]{} U \left( \begin{array}{ccc}
%                                    0 &  0 & 0 \\
%                                    0 & I_{p_k} & 0\\
%                                    0 & 0 & 0
%                                   \end{array} \right) {U}^T = \Pi_{F_{k}}
%\end{align*}

\noindent Conversely assume that for all $k \in \{1, \ldots, r\}$, $F_{k}^{(m)} \xrightarrow[m \to +\infty]{{\rm G}_{k,n}} F_k$. There exists $R^{(m)} \in \xO(I)$ such that $U^{(m)} R^{(m)} \xrightarrow[m \to +\infty]{} U$. Let $\cU \subset \cF_I$ be an open neighbourhood of $V$, then $\pi_I^{-1} (\cU)$ is open in $\xO(n)$ and contains $U$. Consequently $U^{(m)} R^{(m)} \in \pi_I^{-1} (\cU)$ for $m$ large enough and thus. %there exists $\epsilon > 0$ such that  $B(U,\epsilon) \subset \pi_I^{-1} (\cU)$ and then 
%\[
%\pi_I ( B(U,\epsilon) ) \subset \cU \quad \text{and} \quad \bigcup_{R \in \xO(I)} B(U,\epsilon) R\subset \pi_I^{-1}(\cU)  \: .
%\]
%It follows that for $m$ large enough, $U^{(m)} \in \pi_I^{-1}(\cU)$ and thus 
$V^{(m)} = \pi_I (U^{(m)} R^{(m)}) \in \cU$, which proves that $V^{(m)} = \pi_I(U^{(m)})  \xrightarrow[m \to \infty]{\cF_I} V = \pi_I(U)$.
\end{proof}

For each possible type $I$, $\cF_I$ can be provided not only with a topology but with a structure of homogeneous space, we will be more precise in Section~\ref{section:SmoothAndRiemannianF_I} focusing on the Riemannian structure inherited by $\cF_I$.
However, while the structure of flags of a fixed type is thereby provided, it does not provide a structure in the whole set of flags, which motivates our analysis. 

We end this section by defining a natural projection between spaces of flags of different types. Let us give an example, let $(e_1, \ldots e_4)$ be the canonical basis of $\R^4$ and consider the flag of type $(2,1,1)$:
$$\{ 0 \} \subset E_1 = \xspan (e_1, e_2) \subset E_2 = E_1 \oplus \xspan (e_3) \subset E_3 = E_2 \oplus \xspan(e_4) = \R^4 \: .
$$
It is possible to build a canonical flag of type $(2,2)$ from this previous flag by directly adding $\xspan(e_3,e_4)$:
$$\{ 0 \} \subset F_1 = \xspan (e_1, e_2) \subset F_2 = F_1 \oplus \xspan (e_3, e_4) = \R^4 \: .
$$
This operation uses the fact that the type $(2,2)$ is coarser than $(2,1,1)$ in the following sense:

\begin{definition}[Projection between sets of flags] \label{dfn:typeCoarser}
Let $1 \leq r \leq s \leq n$ and $I = (p_1, \ldots, p_r)$ and $J = (q_1, \ldots, q_s)$ be two types (i.e. two compositions of $n$ with non zero integers).
\begin{enumerate}[$\bullet$]
 \item We say that $I$ is coarser than $J$ and we use the notation $J \preccurlyeq  I$ if there exists $1 = i_0 \leq i_1 \leq \ldots \leq i_r = n$ such that
\[
 p_k = \sum_{j = i_{k-1} + 1}^{i_k} q_j \quad \text{for all } k = 1 \ldots r \: .
\]
In other words, $J$ is a subcomposition or a refinement of $I$.
\item In such a case $J \preccurlyeq I$, $\xO(q_1) \times \ldots \times \xO(q_s)$ is a closed subgroup of $\xO(p_1) \times \ldots \times \xO(p_r)$ and there exists a unique continuous application $p_{J \to I} : \cF_J \rightarrow \cF_I$ such that $\pi_I = p_{J \to I} \circ \pi_J$.
\end{enumerate}
\centering
\begin{tikzpicture}[node distance=2.5cm, auto]
\node (A) {$\xO(n)$};
\node(B) [right of=A] {$\faktor{\xO(n)}{\xO(I)}$};
\node (C) [below of=A] {$\faktor{\xO(n)}{\xO(J)}$};
\draw[->](A) to node {$\pi_I$}(B);
\draw[->](A) to node [left] {$\pi_J$}(C);
\draw[->](C) to node [below=1.25ex, right] {$p_{J \to I}$}(B);
\end{tikzpicture}

\end{definition}

Note that if $J \preccurlyeq I$, the canonical projection $\pi_I$ is continuous and constant on the equivalence classes with respect to the action of $\xO(J)$ on $\xO(n)$: for all $U, U^\prime \in \xO(n)$ such that $\pi_J(U) = \pi_J(U^\prime)$, $\pi_I(U) = \pi_I(U^\prime)$. Existence, uniqueness and continuity of $p_{J \to  I}$ follow from the universal property of the quotient $\cF_J$.

\subsection{Flags and Principal Analysis Component (PCA)}
\label{subsec:FlagPCA}

%There is a context allowing to give a natural sense to the fact that two flags of different types are close: principal analysis component
We are interested in the structure of flags because they naturally arise as the product of Principal Analysis Component. In short, when performing a PCA on $n$ random variables $X_1 , \ldots, X_n$ we first compute the ($n$ by $n$) covariance matrix $R$ of the data, which is symmetric positive semi-definite. In this paper, we renormalize the trace $\tr R$ to $1$,  which only changes the eigenvalues by a global factor equal to $\tr R$ and does not change the eigenspaces. Then, we compute the eigenvalues $\lambda_1 \geq \lambda_2 \geq \ldots \geq \lambda_n \geq 0$ and associated eigenvectors $V=(v_1, \ldots , v_n) \in \xO(n)$ of $R$. Let us recall that the application
\begin{equation} \label{eqEigenDec}
\left\lbrace \begin{array}{lll}
\xO(n) \times \R_+^n & \rightarrow & \xSym_+^1(n)  \\
V , (\lambda_1, \ldots, \lambda_n) & \mapsto & V \xdiag ( \lambda_1, \ldots, \lambda_n) V^T
\end{array}
\right.
\end{equation}
is surjective (spectral theorem) but badly not injective, even up to permutation (i.e. ordering of eigenvalues) and $\pm v_1, \ldots \pm v_n$. In the case of a multiple eigenvalue, for instance $\lambda_1 = \lambda_2$, then any orthonormal basis of $\xspan (v_1, v_2)$ could replace $(v_1, v_2)$.

In this setting, comparing two objects through their PCA subspace decomposition amounts to comparing two matrices $A, B \in \xSym_+^1(n)$. An obvious solution is to consider the euclidean distance 
$
\| A - B \| = \sqrt{ \tr ((A-B)^T (A-B) }
$ in $\xSym_+^1(n)$ inherited from $\xM_n (\R)$. However, such a metric does not respect the geometry of the eigen decomposition as evidenced in the following example (see also Figure~\ref{figEuclGeod1}).

\begin{xmpl} \label{ex:euclDist}
let us consider the following example in $\R^2$, $D_0$ is the horizontal axis and $D_\theta$ is a line making an angle $\theta \in [-\pi/2, \pi/2 ]$ with this axis. Assume that in a first case, all points are aligned along the horizontal axis and in a second case, they are all aligned along $D_\theta$. Our normalized PCA gives two matrices 
\[
A = \begin{pmatrix}
1 & 0 \\ 
0 & 0
\end{pmatrix} 
\quad \text{and} \quad 
B = \begin{pmatrix}
\cos^2 \theta & \cos \theta \: \sin \theta \\ 
\cos \theta \: \sin \theta & \sin^2 \theta
\end{pmatrix} \: .
\]
One can check that $\| A - B \| = \sqrt{2} | \sin \theta |$ and the geodesic between $A$ and $B$ associated with the euclidean metric is 
\[
\gamma : t \in [0, 1] \mapsto (1-t) A + t B \: . 
\]
At time $t = 1/2$, one gets
\[
\gamma (t) = \frac{1}{2} \begin{pmatrix}
1 + \cos^2 \theta & \cos \theta \: \sin \theta \\ 
\cos \theta \: \sin \theta & \sin^2 \theta
\end{pmatrix} \: ,
\]
whose eigenvalues are $\frac{1}{2} (1 + |\cos \theta|)$ and $\frac{1}{2} (1 - |\cos \theta|)$. In particular, when $\theta$ is non-zero, this matrix is not anymore the covariance matrix of aligned points. From a geometric perspective, we expect that the geodesic rotates the axis from horizontal to $D_\theta$, keeping eigenvalues $1$ and $0$. 
When $A, B \in \xSym_+ (n)$, $\tr A = \tr B = 1$ encode points that are aligned along two vector spaces of same dimension $1 \leq d \leq n$, that is $A$ and $B$ are orthogonal projectors of rank $d$, 
this can be naturally achieved by taking the Riemannian metric in the Grassmannian $\G$. However, if we think of a geometric object, points are generally spread around a $d$-plane, because of both noise and curvature and the covariance matrix is not anymore proportional to an orthogonal projector.
\end{xmpl}

Geometric information reads more easily from the eigen decomposition than directly from the PCA matrix itself: for instance $d$--planes correspond with $d$ eigenvalues all equal to $1/d$ (after renormalizing the trace to $1$) the other $n-d$ eigenvalues being $0$, and the eigen space corresponding to $1/d$ gives the direction of the plane. More generally, the sequence of increasing eigenspaces, starting with the largest eigenvalue, gives a flag of relevant subspaces whose type is determined by the eigenvalues multiplicities.

In the next section, we rewrite the set of symmetric positive semi-definite matrices of trace $1$ as the set of possible eigen decompositions through the homeomorphism given in Proposition~\eqref{prop:homeo}. With this identification at hand, we will be able to work directly on the geometric information contained in the eigen decomposition and to build distances (see Section~\ref{section:distancesWF}) and metrics (see Section~\ref{section:RiemannianStructureFlags}) better suited than the euclidean distance in $\xM_n(\R)$.
Notice that such a correspondence actually involves a quotient set: ensuring the injectivity of the map \eqref{eqEigenDec} requires to identify eigen decompositions of the same matrix.

\subsection{Weighted flags as a quotient space}
\label{subsec:weigthedFlagTopo}

Given a matrix $A \in \xSym_+^1(n)$, we denote by $\lambda_1 (A) \geq \lambda_2 (A) \geq \ldots \geq \lambda_n(A)$ its eigenvalues and $F_\lambda(A)$ the eigenspace corresponding to eigenvalue $\lambda$.
We first define the set in which eigenvalues vary:
\[
\cW(n) = \left\lbrace (\lambda_1, \ldots, \lambda_n) \in [0,1]^n \: : \: \lambda_1 \geq \lambda_2 \geq \ldots \geq \lambda_n, \: \sum_{i=1}^n \lambda_i = 1 \right\rbrace \: .
\]
However it will be easier to work with weights in the $n$--simplex
\[
\Delta(n) = \left\lbrace (\mu_1, \ldots, \mu_n) \in [0,1]^n \: : \: \sum_{i=1}^n \mu_i = 1 \right\rbrace
\]
thanks to the linear correspondence
\begin{equation} \label{eq:lambdaToMu}
\begin{array}{ccl}
\cW(n) & \rightarrow & \Delta(n) \\
(\lambda_1, \ldots, \lambda_n) & \mapsto &\displaystyle (\mu_1, \ldots, \mu_n) = \left(\lambda_1 - \lambda_2, 2(\lambda_2 - \lambda_3), \ldots, k(\lambda_k - \lambda_{k+1}), \ldots, n(\lambda_n-0) \right)
\end{array}
\end{equation}
of inverse
\begin{equation} \label{eq:muToLambda}
\begin{array}{ccl}
\Delta(n) & \rightarrow & \cW(n) \\
(\mu_1, \ldots, \mu_n) & \mapsto &\displaystyle (\lambda_1, \ldots , \lambda_n) =  \left(\sum_{i=1}^n \frac{\mu_i}{i} , \ldots, \sum_{i=k}^n \frac{\mu_i}{i}, \ldots, \frac{\mu_n}{n} \right)
\end{array} \: .
\end{equation}
We consistently define for $A \in \xSym_+^1(n)$,
\begin{equation} \label{eqMuk}
 \mu_k (A) = k \left( \lambda_k(A) - \lambda_{k+1}(A) \right) \, \text{for } k = 1 \ldots n-1 \quad \text{and} \quad \mu_n(A) = n \lambda_n(A) \: .
\end{equation}
Using \eqref{eqMuk}, it is straightforward to compute the different dimensions of the eigenspaces. For instance, $A \in \xSym_+^1(n)$ is an orthogonal projector or rank $d$ (up to the factor $1/d$ due to its trace renormalized to $1$) if and only if $\mu_d(A) = 1$ and $\mu_i(A) = 0$ for all $i\neq d$. The flag of eigenspaces is then $\{ 0 \} \subset F_{1/d}(A) \subset \R^n = F_{1/d}(A) \oplus F_0(A)$ with type $(d, n-d)$. More generally, the type of the eigenspaces flag reads directly on the sequence $(\mu_1 (A), \ldots, \mu_n(A))$ by locating the non zero values of $\mu_i(A)$. For instance if $(\mu_1, \ldots, \mu_7) = (\underbrace{0, 0, 1/6}_{3}, \underbrace{1/2}_{1}, \underbrace{0, 1/3}_{2}, \underbrace{0}_{1})$, the type of the flag is $(3,1,2,1)$.

\begin{definition}[type in $\Delta(n)$] \label{dfn:typeDelta}
Let $\alpha = (\alpha_1, \ldots, \alpha_n) \in \Delta(n)$, we define its \emph{type} $\tau(\alpha)$ as follows:
\begin{itemize}
\item if $\alpha_n > 0$, let
$
r = \# \{ i \: : \: \alpha_i > 0 \}$ be the number of non zero $\alpha_i$ and $1 \leq d_1 < d_2 < \ldots < d_r = n$ all the corresponding indices. Then $\tau(\alpha) = (d_1, d_2-d_1, \ldots, d_r - d_{r-1})$ ;
\item if $\alpha_n = 0$, let $m = \max \{ i \: : \:  \alpha_i > 0 \}$ and let $(p_1, \ldots, p_{r-1})$ be the type of $(\alpha_1, \ldots, \alpha_m) \in \Delta(m)$ (in the previous sense, so that $p_1 + \ldots + p_{r-1} = m$), then $\tau(\alpha) = (p_1, \ldots, p_{r-1}, n-m)$.
\end{itemize}
For $A \in \xSym_+^1(n)$, we denote by $\tau (A)$ the type of $(\mu_1(A), \ldots, \mu_n(A))$.
%, notice that the type $\tau (A)$ only depends on the sequence of eigenvalues of $A$, or equivalently on $(\alpha_1, \ldots, \alpha_n)$, hence is well defined in the qutient $\cW \cF(n)$. In addition, it si not difficult to check that the type only depends on the indices $1 \leq d_1 < d_2 < \ldots < d_q \leq n$ of zero values of $\alpha_i$, allowing to use the notation $\tau (\alpha) = \tau (\alpha_1, \ldots, \alpha_n)$ for the type in the sequel.
%
\end{definition}

\begin{remk}
\label{remk:typeNonzeroValues}
The type of $\alpha$ is exactly determined by the indices of the nonzero values $\{ i \: : \: \alpha_i > 0 \}$. Moreover, given $\alpha, \beta \in \Delta(n)$,
\[
\{ i \: : \: \beta_i > 0 \} \subset \{ i \: : \: \alpha_i > 0 \} \quad \Longrightarrow \quad \tau(\alpha) \preccurlyeq \tau(\beta) \: .
\]
\end{remk}

Notice that for $A \in \xSym_+^1(n)$, the type $\tau(A)$ in the sense of Definition~\ref{dfn:typeDelta} is exactly the type of the eigenspaces flag in the sense of Definition~\ref{dfn:flag}.

\begin{definition}[Weighted flags] \label{dfn:weightedFlags}
We define the following equivalence relation $\sim$ in $\Delta(n) \times \xO(n)$: let $(\alpha, U) = \left( (\alpha_1, \ldots, \alpha_n), (u_1, \ldots, u_n) \right)$, $(\beta, V) = \left( (\beta_1, \ldots, \beta_n), (v_1, \ldots, v_n) \right) \in \Delta(n) \times \xO(n)$,
\medskip
\begin{align*} 
%\label{eqRelEquiv}
(\alpha, U) \, \sim & \, (\beta, V) \\
& \quad \Longleftrightarrow \quad \alpha = \beta \text{ and } \forall k = 1 \ldots n, \, \left\lbrace \begin{array}{l}
\xspan (u_1, \ldots, u_k) = \xspan (v_1, \ldots, v_k) \\
\text{or } \alpha_k = \beta_k = 0
\end{array} \right.\\
& \quad \Longleftrightarrow \quad \alpha = \beta \text{ and } \pi_I( U ) = \pi_I( V ) \text{ with } I = \tau(\alpha) = \tau(\beta)\\ 
& \quad \Longleftrightarrow \quad (\alpha, U) \text{ and } (\beta, V) \text{ correspond to the eigen decomposition of the same matrix.}
\end{align*}
We denote by $\cW \cF (n) = \faktor{\Delta(n) \times \xO(n)}{ \sim }$ the resulting quotient and $\pi_{\cW\cF} : \Delta(n) \times \xO(n) \rightarrow \cW\cF(n)$ the canonical projection. We refer to the equivalence classes of $\sim$ as \emph{weighted flags}.\\
Given $(\alpha, U) \in \Delta(n) \times \xO(n)$, we denote by $\left(\alpha, \pi_{\tau(\alpha)} ( U )\right)$ or simply $\left(\alpha, \pi(U) \right)$ its class in $\cW \cF (n)$.\\
We will often directly write a weighted flag as $(\alpha, W) \in \cW\cF(n)$ with $\alpha \in \Delta(n)$ and $W \in \cF_{\tau(\alpha)}$.
\end{definition}

\subsubsection*{Topology in weighted flags}
Providing $\Delta(n) \times \xO(n)$ with the topology induced by any norm in $\R^n \times \xM_n(\R)$ endows $\cW\cF(n)$ with the associated quotient topology: 

\begin{proposition}[Topology in weighted flags]
\label{prop:weightedFlagsCV}
Let $\left(\alpha^{(m)} , W^{(m)} \right)_{m \in \N}$, $(\alpha,W) \in \cW\cF(n)$. We introduce
\[
\begin{array}{lll}
I = \tau(\alpha), & U \in \xO(n) \text{ such that } \pi_I(U) = W & \text{ and} \quad (u_1, \ldots, u_n) = U \\
J^{(m)} = \tau(\alpha^{(m)} ), & U^{(m)} \in \xO(n) \text{ such that } \pi_{J^{(m)}} (U^{(m)}) = W^{(m)} & \text{ and} \quad (u_1^{(m)}, \ldots, u_n^{(m)}) = {U^{(m)} } \: .
\end{array}
\]
Then $\left(\alpha^{(m)} , W^{(m)} \right)$ converges to $(\alpha,W) \in \cW\cF(n)$,
denoted by $\displaystyle (\alpha^{(m)} , W^{(m)} ) \xrightarrow[m \to \infty]{\cW\cF(n)} (\alpha, W)$, if and only if 
\begin{equation}
\alpha^{(m)} \xrightarrow[m \to \infty]{} \alpha \text{ and } \forall k = 1 \ldots n, \text{ if } \alpha_k > 0, \, \xspan (u_1^{(m)}, \ldots, u_k^{(m)}) \xrightarrow[m \to \infty]{ {\rm G}_{k,n}} \xspan (u_1, \ldots, u_k) 
\label{eq:flagTopo}
\end{equation}
if and only if
\begin{equation*}
\alpha^{(m)} \xrightarrow[m \to \infty]{} \alpha \text{ and } p_{J^{(m)} \to I} (W^{(m)} ) \xrightarrow[m \to \infty]{\cF_I}  W \: .
\end{equation*}
\end{proposition}

Convergence in $\cF_I$ has been characterized in Proposition~\ref{prop:flagCV}. We have to clarify the meaning of $p_{J^{(m)} \to I} (W^{(m)} )$, since there is no reason why $J^{(m)}  = \tau(\alpha^{(m)})$ and $I = \tau(\alpha)$ would coincide. Let us investigate the connections between $\tau(\alpha^{(m)})$ and $ \tau(\alpha)$. Let $k \in \{1, \ldots,n\}$.
\begin{itemize}
 \item if $\alpha_k > 0$ then for $m$ large enough, $\alpha_k^{(m)} > 0$.
 \item if $\alpha_k = 0$ then $\alpha_k^{(m)}$ could be $0$ or $> 0$ even for $m$ large.
\end{itemize}
Consequently, for $m$ large enough, $\{ k \: : \: \alpha_k > 0 \} \subset \{ k \: : \: \alpha_k^{(m)} > 0 \}$ and thanks to Remark~\ref{remk:typeNonzeroValues} we have $J^{(m)} \preccurlyeq I$ whence the continuous projection $p_{J^{(m)} \to I}$ is well-defined (see Definition~\ref{dfn:typeCoarser}).

\begin{proposition} \label{prop:homeo}
 The following application is an homeomorphism
 \begin{equation*}
\left\lbrace \begin{array}{llll}
h : & \cW\cF(n) & \rightarrow & \xSym_+^1(n)  \\
& (\mu, \pi (U)) & \mapsto & U \xdiag ( \lambda_1, \ldots, \lambda_n) U^T
\end{array}
\right.
\end{equation*}
where $( \lambda_1, \ldots, \lambda_n) \in \cW(n)$ is defined from $\mu \in \Delta(n)$ by \eqref{eq:muToLambda}.
\end{proposition}

\begin{proof}
The application $h$ is defined from the continuous and surjective application 
\begin{equation*}
\left\lbrace \begin{array}{rlll}
g : & \Delta(n) \times \xO(n) & \rightarrow & \xSym_+^1(n)  \\
& (\mu,  U) & \mapsto & U \xdiag ( \lambda_1, \ldots, \lambda_n) U^T
\end{array}
\right.
\end{equation*}
that is constant with respect to the equivalence classes of $\sim$ so that $h$ is the unique application application from $\cW\cF(n)$ to $\xSym_+^1(n)$ such that $g = h \circ \pi_{\cW\cF}$ and $h$ is continuous by universal property of quotient topology. The injectivity of $h$ is precisely enforced by identifying all possible eigen decompositions of the same matrix that is exactly quotienting with respect to the equivalence relation $\sim$ (see Definition~\ref{dfn:weightedFlags}).

Let us check the continuity of $h^{-1}$. Let $(A^{(m)} )_{m\in \N}$ be a sequence in $\xSym_+^1(n)$ converging to $A \in \xSym_+^1(n)$ (for the topology induced by any matrix norm).
We define $\left( \mu^{(m)}, W^{(m)} \right) = h^{-1} \left( A^{(m)} \right)$ and $\left( \mu, W \right) = h^{-1} \left( A \right)$. By definition of $h$, $(\mu, W)$ is the class of possible eigen decomposition of $A$ (and similarly for $A^{(m)}$), the question is thus to check the continuity of the eigen decomposition near $A$ in $\xSym_+^1(n)$.
We first recall the continuity of the eigenvalues of real symmetric matrices\footnote{In the case of real symmetric matrices, the Lipschitz continuity holds for each ordered eigenvalue $A \mapsto \lambda_k(A)$ and follows from Weyl inequality: for $A, B \in \xSym(n)$, $|\lambda_k(A+B) - \lambda_k(A) | \leq \|B\|_{op}$.}:
%see also https://terrytao.wordpress.com/2010/01/12/254a-notes-3a-eigenvalues-and-sums-of-hermitian-matrices/#eigenperturb
for each $k = 1, \ldots, n$, $\lambda_k \left( A^{(m)} \right) \xrightarrow[m \to +\infty]{} \lambda_k(A)$ and thanks to \eqref{eq:lambdaToMu} $\mu_k^{(m)} \xrightarrow[m \to \infty]{} \mu$.
Let us introduce
\[
\begin{array}{lll}
I = \tau(\alpha), & U \in \xO(n) \text{ such that } \pi_I(U) = W & \text{ and} \quad (u_1, \ldots, u_n) = U \\
J^{(m)} = \tau(\alpha^{(m)} ), & U^{(m)} \in \xO(n) \text{ such that } \pi_{J^{(m)}} (U^{(m)}) = W^{(m)} & \text{ and} \quad (u_1^{(m)}, \ldots, u_n^{(m)}) = {U^{(m)} } \: .
\end{array}
\]
Let $k$ such that $\mu_k > 0$, according to \eqref{eq:flagTopo} it remains to prove that 
\[
 \xspan (u_1^{(m)}, \ldots, u_k^{(m)}) \xrightarrow[m \to \infty]{ {\rm G}_{k,n}} \xspan (u_1, \ldots, u_k) \: .
\]
We apply a result from \cite{Samworth}  (Theorem~2) originally proved in \cite{DavisKahan}.
% see also http://yueqicao.top/2021/01/12/Davis-Kahan-s-Theorem/
Let $1 \leq t < s \leq n$, 
\[
E = \xspan (u_t, u_{t+1}, \ldots , u_s) \quad \text{and} \quad E^{(m)} = \xspan (u_t^{(m)}, u_{t+1}^{(m)}, \ldots , u_s^{(m)}) \: .
\]
Assume that $\delta = \min \left( |\lambda_{t-1}(A) - \lambda_t(A) | , |\lambda_{s}(A) - \lambda_{s+1}(A) | \right) > 0$ (with the convention $\lambda_0 = +\infty$ and $\lambda_{n+1} = - \infty$) then, 
\[
\| E^{(m)} - E \|_F \leq 2\sqrt{2} \frac{\| A^{(m)} - A\|_F}{\delta}
\]
where $\| E^{(m)} - E \|_F$ stands for the Frobenius norm between the orthogonal projectors onto $E^{(m)}$ and $E$. We can apply this result with $t = 1$ and $s = k$. As such choices for $t$ and $s$ correspond to $\delta = | \lambda_k - \lambda_{k+1} | = \frac{\mu_k}{k} > 0$ (see \eqref{eqMuk}), we can conclude that
\[
\left( \mu^{(m)}, W^{(m)} \right) = h^{-1} \left( A^{(m)} \right) \xrightarrow[m \to +\infty]{\cW\cF(n)} \left( \mu, W \right) = h^{-1} \left( A \right) \: .
\]
Note that $\cW\cF(n)$ is a sequential topological space as the quotient of a metric space so that sequential continuity of $h^{-1}$ is equivalent to continuity.
% see https://en.wikipedia.org/wiki/Sequential_space
\end{proof}

Proposition~\ref{prop:homeo} states that so far $\cW\cF(n)$ and $\xSym_+^1(n)$ are equivalent as topological spaces. %Indeed, it is equivalent to pick up an element in $\cW\cF(n)$ or a matrix in $\xSym_+^1(n)$, 
Nevertheless, the presentation in $\cW\cF(n)$ directly gives access to the eigen decomposition and we will use this fact to build distances and even metrics strongly relying on the eigen elements and differing from the euclidean metric in $\xSym_+^1(n)$, keeping in mind Example~\ref{ex:euclDist}.

\section{Distances in weighted flags}
\label{section:distancesWF}

As a subset of $\xM_n (\R)$, $\xSym_+^1 (n)$ inherits the usual euclidean distance 
\[
d^{eucl} (A,B) = \sqrt{\tr ((A-B)(A-B)^T)} = \sqrt{\tr ((A-B)^2)} , \quad A,B \in \xSym_+^1 (n) \: .
\]
We could try transferring this structure to $\cW\cF(n)$, however,
the euclidean geodesic in $\xM_n (\R)$ between $A$ and $B$ is simply parametrized by
\[
\gamma(t) = (1-t) A + t B \in \xSym_+^1 (n), \quad t \in [0,1] \: .
\]
and as illustrated in Example~\ref{ex:euclDist}, even though the geodesic stays in $\xSym_+^1 (n)$, it does not respect the type of the flag of eigenspaces which can completely change along the geodesic. 
Building upon the homeomorphism evidenced in Proposition~\ref{prop:homeo} between the space of weighted flags $\cW\cF(n)$ and $\xSym_1^+(n)$, we intend to follow a converse path and transfer structure from $\cW\cF(n)$ to $\xSym_1^+(n)$, which hence requires to
push further the structure of $\cW\cF(n)$. Section~\ref{subsec:distancesWF} is a rather straightforward attempt to define a length space that however fails and is therefore unnecessary to the understanding of the rest of the paper. More precisely, we investigate in Section~\ref{subsec:distancesWF} two distances, compatible with the quotient topology of $\cW\cF(n)$ and based on the injection \eqref{eq:injGrassmannProduct} of $\cW\cF(n)$ into a larger quotient space, more precisely a product of cones, that can be turned into a length space. A distance arising from a length structure allows to consider shortest paths between points, unfortunately, the length structure does not restrict to $\cW\cF(n)$: the quotient defining $\cW\cF(n)$ is more complex than a product of cones. In Section~\ref{section:stratification}, we investigate further the stratification of $\cW\cF(n)$ with respect to the type of the flag: $\cW\cF(n)$ consists in gluing together smooth manifolds of the form $\mathring{\Delta}(r) \times \cF_{I}$, for $r = 1, \ldots, n$ (see \eqref{eq:diffeoCell}) along the $r$--skeleton of the simplex $\Delta(n)$.

\subsection{Distances in weighted flags}
\label{subsec:distancesWF}

As previously mentioned, though not necessary later in the paper, Section~\ref{subsec:distancesWF} illustrates how natural attempts to provide $\cW\cF(n)$ with a length space structure would fail.
A first attempt relies on the following observation. We recall that for a topological space $Y$, the quotient $\faktor{\left( [0,1] \times Y \right) }{ \left( \{0\} \times Y \right) }$ is the topological cone obtain by identifying all points in $\{0\} \times Y$.
Let us consider the continuous application
\[
\begin{array}{cll}
 \Delta(n) \times \xO(n) & \rightarrow &  \displaystyle \prod_{i=1}^n \left( \faktor{\left( [0,1] \times {\rm G}_{i,n} \right) }{ \left( \{0\} \times {\rm G}_{i,n} \right) } \right) \\
 (\alpha, U) & \mapsto &  \left(\alpha_i, \: \xspan (u_1, \ldots, u_i) \right)_{i = 1 \ldots n} 
\end{array} \: .
\]
that is constant on the equivalence classes of $\sim$ by Definition~\ref{dfn:weightedFlags}). It consequently induces the following continuous application on the quotient
%\begin{enumerate}[$\bullet$]
% \item $\{ \alpha_{d_k} \: : \: k = 1 \ldots r \} = \{ \alpha_i > 0 \}$ are the only non zero entries in $\alpha$, we recall that $d_k = p_1 + \ldots + p_k$.
% \item $(\alpha, [| U |])$ uniquely defines the vector spaces $\xspan (u_{d_k-p_k}, u_{d_k-p_k +1 }, \ldots , u_{d_k})$ for $k = = 1 \ldots r$.
% \item similarly, $(\alpha, [| U |])$ uniquely defines the vector spaces $E_{d_k}(U) := \xspan (u_1, u_2, \ldots , u_{d_k})$ for $k = 1 \ldots r$ by progressively taking the direct sum.
%\end{enumerate}
%In particular, it is possible to associate with  $(\alpha, [| U |])$ the family of vector spaces $(E_{d_k}(U))_{k = 1 \ldots r}$ 
%allowing to define the injection
\begin{equation} \label{eq:injGrassmannProduct}
\begin{array}{rcll}
f : & \cW\cF(n) & \rightarrow &  \displaystyle \prod_{i=1}^n \left( \faktor{\left( [0,1] \times {\rm G}_{i,n} \right) }{ \left( \{0\} \times {\rm G}_{i,n} \right) } \right)  \\
& (\alpha, W = \pi(U) ) & \mapsto &  \left(\alpha_i, \: \xspan (u_1, \ldots, u_i) \right)_{i = 1 \ldots n}
\end{array}
\end{equation}
which is injective, again by Definition~\ref{dfn:weightedFlags}. It is then natural to define a distance directly on the larger topological space $\displaystyle \prod_{i=1}^n \left( \faktor{\left( [0,1] \times {\rm G}_{i,n} \right) }{ \left( \{0\} \times {\rm G}_{i,n} \right) } \right)$ and then consider the induced distance on $\cW\cF$ through $f$. To this end, we can make use of known distances on a quotient of the form $\displaystyle \left( [0,1] \times X \right) / \left( \{0\} \times X \right)$ for a given metric space $(X,d)$ such as Krakus distance and the conic distance. In our case we consider $X = {\rm G}_{i,n}$ and we recall that 
the Riemannian distance $d_{{\rm G}_{i,n}}$ on ${\rm G}_{i,n}$ satisfies for $A, B \in {\rm G}_{i,n}$,
\[
d_{{\rm G}_{i,n}} (A,B) = \sqrt{\theta_1^2 + \ldots + \theta_i^2} \leq \frac{\pi}{2} \sqrt{i}
\]
where $\theta_1 \geq \ldots \geq \theta_i \in \left[ 0, \frac{\pi}{2} \right]$ are the principal angles between $A$ and $B$. We then fix ${d}_{i} = \frac{1}{\sqrt{i}} d_{{\rm G}_{i,n}}$ so that $\diam {\rm G}_{i,n} \leq 2$. In the following paragraphs, we give the expression of two distances on $\cW\cF(n)$ obtained by restriction of a distance on $\displaystyle \prod_{i=1}^n \left( \faktor{\left( [0,1] \times {\rm G}_{i,n} \right) }{ \left( \{0\} \times {\rm G}_{i,n} \right) } \right)$.

\subsubsection*{Krakus distance}
Let $(X, d)$ be a metric space such that ${\rm diam} \: X \leq 2$, the following application defines a distance on the quotient $\displaystyle Y = [0,1] \times X / \{0\} \times X$ (see Section~4 in \cite{Krakus})
\[
\begin{array}{lccll}
\widetilde{d} & : & ([0,1] \times X) \times ([0,1] \times X) & \rightarrow & \R_+ \\
           &   & (r,x),(s,y) & \mapsto &   |r-s| + \min(r,s) d(x,y) \: .
\end{array}
\]
Therefore, it is natural to consider the following distance $d_i^{kr}$ on $\displaystyle  \faktor{\left( [0,1] \times {\rm G}_{i,n} \right) }{ \left( \{0\} \times {\rm G}_{i,n} \right) } $: for $E, F \in {\rm G}_{i,n}$ and $r,s \in [0,1]$,
\[
d_i^{kr} ( (r,E), (s,F) ) = |r-s| + \min(r,s) d_i(E,F)
\]
and then $d^{kr} =  \sum_{i=1}^n d_i^{kr}$ provides a distance on $\displaystyle \prod_{i=1}^n \left( \faktor{\left( [0,1] \times {\rm G}_{i,n} \right) }{ \left( \{0\} \times {\rm G}_{i,n} \right) } \right)$. It is then possible to give the expression of $d^{kr}$ when restricted to $\cW \cF(n)$,
\begin{multline}
d^{kr} \left( (\alpha, \pi(U)) , (\beta, \pi(V)) \right)  =   
 \displaystyle \sum_{i=1}^n |\alpha_i - \beta_i | 
         + \min ( \alpha_i , \beta_i )  d_i \left(\xspan (u_1, \ldots, u_i) , \xspan (v_1, \ldots, v_i) \right)          
           \: .
\end{multline}

However, in order to work with shortest paths in $\cW \cF(n)$, we aim at defining an ``intrinsic'' distance, that is, associated with a length structure (see Definition~2.1.6 in \cite{Burago}). The following and so called conic distance is a natural distance transferring the intrinsic property of $(X,d)$ to the quotient $[0,1] \times X / \{0\} \times X$.

\subsubsection*{Euclidean cone over a length space}
Let $(X, d)$ be a metric space such that ${\rm diam} \: X \leq \pi$, the following application defines a distance on the quotient $\displaystyle Y = [0,1] \times X / \{0\} \times X$ (see Definition~3.6.12 in \cite{Burago})
\[
\begin{array}{lccll}
\widetilde{d} & : & ([0,1] \times X) \times ([0,1] \times X) & \rightarrow & \R_+ \\
          &   & (r,x),(s,y) & \mapsto & \left( r^2 + s^2 -2rs \cos d(x,y) \right)^{1/2} \\
          &   &             &         & = \left( |r-s|^2 + 2rs (1 - \cos d(x,y)) \right)^{1/2} \: .
\end{array}
\]
Moreover, if $(X, d)$ is a length space, then $(Y, \widetilde{d})$ is also a length space (see Theorem~3.6.17 in \cite{Burago}.
% intrisic metric is a metric associated with a length space
Therefore, it is natural to consider the following distance $d_i^{con}$ on $ \left( [0,1] \times {\rm G}_{i,n} \right) / \left( \{0\} \times {\rm G}_{i,n} \right)$: for $E, F \in {\rm G}_{i,n}$ and $r,s \in [0,1]$,
\[
d_i^{con} ( (r,E), (s,F) ) = \left( r^2 + s^2 -2rs \cos d_i (E,F) \right)^{1/2}
\]
and then $d^{con} = \left( \sum_{i=1}^n (d_i^{con})^2 \right)^\frac{1}{2}$ provides a distance that turns $ \prod_{i=1}^n \left( [0,1] \times {\rm G}_{i,n} \right) / \left( \{0\} \times {\rm G}_{i,n} \right)$ into a length space. It is possible to give the expression of $d^{con}$ when restricted to $\cW \cF(n)$,
\begin{multline}
d^{con} \left( (\alpha, \pi(U)) , (\beta, \pi(V)) \right)^2  =   
 \displaystyle \sum_{i=1}^n \alpha_i^2 + \beta_i^2 
         - 2 \alpha_i \beta_i \cos d_i (\xspan (u_1, \ldots, u_i) , \xspan (v_1, \ldots, v_i) )          
           \: .
\end{multline}
Unfortunately, the length space property does not transfer to $\cW\cF(n)$: given $(\alpha, \pi(U))$ and $(\beta, \pi(V))$ in $\cW\cF(n)$, even if one can define a shortest path between them in $\displaystyle \prod_{i=1}^n \left( [0,1] \times {\rm G}_{i,n} \right) / \left( \{0\} \times {\rm G}_{i,n} \right)$, there is no guarantee that the path would stay in $\cW\cF(n)$: the family of vector spaces in $({\rm G}_{i,n})_{i=1 \ldots n}$ does not necessarily stay a nested family (with respect to the type of the element). 

\begin{remk}
Notice that for both Krakus and conic distance, a sequence $(\mu^{(m)} , W^{(m)})_{m \in \N} \in \cW \cF(n)$ of weighted flags converges to $(\mu , W) \in \cW \cF(n)$ if and only if
\begin{multline*}
| \mu^{(m)} - \mu | \xrightarrow[m \to \infty]{} 0 \quad \text{and for all } i \in \{1, \ldots, n\} \text{ s.t. } \mu_i > 0, \\ d_{{\rm G}_{i,n}} \left( \xspan(w_1^{(m)}, \ldots, w_i^{(m)}), \xspan( w_1, \ldots, w_i)  \right) \xrightarrow[i \to \infty]{} 0  \\
\text{if and only if} \quad (\mu^{(m)} , W^{(m)}) \xrightarrow[m \to \infty]{\cW\cF(n)} (\mu , W)  \text{ by } \eqref{eq:flagTopo} \: .
\end{multline*}
In other words, both distances induce the same topology, which was the quotient topology, in $\cW\cF(n)$.
\end{remk}

\subsection{Stratification of weighted flags} \label{section:stratification}

A key issue with the conic distance introduced in the previous section stems from the loss of nestedness of the subspaces through the embedding \eqref{eq:injGrassmannProduct} of $\cW\cF(n)$  into $\prod_{i=1}^n \left( [0,1] \times {\rm G}_{i,n} \right) / \left( \{0\} \times {\rm G}_{i,n} \right)$. Loosely speaking, $\cW\cF(n)$ is not a ``simple'' product of cones but its structure is more involved. 
Indeed, $\cW\cF(n)$ glues together all possible flag spaces $\cF_I$ according to the structure of the $n$--simplex $\Delta(n)$, the type $I$ determines where to glue $\cF_I$, it can be on a face, edge, vertex and their general $n$--dimensional counterparts. Let us define the different elementary cells in the structure.

Given $r \in \{1, \ldots, n\}$, let $K = \{k_1, \ldots, k_r\} \subset \{1, \ldots, n\}$ be the set of indices of the positive $\mu_i$'s. We introduce
\begin{multline} \label{eq:stratifCell}
M(r \: ; K) := \left\lbrace ( \mu, W) \in \cW\cF(n) \: \left| \begin{array}{l}
\mu_j > 0 \text{ for } j \in K   \\
\mu_j = 0 \text{ for } j \notin K  
\end{array} \right. \right\rbrace \quad \text{and} \\ \quad \mathring{\Delta} ( n \: ; K ) = \left\lbrace \mu \in \Delta(n) \: \left| \begin{array}{l}
\mu_j > 0 \text{ for } j \in K   \\
\mu_j = 0 \text{ for } j \notin K  
\end{array}  \right. \right\rbrace \: .
\end{multline}
Let $(\mu,V) \in M(r \: ; K)$, we recall that the type of $\mu$ is entirely determined by the set $\{ j \: : \: \mu_j > 0 \} = K$ (see Remark~\ref{remk:typeNonzeroValues}) so that the type $\tau(\mu)$ is constant in $M(r ; K)$ and only depends on $K$. We denote by $I(K)$ this type, then $M(r ; K)$
is a smooth manifold and more precisely, we have the diffeomorphism 
\begin{multline} \label{eq:diffeoCell}
M(r \: ; K) = \mathring{\Delta} ( n \: ; K ) \times \cF_{I(K)} \simeq \left\lbrace \begin{array}{lcl}
\mathring{\Delta}(r) \times \cF_{I(K)} & \text{if} & r \geq 2, \\
\cF_{I(K)} & \text{if} & r = 1 
\end{array} \right. \: , \\
\quad \text{with } \mathring{\Delta}(r) = \mathring{\Delta}(r \: ; \{1, \ldots, r \}) = \Delta(r) \cap ]0,1[^r \: .
\end{multline}
The closure in $\cW\cF(n)$ of this set is then
\begin{align*}
\overline{ M(r \: ; K) } & =  \left\lbrace ( \mu, W) \in \cW\cF(n) \: | \:
\mu_j = 0  \text{ for } j \notin K  \right\rbrace \\
& = \bigsqcup_{\substack{K^\prime \in \mathcal{P}(\{1,\ldots,n\}) \\ K^\prime \subset K}} M(|K^\prime| \: ;  K^\prime ) \\
& = \bigsqcup_{r^\prime = 1}^r \, \bigsqcup_{\substack{K^\prime \subset K \\  |K^\prime| = r^\prime}} M(r^\prime \: ;  K^\prime ) 
\end{align*}
where $|K^\prime|$ is the cardinality of $K^\prime$. We can finally define the different strata of $\cW\cF(n)$ as follows: for $r \in \{1, \ldots, n\}$,
\[
 X_{r} = \bigcup_{\substack{K \in \mathcal{P}(\{1,\ldots,n\}) \\ |K| = r}} \overline{ M(r \: ; K) } = \bigsqcup_{\substack{K \in \mathcal{P}(\{1,\ldots,n\}) \\ |K| \leq r}} M(r \: ;  K ) 
\]
and we obtain the filtration
\[
\cW\cF(n) = X_{n} = \overline{M(n \: ; \: \{1, \ldots, n\})} \supset X_{n-1} \supset \ldots \supset X_2 \supset X_1 \: .
\]
As a conclusion, $\cW\cF(n)$ is a stratified space whose structure is given by the simplex $\Delta(n)$ and we note that for $r \in \{1, \ldots, n\}$,
\begin{enumerate}[$\bullet$]
 \item $X_{r}$ corresponds to weighted flags $(\mu, W) \in \cW\cF(n)$ with weights $\mu$ in the $r$--skeleton of the simplex $\Delta(n)$,
 \item similarly, $\displaystyle X_{r} \setminus X_{r-1}$ exactly corresponds to the $r$--dimensional cells of the simplex $\Delta(n)$, 
 \item $M(n \: ; \: \{1, \ldots, n\}) = X_{n} \setminus X_{n-1}$ is dense in $X_{n} = \cW\cF(n)$, we will shorten the notation to $M(n) = M(n \: ; \: \{1, \ldots, n\})$ hereafter.
 \item $M (r^\prime \: ; \: K^\prime ) \subset \overline{ M(r ; \: K) }$ if and only if $r^\prime \leq r$ and $K^\prime \subset K$ if and only if the type $I(K^\prime)$ associated with $K^\prime$ is coarser than the type $I(K)$ associated with $K$: $I(K) \preccurlyeq I(K^\prime)$ (see Definition~\ref{dfn:typeCoarser}) whence there is in this case a natural projection from $ M(r \: ; \: K)$ to $M (r^\prime \: ; \: K^\prime )$ induced by ${\rm id}_{\R^n} \times p_{I(K) \to I(K^\prime)}$.
\end{enumerate}

\section{Smooth and Riemannian structure of flags of a fixed type}
\label{section:SmoothAndRiemannianF_I}

We know from \eqref{eq:diffeoCell} that elementary cells in the structure of $\cW\cF(n)$ are diffeomorphic to a product of the form $\mathring{\Delta}(r) \times \cF_{I}$ that can hence be endowed with a Riemannian metric on the Cartesian product induced by Riemannian metrics on $\mathring{\Delta}(r)$ and $\cF_I$. The purpose of this section is to prepare the definition of a Riemannian metric in each of these cells $\mathring{\Delta}(r) \times \cF_{I}$ (in Proposition~\ref{prop:WFriemannianMetrics}), also based on the Riemannian metric on $\cF_I$ though different from the aforementioned product metric to be more consistent with the global structure of $\cW\cF(n)$. For this reasons, we recall useful facts concerning the smooth manifold structure (in Section~\ref{section:SmoothStructure}) and then the Riemannian structure of $\cF_I$ (see Section~\ref{subsec:RiemaniannF_I}), where the type $I = (p_1, \ldots, p_r)$ is fixed hereafter. 
As $\cF_I = \faktor{\xO(n)}{\xO(I)}$, its (smooth) Riemannian structure is inherited from the (smooth) Riemannian structure of $\xO(n)$.
More precisely, we recall that the tangent space to $\cF_I$ is isomorphic to a subspace of the tangent space to $\xO(n)$ called the horizontal space (see \eqref{eq:flagTgtSpace} and Remark~\ref{remk:horizontalSpace}) and that $\xC^1$ paths in $\cF_I$ can be lifted in $\xO(n)$ with the additional requirement that the tangent to the lifted path is horizontal (See Proposition~\ref{prop:horizontalLift}). We then recall  how a Riemannian metric in $\xO(n)$ induces a Riemannian metric $g^I$ in $\cF_I$ (Proposition~\ref{prop:flagRiemannianStruct}) and how both associated lengths $L_{\xO(n)}$ and $L_I$ compare (Proposition~\ref{prop:submersionLength}).
The Riemannian manifold structure of the set of flags of a given type $\cF_I$ has already been carefully investigated and we refer to \cite{Lim2019} for additional details completing this section. 

\subsection{Smooth structure}
\label{section:SmoothStructure}

% Thm 4.2 in Helgason and more generally chap II par 4 p.123
% p. 148-149 Mneimé Testard : espace tangent supplementaire de ker d\pi
% Thm 9.2 Boothby
% Thm 21.17 Lee
% Thm 1.97 GHL
We recall the following properties concerning the differentiable structure of the coset $\cF_I$ without proofs. They are consequence of the fact that $\cF_I$ is the coset of a compact Lie group by a closed subgroup. We refer to \cite{GHL}, \cite{Lee} and \cite{Helgason} for general results on such cosets. 
\begin{enumerate}[$\bullet$]
\item $\xO(n)$ is a Lie group and the tangent space to $\xO(n)$ at some $U \in \xO(n)$ is $T_U \xO(n) = U \xSkew(n)$.
\item $\xO(I)$ is a Lie group and the tangent space to $\xO(I)$ at some $W \in \xO(I)$ is 
\begin{align*}
T_U \xO(I) = U \xSkew(I) \text{ where }
 \xSkew(I) & = \left\lbrace {\rm diag}(A_1, \ldots, A_r) \: : \: \forall i = 1, \ldots, r, \, A_i \in \xSkew(p_i) \right\rbrace \\
 & \simeq \xSkew(p_1) \times \ldots \times \xSkew(p_r) \: .
\end{align*}
\item There exists a unique smooth structure in $\cF_I$ such that the quotient map $\pi_I : \xO(n) \rightarrow \cF_I$ is a smooth submersion (see Thm. 1.95 in \cite{GHL} or Thm. 21.17 in \cite{Lee}). Moreover, $\pi_I$ is a smooth fibration with fiber $\xO(I)$ and
\[
\ker T_U \pi_I = U \xSkew(I) \quad \text{and} \quad T_{\pi_I (U)} \cF_I \simeq \faktor{U \xSkew(n)}{U \xSkew(I)} \: .
\]
It is possible to decompose $\xSkew(n)$ into the direct sum $\xSkew(n) = \mathfrak{m}_{I} \: \oplus \:  \xSkew(I)$ with
\begin{align} \label{eq:mI}
%& = \left\lbrace \begin{array}{l}
%B \in \xSkew(n) \\
%B = (b_{ij})_{ij}
%\end{array}  : \: b_{ij} = 0 \text{ for } \left| \begin{array}{l}
% 1 \leq i,j \leq p_1 \text{ or}\\
% p_1 < i,j \leq p_1 + p_2 \\
% \text{ or } \ldots \text{ or}\\
%n-p_r < i,j \leq n
%\end{array} \right. \right\rbrace \oplus \left\lbrace {\rm diag}(A_1, \ldots, A_r) \: : \: A_k \in \xSkew(p_k) \right\rbrace \\
 \mathfrak{m}_{I}  = \left\lbrace \left[ \begin{array}{cccc}
0_{p_1} & B_{1,2} & \ldots & B_{1,r} \\ 
-B_{1,2}^T & 0_{p_2} & \ldots & B_{2,r} \\ 
\vdots & \vdots  & 0_{\ddots} & \vdots \\ 
- B_{1,r}^T & - B_{2,r}^T & \ldots & 0_{p_r}
\end{array} \right] : B_{i,j} \in \xM_{p_i,p_j} (\R)\right\rbrace  \text{ so that } T_{\pi_I(Id)} \cF_I \simeq \mathfrak{m}_{I} \: .
%\left\lbrace \left[ \begin{array}{cccc}
%A_1 & 0 & \ldots & 0 \\ 
%0 & A_2 & \ldots & 0 \\ 
%\vdots & \vdots & \ddots & \vdots \\ 
%0 & 0 & \ldots & A_r
%\end{array} \right] : A_k \in \xSkew(p_k) \right\rbrace
\end{align}
We then have $U \xSkew(n) = U \mathfrak{m}_{I} \: \oplus \:  U \xSkew(I)$ and consequently, for $W = \pi_I(U) \in \cF_I$,
\begin{equation} \label{eq:flagTgtSpace}
T_{\pi_I(U)} \cF_I \simeq U \:  \mathfrak{m}_{I} \: .
\end{equation}
Note that for $R = \xdiag (R_1, \ldots, R_r) \in \xO(I)$ and $B \in \mathfrak{m}_I$ the block matrix $B = ( B_{i,j} )_{i,j = 1 \ldots r}$ with $B_{i,j} \in \xM_{p_i,p_j}(\R)$ and $B_{j,i} = - B_{i,j}^T$, we have $R^T B R$ is the block matrix
\begin{equation} \label{eq:stabilitymI}
R^T B R = \left( R_i^T B_{i,j} R_j \right)_{i,j = 1 \ldots r} \quad \text{in particular} \quad R^T \: \mathfrak{m}_I \: R = \mathfrak{m}_I \: .
\end{equation}

\end{enumerate}

%
%As $G = \xO(n)$ is a Lie group and $H = \xO(p_1) \times \ldots \xO(p_r)$ is a closed subgroup of $G$, then $H$ acts smoothly, freely and properly on $G$ and the Quotient Manifold Theorem (see \cite{Lee}) asserts that there is a unique  structure of smooth manifold on  which is compatible with the quotient topological structure. 
%For $W \in \xO(n)$, we recall that the class $\left[| W \right|]$ is the class of all $W {\rm diag} (R_1, \ldots, R_r)$ for $R_i \in \xO(p_i)$. 
%

\begin{remk}[horizontal space] \label{remk:horizontalSpace}
Anticipating on the Riemannian structure, note that $H_U^I = U \:  \mathfrak{m}_{I}$ is not any complement but it is the orthogonal complement of $\ker T_U \pi_I = U \xSkew(I)$ in $T_U \xO(n)$: it is called the \emph{horizontal subspace} of $T_U \xO(n)$ and denoted by $H_U^I$ or simply $H_U$, while $V_U^I = \ker T_U \pi_I$ is called the \emph{vertical subspace} of $T_U \xO(n)$ hereafter. Note that if $J \preccurlyeq I$, we have $\mathfrak{m}_I \subset \mathfrak{m}_J$ and thus for all $U \in \xO(n)$, $H_U^I \subset H_U^J$.
\end{remk}
The coset $\cF_I = \faktor{\xO(n)}{\xO(I)}$ is a principal $\xO(I)$--bundle (see \cite{KobayashiNomizu} Example 5.1) which allows to lift $\xC^1$ paths in $\cF_I$ to $\xC^1$ paths in $\xO(n)$ that are \emph{horizontal paths} meaning that the tangent vector along the lifted path belongs to the horizontal space. Let us state this important fact, we refer to Prop. 3.1 in \cite{KobayashiNomizu}.

\begin{proposition}[Horizontal Lift] \label{prop:horizontalLift}
Let $\gamma : [a,b] \rightarrow \cF_I$ be a $\xC^1$ (resp. piecewise $\xC^1$) path and let $U \in \xO(n)$ such that $\pi_I(U) = \gamma(a)$. Then, there exists a unique $\xC^1$ (resp. piecewise $\xC^1$) path $\widetilde{\gamma} : [a,b] \rightarrow \xO(n)$ satisfying $\widetilde{\gamma}(a) = U$, $\pi_I \circ \widetilde{\gamma} = \gamma$ and
$\forall t \in [a,b]$ (resp. $\forall t \in [a,b] \setminus \{a_0, a_1, \ldots, a_N\}$), $\dt{\widetilde{\gamma}}(t) \in H^I_{\widetilde{\gamma}(t)} = \widetilde{\gamma}(t) \: \mathfrak{m}_{I}$. We will refer to $\widetilde{\gamma}$ as the \emph{horizontal lift of $\gamma$ starting at $U$}.
\end{proposition}

% https://math.stackexchange.com/questions/3524475/lifting-curves-from-riemannian-submersions-reconciling-claims-from-two-books
%
%https://mathoverflow.net/questions/419306/path-lifting-property-for-pim-rightarrow-m-g-for-g-compact-lie-acting-smoo
%
% 5.1 in Ortega-Ratiu and ShapeHMpaper.pdf p.4
%
% Mneimé-Testard 4.A.5

We refer to \cite{KobayashiNomizu} for the complete proof and just say a few words about the existence of a $\xC^1$ lift. First note that $\pi_I$ being a fibration, it admits local section which allows to locally lift a path and the compactness of $[a,b]$ allows to have finitely many $\xC^1$ pieces that may differ on overlaps. However, consider two such pieces $\gamma_1, \gamma_2 : ]t - \delta, t + \delta[ \rightarrow \xO(n)$ on the time overlap $]t - \delta, t + \delta[$, then $U_1 = \gamma_1(t)$ and $U_2 = \gamma_2(t)$ satisfy $\pi_I(U_1) = \pi_I(U_2)$ and thus, there exists $R \in \xO(I)$ such that $U_1 = U_2 R$, now changing $\gamma_2$ for $\gamma_2 R$ allows to obtain a $\xC^1$ path connecting both pieces. Iterating this process finitely many times produces a $\xC^1$ lift and it remains to deal with the differential system conveying the horizontality condition.

\subsection{Riemannian structure}
\label{subsec:RiemaniannF_I}

We recall that $\xO(n)$ can be provided with the Riemannian metric $g^{\xO(n)}$ induced by the euclidean metric in $\xM_n(\R)$ that is: $\forall U \in \xO(n)$, $\forall X, Y \in T_U \xO(n) = U \xSkew(n)$, there exist $B, C \in \xSkew(n)$ such that $X = UB$, $Y = UC$ and
\[
g_U^{\xO(n)} (X,Y) = \tr (X Y^T) = \tr (B C^T) = \sum_{i,j=1}^n b_{ij} c_{ij} \: .
\]
It is also possible to define another Riemannian metric on $\xO(n)$ by changing the initial euclidean one. In Section~\ref{section:RiemannianStructureFlags}, we will consider the following cases: given nonzero $(\Delta_{ij})_{i,j=1 \ldots n}$, one can define
\begin{equation} \label{eq:RiemannianMetricOn}
g_U^{\xO(n)} (X,Y) = \sum_{i,j=1}^n \Delta_{ij}^2 b_{ij} c_{ij} = \sum_{i,j=1}^n \Delta_{ij}^2 (U^T X)_{ij} (U^T Y)_{ij}  \: .
\end{equation}
The Riemannian structure on $\xO(n)$ can be transferred to the coset $\cF_I$ in the following sense (see 2.28 in \cite{GHL}):
\begin{proposition} \label{prop:flagRiemannianStruct}
Given a Riemannian metric $g^{\xO(n)}$ on $\xO(n)$ as in \eqref{eq:RiemannianMetricOn}, there exists on $\cF_I$ a unique Riemannian metric $g^I$ such that $\pi_I$ is a Riemannian submersion, i.e.
\begin{itemize}
 \item $\pi_I$ is a smooth submersion,
 \item for any $U \in \xO(n)$, $T_U \pi_I$ is an isometry between the horizontal space $H_U^I = U \:  \mathfrak{m}_{I}$ (see Remark~\ref{remk:horizontalSpace}) and $T_{\pi_I(U)} \cF_I$.
\end{itemize}
Moreover, the metric $g^I$ can be defined as follows: let $V \in \cF_I$ and $X, Y \in T_{V} \cF_I$. Given $U \in \xO(n)$ such that $V = \pi_I(U)$, there exist unique $\widetilde{X}, \widetilde{Y} \in U \:  \mathfrak{m}_{I}$ such that $X = T_U \pi_I \cdot \widetilde{X}$ and $Y = T_U \pi_I \cdot \widetilde{Y}$. Then
\begin{equation} \label{eq:flagMetric}
g^I_V(X,Y) = g_U^{\xO(n)} (\widetilde{X}, \widetilde{Y})
%= \tr (\widetilde{X} \widetilde{Y}^T) \: .
\end{equation}
\end{proposition}

\begin{remk}
Note that in the case where $g^{\xO(n)}$ is the usual Riemannian metric in $\xO(n)$ (i.e. $\Delta_{ij}$ are all equal to $1$ in \eqref{eq:RiemannianMetricOn}), we recover from \eqref{eq:flagMetric} the following expression for $g^I$: $g^I_V(X,Y) = \tr (\widetilde{X} \widetilde{Y}^T)$.
\end{remk}

\noindent We recall that $(\xO(n), g^{\xO(n)})$ and $(\cF_{I}, g^{I})$ being compact Riemannian manifolds, they induce complete metric spaces $(\xO(n), d^{\xO(n)})$ and $(\cF_{I}, d^{I})$. We start with investigating the length structure associated with $g^I$. Given a piecewise $\xC^1$ path $\gamma : [a,b] \rightarrow \cF_I$ (that is $\gamma$ is continuous and there exists $a = a_0 < a_1 < \ldots < a_N = b$ such that $\gamma_{| [a_{i-1},a_{i}]}$ is $\xC^1$ for $i = 1, \ldots, N$), the length of $\gamma$ is 
\begin{equation} \label{eq:lengthFI}
L_I(\gamma) = \int_a^b \sqrt{ g^I_{\gamma(t)} (\dt{\gamma}(t), \dt{\gamma}(t)) } \: dt \: .
\end{equation}
We recall that the distance $d_I$ in $\cF_I$ induced by the Riemannian metric $g^I$ through the length $L_I$ is then defined as follows: for $V_1, V_2 \in \cF_I$,
\[
d_I (V_1, V_2) = \inf \left\lbrace L_I (\gamma) \: : \: \gamma \text{ is a piecewise } \xC^1 \text{ path from } V_1 \text{ to } V_2 \right\rbrace \: .
\]

\begin{remk} \label{remk:dIquotientTopo}
Note that $d_I$ induced by the Riemannian metric $g^I$ induces the manifold topology, that is the quotient topology in $\cF_I$.
\end{remk}
% 13.29 lee_smooth_manifolds: Riemannian distance induces manifold topolgy.
%
\noindent We could similarly define the length $L_{\xO(n)}$ of a piecewise $\xC^1$ path in $(\xO(n), g^{\xO(n)})$ as well as the induced distance function $d_{\xO(n)}$.
We collect several useful connections between lengths $L_{\xO(n)}$ and $L_I$ in Proposition~\ref{prop:submersionLength} and distances $d_{\xO(n)}$ and $d_I$ in Proposition~\ref{prop:submersionContraction}.

\begin{proposition} \label{prop:submersionLength}
Let $g^{\xO(n)}$ be a Riemannian metric in $\xO(n)$ as in \eqref{eq:RiemannianMetricOn} and $g^I$ the induced metric in $\cF_I$ defined by Proposition~\ref{prop:flagRiemannianStruct}.
Let $\widetilde{\gamma} : [a,b] \rightarrow \xO(n)$ be a piecewise $\xC^1$ path and $\gamma_I := \pi_I \circ \widetilde{\gamma}$, then 
\begin{enumerate}[$(i)$]
\item $L_I (\pi_I \circ \widetilde{\gamma}) \leq L_{\xO(n)} (\widetilde{\gamma})$ and for all $t \neq a_i$, $g_{\gamma_I(t)}^{I}  \left( \dt{\gamma}_I (t), \dt{\gamma}_I (t) \right) \leq g_{\widetilde{\gamma}(t)}^{\xO(n)} (\dt{\widetilde{\gamma}}(t),\dt{\widetilde{\gamma}}(t))$.
\end{enumerate}
We additionally assume that $\widetilde{\gamma}$ is $I$--horizontal, i.e. $\forall t \neq a_i$, $\dt{\widetilde{\gamma}}(t) \in H_{\widetilde{\gamma}(t)}^I$. Then,
\begin{enumerate}[$(i)$]
\setcounter{enumi}{1}
 \item $L_I (\pi_I \circ \widetilde{\gamma}) = L_{\xO(n)} (\widetilde{\gamma})$,
 \item if $J \preccurlyeq I$ then $L_I (\pi_I \circ \widetilde{\gamma}) = L_J (\pi_J \circ \widetilde{\gamma})$,
 \item if $I \preccurlyeq J$ then $L_J (\pi_J \circ \widetilde{\gamma}) \leq L_I (\pi_I \circ \widetilde{\gamma})$,
\end{enumerate}
Let $\gamma : [a,b] \rightarrow \cF_I$ be a piecewise $\xC^1$--path, then 
\begin{enumerate}[$(i)$]
\setcounter{enumi}{4}
 \item there exists a piecewise $\xC^1$ $I$--horizontal path $\widetilde{\gamma} : [a,b] \rightarrow \xO(n)$ such that $\pi_I \circ \widetilde{\gamma} = \gamma$ and $L_I (\gamma) = L_{\xO(n)} (\widetilde{\gamma})$,
 \item if moreover $J \preccurlyeq I$, there exists a piecewise $\xC^1$ $I$--horizontal path $\widetilde{\gamma} : [a,b] \rightarrow \xO(n)$ such that $L_J (\pi_J \circ \widetilde{\gamma}) = L_I (\gamma)$ and $p_{J \to I} ( \pi_J \circ \widetilde{\gamma} ) = \gamma$.
\end{enumerate}
\end{proposition}

\begin{proof}
We first prove $(i)$. Take a piecewise $\xC^1$ path $\widetilde{\gamma} : [a,b] \rightarrow \xO(n)$, not necessarily horizontal.
We denote $\gamma_I = \pi_I \circ \widetilde{\gamma}$. For a fixed $t \in [a,b]$, $t \neq a_i$, decompose 
\[
\dt{\widetilde{\gamma}}(t) = U_H + U_V \text{ where } U_H \in H_{\widetilde{\gamma}(t)}^I \text{ and } U_V \in V_U^I = \ker T_{\widetilde{\gamma}(t)} \pi_I  \text{ are orthogonal}\: .
\]
We then have $\dt{\gamma}_I (t) = T_{\widetilde{\gamma}(t)} \pi_I \cdot \dt{\widetilde{\gamma}}(t) = T_{\widetilde{\gamma}(t)} \pi_I \cdot U_H$ which implies by \eqref{eq:flagMetric} that
\begin{align*}
g_{\gamma_I(t)}^{I} & \left( \dt{\gamma}_I (t), \dt{\gamma}_I (t) \right)  = g_{\widetilde{\gamma}(t)}^{\xO(n)} (U_H,U_H) \\
& \leq g_{\widetilde{\gamma}(t)}^{\xO(n)} (U_H,U_H) + g_{\widetilde{\gamma}(t)}^{\xO(n)} (U_V,U_V)  = g_{\widetilde{\gamma}(t)}^{\xO(n)} (\dt{\widetilde{\gamma}}(t), \dt{\widetilde{\gamma}}(t)) \text{ by orthogonality of } U_H \text{ and } U_V \: , 
\end{align*}
which implies
$L_I (\gamma_I) \leq L_{\xO(n)} (\widetilde{\gamma})$ and then $(i)$.

\noindent We now additionally assume that $\widetilde{\gamma}$ is $I$--horizontal, then, with the previous notations, $U_H = \dt{\widetilde{\gamma}} \in H_{\widetilde{\gamma}(t)}^I$ and $U_V = 0$ so that 
\[
g_{\gamma_I(t)}^{I} \left( \dt{\gamma}_I (t), \dt{\gamma}_I (t) \right) = g_{\widetilde{\gamma}(t)}^{\xO(n)} (\dt{\widetilde{\gamma}}(t), \dt{\widetilde{\gamma}}(t)) \quad \Longrightarrow \quad L_I (\gamma_I) = L_{\xO(n)} (\widetilde{\gamma}) \: .
\]
If $J \preccurlyeq I$, then
$\forall t \in [a,b] \setminus \{a_0, a_1, \ldots , a_N\}$ we have $\dt{\widetilde{\gamma}}(t) \in H_{\widetilde{\gamma}(t)}^I \subset H_{\widetilde{\gamma}(t)}^J$ and therefore $\widetilde{\gamma}$ is $J$--horizontal and we infer from $(i)$:
\[
L_I (\pi_I \circ \widetilde{\gamma}) = L_{\xO(n)} (\widetilde{\gamma}) = L_J (\pi_J \circ \widetilde{\gamma}) 
\]
If $I \preccurlyeq J$, similarly as for proving $(i)$, given a fixed $t \in [a,b]$, $t \neq a_i$, we have $H_U^J \subset H_U^I$ and the following orthogonal decompositions:
\begin{align*}
T_U \xO(n) & = H_U^I \oplus \ker T_{\widetilde{\gamma}(t)} \pi_I = H_U^J \oplus \ker T_{\widetilde{\gamma}(t)} \pi_J \\
H_U^I & = H_U^J \oplus \left( H_U^I \cap \ker T_{\widetilde{\gamma}(t)} \pi_J \right) \: .
\end{align*}
We can decompose 
\[
\dt{\widetilde{\gamma}}(t) = U_H + U_V \text{ where } U_H \in H_{\widetilde{\gamma}(t)}^J \text{ and } U_V \in \left( H_U^I \cap \ker T_{\widetilde{\gamma}(t)} \pi_J \right) \text{ are orthogonal}\: .
\]
Let $\gamma_J = \pi_J \circ \widetilde{\gamma}$. We then have $\dt{\gamma}_J (t) = T_{\widetilde{\gamma}(t)} \pi_J \cdot \dt{\widetilde{\gamma}}(t) = T_{\widetilde{\gamma}(t)} \pi_J \cdot U_H$ which implies by \eqref{eq:flagMetric} that
\[
g_{\gamma_J(t)}^{I} \left( \dt{\gamma}_J (t), \dt{\gamma}_J (t) \right) = g_{\widetilde{\gamma}(t)}^{\xO(n)} (U_H,U_H) \leq g_{\widetilde{\gamma}(t)}^{\xO(n)} (\dt{\widetilde{\gamma}}(t), \dt{\widetilde{\gamma}}(t)) \quad \Longrightarrow \quad L_J (\gamma_J) \leq L_{\xO(n)} (\widetilde{\gamma}) \: .
\]
We conclude the proof of $(iv)$ using $\widetilde{\gamma}$ $I$--horizontal and $(ii)$: $L_J (\gamma_J) \leq L_{\xO(n)} (\widetilde{\gamma}) = L_I (\gamma_I)$.

\noindent The last points $(v),(vi)$ directly follows from Proposition~\ref{prop:horizontalLift} taking an horizontal lift $\widetilde{\gamma}$ of $\gamma$ and applying $(ii)$ and $(iii)$, and recalling that $\pi_I = p_{I \to J} \circ \pi_J$ for $J \preccurlyeq I$.
\end{proof}

\begin{remk}[Riemannian submersion between $\cF_J$ and $\cF_I$]
We assume that $J \preccurlyeq I$ and we recall (see Definition~\ref{dfn:typeCoarser}) that $p_{J \to I}$ is the continuous application satisfying $\pi_I = p_{J \to I} \circ \pi_J$. Now adding that both $\pi_I$ and $\pi_J$ are smooth submersions (see Proposition~\ref{prop:flagRiemannianStruct}) then $p_{J \to I}$ is also a smooth submersion. Moreover given  $U \in \xO(n)$, we have
\[
H_U^I \subset H_U^J, \quad T_U \xO(n) = H_U^J \oplus V_U^J \quad \text{and} \quad H_U^J = H_U^I \oplus \left( H_U^J \cap V_U^I \right) \: .
\]
One can check that the following orthogonal decomposition holds
\[
T_{\pi_J(U)} \cF_J = T_U \pi_J (H_U^J) =  T_U \pi_J (H_U^I) \oplus \ker T_{\pi_J(U)} p_{J \to I}
\]
and $p_{J \to I}$ is a Riemannian submersion with horizontal space at $\pi_J(U)$ being $T_U \pi_J (H_U^I)$. The path $\pi_J \circ \widetilde{\gamma}$ given by Proposition~\ref{prop:submersionLength} $(vi)$ is thus a horizontal lift of $\gamma$ in $\cF_J$.

\centering
\begin{tikzpicture}[node distance=2.5cm, auto]
\node (A) {$T_U \xO(n) = H_U^I \oplus \left( H_U^J \cap V_U^I \right) \oplus V_U^J$};
\node(B) [right of=A, node distance=7cm] {$T_{\pi_I(U)} \cF_I = T_U \pi_I (H_U^I)$};
\node (C) [below of=A] {$T_{\pi_J(U)} \cF_J = T_U \pi_J (H_U^I) \oplus T_U \pi_J (H_U^J \cap V_U^I)$};
\draw[->](A) to node {$T_U \pi_I$}(B);
\draw[->](A) to node [left] {$T_U \pi_J$}(C);
\draw[->](C) to node [below=1.25ex, right] {$T_{\pi_J(U)} p_{J \to I}$}(B);
\end{tikzpicture}
\end{remk}

\begin{proposition} \label{prop:submersionContraction}
Let $g^{\xO(n)}$ be a Riemannian metric in $\xO(n)$ as in \eqref{eq:RiemannianMetricOn} and $g^I$ be the induced metric in $\cF_I$ defined in Proposition~\ref{prop:flagRiemannianStruct}.
For all $U_1, U_2 \in \xO(n)$,
\[
d_I (\pi_I(U_1), \pi_I(U_2)) \leq d_{\xO(n)} (U_1, U_2) \: .
\]
We assume $J \preccurlyeq I$ then for all $V_1, V_2 \in \cF_J$,
\[
d_I (p_{J \to I}(V_1), p_{J \to I} (V_2)) \leq d_{J} (V_1, V_2)  \: .
\]
\end{proposition}

In other words, the applications $\pi_I$ and $p_{J \to I}$ shorten distances, which is more generally true for Riemannian submersions.

\begin{proof}
Let $U_1, U_2 \in \xO(n)$ and $\widetilde{\gamma}$ be a piecewise $\xC^1$ path from $U_1$ to $U_2$ then $\gamma_I = \pi_I \circ \widetilde{\gamma}$ is a piecewise $\xC^1$ path from $\pi_I(U_1)$ to $\pi_I(U_2)$ so that by Proposition~\ref{prop:submersionLength} $(i)$,
\[
d_I (\pi_I(U_1), \pi_I(U_2)) \leq L_I (\gamma_I) \leq L_{\xO(n)} (\widetilde{\gamma}) \: ,
\]
and we conclude taking the infimum over all such paths $\widetilde{\gamma}$.

We assume $J \preccurlyeq I$. Let $V_1, V_2 \in \cF_J$ and $\gamma$ be a piecewise $\xC^1$ path from $V_1$ to $V_2$.
Applying Proposition~\ref{prop:submersionLength} $(v)$, take a piecewise $\xC^1$ $J$--horizontal lift $\widetilde{\gamma}$ between $U_1$ and $U_2 \in \xO(n)$, $\gamma = \pi_J \circ \widetilde{\gamma}$ and $V_1 = \pi_J(U_1)$, $V_2 = \pi_J(U_2)$. We then have that $\gamma_I = \pi_I \circ \widetilde{\gamma}$ is a piecewise $\xC^1$ path from $\pi_I(U_1)$ to $\pi_I(U_2)$ so that by Proposition~\ref{prop:submersionLength} $(iv)$,
\[
d_I(\pi_I(U_1), \pi_I(U_2)) \leq L_I(\gamma_I) \leq L_J(\gamma)
\]
As $\pi_I = p_{J \to I} \circ \pi_J$, we have for $i = 1,2$, $\pi_I(U_i) = p_{J \to I}(V_i)$ and we conclude taking the infimum over all such paths $\gamma$.

\end{proof}

\section{Weighted flags arising as the completion of a Riemannian structure}
\label{section:RiemannianStructureFlags}

We recall that we introduced, in Section~\ref{section:distancesWF}, distances compatible with the topology of weighted flags and nevertheless missing an associated length structure allowing to define shortest paths between weighted flags. We introduce in the present section a family of Riemannian metrics, well-defined in each cell  $M(r \: ;  K)$ of the stratification of $\cW\cF(n)$ and such that the metric completion of the elementary cell $M(n)$ exactly recovers $\cW\cF(n)$. More precisely, Section~\ref{subsec:RiemaniannStrata} starts with investigating $\cW\cF(2)$ and observing that it is a topological cone over $\mathbb{S}^1$ and it can hence be provided with a conical metric. The careful study of $\cW\cF(3)$ then leads to a better understanding of the stratification and more specifically of the 
relative inclusions of the horizontal spaces to the elementary cells (see Remark~\ref{rk:linkMuBij}, resulting in the definition of the Riemannian metrics $g$ in Proposition~\ref{prop:WFriemannianMetrics}. In Section~\ref{subsec:lenghtM}, we define a (global) length structure $\xLM$ in $\cW\cF(n)$ (see \eqref{eq:lengthStructure}) that extends the Riemannian length $L_g$ (see Proposition~\ref{prop:WFriemannianMetrics}$(iii)$) defined in $M(n)$. We prove that the respective metric space $(M(n), d_g)$ and $(\cW\cF(n) , \xdM)$ induced by the respective lengths $L_g$ and $\xLM$ agree as well in the sense of Theorem~\ref{thm:completion}: $(\cW\cF(n) , \xdM)$ is the metric completion of $(M(n), d_g)$. The proof consists in two main steps Proposition~\ref{prop:dMtopoEquivalent} and Proposition~\ref{prop:innerPathApproximation}. Proposition~\ref{prop:dMtopoEquivalent} first check that $\xdM$ induces the weighted flag topology. When comparing $(M(n), d_g)$ and $(M(n), \xdM)$, it is easy to see that $\xdM \leq d_g$ since the length structure $\xLM$ extends $L_g$, only allowing to consider more paths: it includes paths crossing the other strata $M( r \: ; K)$ of $\cW\cF(n)$.
Proposition~\ref{prop:innerPathApproximation} then shows that the converse inequality between $d_g$ and $\xdM$ is almost true: it is possible to approximate paths escaping $M(n)$ pointwisely with path staying in $M(n)$.
Section~\ref{subsec:lengthWFalternative} points out that it may be more natural to define a similar though different length structure $\xLWF$ (Definition~\ref{dfn:piecewiseC1pathII}) by considering concatenations of piecewise $\xC^1$ paths lying in the different elementary cells $M (r \: ; K)$, and not only in $M(n)$: such a distance allows to connect weighted flags belonging to the same elementary cell with a path itself remaining it the common cell, as a very specific consequence, it is possible to connect weighted flags in the same Grassmannian $\G \simeq M(1, \{d\})$ with a path of such weighted flags in $\G$. We show that such a length induces a distance $\xdWF$ in $\cW\cF(n)$ that coincides with $\xdM$ (Proposition~\ref{prop:xdMxdWF}).
\subsection{Riemannian metrics in the strata}
\label{subsec:RiemaniannStrata}

In this section, we define a Riemannian metric in each cell $M(r \: ;  K)$ (see  Eq. \eqref{eq:stratifCell}) of the stratification of $\cW\cF(n)$.
Recall that $K$ is the set of indices of the positive $\mu_i$'s, that is the indices where we change from one subspace to the next in the flag $\cF_{I(K)}$. 
Indeed, for $r = |K| \in \{2, \ldots, n \}$, though the strata $X_{r} \setminus X_{r-1}$ are not manifolds, as we have already seen in \eqref{eq:diffeoCell}, for an elementary cell we have the diffeomorphism
\begin{align*}
M(r \: ;  K)   =  \mathring{\Delta}(n \: ; K) \times \cF_{I(K)}  \simeq \mathring{\Delta}(r) \times \cF_{I(K)} \: .
\end{align*}
The inner $r$--simplex $\mathring{\Delta}(r) = \{ (x_1, \ldots, x_{r}) \in ]0,1[^{r} \: : \: x_1 + \ldots + x_{r} = 1 \}$ is a smooth manifold of dimension $r-1$.
The tangent space to $\mathring{\Delta}(n \: ; K)$ does not depend on the point $\mu \in \mathring{\Delta}(n \: ; K)$ and we have
\[
T_\mu \mathring{\Delta}(n \: ; K) = 
\left\lbrace \alpha \in \R^n \: \left| \begin{array}{l}
\alpha_{j}  =  0 \text{ for } j \notin K   \\
\alpha_1 + \ldots + \alpha_n = 0 
\end{array} \right. \right\rbrace \: .
\]
The Riemannian metric on $\mathring{\Delta}(n \: ; K)$ is then the one induced by the euclidean metric in $\R^n$:
\begin{equation} \label{eq:euclSimplex}
g^{\R^n}_\mu (\alpha,\beta) = g^{\R^n} (\alpha,\beta) = \sum_{i=1}^n \alpha_i \beta_i = \sum_{j \in K} \alpha_j \beta_j \: .
\end{equation}
Note that for $r = 1$, $K = \{k_1\}$ and $\mathring{\Delta}(n \: ; K) = \{ (0, \ldots, 0, 1, 0, \ldots, 0)\}$ (with $1$ corresponding to index $k_1$) is reduced to a point, in this case $T_\mu \mathring{\Delta}(n \: ; K) = \{ 0 \}$.
%
%%%%

\noindent It is then possible to define the product Riemannian metric on $ M(r \: ; K)$ as follows: let $(\mu, W ) \in M(r \: ; K) = \mathring{\Delta}(n \: ; K) \times \cF_{I(K)}$ and $( \alpha, X)$, $(\beta, Y) \in T_\mu {\mathring{\Delta}(n \: ; K)} \times T_{W} \cF_{I(K)}$, then, using \eqref{eq:euclSimplex} and \eqref{eq:flagMetric}, we could define
\begin{align*} 
g_{\mu , W}^{\text{product}} \left( ( \alpha, X), (\beta, Y) \right) & = g^{\R^n} (\alpha, \beta) + g_{W}^{I(K)} (X,Y)  \: .
\end{align*}
For $I = I(K)$, let $U \in \xO(n)$ such that $W = \pi_I(U)$, there exist unique $\widetilde{X}, \widetilde{Y} \in H_U^I = U \:  \mathfrak{m}_{I}$ such that $X = T_U \pi_I \cdot \widetilde{X}$ and $Y = T_U \pi_I \cdot \widetilde{Y}$, then $B = U^T \widetilde{X}, C = U^T \widetilde{Y} \in \mathfrak{m}_{I} \subset \xSkew(n)$ and (recalling \eqref{eq:RiemannianMetricOn})
\begin{align} 
g_{\mu , W}^{\text{product}} \left( ( \alpha, X), (\beta, Y) \right)  & = g^{\R^n} (\alpha, \beta) + g_W^{I} (X, Y) = g^{\R^n} (\alpha, \beta) + g_U^{\xO(n)} (\widetilde{X}, \widetilde{Y}) \nonumber \\
&=  \sum_{i=1}^{n} \alpha_i \beta_i + 2 \sum_{1 \leq i < j \leq n} \Delta_{ij}^2 b_{ij} c_{ij} \: . \label{eq:productMetric}
\end{align}
In order to continuously glue those metrics together according to the stratification, we are going to define a modified Riemannian metric that collapses consistently near the boundary $\overline{M(r \: ; K)} \setminus M(r \: ; K)$ of each cell. To this aim, we will adapt the coefficients $\Delta_{ij}$ of the metric in $\xO(n)$ with respect to the type of $(\mu, W)$. Loosely speaking $(\mu, W)$ tend to the boundary of a cell $M(r \: ; K)$ means that the type of $\mu$ changes and more precisely, some of the $\mu_k$ tend to $0$: we will observe that it is possible to let some $\Delta_{ij}$ tend to $0$ consistently. 
We investigate the cases $n = 2$ and then $n=3$ before coming to the general gluing, we keep the same notations as in \eqref{eq:productMetric} all along.
\subsubsection{The case $n=2$}
We have $3$ pieces: $M(2 \:  ; \{1,2\}) = \mathring{\Delta}(2 \: ; \{1,2\}) \times \cF_{(1,1)}$ is the principal stratum that can be identified with $\mathring{\Delta}(2) \times \mathbb{S}^1$. 
The two pieces are of lower dimension $M(1 \: ; \{1\}) = \{(1,0)\} \times \cF_{(1,1)} \simeq \{(1,0)\} \times \mathbb{S}^1$ and $M(1 \: ;\{2 \} ) = \{ (0,1)\} \times \cF_{(2)} = \{(0,1)\} \times \{ \R^2 \}$.
\begin{figure}[!htp]
\centering
\includegraphics[scale=0.55]{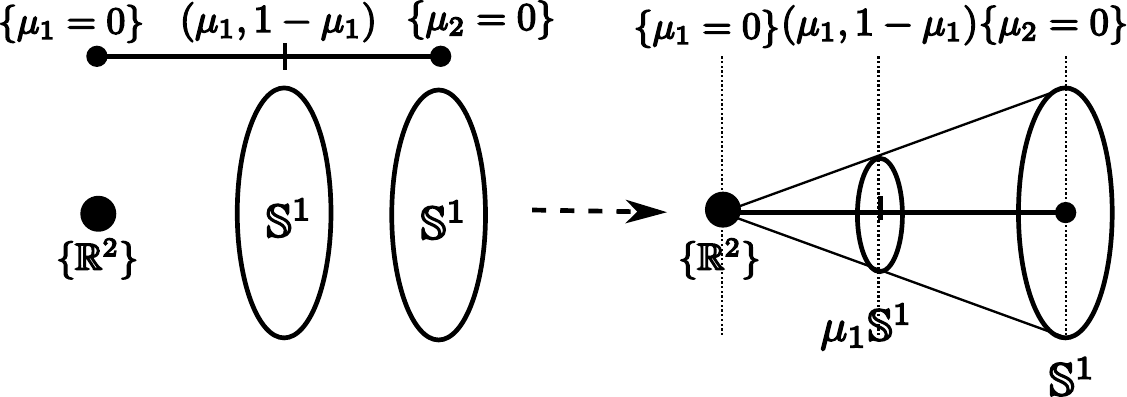}
\caption{Weighted flags in $\R^2$\label{weightedFlagsR2}}
\end{figure}
In order to respect the quotient structure (topology) of $\cW\cF(2)$, a conical metric sounds relevant (see Fig.\ref{weightedFlagsR2}):
\[
g_{\mu , W} \left( ( \alpha, X), (\beta, Y) \right) = \sum_{i=1}^{2} \alpha_i \beta_i + (\mu_1)^2 \sum_{1 \leq i < j \leq 2} b_{ij} c_{ij} =  2 \alpha_1 \beta_1 + 2 (\mu_1)^2 b_{12} c_{12} \: ,
\]
where we kept the same notations as in Eq. \eqref{eq:productMetric}, we recall that $\alpha_2 = -\alpha_1$ and $\beta_2 = - \beta_1$.

%Let us then consider the {\bf case $n=3$}. 
\subsubsection{The case $n=3$}
\begin{figure}[!htp]
\centering
\includegraphics[scale=0.65]{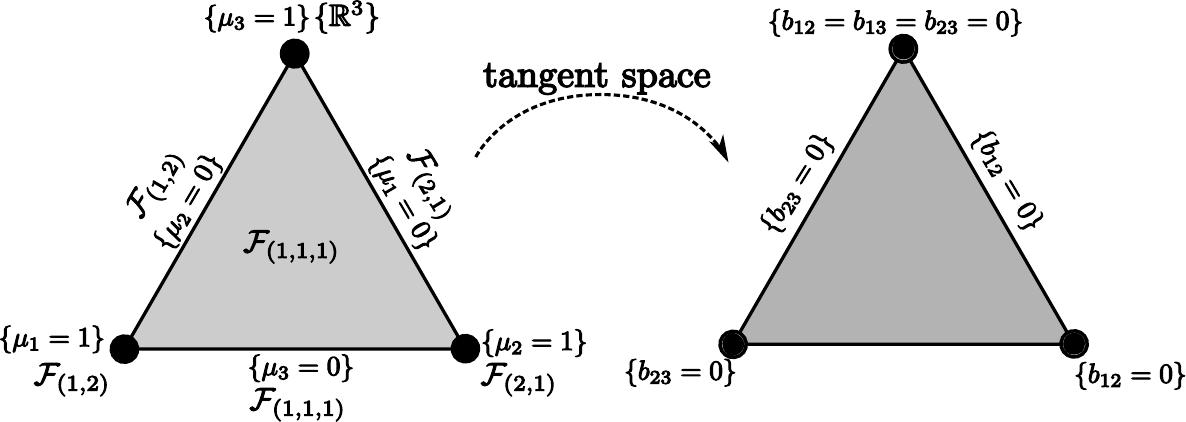}
\caption{Weighted flags in $\R^3$. Notice that the bottom line of the simplex corresponds to rank 2 matrices, while all points above encode rank 3 ones.\label{fig:n3}}
\end{figure}

.
We now have $7$ pieces (see Fig.~\ref{fig:n3}): the principal stratum with all non-zero weights $\mu_i$ is  $M(3 \: ; \{1,2,3\}) = \mathring{\Delta}(3 \: ; \{1,2,3\}) \times \cF_{(1,1,1)} \simeq \mathring{\Delta}(3) \times \cF_{(1,1,1)}$. Three other pieces correspond to exactly one of the $\mu_i$ vanishing: %is zero:
\begin{multline*}
M(2 \: ; \{2,3\}) = \mathring{\Delta}(2 \: ; \{2,3\}) \times \cF_{(2,1)}, \quad M(2 \: ; \{1,3\}) = \mathring{\Delta}(2 \: ; \{1,3\}) \times \cF_{(1,2)} \\ \quad \text{and} \quad M(2 \: ; \{1,2\}) = \mathring{\Delta}(2 \: ; \{1,2\}) \times \cF_{(1,1,1)} \: .
\end{multline*}
Notice that there is a special role for $\mu_3=0$ in this decomposition that differ from the vanishing of the other $mu_i$'s: this means than $\lambda_3=0$, that is a rank deficiency. However, the eigenspace is still the full flag $\cF(1,1,1)$ because all the eigenvalues are different.

To finish, we have three pieces corresponding to exactly two vanishing $\mu_i$'s, that is exactly one of the $\mu_i$ is equal to one, which precisely encodes for us  the corresponding Grassmannian ${\rm G}_{i,n} \simeq \cF_{(i,n-i)}$:
\begin{multline*}
M(1 \: ; \{1\}) = \{ (1,0,0) \} \times \cF_{(1,2)} \simeq {\rm G}_{1,3}, \quad M(1 \: ; \{2\}) = \{ (0,1,0) \} \times \cF_{(2,1)} \simeq {\rm G}_{2,3} \\ \quad \text{and} \quad M(1 \: ; \{3\}) = \{ (0,0,1) \} \times \cF_{(3)} \simeq {\rm G}_{3,3} = \{ \R^3 \} \: .
\end{multline*}

We can also explicit the tangent space at identity for each type of flag:
\begin{multline*}
\mathfrak{m}_{(1,1,1)} = \left\lbrace \left[ \begin{array}{ccc}
0 & b_{12} & b_{13} \\ 
-b_{12} & 0 & b_{23} \\ 
-b_{13} & -b_{23} & 0
\end{array} \right] : (b_{12}, b_{13}, b_{23}) \in \R^3 \right\rbrace, \quad \mathfrak{m}_{(3)} = \left\lbrace 0 \right\rbrace, \\
\mathfrak{m}_{(1,2)} = \left\lbrace \left[ \begin{array}{ccc}
0 & b_{12} & b_{13} \\ 
-b_{12} & 0 & 0 \\ 
-b_{13} & 0 & 0
\end{array} \right] : (b_{12}, b_{13}) \in \R^2 \right\rbrace, \, 
\mathfrak{m}_{(2,1)} = \left\lbrace \left[ \begin{array}{ccc}
0 & 0 & b_{13} \\ 
0 & 0 & b_{23} \\ 
-b_{13} & -b_{23} & 0
\end{array} \right] : (b_{13}, b_{23}) \in \R^2 \right\rbrace
\end{multline*}
We summarize those information in Figure~\ref{fig:n3}.
We observe that $b_{12} = 0 \Leftrightarrow \mu_1 = 0$, $b_{23} = 0 \Leftrightarrow \mu_2 = 0$ and $b_{13} = 0 \Leftrightarrow \mu_1 = \mu_2 = 0 \Leftrightarrow \mu_1 + \mu_2 = 0$. Gathering these information, this suggests to pinch the metric as follows:
\[
g_{\mu , W } \left( ( \alpha, X), (\beta, Y) \right) = \sum_{i=1}^{3} \alpha_i \beta_i + 2 (\mu_1)^2 b_{12} c_{12} + 2 (\mu_2)^2 b_{23} c_{23} + 2 (\mu_1 + \mu_2)^2 b_{13} c_{13}  \: .
\]
\\
\\
\\
\\

\subsubsection{The general case}\

\noindent
\begin{minipage}{0.68\textwidth}
%Coming back to the general case, 
\begin{remk}[Structural property of $\mathfrak{m}_I$]
\label{rk:linkMuBij}
Let $(\mu, W) \in \cW\cF(n)$ and let $I$ such that $(\mu, W) \in \Delta(n) \times \cF_I$. Let $B = (b_{ij})_{1 \leq i,j \leq n} \in \mathfrak{m}_I$,
we similarly observe that the whole diagonal block containing $b_{ii}$ up to $b_{jj}$ is zero if and only if
$\mu_i = \ldots = \mu_{j-1} = 0$,
if and only if $\mu_i + \ldots + \mu_{j-1} = 0$ for instance.
This suggest to modify $\eqref{eq:productMetric}$ by setting $\Delta_{ij} = \Delta_{ij} (\mu) = \mu_i + \ldots + \mu_{j-1}$:
\begin{equation*} 
g_{\mu , W} \left( ( \alpha, X), (\beta, Y) \right) = \sum_{i=1}^{n} \alpha_i \beta_i + 2 \sum_{1 \leq i < j \leq n} \left( \sum_{l = i}^{j-1} \mu_l \right)^2 b_{ij} c_{ij} \: .
\end{equation*}
\end{remk}
\end{minipage}
\begin{minipage}{0.29\textwidth}
\centering
\includegraphics[scale=0.45]{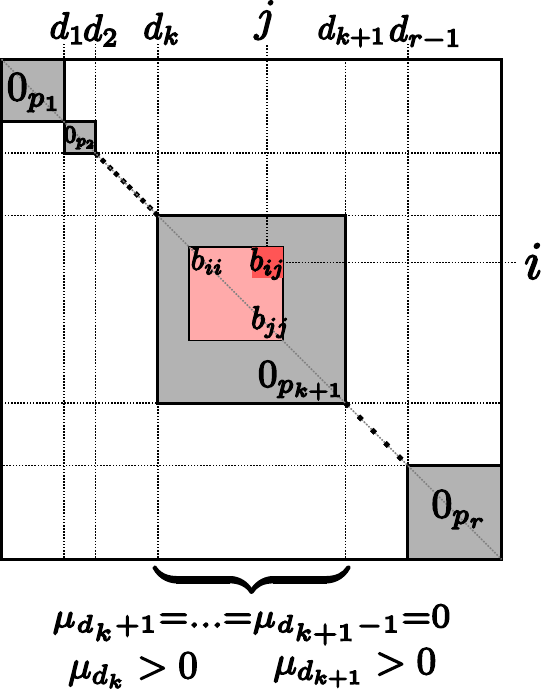}
\end{minipage}

\noindent We introduce a notation to access directly to the indices $(i,j)$ in the diagonal blocks in grey in the picture: given $I = (p_1, \ldots, p_r)$, the set of block diagonal indices is denoted by
\begin{equation} \label{eq:diagonalIndicesI}
X_I = \{ (i,j) \: : \: p_1 + \ldots + p_{k} + 1 \leq i < j \leq p_1 + \ldots + p_{k+1}  \text{ for some } k \in \{ 1, \ldots, r-1 \} \} \: .
\end{equation}

We provide in Proposition~\ref{prop:WFriemannianMetrics} a family of explicit metrics adapted to the stratification of $\cW\cF(n)$. A much more general construction, called \emph{iterated edge metric}, for abstract stratified spaces is investigated in \cite{Mondello}. Nevertheless, as we aim at performing numerical computations, we exhibit the explicit metric \eqref{eq:metric}.
\begin{proposition} \label{prop:WFriemannianMetrics}
Let $(\mu, W ) \in \cW\cF(n)$ and let $K \subset \{1, \ldots, n\}$ such that $(\mu, W ) \in M(|K| \: ; K) = \mathring{\Delta}(n \: ; K) \times \cF_{I(K)}$. We denote $I = I(K)$ and let $( \alpha, X)$, $(\beta, Y) \in T_\mu {\mathring{\Delta}(n \: ; K)} \times T_{W} \cF_{I}$. Then, for $U \in \xO(n)$ such that $W = \pi_I(U)$, there exist unique $\widetilde{X}, \widetilde{Y} \in H_U^I = U \:  \mathfrak{m}_{I}$ such that $X = T_U \pi_I \cdot \widetilde{X}$ and $Y = T_U \pi_I \cdot \widetilde{Y}$, we introduce $B = U^T \widetilde{X}, C = U^T \widetilde{Y} \in \mathfrak{m}_{I} \subset \xSkew(n)$ and we define
\begin{equation} \label{eq:metric}
g_{\mu , W} \left( ( \alpha, X), (\beta, Y) \right) = \sum_{i=1}^{n} \alpha_i \beta_i +  \sum_{1 \leq i < j \leq n} \left( f (\mu_{i \to j}) \right)^2 b_{ij} c_{ij} \: ,
\end{equation}
where $\mu_{i \to j} = (\underbrace{0, \ldots, 0}_{\in \R^{i-1}}, \underbrace{ \mu_i, \mu_{i+1}, \ldots, \mu_{j-1} }_{\in \R^{j-i}} , \underbrace{0, \ldots, 0}_{\in \R^{n-j+1}})$ and $f : [0,1]^n \rightarrow \R_+$ is Lipschitz continuous in $[0,1]^n$, $\xC^2$ in $[0,1]^n \setminus \{0\}$, and satisfies $f (\mu) = 0$ if and only if $\mu =0$. Then,
\begin{enumerate}[$(i)$]
 \item $g$ is well-defined,
 %and $(\mu,W) \mapsto g_{\mu,W}$ is continuous,
 %
 \item $g$ induces a Riemannian metric on each cell $M(r \: ; K)$ in the stratification of $\cW\cF(n)$.
 \item if $\gamma = (\mu,W) : [a,b] \rightarrow M(r \: ; K)$ is a piecewise $\xC^1$ path and $U : [a,b] \rightarrow \xO(n)$ is a piecewise $\xC^1$ lift of $W$, the length $L_g$ induced by $g$ satisfies
 \[
L_g(\gamma) = \int_a^b \sqrt{g_{\gamma(t)} (\dt{\gamma}(t), \dt{\gamma}(t))} \: dt = \int_a^b \sqrt{ |\dt{\mu}(t)|^2 + \sum_{1\leq i < j \leq n} f(\mu_{i \to j}(t))^2 \left( U(t)^T \dt{U}(t) \right)_{ij}^2 } \: dt \: .
 \]
\end{enumerate}
\end{proposition}

\begin{proof}
We first prove that $g$ is well-defined. Let $K \subset \{1, \ldots, n\}$, $|K|=r$, and $(\mu , W) \in M(r \: ; K) \subset \cW\cF(n)$. Let $I = I(K)$ such that $W \in \cF_I$, we have to prove that the second term in the right hand side of \eqref{eq:metric} does not depend on the choice of $U \in \xO(n)$ such that $\pi_I(U) = W$.

\noindent Let $X, Y \in T_{W} \cF_{I}$ and let $U, U^\prime \in \xO(n)$, $R \in \xO(I)$ satisfy $\pi_I(U) = \pi_I(U^\prime) = W$ and $U^\prime = UR$. Then, there exist unique $\widetilde{X}, \widetilde{Y} \in H_U^I = U \:  \mathfrak{m}_{I}$ and $\widetilde{X^\prime}, \widetilde{Y^\prime} \in H_{U^\prime}^I = U^\prime \:  \mathfrak{m}_{I}$  such that 
\[
X = T_{U} \pi_I \cdot \widetilde{X} = T_{U^\prime} \pi_I \cdot \widetilde{X^\prime} \quad \text{and} \quad Y = T_{U} \pi_I \cdot \widetilde{Y} = T_{U^\prime} \pi_I \cdot \widetilde{Y^\prime}  \: .
\]
We finally introduce $B = U^T \widetilde{X}$, $C = U^T \widetilde{Y} \in \mathfrak{m}_{I}$ and $B^\prime = {U^\prime}^T \widetilde{X^\prime}$, $C^\prime = {U^\prime}^T \widetilde{Y^\prime} \in \mathfrak{m}_{I}$.

$\bullet$ Let us check that 
\begin{equation} \label{eq:BprimeRB}
B^\prime = R^T B R \quad \text{and} \quad C^\prime = R^T C R \: . 
\end{equation}
The mapping $f_R : \xO(n) \rightarrow \xO(n)$ is smooth and satisfies for $V \in \xO(n)$, $H \in T_V \xO(n) = V \xSkew(n)$, $T_V f_R \cdot H = HR$. Moreover $\pi_I \circ f_R = \pi_I$ so that $T_U \pi_I = T_{U^\prime} \pi_I \circ T_U f_R$ and thus 
\[
X = T_{U} \pi_I \cdot \widetilde{X} = T_{U^\prime} \pi_I \cdot \widetilde{X} R \quad \text{and} \quad X = T_{U^\prime} \pi_I \cdot \widetilde{X^\prime} \quad \text{with} \quad \widetilde{X} R, \: \widetilde{X^\prime} \in H_{U^\prime}^I
\]
We recall (see \eqref{eq:stabilitymI}) that $m_I R = R \: m_I$ so that $\widetilde{X} R \in U \: \mathfrak{m}_I \: R = UR \mathfrak{m}_I = H_{U^\prime}^I$ as well as $\widetilde{X^\prime} \in H_{U^\prime}^I$ and by uniqueness we have $\widetilde{X^\prime} = \widetilde{X} R$. Therefore $B^\prime = {U^\prime}^T \widetilde{X^\prime} = R^T U^T \widetilde{X} R = R^T B R$ and the similar relation holds for $C^\prime$ and $C$.

$\bullet$ We write $I = (p_1, \ldots, p_r)$ and for all $k = 1 \ldots r$, $d_0 = 0$ and $d_k = d_{k-1} + p_k$. As $(\mu , W) \in \Delta(n) \times \cF_I$, we have for $k = 1 \ldots r-1$,
$\mu_{d_k} > 0$ while for $i \notin \{d_1, \ldots d_r\}$, $\mu_i = 0$ (we have no information on $\mu_{d_r} = \mu_n$). Then, given $1 \leq k \leq l \leq r$ then for all $(i,j)$ satisfying $d_{k-1} < i \leq d_k$ and $d_{l-1} < j \leq d_l$ the $n$--uplet $\mu_{i \to j}$ is constant and the non zero coefficients are $\mu_{d_k}$, $\mu_{d_{k+1}}$ up to $\mu_{d_{l-1}}$:
\begin{align}
 \mu_{i \to j} & =  (0, \ldots, 0,  \mu_i, \mu_{i+1}, \ldots, \mu_{j-1}  , 0, \ldots, 0) =  (0, \ldots, 0, \mu_{d_k}, 0, \ldots, 0 , \mu_{d_{k+1}},  \ldots, \mu_{d_{l-1}}, 0, \ldots, 0) \nonumber \\
 & = \left\lbrace \begin{array}{cl}
                    \mu_{d_k \to d_l} \neq 0 & \text{if } k < l \\
                    0 & \text{if } k = l
                  \end{array} \right. \: . \label{eq:muitoj}
\end{align}

$\bullet$ Let us write our matrices $B, B^\prime, C, C^\prime \in \mathfrak{m}_I$ blockwise with respect to the type $I$, for instance we write $B$ as
\[
 B = \left( B_{k,l} \right)_{k,l = 1 \ldots r} \quad \text{with} \quad B_{k,l} \in \xM_{p_k,p_l}(\R) \: .
\]
We also write $R = \xdiag (R_1, \ldots, R_r)$ and from \eqref{eq:BprimeRB}, we infer that for all $k,l = 1 \ldots r$,
\begin{equation} \label{eq:blockBprimeB}
B_{k,l}^\prime = R_k^T B_{k,l} R_l \quad \text{and} \quad C_{k,l}^\prime = R_k^T C_{k,l} R_l \quad \Longrightarrow \quad
B_{k,l}^\prime {C_{k,l}^\prime}^T = R_k^T B_{k,l} C_{k,l}^T R_k \: .
\end{equation}
Then, using \eqref{eq:muitoj} and writing the sum in \eqref{eq:metric} with respect to the blockwise aforementioned decomposition, we obtain
\begin{align}
\sum_{1 \leq i < j \leq n} \left( f (\mu_{i \to j}) \right)^2 b_{ij} c_{ij} & = \sum_{1 \leq k < l \leq r} \left( f (\mu_{d_k \to d_l}) \right)^2 \sum_{\substack{i = d_{k-1} + 1, \ldots, d_k \\ j = d_{l-1} + 1, \ldots, d_l}} b_{ij} c_{ij} \nonumber \\
& = \sum_{1 \leq k < l \leq r} \left( f (\mu_{d_k \to d_l}) \right)^2 \underbrace{ \tr (B_{k,l} C_{k,l}^T) }_{= \tr (B_{k,l}^\prime {C_{k,l}^\prime}^T) \text{ by \eqref{eq:blockBprimeB}} } \label{eq:gScalarProduct} \\
& = \sum_{1 \leq i < j \leq n} \left( f (\mu_{i \to j}) \right)^2 b_{ij}^\prime c_{ij}^\prime \nonumber
\end{align}
We can eventually conclude that $g$ is well-defined. Furthermore, it is not difficult to check that $g_{\mu,W}$ defines a scalar product on $T_{\mu, W} M(r \: ; K)$, using for instance the above expression \eqref{eq:gScalarProduct}  together with $\mu_{d_k \to d_l} \neq 0$ for $k<l$ (from \ref{eq:muitoj}) that implies $(f (\mu_{d_k \to d_l}))^2 > 0$.

$\bullet$ It remains to check that given two vector fields $X = X_{\mu,W}$, $Y = Y_{\mu,W}$, the application $(\mu, W) \in M(r \: ; K) \mapsto g_{\mu,W} (X,Y)$ is smooth. As $\pi_I : \xO(n) \rightarrow \cF_I$ is a fibration (see Section~\ref{section:SmoothStructure}), there is a local smooth section $S : \cV \rightarrow \xO(n)$ defined on an open set $W \in \cV \subset \cF_I$. Then, $V \in \cV \mapsto p_{H_{S(V)}^I}$ mapping $V$ to the orthogonal projector onto the horizontal space $H_{S(V)}^I$ is smooth as well as $V \mapsto \widetilde{X} = p_{H_{S(V)}^I} \:  T_V S \cdot X $ and $V \mapsto B = S(V)^T \: \widetilde{X}$. The smoothness of $g$ then follows from expression \eqref{eq:metric}.
% See end of proof of 2.28 p 61 in  \cite{GHL}.

$\bullet$ Concerning the expression of $L_g$ given in $(iii)$, from the definition of $g$ we have
 \begin{equation} \label{eq:tildeXUH}
L_g(\gamma) = \int_a^b \sqrt{ |\dt{\mu}(t)|^2 + \sum_{1\leq i < j \leq n} f(\mu_{i \to j}(t))^2 \left( U(t)^T \widetilde{X}(t) \right)_{ij}^2 } \: dt \: .
 \end{equation}
and we are left with the identification of $\widetilde X \in H_U^I$ such that $\dt{W} = T_U \pi_I \cdot \widetilde{X}$. Notice that $W = \pi_I(U)$ in $[a,b]$ so that $\dt{W} = T_U \pi_I \cdot \dt{U}$ but in general we do not have $\dt{U} \in H_U^I$. However, writing $\dt{U} = (\dt{U})_H + (\dt{U})_V$ with $(\dt{U})_H \in H_U^I = U \mathfrak{m}_I$ and $(\dt{U})_V \in \ker T_U \pi_I = U \xSkew(I)$, we have $\dt{W} = T_U \pi_I \cdot (\dt{U})_H$, so that $\widetilde{X} = (\dt{U})_H$ and
\begin{equation} \label{eq:lengthLiftedUij}
\left( U^T \dt{U} \right)_{ij} = \left( U^T (\dt{U})_V \right)_{ij} + \left( U^T (\dt{U})_H \right)_{ij} = \left\lbrace  \begin{array}{ccl}
\left( U^T (\dt{U})_V \right)_{ij} & \text{if} & (i,j) \in X_I \\
\left( U^T (\dt{U})_H \right)_{ij} & \text{if} & (i,j) \notin X_I 
\end{array} \right. \: .
\end{equation}
We can conclude observing that $(i,j) \in X_I \Leftrightarrow f(\mu_{i \to j}) = 0$ and then we can equivalently replace $\widetilde{X}$ with $(\dt{U})_H$ or simply $\dt{U}$ in \eqref{eq:tildeXUH}.
\end{proof}

\subsubsection*{Matrix metric tensor in a smooth frame}

Let $K \subset \{1, \ldots, n \}$, $r = |K| \geq 1$, $k_0 \in K$ and $I = I(K)$ such that $M(r \: ; K) = \mathring{\Delta}(r \: ; K ) \times \cF_I$. Let $(\mu, W) \in M(r \: ; K)$.

\noindent Recalling that
\[T_\nu \mathring{\Delta}(n,K) = \{ \alpha \in \R^n \: : \: \alpha_j = 0 \text{ for } j \notin K \text{ and } \alpha_1 + \ldots + \alpha_n = 0 \} 
\]
is a subspace of dimension $r-1$, one can fix an orthonormal basis $(f_1, \ldots, f_{r-1})$ of $T_\nu \mathring{\Delta}(n,K)$.

\noindent We introduce the orthogonal basis $\left(\Omega^{ij} \right)_{1 \leq i  < j \leq n}$ of $\xSkew(n)$ where the matrix $\Omega^{ij}$ has two nonzero coefficients: its coefficient $(i,j)$ is $1$ and its coefficient $(j,i)$ is $-1$.
As $\pi_I : \xO(n) \rightarrow \cF_I$ is a fibration (see Section~\ref{section:SmoothStructure}), there is a local smooth section $S : \cV \subset \cF_I \rightarrow \xO(n)$ defined on an open neighbourhood $\cV$ of $W$. Moreover, $T_{S(V)} \pi_I : H_{S(V)}^I = S(V) \mathfrak{m}_I \rightarrow T_V \cF_I$ is an isometry and thus, defining for $1 \leq i < j \leq n$, $X^{ij} : V \in \cV \subset \cF_I \mapsto T_{S(V)} \pi_I \cdot S(V) \Omega^{ij}$, we obtain that 
\begin{equation} \label{eq:smoothFrame}
\left( ( f_k, 0 )_{1 \leq k \leq r-1} , (0, X^{ij} )_{\substack{1 \leq i < j \leq n \\ (i,j) \notin X_I}} \right)
\end{equation} is a smooth frame of $M(r \: ; K)$ satisfying for $(\nu , V) \in \mathring{\Delta}(r \: ; K ) \times \cV \subset M(r \: ; K)$:
\[
g_{\nu, V} \left( ( 0, X^{ij} ) , ( 0 , X^{i'j'}) \right) = \left\lbrace \begin{array}{lcl}
0                   & \text{if} & (i,j) \neq (i', j') \\                                                               
f( \nu_{i \to j})^2 & \text{if} & (i,j) = (i',j')                                                                       \end{array} \right. \: .
\]
Furthermore, for $1 \leq k, l \leq r-1$,
\[
g_{\nu, V} \left( ( f_k,0) , ( f_l , 0) \right) = f_k \cdot f_l = \left\lbrace \begin{array}{lcl}
0 & \text{if} & k \neq l \\                                                               
1  & \text{if} & k = l                                                                       \end{array} \right. \: .
\]

\noindent In particular, the metric tensor of the metric $g$ in the smooth frame \eqref{eq:smoothFrame} is the following diagonal matrix:
{\tiny{\[ \begin{array}{c|ccc|ccc|}
        & f_1 &  \ldots & f_{r-1} & & (X^{ij} )_{\substack{1 \leq i < j \leq n \\ (i,j) \notin X_I}} & \\ 
 \hline 
 f_1    &       &        &       &   & & \\ 
 \vdots &    & I_{r-1}  &  & & 0 & \\ 
 f_{r-1}     &      &        &      &  & & \\ 
\hline
   &          &      &      &   \ddots   & &    0 \\ 
    &          &      &      &      & &    \\ 
 X^{ij}    &          &   0   &      & & f(\nu_{i \to j})^2    &   \\ 
     &         &     &     &     &  &    \\ 
        &          &      &      &  0   & &  \ddots  \\ 
 \hline
 \end{array} \]}
}
%
%As a consequence, the Riemannian measure on $M(n)$ is 
%\[
%\sqrt{n} \prod_{1 \leq i < j \leq n} \left(\mu_i + \ldots + \mu_{j-1} \right)
%\]
%$\displaystyle \frac{(n-1)(n+2)}{2}$--dimensional.

\subsection{Length structure in weighted flags}
\label{subsec:lenghtM}

The goal of next section is to show that the Riemannian metric $g$ defined in \eqref{eq:metric} induces a metric space $(M(n), d_g)$ whose completion is exactly $(\cW \cF(n) , \xdM)$ (see \eqref{eq:RiemannianDistance} below for the definition of $\xdM$ and see Theorem~\ref{thm:completion} for the statement). The metric space $(\cW \cF(n) ,  \xdM)$ is a length space associated with the length structure \eqref{eq:lengthStructure}.

We recall that
\begin{itemize}
 \item on one hand, $M(n) = \mathring{\Delta}(n) \times \cF_{(1,\ldots,1)}$ is dense in $\cW\cF(n)$ (with respect to the quotient topology in $\cW\cF(n)$),
 \item on the other hand, it can be completed into the manifold with boundary $\widetilde{M(n)} = \Delta(n) \times \cF_{(1,\ldots,1)}$, the topology in this case being the product one in $\Delta(n) \times \cF_{(1,\ldots,1)}$ and not the weighted flag topology.
\end{itemize}
Loosely speaking, one can picture the difference between $\overline{M(n)}$ and $\widetilde{M(n)}$ by starting from $]0,1] \times \mathbb{S}^1$ and closing it as a cone $\faktor{[0,1] \times \mathbb{S}^1}{\{0\} \times \mathbb{S}^1}$ or as a cylinder $[0,1] \times \mathbb{S}^1$. We can then define piecewise $\xC^1$ paths in $\cW\cF(n)$.

\begin{definition}[piecewise $\xC^1$ path - I] \label{dfn:piecewiseC1pathI}
Let $\gamma : [a,b] \subset \R \rightarrow \cW\cF(n)$ be continuous.
We first say that $\gamma$ is a piecewise $\xC^1$ path in $\overline{M(n)} = \cW\cF(n)$ if there exist $a = a_0 < a_1 < \ldots < a_{N-1} < a_N = b$ such that for all $i$, $\gamma(]a_{i-1}, a_i[) \subset M(n)$ and $\gamma_{| ]a_{i-1}, a_i[}$ extends on $[a_{i-1}, a_i]$ into a $\xC^1$ path in the smooth manifold with boundary $\widetilde{M(n)} = \Delta(n) \times \cF_{(1, \ldots, 1)}$.
\end{definition}

\noindent We define the length $\xLM$ of a piecewise $\xC^1$ path $\gamma : [a,b] \rightarrow \cW\cF(n)$ as
% \xLM
\begin{equation} \label{eq:lengthStructure}
\xLM(\gamma) = \int_a^b \sqrt{ g_{\gamma(t)} (\dt{\gamma}(t), \dt{\gamma}(t)) } \: dt
\end{equation}
as well as the associated distance
\begin{equation} \label{eq:RiemannianDistance}
\xdM \left((\mu, W  ), (\nu, Q ) \right) = \inf \left\lbrace \xLM (\gamma) \: \left| \begin{array}{l}
\gamma : [0,1] \rightarrow \cW\cF(n) \text{ piecewise } \xC^1 \text{ path} \\
\gamma(0) = (\mu,  W  ) \text{ and } \gamma(1) = (\nu,  Q  )
\end{array} \right. \right\rbrace \: .
\end{equation}

\begin{remk} \label{remk:lengthDistanceMg}
Considering the Riemannian structure $(M(n), g)$ defined in Proposition~\ref{prop:WFriemannianMetrics}, as well as the induced length $L_g$ and distance $d_g$, the following observations are consequences of the definitions:
\begin{itemize}
\item for piecewise $\xC^1$ paths contained in $M(n)$, $L_g$ and $\xLM$ coincide,
\item given $(\mu,W), \: (\nu,Q) \in M(n)$,
$
\xdM ((\mu,W), (\nu,Q)) \leq  d_g ((\mu,W), (\nu,Q)) 
$
while the converse inequality is not obvious and will be proven in Proposition~\ref{prop:innerPathApproximation}. The obstruction comes from the fact that it could be more economic in terms of length to use a path not entirely contained in $M(n)$ even though we try to connect points in $M(n)$.
\end{itemize}
\end{remk}

It is not difficult to check that $\xLM$ defines a length structure in $\cW\cF(n)$ in the sense of \cite{Burago} 2.1.1. Indeed, $\xLM$ is additive, $t \mapsto \xLM (\gamma_{| [a,t]})$ is continuous, $\xLM$ is invariant under reparametrization (change of variables). The last requirement is that $\xLM$ agrees with the (quotient) topology of $\cW\cF(n)$, that is if $(\mu, W) \in \cW\cF(n)$ and $\cV$ is an open set containing $(\mu, W)$, then the length from $(\mu, W)$ to $\cW\cF(n) \setminus \cV$ is strictly positive, more precisely:
\begin{equation} \label{eq:topoAgreement}
\inf \left\lbrace \xLM (\gamma) \: : \: \gamma(a) = (\mu, W) \text{ and } \gamma(b) \notin \cV \right\rbrace > 0 \: .
\end{equation}
In Proposition~\ref{prop:dMtopoEquivalent}, we prove that $\xdM$ induces the quotient topology in $\cW\cF(n)$ which is even stronger than \eqref{eq:topoAgreement}.
% \eqref{eq:topoAgreement} is indeed weaker than same topology, see \cite{Burago} ex 2.1.5 p28.

\begin{proposition} \label{prop:dMtopoEquivalent}
Let $\left( \mu^{(m)} ,  W^{(m)}  \right)_{m \in \N}$ be a sequence in $\cW\cF(n)$ and let $(\mu, W ) \in \cW\cF(n)$. Then,
\[
\xdM \left(( \mu^{(m)} ,  W^{(m)}   ), (\mu, W ) \right) \xrightarrow[m \to +\infty]{} 0 \quad \Longleftrightarrow \quad ( \mu^{(m)} ,  W^{(m)}   ) \xrightarrow[m \to +\infty]{\cW\cF(n)}  (\mu, W ) \: .
\]
\end{proposition}
%
%In other words, the topology induced by $\xdM$ on $\cW\cF(n)$ agrees with the quotient topology.

\begin{proof}
Let us first fix some notations.
Let $\left( \mu^{(m)} ,  W^{(m)}  \right)_{m \in \N}$ be a sequence in $\cW\cF(n)$ and let $(\mu, W ) \in \cW\cF(n)$. Let $I = (p_1, \ldots, p_r)$ be the type of $\mu$.
Let $m \in \N$ and fix a piecewise $\xC^1$ path
\[ \begin{array}{lcccl}
\gamma & : & [0,1] & \rightarrow & \cW\cF(n) \\
      &   &   t   & \mapsto     & \gamma(t) = \left( \nu(t), Q(t) \right)
\end{array} \quad \text{such that } \left\lbrace \begin{array}{lcl}
\gamma(0) & = & \left( \mu^{(m)} ,  W^{(m)}  \right) \\
\gamma(1) & = & \left( \mu ,  W  \right)
\end{array} \right. \: .
\]
Let $U : [0,1] \rightarrow \xO(n)$ be a piecewise $\xC^1$ lift of $Q$ so that the length of $\gamma$ is
\begin{equation} \label{eq:lengthGamma}
\xLM (\gamma) = \int_0^1 \left( \sum_{k=1}^n \dt{\nu}_k(t)^2 + \sum_{1 \leq i < j \leq n} f \left( \nu_{i \to j}(t) \right)^2 \: \left( U(t)^T \dt{U}(t) \right)_{ij}^2 \right)^{\frac{1}{2}} \: dt \: .
\end{equation}

{\bf Step $1$:} We first assume that $\xdM \left(( \mu^{(m)} ,  W^{(m)}   ), (\mu, W ) \right) \xrightarrow[m \to +\infty]{} 0$.

%\noindent We will prove that $\mu^{(m)} \xrightarrow[m \to +\infty]{} \mu$ in \eqref{eq:muToZero}, then we will prove that 
\noindent In particular $\displaystyle \xLM(\gamma) \geq \int_0^1 \left( \sum_{k=1}^n \dt{\nu}_k(t)^2 \right)^{\frac{1}{2}} \: dt = \int_0^1 | \dt{\nu}(t) | \: dt  \geq | \nu (0) - \nu(1) | = | \mu^{(m)} - \mu |$ and taking the infimum over all such piecewise $\xC^1$ paths $\gamma$, it follows that 
\begin{equation} \label{eq:muToZero}
| \mu^{(m)} - \mu | \leq \xdM \left( \left( \mu^{(m)} ,  W^{(m)}  \right) , \left( \mu ,  W  \right) \right) \xrightarrow[m\to \infty]{} 0 \: .
\end{equation}
Let $K = \{ j \in \{1, \ldots, n\} \: : \: \mu_j > 0 \}$. By \eqref{eq:muToZero}, we have $W^{(m)} \in \cF_{J(m)}$ with $J^{(m)} \preccurlyeq I$ for $m$ large enough.
Let $0 < \delta < 1$ such that $\displaystyle \min_{j \in K} \mu_j > 2 \delta$, and by \eqref{eq:muToZero} we can fix $N \in \N$ such that for all $m \geq N$, $\displaystyle \min_{j \in K} (\mu_j^{(m)} - \delta) > \delta$ and moreover
\begin{equation} \label{eq:leqDeltaSquare}
{ \xdM \left( \left( \mu^{(m)} ,  W^{(m)} \right) , \left( \mu ,  W \right) \right) }^2 - | \mu^{(m)} - \mu |^2 + 2^{-m} < \delta^2 < \min \left\lbrace \left. (\mu_{j}^{(m)} - \delta ) (\mu_{j} - \delta) \: \right| \: j \in K \right\rbrace \: . 
\end{equation}
We fix $m \geq N$ and a piecewise $\xC^1$ path $\gamma = (\nu,Q) : [0,1] \rightarrow \cW\cF(n)$ such that
$\gamma(0)  =  \left( \mu^{(m)} ,  W^{(m)}  \right)$ and $
\gamma(1) =  \left( \mu ,  W  \right)$ and $\gamma$ satisfies
\begin{equation} \label{eq:gammaAlmostMin}
{\xLM (\gamma)}^2 \leq {\xdM \left( ( \mu^{(m)} ,  W^{(m)} ), ( \mu,  W ) \right) }^2 + 2^{-m}
\end{equation}
(i.e. $\gamma$ reaches the infimum up to $2^{-m}$ in \eqref{eq:RiemannianDistance}). Assume by contradiction that there exists $i_0 \in K$ and $c \in ]0,1[$ such that $\nu_{i_0}(c) \leq  \delta$. We can now estimate the gap between the length of $\nu$ and the length $|\mu^{(m)} - \mu|$ of the segment connecting $\nu(0) = \mu^{(m)}$ with $\nu(1) = \mu$ (see \eqref{eq:gammaThrough0}) to obtain a contradiction.
Indeed, we then have $\nu_{i_0}(0) - \nu_{i_0}(c) \geq \nu_{i_0}(0) -  \delta = \mu_{i_0}^{(m)} -  \delta > 0$ and similarly $\nu_{i_0}(1) - \nu_{i_0}(c) \geq \nu_{i_0}(1) -  \delta = \mu_{i_0} -  \delta > 0$ so that 
\begin{equation} \label{eq:NuI0}
|2 \nu_{i_0}(c) - \nu_{i_0}(0) - \nu_{i_0}(1) | = (\nu_{i_0}(0) - \nu_{i_0}(c)) + (\nu_{i_0}(1) - \nu_{i_0}(c)) \geq \mu_{i_0}^{(m)} + \mu_{i_0} - 2 \delta \: .
\end{equation}

\begin{figure}[!htp]
\centering
\includegraphics[scale=0.50]{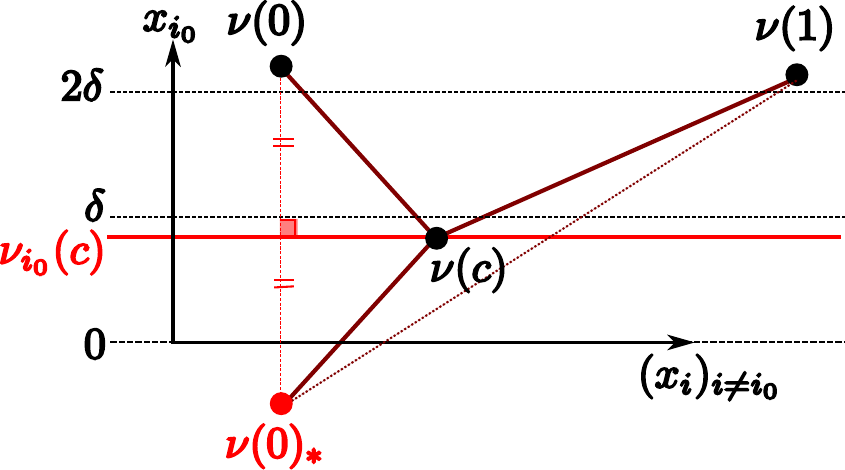}
\caption{Definition of $\nu(0)_\ast$. \label{fig:symPath}}
\end{figure}
\noindent Moreover, note that $\nu(0)_\ast$ defined as $\nu_{i_0}(0)_\ast = 2 \nu_{i_0}(c) - \nu_{i_0}(0)$ and $\nu_i(0)_\ast = \nu_i(0)$ for $i \neq i_0$, is the image of $\nu(0)$ by the symmetry with respect to the hyperplane $x_{i_0} = \nu_{i_0}(c)$ in $\R^n$, see Figure~\ref{fig:symPath}. We then have $|\nu(0) - \nu(c)| = |\nu(0)_\ast - \nu(c)|$ and we can infer that
\begin{align*}
\xLM(\gamma) & = \xLM(\gamma_{| [0,c]}) +  \xLM(\gamma_{| [c,1]}) \geq |\nu(0) - \nu(c)| + |\nu(c) - \nu(1)| = |\nu(0)_\ast - \nu(c)| + |\nu(c) - \nu(1)|\\
&  \geq |\nu(0)_\ast - \nu(1)| = \sqrt{ \sum_{i \neq i_0} | \nu_i(0) - \nu_i(1) |^2 + |2 \nu_{i_0}(c) - \nu_{i_0}(0) - \nu_{i_0}(1) |^2 } \: .
\end{align*}
Thanks to \eqref{eq:NuI0}, we then have
\begin{align}
\xLM(\gamma)^2 - | \mu^{(m)} - \mu |^2 & \geq  \sum_{i \neq i_0} | \mu_i^{(m)} - \mu_i |^2 + (\mu_{i_0}^{(m)} + \mu_{i_0} - 2 \delta)^2  -  | \mu^{(m)}  - \mu |^2  \nonumber \\
& = (\mu_{i_0}^{(m)} + \mu_{i_0} - 2 \delta)^2 - |\mu_{i_0}^{(m)} - \mu_{i_0}|^2  = 4 (\mu_{i_0}^{(m)} -  \delta) (\mu_{i_0} -  \delta) \: . \label{eq:gammaThrough0}
\end{align}
From \eqref{eq:leqDeltaSquare}, \eqref{eq:gammaAlmostMin} and \eqref{eq:gammaThrough0} we obtain the contradiction
\[
4 \delta^2 < 4 (\mu_{i_0}^{(m)} -  \delta) (\mu_{i_0} -  \delta) \leq { \xdM \left( \left( \mu^{(m)} ,  W^{(m)} \right) , \left( \mu ,  W \right) \right) }^2 - | \mu^{(m)} - \mu |^2 + 2^{-m} < \delta^2 \: .
\]
We conclude that for all $j \in K$ and for all $t \in [0,1]$, $\nu_j(t) > \delta$.
The application $x \mapsto f(x)^2$ is continuous and positive (nonzero) on the compact set $S_\delta = \cup_{i = 1}^n \{ x \in [0,1]^n \: : \: x_i \geq \delta \}$ and we can define $f_\delta > 0$ such that $f_\delta^2 = \min_{S_\delta} f^2$.
Let $(i,j) \in \{1, \ldots,n\}^2$, $i<j$, if there exists $i_0 \in K$ such that $ i \leq i_0 \leq j-1$ then for all $t \in [0,1]$, $\nu_{i_0}(t) \geq \delta$ and thus $\nu_{i \to j}(t) \in S_\delta$ implying
\[
\min_{t \in [0,1]} f(\nu_{i \to j}(t))^2 \geq f_\delta^2 \: .
\]
%
%Wrinting $I = (p_1, \ldots, p_r)$, we denote by
%\begin{align*}
%X & = \{ (i,j) \: : \: p_1 + \ldots + p_{k} + 1 \leq i , j \leq p_1 + \ldots + p_{k+1}  \text{ for some } k \in \{ 1, \ldots, r-1 \} \} \: .
%\end{align*}
Consequently, for all $(i,j) \notin X_I$ (we recall that $X_I$ are indices corresponding to the diagonal blocks with respect to the type $I$, see \eqref{eq:diagonalIndicesI}), we have $\displaystyle \min_{t \in [0,1]} f(\nu_{i \to j}(t))^2 \geq f_\delta^2$ and inserting this lower bound in \eqref{eq:lengthGamma} we obtain 
\begin{align} \label{eq:flagCauchySeq}
\xLM (\gamma) & \geq f_\delta \int_0^1 \sqrt{ \sum_{\substack{1 \leq i < j \leq n\\ (i,j) \notin X_I }} \left( U(t)^T \dt{U}(t) \right)_{ij}^2 } \: dt 
%& = f_\delta \int_a^b \sqrt{ \sum_{\substack{1 \leq i < j \leq n\\ (i,j) \notin X }} \left( U(t)^T (U^\dt(t))_H \right)_{ij}^2 } \: dt
%= \delta  \int_0^1 \sqrt{g_{Q(t)} (Q^\prime(t),Q^\prime(t))} \: dt  \geq \delta d_{\cF_I} \left( W^{(p)} , W^{(q)} \right) 
\end{align}

\noindent Let $U_I = \pi_I \circ U : [0,1] \rightarrow \cF_I$, then for any $t \in [0,1]$, $\dt{U}(t) = (\dt{U}(t))_H + (\dt{U}(t))_V \in U(t) \xSkew(n)$ where $(\dt{U}(t))_H \in H^I_{U(t)} = U(t) \mathfrak{m}_I$ and $(\dt{U}(t))_V \in \ker T_{U(t)} \pi_I$ are orthogonal (w.r.t. the usual euclidean structure in $\xM_n(\R)$).
We then have $\dt{U}_I (t) = T_{U(t)} \pi_I \cdot \dt{U}(t) = T_{U(t)} \pi_I \cdot (\dt{U}(t))_H$ which implies by \eqref{eq:flagMetric} that
\begin{align}
g_{U_I(t)}^{I} \left( \dt{U}_I (t), \dt{U}_I (t) \right) & = g_{U(t)}^{\xO(n)} ((\dt{U}(t))_H, (\dt{U}(t))_H) \nonumber \\
& = \sqrt{ \sum_{\substack{1 \leq i < j \leq n\\ (i,j) \notin X_I }} \left( U(t)^T (\dt{U}(t))_H \right)_{ij}^2 } \text{ using } U(t)^T (\dt{U}(t))_H \in \mathfrak{m}_I \nonumber \\
& = \sqrt{ \sum_{\substack{1 \leq i < j \leq n\\ (i,j) \notin X_I }} \left( U(t)^T \dt{U}(t) \right)_{ij}^2 } \: ,
\label{eq:flagCauchySeq2}
\end{align}
where we used that $U(t)^T (\dt{U}(t))_V \in \xSkew(I)$ and $U(t)^T (\dt{U}(t))_H \in \mathfrak{m}_I$ and thus (see \eqref{eq:lengthLiftedUij}):
\[\left( U(t)^T \dt{U}(t) \right)_{ij} =
\left( U(t)^T (\dt{U}(t))_H \right)_{ij} \text{ if } (i,j) \notin X_I  \: .
\]

\noindent From \eqref{eq:gammaAlmostMin}, \eqref{eq:flagCauchySeq} and \eqref{eq:flagCauchySeq2} we infer that
\begin{align}
{ \xdM  \left( ( \mu^{(m)} ,  W^{(m)} ), ( \mu,  W ) \right) }^2 + 2^{-m} & \geq {\xLM (\gamma)}^2 \nonumber \\
& \geq f_\delta^2 \left( \int_0^1 g_{U_I(t)}^{I} \left( \dt{U}_I (t), \dt{U}_I (t) \right) \: dt \right)^2 \geq f_\delta^2 L_I(\pi_I \circ U)^2 \nonumber \\
& \geq f_\delta^2 d_I (\pi_I(U(0)), \pi_I(U(1)))^2 = f_\delta^2 d_I \left(p_{J^{(m)} \to I}(W^{(m)}) , W \right)^2 \: . \label{eq:flagCv}
\end{align}
Letting $m$ tend to $+\infty$, \eqref{eq:muToZero} and \eqref{eq:flagCv} exactly mean that $\displaystyle (\mu^{(m)}, W^{(m)}) \xrightarrow[m \to \infty]{\cW\cF(n)} (\mu, W)$.

{\bf Step $2$:} We conversely assume that $\displaystyle (\mu^{(m)}, W^{(m)}) \xrightarrow[m \to \infty]{\cW\cF(n)} (\mu, W)$.

\noindent We can fix $N \in \N$ such that for all $m \geq N$, $W^{(m)} \in \cF_{J^{(m)}}$ with $J^{(m)} \preccurlyeq I$. Let $m \geq N$ and let $Q : [0,1] \rightarrow \cF_I$ be a piecewise $\xC^1$ path between $p_{J^{(m)} \to I} (W^{(m)})$ (we simply write $W^{(m)}$ hereafter) and $W$ in $\cF_I$ satisfying
\begin{equation} \label{eq:WmWQ}
d_{I} ( W^{(m)}, W) \geq L_{I} (Q) - 2^{-m} \: .
\end{equation}
Applying Proposition~\ref{prop:submersionLength}, let $U : [0,1] \rightarrow \xO(n)$ be a $I$--horizontal piecewise $\xC^1$ lift of $Q$, so that $L_I(Q) = L_{\xO(n)} (U)$.
Let us introduce $\displaystyle \widetilde{\mu^{(m)}} := \frac{1}{1 + n 2^{-m}} \left(\mu^{(m)} + 2^{-m} \right)$ so that for all $k = 1 , \ldots, n$, $\widetilde{\mu_k^{(m)}} > 0$ and $\displaystyle \sum_{k=1}^n \widetilde{\mu_k^{(m)}} = \frac{1}{1 + n 2^{-m}} \left(n 2^{-m} + \sum_{k=1}^n \mu_k^{(m)} \right) = 1$ and thus $\widetilde{\mu^{(m)}} \in \mathring{\Delta}(n)$.
We then define the piecewise $\xC^1$ path $\gamma : [0,1] \rightarrow \cW\cF(n)$ by $\gamma(t) = (\nu(t),\pi(U(t)))$ with $\nu(t) = (1-t) \widetilde{\mu^{(m)}} + t \mu$ for all $t \in [0,1]$, note that for all $t \neq 1$, the type of $\nu(t)$ satisfies $\tau (\nu(t)) = (1, \ldots, 1)$ so that $\pi(U(t)) = \pi_{(1, \ldots, 1)}(U(t))$ and $\gamma(t) \in M(n)$. As $\pi_{J^{(m)}}(U(0)) = W^{(m)}$, we can then connect $(\mu^{(m)} , W^{(m)})$ to $(\widetilde{\mu^{(m)}} , \pi_{(1,\ldots,1)} (U(0) )$ by a path $t \mapsto \left( (1-t)  \mu^{(m)} + t \widetilde{\mu^{(m)}} , \pi(U(0)) \right)$ whose length is 
\begin{align}
\left|\mu^{(m)} - \widetilde{\mu^{(m)}} \right| & = \frac{1}{1 + n 2^{-m}} \left|(1 + n 2^{-m}) \mu^{(m)} - (\mu^{(m)} + 2^{-m}) \right| = \frac{2^{-m}}{1 + n 2^{-m}} \left|n \mu^{(m)} - 1 \right| \nonumber \\
& \leq C_n 2^{-m} \label{eq:mumTildemum}
\end{align}
where $C_n > 0$ only depends on $n$ (we used that $\mu^{(m)} \in \Delta(n)$ thus $| \mu^{(m)} | \leq 1$).
By concatenation of both paths and using \eqref{eq:WmWQ} and \eqref{eq:mumTildemum}, we obtain
\begin{align*}
\xdM & \left( \left( \mu^{(m)} ,  W^{(m)} \right) , \left( \mu , W \right) \right) \leq
\xLM (\gamma) + \left|\mu^{(m)} - \widetilde{\mu^{(m)}} \right| \\
 & \leq \int_0^1 \sqrt{ \sum_{k=1}^n (\mu_k - \widetilde{\mu_k^{(m)}})^2 + \sum_{1 \leq i < j \leq n}  \underbrace{f(\nu_{i \to j}(t))^2}_{\leq C^2 = (\max f)^2}  \: \left( U(t)^T \dt{U}(t) \right)_{ij}^2 } \: dt + C_n 2^{-m}\\
& \leq \int_0^1 \sqrt{ \sum_{k=1}^n (\mu_k - \widetilde{\mu_k^{(m)}})^2 }\: dt + C \int_0^1 \sqrt{ \sum_{1 \leq i < j \leq n}  \left( U(t)^T \dt{U}(t) \right)_{ij}^2 } \: dt +  C_n 2^{-m}  \\
& \leq | \mu^{(m)} - \mu | + C \underbrace{L_{\xO(n)} (U)}_{L_{I} (Q)} + 2 C_n 2^{-m} \\
& \leq | \mu^{(m)} - \mu | + C d_{I} (W, W^{(m)}) + (2 C_n + C) 2^{-m} \xrightarrow[m \to \infty]{} 0 \: .
\end{align*}
\end{proof}

\begin{remk}
 Note that in the proof of Proposition~\ref{prop:dMtopoEquivalent}, more precisely in Step $1$, we observe that given two points $(\mu, W)$, $(\overline{\mu}, \overline{W})$ close enough, almost shortest paths do not cross lower strata. Unfortunately, such property is only local at this stage since it is proven for close enough points.
\end{remk}

´
\begin{proposition} \label{prop:innerPathApproximation}
Let $\gamma : [a,b] \rightarrow \cW\cF(n)$ be a piecewise $\xC^1$ path and let $\eta > 0$. Then, there exists a piecewise $\xC^1$ path $\gamma_\eta : [a,b] \rightarrow \cW\cF(n)$ such that $\gamma_\eta(a) = \gamma(a)$, $\gamma_\eta(b) = \gamma(b)$, for all $t \in (a,b)$, $\gamma_\eta(t) \in M(n)$ and
\[
L_g (\gamma_\eta) \leq \xLM(\gamma) + \eta \: .
\]
In particular $d_g$ and $\xdM$ coincide in $M(n)$.
\end{proposition}

\begin{proof}
Let $\eta > 0$ and let $\epsilon = \epsilon_\eta$ to be set at the end of Step $5$. Let $\gamma  = (\nu,Q) : [a,b] \rightarrow \cW\cF(n)$ be a piecewise $\xC^1$ path.

{\bf Step $1$:} 
Given Definition~\ref{dfn:piecewiseC1pathI} of a piecewise $\xC^1$ path, we can reduce the problem to the following one, up to considering a finite number of restrictions of $\gamma$.
Let $c \in ]a,b[$, we assume that 
\begin{itemize}
 \item there exists $K \subset \{1, \ldots, n\}$ with $1 \leq |K| \leq n-1$ such that $\gamma(c) \in M(r \: ; K)= \mathring{\Delta}(r \: ; K) \times \cF_I$,
 \item $\gamma$ is $\xC^1$ in $[a,b] \setminus \{c\}$,
 \item for all $t \in [a,b]$, $t \neq c$, $\gamma(t) \in M(n)$.
\end{itemize}
By Definition~\ref{dfn:piecewiseC1pathI} of piecewise $\xC^1$ path, $\gamma_{[a,c[}$ extends into a $\xC^1 $ path $\gamma^- = (\nu^-, Q^-) : [a,c] \rightarrow \Delta(n) \times \cF_{(1, \ldots, 1)}$ and $\gamma_{]c,b]}$ into a $\xC^1 $ path $\gamma^+ = (\nu^+, Q^+) : [c,b] \rightarrow \Delta(n) \times \cF_{(1, \ldots, 1)}$.

{\bf Step $2$:} As $\gamma$ is continuous with respect to the weighted flag topology, then $\nu : [a,b] \rightarrow (\Delta(n), |\cdot |)$ is continuous and $\lim_{t \to c} \nu(t) = \nu(c) =: \mu$ and by assumption the type of $\mu$ is $I$.
By definition of $X_I$, $\forall (i,j) \in X_I$, $\mu_{i \to j} = 0$ and $f(\mu_{i \to j}) = 0$ and by continuity of $f$ and $\gamma$, $\lim_{t \to c} f(\nu_{i \to j}(t)) = 0$.
Furthermore, $Q^-$, $Q^+$ are $\xC^1$ path in $\cF_{(1,\dots, 1})$ and therefore their length defined in \eqref{eq:lengthFI} satisfies
\[
\lim_{\substack{t \to c\\ t < c}} L_{(1,\ldots,1)} \left(Q^-_{| [t,c]} \right) = 0 \quad \text{and} \quad \lim_{\substack{t \to c\\ t > c}} L_{(1,\ldots,1)} \left(Q^+_{| [c,t]} \right) = 0 \: .
\]
Therefore we can fix $c^- \in ]a,c[$ and $c^+ \in ]c,b[$ such that 
\begin{enumerate}[(i)]
 \item for all $t \in [c^-, c^+]$, $|\nu (t) - \mu | \leq \epsilon$,
 \item for all $(i,j) \in X_I$, $f( \nu_{i \to j} (c^-) ) \leq \epsilon$ and $f( \nu_{i \to j} (c^+) ) \leq \epsilon$,
 \item $L_g \left( \gamma_{| ]c^-, c^+[} \right) \leq \epsilon$.
 \item $\displaystyle L_{(1,\ldots,1)} \left(Q^-_{| [c^-,c]} \right) \leq  \epsilon \quad \text{and} \quad L_{(1,\ldots,1)} \left(Q^+_{| [c,c^+]} \right) \leq \epsilon$.
\end{enumerate}
We define $\gamma_\nu : [0,1] \rightarrow M(n)$ by $\gamma_\nu(t) = ( (1-t) \nu(c^-) + t \nu(c^+), Q^-(c) )$ whose length is exactly
\[
 L_g(\gamma_\nu) = |\nu (c^+) - \nu (c^-) | \leq 2 \epsilon \: .
\]
The next steps consist in connecting $\gamma(c^-) = (\nu(c^-), Q(c^-))$ with $\gamma(c^+) = (\nu(c^+), Q(c^+))$ by an alternative path passing through $(\nu(c^-), Q^-(c))$, $ (\nu(c^+), Q^-(c))$, $(\nu(c^+), Q^+(c))$, and remaining in $M(n)$, see Figure~\ref{fig:approxPath}.

\begin{figure}[!htp]
\centering
\includegraphics[scale=0.80]{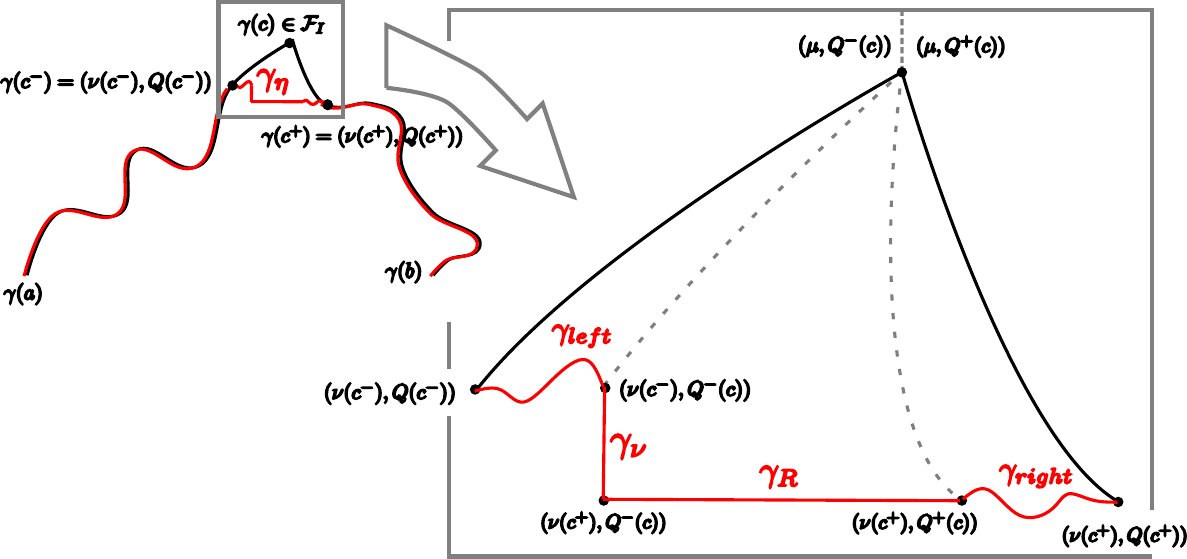}
\caption{Definition of $\gamma_\eta$. \label{fig:approxPath}}
\end{figure}

{\bf Step $3$:}
By continuity of $\gamma$, $Q^-(c)$ and $Q^+(c)$ both project on $Q(c)$ in $\cF_I$ in the sense that $p(Q^{-}(c)) = p (Q^{+}(c)) = Q(c)$ with $p = p_{(1,\ldots,1) \to I}$.
Consequently, we can choose $V^{-}, V^{+} \in \xO(n)$ such that $\pi_{(1,\ldots,1)} (V^-) = Q^-(c)$ and $ \pi_{(1, \ldots, 1)}(V^{+}) = Q^+(c)$ and $V^{+} = V^{-} R_0$ with $R_0 \in \xO_0(I)$ ($\xO_0(I)$ is the connected component of $\xO(I)$ containing $I_n$).

%As $\gamma$ is piecwise $\xC^1$, there exists $(\mu, W^{-}) = \lim_{t \to c_-} \gamma(t)$ and $(\mu, W^{+}) = \lim_{t \to c_+} \gamma(t)$ both belonging to $\Delta(n) \times \cF_{(1,\ldots, 1)}$ and such that $p(W^{-}) = p (W^{+}) =: W$ with $p = p_{(1,\ldots,1) \to I}$. Choose $U^{-}, U^{+} \in \xO(n)$ such that $\pi_{(1,\ldots,1)} (U^-) = W^-$ and $ \pi_{(1, \ldots, 1)}(U^{+}) = W^+$ and $U^{+} = U^{-} R_0$ with $R_0 \in \xO_0(I)$ ($\xO_0(I)$ is the connected component of $\xO(I)$ containing $I_n$). 
\noindent Consider a shortest $\xC^1$ path $R : [0,1] \rightarrow \xO_0(I) \subset \xO(I)$ such that $R(0) = I_n$ and $R(1) = R_0$ and then $\gamma_R : [0,1] \rightarrow M(n)$ defined by $\gamma_R(t) = \left(\nu(c^+), \pi_{(1,\ldots,1)}(V^- R(t)) \right)$. Denoting $V = V^- R$, we have $V^T \dt{V} = R^T (V^-)^T V^- \dt{R} = R^T \dt{R}$, as $R$ is a path in $\xO(I)$, $\dt{R}(t) \in T_{R(t)} \xO(I) = R(t) \xSkew(I)$ i.e. $R^T \dt{R} \in \xSkew(I)$. We then have $\left( R^T \dt{R} \right)_{ij}^2 = 0$ for $1 \leq i < j \leq n$, $(i,j) \notin X_I$ and thus
\begin{align*}
L_g(\gamma_R) & = \int_0^1 \sqrt{\sum_{(i,j) \in X_I} f(\nu(c^+)_{i \to j})^2 \left( R(t)^T \dt{R}(t) \right)_{ij}^2}  \: dt \leq \epsilon  \int_0^1 \sqrt{\sum_{(i,j) \in X_I} \left( R(t)^T \dt{R}(t) \right)_{ij}^2}  \: dt \\
& \leq \epsilon \: L_{\xO(I)}(R) = \epsilon \: d_{\xO(I)} (I_n, R_0) \leq \epsilon \: \xdiam(\xO_0(I)) \: .
\end{align*}

{\bf Step $4$:} It remains to connect $(\nu(c^-), Q(c^-))$ to $(\nu(c^-), Q^-(c))$ (as well as $(\nu(c^+), Q(c^+))$ to $(\nu(c^+), Q^+(c))$). 
Let $U : [a,c] \rightarrow \xO(n)$ be a $\xC^1$ lift of $Q^-$ such that $U(c) = V^-$. 
We then define $\gamma_{left} : [c^-,c] \rightarrow M(n)$ by $\gamma_{left}(t) = (\nu(c^-), \pi_{(1,\ldots,1)} U(t))$ and compare $L_g(\gamma_{left})$ with $L_g(\gamma_{]c^-, c[})$. 
Note that for all $c^- \leq t < c$, $\pi_{(1,\ldots,1)} U(t)) = Q(t)$ (the ``flag component'' of $\gamma$).
We use that $f$ is Lipschitz and bounded in $[0,1]^n$ and $f^2$ is thus Lipschitz with constant denoted by $\xLip(f^2)$ so that \[
| f(\nu(t)_{i \to j})^2 - f(\nu(c^-)_{i \to j})^2 | \leq \xLip(f^2) | \nu(t)_{i \to j} - \mu_{i \to j}| \leq \xLip(f^2) | \nu(t) - \mu| \leq  \xLip(f^2) \epsilon \: .                                      \]
We obtain
\begin{align*}
  L_g(\gamma_{left}) & = \int_{c^-}^c \sqrt{\sum_{1 \leq i < j \leq n} f(\nu(c^-)_{i \to j})^2 \left( U(t)^T \dt{U}(t) \right)_{ij}^2}  \: dt \\
  & \leq \int_{c^-}^c \sqrt{\sum_{1 \leq i < j \leq n} \left(f(\nu(t)_{i \to j})^2 + \xLip(f^2) \epsilon \right) \left( U(t)^T \dt{U}(t) \right)_{ij}^2}  \: dt \\
  & \leq \int_{c^-}^c \sqrt{\sum_{1 \leq i < j \leq n} f(\nu(t)_{i \to j})^2 \left( U(t)^T \dt{U}(t) \right)_{ij}^2}  \: dt + \sqrt{\xLip(f^2) \epsilon} \int_{c^-}^c \sqrt{\sum_{1 \leq i < j \leq n}  \left( U(t)^T \dt{U}(t) \right)_{ij}^2}  \: dt \\
  & \leq \xLM \left(\gamma_{| [c^-, c[} \right) + \sqrt{\xLip(f^2) \:  \epsilon} \: L_{(1,\ldots,1)} \left(Q^-_{| [c^-,c]} \right) \\
  & \leq \xLM \left(\gamma_{| [c^-, c[} \right) + \epsilon \sqrt{\xLip(f^2) \:  \epsilon}  \text{ by Step $2 (iv)$.}
\end{align*}
We can similarly define $\gamma_{right}$ such that $L_g(\gamma_{right}) \leq \xLM \left(\gamma_{| ]c, c^+]} \right) + \epsilon \sqrt{\xLip(f^2) \:  \epsilon}$.

{\bf Step $5$:} We can concatenate together our different paths in the following order: $\gamma_{| [a,c^-]}$, $\gamma_{left}$, $\gamma_\nu$, $\gamma_R$, $\gamma_{right}$ and finally $\gamma_{| [c^+,b]}$ to construct a piecewise $\xC^1$ path $\gamma_\eta$. We recall (Remark~\ref{remk:lengthDistanceMg}) that $L_g$ and $\xLM$ agree on paths in $M(n)$, therefore $\gamma_\eta$ satisfies
\[
\begin{array}{cclccccccc}
L_g(\gamma_\eta) - \xLM(\gamma) & =   & - \xLM \left(\gamma_{| [c^-, c^+]} \right)  + L_g(\gamma_{left}) +  L_g(\gamma_{right}) & + & L_g(\gamma_\nu) & + & L_g(\gamma_R)  \\
& \leq &  2 \epsilon \sqrt{\xLip(f^2) \:  \epsilon}  & + & 2 \epsilon & + & \epsilon \: \xdiam ( \xO_0(I)) \\
& \leq & \eta \: ,
\end{array}
\]
initially setting $\epsilon = \epsilon_\eta > 0$ so that the last inequality holds.

{\bf Step $6$:}
Thanks to Remark~\ref{remk:lengthDistanceMg}, we are left with the proof of $d_g \leq \xdM$ in $M(n)$. Let $(\mu,W)$, $(\nu, Q) \in M(n)$ and $\eta > 0$. Let $\gamma : [a,b] \rightarrow \cW\cF(n)$ be a piecewise $\xC^1$ path such that $\xLM(\gamma) \leq \xdM ( (\mu,W), (\nu, Q) ) + \eta$. Let us now apply the first part of Proposition~\ref{prop:innerPathApproximation} and let $\gamma_\eta : [a,b] \rightarrow M(n)$ be a piecewise $\xC^1$ path between $(\mu,W)$ and $(\nu, Q)$ satisfying $L_g (\gamma_\eta) \leq \xLM(\gamma) + \eta$. We consequently obtain
\[
d_g ( (\mu,W), (\nu, Q) ) \leq L_g (\gamma_\eta) \leq \xLM(\gamma) + \eta \leq \xdM ( (\mu,W), (\nu, Q) ) + 2\eta
\]
and we conclude letting $\eta$ tend to $0$.
\end{proof}

\begin{theorem} \label{thm:completion}
The metric space $(\cW\cF(n), \xdM)$ is the completion of the Riemannian metric space $(M(n), d_g)$.
\end{theorem}

\begin{proof}
We know from Proposition~\ref{prop:dMtopoEquivalent} that  the topology induced by $\xdM$ coincide with the quotient topology in $\cW\cF(n)$ so that $(\cW\cF(n), \xdM)$ is a compact metric space and thus it is complete. We additionally know from Proposition~\ref{prop:innerPathApproximation} that $(M(n) , \xdM)$ and $(M(n), d_g)$ coincide. It remains to prove that $\overline{M(n)} = \cW\cF(n)$ where the closure $\overline{M(n)}$ is taken with respect to $\xdM$, which is equivalent to taking the closure with respect to the quotient topology since both topology coincide and we already know that $M(n)$ is dense in $\cW\cF(n)$ (for the quotient topology, see Section~\ref{section:stratification}).
\end{proof}

\subsection{Another possible length structure in weighted flags}
\label{subsec:lengthWFalternative}

It is also possible to define a length structure in $\cW\cF(n)$ using the fact that $g$ defined in \eqref{eq:metric} provides a Riemannian metric in each elementary cell in $\cW\cF(n)$. More precisely, we recall that for a non emptyset $K \subset \{1, \ldots, n\}$, $M(|K| \: ; K) = \mathring{\Delta}(|K| \: ; K) \times \cF_{I(K)}$ that can be completed (with respect to the product topology) into the smooth manifold with boundary 
\[
\widetilde{M} = \Delta(|K| \: ; K) \times \cF_{I(K)} \quad \text{with} \quad \Delta(|K| \: ; K) = \{ \mu \in [0,1] \: : \: \mu_j = 0 \text{ for } j \notin K \} \: .
\]
We can then extend the notion of piecewise $\xC^1$ paths in $\cW\cF(n)$ given in Definition~\ref{dfn:piecewiseC1pathI}.

\begin{definition}[piecewise $\xC^1$ path - II] \label{dfn:piecewiseC1pathII}
Let $\gamma : [a,b] \subset \R \rightarrow \cW\cF(n)$ be continuous.
\begin{itemize}
 \item Let $M = M(|K| \: ; K)$. We first say that $\gamma$ is a piecewise $\xC^1$ path in $\overline{M}$ if there exist $a = a_0 < a_1 < \ldots < a_{N-1} < a_N = b$ such that for all $i$, $\gamma(]a_{i-1}, a_i[) \subset M$ and $\gamma_{| ]a_{i-1}, a_i[}$ extends to a $\xC^1$ path in the smooth manifold with boundary $\widetilde{M}$.
 \item The restriction of $g$ (defined in \eqref{eq:metric}) to $M(|K| \; : K)$ is a Riemannian metric and we define 
 \[
L_{\overline{M(|K| \; : K)}} (\gamma) = \sum_{i=1}^N \int_{a_{i-1}}^{a_i} \sqrt{ g_{\gamma(t)} (\dt{\gamma}(t), \dt{\gamma}(t)) } \: dt
 \]

 \item We then say that $\gamma$ is a piecewise $\xC^1$ path if there exist $a = a_0 < a_1 < \ldots < a_{N-1} < a_N = b$ such that for all $i$, there exists $K_i$ such that $\gamma_{| [a_{i-1}, a_i]}$ is a piecewise $\xC^1$ path in $\overline{M(|K_i| \: ; K_i)}$.

 \item We define the length $\xLWF$ of a piecewise $\xC^1$ path $\gamma : [a,b] \rightarrow \cW\cF(n)$ as
\[
\xLWF (\gamma) = \sum_{i=1}^N L_{\overline{M(|K_i| \; : K_i)}} \left(\gamma_{| [a_{i-1},a_i]} \right)
\]
as well as the associated distance
\begin{equation} \label{eq:RiemannianDistanceII}
\xdWF \left((\mu, W  ), (\nu, Q ) \right) = \inf \left\lbrace \xLWF (\gamma) \: \left| \begin{array}{l}
\gamma : [a,b] \rightarrow \cW\cF(n) \text{ piecewise } \xC^1 \text{ path} \\
\gamma(a) = (\mu,  W  ) \text{ and } \gamma(b) = (\nu,  Q  )
\end{array} \right. \right\rbrace \: .
\end{equation}
\end{itemize}
\end{definition}

\begin{proposition}
\label{prop:xdMxdWF}
The length structure $\xLWF$ extends $\xLM$ and the induced distances $\xdWF$ and $\xdM$ coincide in $\cW\cF(n)$.
\end{proposition}

\begin{proof}
First note that if $\gamma : [a,b] \rightarrow \cW\cF(n)$ is a piecewise $\xC^1$ path in the sense $I$ of Definition~\ref{dfn:piecewiseC1pathI} then it is in particular a piecewise $\xC^1$ path in the sense $II$ of Definition~\ref{dfn:piecewiseC1pathII} and moreover $\xLM(\gamma) = \xLWF (\gamma)$. Consequently, the length structure $\xLWF$ extends $\xLM$ to a larger class of admissible paths (and in particular paths that are entirely contained in some elementary cell $M( |K| \, ; K) \subset \cW\cF(n) \setminus M(n)$, which is natural to connect points both belonging to $M( |K| \, ; K)$).
 
We now check that the induced distances $\xdWF$ and $\xdM$ coincide in $\cW\cF(n)$.
As $\xLWF$ extends $\xLM$, we already have the inequality $\xdWF \leq \xdM$. Conversely, given $\eta > 0$ and a $\xC^1$ path $\gamma = (\nu, Q) : [a,b] \rightarrow \overline{M(|K| \; : K)} \subset \cW\cF(n)$ satisfying $\gamma(]a,b[) \subset M(|K| \; : K)$, there exists a piecewise $\xC^1$ path $\gamma_\eta : [a, b] \rightarrow \cW\cF(n)$ connecting $\gamma(a)$ and $\gamma(b)$ and satisfying
 \[
\gamma_\eta(]a,b[) \subset M(n) \quad \text{and} \quad \xLM(\gamma_\eta) \leq \xLWF(\gamma) +\eta \: .
 \]
 Indeed, take a $\xC^1$ lift $U : [a,b] \rightarrow \xO(n)$ of $Q$ and for $\epsilon > 0$ (to be chosen later with respect to $\eta$), define $\widetilde{\nu}_\epsilon : t  \mapsto \frac{1}{1 + n \epsilon} \left(\nu(t) + \epsilon \right)$ and $\widetilde{\gamma}_\epsilon = \left( \widetilde{\nu}_\epsilon , \pi_{(1, \ldots, 1)} (U) \right)$ whose image is included in $M(n)$. Then, $\widetilde{\nu}_\epsilon$ converges to $\nu$ uniformly in $[a,b]$ and
 \begin{align*}
 \xLM (\widetilde{\gamma}_\epsilon)  = & \int_a^b \left( \frac{1}{(1 + n \epsilon)^2} \sum_{k=1}^n \dt{\nu}_k(t)^2 + \sum_{1 \leq i < j \leq n} f \left( (\widetilde{\nu}_\epsilon)_{i \to j}(t) \right)^2 \: \left( U(t)^T \dt{U}(t) \right)_{ij}^2 \right)^{\frac{1}{2}} \: dt \\
 \xrightarrow[\epsilon \to 0]{} & \int_a^b \left(  \sum_{k=1}^n \dt{\nu}_k(t)^2 + \sum_{1 \leq i < j \leq n} f \left( \nu_{i \to j}(t) \right)^2 \: \left( U(t)^T \dt{U}(t) \right)_{ij}^2 \right)^{\frac{1}{2}} \: dt \\
 = & \, L_{\overline{M(|K| \; : K)}} (\gamma) = \xLWF(\gamma) \: .
 \end{align*}
It is then possible to fix $\epsilon > 0$ small enough (depending on $\eta$) so that $\xLM (\widetilde{\gamma}_\epsilon) \leq \xLWF(\gamma) + \frac{\eta}{2}$. To define $\gamma_\eta$, it remains to concatenate two segments with $\widetilde{\gamma}_\epsilon$: 
\begin{align*}
t \in [0,1] & \mapsto \left( (1-t) \nu(a) + t \widetilde{\nu}_\epsilon (a), \pi(U(a)) \right) & \text{connecting } & (\nu(a),\pi(U(a))) & \text{with } & (\widetilde{\nu}_\epsilon (a),\pi(U(a))) \\
t \in [0,1]  & \mapsto \left( (1-t) \nu(b) + t \widetilde{\nu}_\epsilon (b), \pi(U(b))\right) & \text{connecting } & (\nu(b),\pi(U(b))) & \text{with } & (\widetilde{\nu}_\epsilon (b),\pi(U(b)))
\end{align*}
whose lengths are bounded by $C_n \epsilon$ (see \eqref{eq:mumTildemum}). We finally obtain, possibly decreasing $\epsilon$,
\[
\xLM (\gamma_\eta) \leq \xLWF(\gamma) + \frac{\eta}{2} + 2 C_n \epsilon \leq \xLWF(\gamma) + \eta \: .
\]
Iterating finitely many times the construction of $\gamma_\eta$, given any piecewise $\xC^1$ path $\gamma$ in the sense of Definition~\ref{dfn:piecewiseC1pathII}, there exists a piecewise $\xC^1$ path $\gamma_\eta$ in the sense of Definition~\ref{dfn:piecewiseC1pathI}, connecting the same endpoints and such that $\xLM (\gamma_\eta) \leq \xLWF(\gamma) + \eta$, eventually proving that $\xdM$ and $\xdWF$ define the same distance in $\cW\cF(n)$.
\end{proof}

\begin{remk}
As the metric space $(\cW\cF(n), \xdWF)$ is complete and compact, then every two points connected by a rectifiable curve are connected by a shortest path (see \cite{Burago}, section 2.5), nevertheless, this shortest path refers to the length structure $L$ induced by $\xdWF$ (see Definition 2.3.1 in \cite{Burago}) that always satisfy $L \leq \xLWF$. It could be interesting to better understand the connection between $L$ and $\xLWF$. Though the equality does not hold in general, if $L_g$ and $d_g$ are associated to a Riemannian metric in a connected smooth manifold, then the length structure $L$ induced by $d_g$ is exactly the Riemannian length $L_g$ (when restricted to piecewise $\xC^1$ paths). 
% see 2.1.5 in lengthStructRiemManifold
%Note that unlike the admissible paths for $\xLWF$ (see Definition~\ref{dfn:piecewiseC1pathII}), the admissible paths for $\xLM$ (see Definition~\ref{dfn:piecewiseC1pathI}) do not allow paths staying in $\cW\cF(n) \setminus M(n)$ while it seems reasonnable to look for shortest paths inside the closure $\overline{M(|K| \: ; K)}$ of some elementary cell.
 %
\end{remk}

\section{Geodesic equations, numerical examples}
\label{section:geodesic}

The purpose of this section is twofold. Using classical tools of calculus of variations, we first compute in Section~\ref{subsec:geodesicEq} the geodesic equations (Proposition~\ref{prop:shortestDeltanOn}$(i)$--$(ii)$) characterizing critical points of the energy $\cE$ associated with the relevant length in $\mathring{\Delta}(n) \times \xO(n)$. In Section~\ref{subsec:geodNum}, we then propose a numerical implementation for computing solutions of the aforementioned geodesic equations and present some examples of geodesic between weighted flags.
\subsection{Optimality conditions for smooth paths: geodesic equations}
\label{subsec:geodesicEq}

In this section, we explicit necessary optimality conditions for smooth shortest paths in the Riemannian manifold $(M(n) , g)$. To this end, we investigate optimality conditions for smooth shortest paths in $\mathring{\Delta}(n) \times \xO(n)$ endowed with the ``lifted'' metric $\widetilde{g}$ defined at $(\mu,U) \in \mathring{\Delta}(n) \times \xO(n)$ by
\begin{multline}
\forall \alpha, \beta \in T_\mu \mathring{\Delta}(n), \, \forall \widetilde{X}, \widetilde{Y} \in U \xSkew(n), \\
\widetilde{g}_{\mu,U} ((\alpha,\widetilde{X}), (\beta,\widetilde{Y})) = \alpha \cdot \beta + \sum_{1 \leq i < j \leq n} f(\mu_{i \to j} )^2  (U^T \widetilde{X})_{ij} (U^T \widetilde{Y})_{ij} \: ,
\end{multline}
and we transfer those conditions to smooth paths in $(M(n) , g)$. We introduce the energy of a curve $\widetilde{\gamma} = (\mu, U)  : [0, L] \rightarrow \mathring{\Delta} (n) \times \xO(n)$:
\begin{equation}
{\cE}_{\widetilde{g}} (\widetilde{\gamma}) = \frac{1}{2} \int_0^L {\widetilde{g}}_{\widetilde{\gamma}(t)} \left( \dt{\widetilde{\gamma}}(t), \dt{\widetilde{\gamma}}(t) \right) \: dt = \frac{1}{2} \int_0^L \left( | \dt{\mu}(t) |^2 + \sum_{1 \leq i < j \leq n} f(\mu_{i \to j}(t))^2 (U(t)^T \dt{U}(t))_{ij}^2 \right) \: dt
\end{equation}
and we recall that $\widetilde{\gamma}$ is a geodesic if and only if $\widetilde{\gamma}$ is a critical point of ${\cE}_{\widetilde g}$. We add that minimizers of ${\cE}_{\widetilde g}$ are minimizers of the length $L_{\widetilde g}$ that are arclength parametrized so that smooth shortest paths are in particular geodesics (of course there are geodesics that are only critical points and not energy/length minimizing).

\begin{proposition} \label{prop:shortestDeltanOn}
Assume $\widetilde{\gamma} = (\mu, U)  : [0, L] \rightarrow \mathring{\Delta} (n) \times \xO(n)$ is a minimizer of  ${\cE}_{\widetilde{g}}$ among smooth paths with same endpoints, then $\widetilde{\gamma}$ is a geodesic and using the notation $B = (b_{ij})_{ij}$, $b_{ij} = (U^T \dt{U})_{ij}$, $\widetilde{\gamma}$ satisfies the following differential equations:
\begin{enumerate}[$(i)$]
 \item for all $k \in \{1, \ldots, n\}$, $\mu_k$ satisfies
 \[
\ddt{\mu}_k = \frac{1}{n}\sum_{\substack{l = 1 \\ l \neq k}}^n \sum_{1 \leq i \leq l < j \leq n} f(\mu_{i \to j}) \: \partial_l f (\mu_{i \to j}) b_{ij}^2  - \frac{n-1}{n} \sum_{1 \leq i \leq k < j \leq n} f(\mu_{i \to j}) \: \partial_k f (\mu_{i \to j}) b_{ij}^2  \: .
\]
\item $U$ satisfies
\begin{equation*}
\dt{U} = U B \quad \text{with} \quad b_{ij} = \frac{1}{f(\mu_{i \to j})^2} \left( U^T U(0) C(0) U(0)^T U \right)_{ij} \quad \text{and} \quad C(0) = \left( f(\mu_{i \to j}(0)^2 b_{ij}(0) \right)_{ij} \: .
\end{equation*}
\item Assume that $\gamma : [0,L] \rightarrow M(n)$ is a smooth arclength parametrized path minimizing $L_g$ among smooth paths with same endpoints, then $\gamma = (\mu , \pi_{(1, \ldots, 1)}(U))$ where $(\mu, U) : [0,L] \rightarrow \mathring{\Delta} (n) \times \xO(n)$ is smooth and satisfies both $(i)$ and $(ii)$.
\end{enumerate}

\end{proposition}

\begin{proof}
Let $\epsilon \neq 0$ and $h : [0,L] \rightarrow \R$ be a $\xC^1$ function satisfying $h(0) = h(L)=0$.

{\bf Step $1$:} We start with performing small variations $\gamma_\epsilon$ of $\gamma = (\mu,U)$ in the $\mu$--part.
Let $k,l \in \{1, \ldots, n\}$ and for $k \neq l$, let $\omega_{kl} = e_k - e_l \in T_\mu \mathring{\Delta}(n)$. Let us consider $\gamma_\epsilon = (\mu + \epsilon h \: \omega_{kl}, U) : [0,L] \rightarrow \mathring{\Delta} (n) \times \xO(n)$ with same endpoints as $\widetilde{\gamma}$ so that for all $\epsilon \neq 0$, ${\cE}_{\widetilde{g}}(\gamma_\epsilon) \geq {\cE}_{\widetilde{g}}(\widetilde{\gamma})$ and consequently
\[
\left. \frac{d}{d\epsilon} \right|_{\epsilon = 0} {\cE}_{\widetilde{g}} (\gamma_\epsilon) = 0 \: .
\]
We use the notation $b_{ij} = (U^T \dt{U})_{ij}$ and we obtain differentiating under the integral:
\begin{align*}
\left.  \frac{d}{d\epsilon} \right|_{\epsilon = 0} & {\cE}_{\widetilde{g}} (\gamma_\epsilon)  = \int_0^L  (\dt{\mu}_k(t) - \dt{\mu}_l(t)) \dt{h}(t) \: dt + \int_0^L h(t) \sum_{1 \leq i \leq k < j \leq n} f(\mu_{i \to j}(t)) \: \partial_k f (\mu_{i \to j}(t)) b_{ij}(t)^2 \: dt \\
&  - \int_0^L h(t) \sum_{1 \leq i \leq l < j \leq n} f(\mu_{i \to j}(t)) \: \partial_l f (\mu_{i \to j}(t)) b_{ij}(t)^2  \: dt \: .
\end{align*}
Integrating by parts, we obtain the following differential equations for all $k \neq l$,
\[
\ddt{\mu}_k - \ddt{\mu}_l + \sum_{1 \leq i \leq k < j \leq n} f(\mu_{i \to j}) \: \partial_k f (\mu_{i \to j}) b_{ij}^2 - \sum_{1 \leq i \leq l < j \leq n} f(\mu_{i \to j}) \: \partial_l f (\mu_{i \to j}) b_{ij}^2 = 0 \: ,
\]
and summing over $l \neq k$, using $\sum_l \ddt{\mu}_l = 0$, we infer for all $k \in \{1, \ldots, n\}$:
\[
n \ddt{\mu}_k  + (n-1) \sum_{1 \leq i \leq k < j \leq n} f(\mu_{i \to j}) \: \partial_k f (\mu_{i \to j}) b_{ij}^2 - \sum_{l \neq k} \sum_{1 \leq i \leq l < j \leq n} f(\mu_{i \to j}) \: \partial_l f (\mu_{i \to j}) b_{ij}^2 = 0 \: .
\]

{\bf Step $2$:} We carry on with small deformations $\gamma_\epsilon$ of $\gamma = (\mu, U)$ in the $U$--part.
Let $\Omega \in \xSkew(n)$ and let us consider $U_\epsilon = \exp(\epsilon h \: \Omega ) U \in \xO(n)$. As for a matrix $A$, $\exp (A)^T = \exp(A^T)$ and $\exp(-A)\exp(A) = I_n$, we have
\[
U_\epsilon^T = U^T \exp( \epsilon h  \Omega )^T = U^T \exp( - \epsilon h  \Omega ) \: .
\]
Moreover, for all $t,s \in [0,L]$, the matrices $\epsilon h(t)  \Omega$ and $\epsilon h(s)  \Omega $ commute and thus
\[
 \dt{U}_\epsilon = \exp( \epsilon h  \Omega) \dt{U} + \epsilon \dt{h} \: \exp( \epsilon h \Omega ) \Omega  U \: .
\]
Consequently $\displaystyle U_\epsilon^T \dt{U}_\epsilon = U^T \dt{U} + \epsilon \dt{h} \: U^T \Omega U \in \xSkew(n)$ and then
\begin{align}
  \frac{d}{d\epsilon}  (U_\epsilon^T \dt{U}_\epsilon)_{ij}^2 & = \frac{d}{d\epsilon} \left( b_{ij} + \epsilon \dt{h} \left[ U^T \Omega U \right]_{ij} \right)^2  = 2 \left( b_{ij} + \epsilon \dt{h} \left[ U^T \Omega U \right]_{ij} \right) \dt{h} \left[ U^T \Omega U \right]_{ij}  \: . \label{eq:depsilonUTUprime}
\end{align}

\noindent Let us consider $\gamma_\epsilon = (\mu, U_\epsilon) : [0,L] \rightarrow \mathring{\Delta}(n) \times \xO(n)$ with same endpoints as $\widetilde{\gamma}$: as in Step $1$, differentiating under the integral and using \eqref{eq:depsilonUTUprime}, we obtain 
\[
0 = \left.  \frac{d}{d\epsilon} \right|_{\epsilon = 0} {\cE}_{\widetilde{g}} (\gamma_\epsilon)  =  \int_0^L \dt{h}(t) \sum_{1 \leq i < j \leq n} f(\mu_{i \to j}(t))^2  b_{ij}(t) \left[ U(t)^T \Omega U(t) \right]_{ij}   \: dt
\]
and integrating by parts, we infer that for all $t \in [0,L]$,
\begin{equation} \label{eq:Omega}
\sum_{1 \leq i < j \leq n} f(\mu_{i \to j}(t))^2  b_{ij}(t) \left[ U(t)^T \Omega U(t) \right]_{ij}  = \sum_{1 \leq i < j \leq n} f(\mu_{i \to j}(0))^2  b_{ij}(0) \left[ U(0)^T \Omega U(0) \right]_{ij} \: .
\end{equation}

{\bf Step $3$:} Let $t \in [0,L]$ and let us denote $C(t) = (c_{ij}(t))_{ij} \in \xSkew(n)$ with $c_{ij}(t) = f(\mu_{i \to j}(t))^2 b_{ij}(t)$ for $1 \leq i < j \leq n$. Then, $C(t)$ and $U(t)^T \Omega U(t)$ being skew-symmetric, \eqref{eq:Omega} rewrites:
\begin{align*}
& \frac{1}{2} \tr \left(C(t) (U(t)^T \Omega U(t))^T \right) = \frac{1}{2} \tr \left(C(0) (U(0)^T \Omega U(0))^T \right) \\
\Longleftrightarrow \quad & \tr \left(C(t) U(t)^T \Omega U(t) \right) = \tr \left(C(0) U(0)^T \Omega U(0) \right)  \\
\Longleftrightarrow \quad & \tr \left(U(t) C(t) U(t)^T \Omega \right) = \tr \left(U(0) C(0) U(0)^T \Omega  \right) \: .
\end{align*}
The previous equality holds for any $\Omega \in \xSkew(n)$ so that 
$
U(t) C(t) U(t)^T = U(0) C(0) U(0)^T
$ and, eventually, for all $t \in [0,L]$,
\[
b_{ij}(t) = \frac{1}{f(\mu_{i \to j}(t))^2} \left( U(t)^T U(0) C(0) U(0)^T U(t) \right)_{ij}
\]
Note that we just recovered the optimality condition in the case where $f = 1$ and for smooth shortest paths in $\xO(n)$: we have $B = C$ and we then obtain the differential equation $U(t)^T \dt{U}(t) = U(t)^T U(0) B(0) U(0)^T U(t)$ that is $\dt{U}(t) = Y U(t)$ and consistently $U(t) = \exp (tY) U(0)$ with $Y = U(0) B(0) U(0)^T $.

{\bf Step $4$:} Let $\widetilde{\gamma} = (\mu , U) : [0,L] \rightarrow \mathring{\Delta}(n) \times \xO(n)$ be smooth and such that $\gamma = (\mu , \pi_{(1, \ldots, 1)} (U))$ then $L_g (\gamma) = L_{\widetilde{g}} (\widetilde{\gamma})$. Assume by contradiction that there exists $\widetilde{\gamma}_\ast = (\mu_\ast , U_\ast)$ with same endpoints as $\widetilde{\gamma}$ and such that $L_{\widetilde{g}} (\widetilde{\gamma}_\ast) < L_{\widetilde{g}} (\widetilde{\gamma})$ then $L_g ( \mu_\ast , \pi_{(1, \ldots, 1)} (U_\ast)) = L_{\widetilde{g}} (\widetilde{\gamma}_\ast) < L_g (\gamma)$ which is impossible. Therefore $\widetilde{\gamma}$ minimizes ${\cE}_{\widetilde g}$ and thus $(\mu , U)$ satisfy both $(i)$ and $(ii)$.
\end{proof}

\begin{remk}[Initial speed $B(0)$ has only two nonzero coefficients.]
Note that in the particular case where there exists $1 \leq i_0 < j_0 \leq n$ such that $b_{ij}(0) = 0$ for all $(i,j) \notin \{ (i_0,j_0) , (j_0, i_0) \}$ then $C(0)$ and $B(0)$ are proportional: $C(0) = f(\mu_{i_0 \to j_0}(0))^2 B(0)$ and $(ii)$ takes the simpler form
\[
\dt{U} = U \frac{f(\mu_{i_0 \to j_0}(0))^2}{f(\mu_{i_0 \to j_0})^2} U^T U(0) B(0) U(0)^T U = \frac{f(\mu_{i_0 \to j_0}(0))^2}{f(\mu_{i_0 \to j_0})^2} U(0) B(0) U(0)^T U
\]
and thus $U(t) = \exp \left( \int_0^t \frac{1}{f(\mu_{i_0 \to j_0}(s))^2} \: ds  \: C(0) \right) U(0)$.
% C(t) constant ?
%and $C(t) = U(0)^T \exp \left( -\int_0^t \frac{1}{f(\mu_{i_0 \to j_0})^2(s)} \: ds  \: C(0) \right) U(0) C(0) U(0)^T \exp \left( -\int_0^t \frac{1}{f(\mu_{i_0 \to j_0})^2(s)} \: ds  \: C(0) \right) U(0)$
\end{remk}

\subsection{Numerical examples of geodesics}
\label{subsec:geodNum}

Proposition~\ref{prop:shortestDeltanOn} provides a differential system to compute smooth geodesics in $(M(n), g)$.
While both the qualitative study and the numerical approximation of Proposition~\ref{prop:shortestDeltanOn} $(i)$--$(ii)$ is a broad issue, we propose hereafter a simple numerical scheme to obtain numerical examples of such geodesics.
Note that we address the issue of computing a numerical geodesic given an initial point and an initial speed, which is of course much easier than computing a numerical geodesic between two points.
While $(i)$ can be numerically tackled by a classical ODE scheme, we approximate a solution of $(ii)$ on an interval $[t , t + \delta t]$ with $U(t + \delta t) = U(t) \exp (\delta t B(t))$ (instead of $U(t + \delta t) = U(t) (I_n + \delta t B(t))$ for instance), in order to ensure $U(t) \in \xO(n) \, \Rightarrow \, U(t+\delta t) \in \xO(n)$. 

\subsubsection*{Implementation} 
The source code is available at \url{https://www.imo.universite-paris-saclay.fr/~blanche.buet/research_buet.html} (in the Notebook section), for a static version, see \url{https://www.imo.universite-paris-saclay.fr/~blanche.buet/research_buet_files/flagfoldsDemo.html}.

Let us be more precise concerning the numerical computations. Given a uniform discretization of $[0,L]$ with $0 = t_0 < t_1 < \ldots < t_N = L$: for all $p$, $t_{p+1} - t_p = h = \frac{L}{N}$, we propose the following numerical scheme to approximate $(\mu(t_p) , U(t_p) )$ with $(\mu^{(p)}, U^{(p)})$ in $\mathring{\Delta}(n) \times {\rm O}(n)$:

$\bullet$ Initial data are $(\mu^{(0)}, U^{(0)}) \in \mathring{\Delta}(n) \times {\rm O}(n)$ and $(\dt{\mu})^{(0)} \in \{(x_1, \ldots, x_n) \in \mathbb{R}^n \: : \: x_1 + \ldots + x_n = 0 \}$, $B^{(0)} \in {\rm Skew(n)}$, from which we can compute $ C^{(0)} = \left( f\left(\mu_{i \to j}^{(0)}\right)^2 b_{ij}^{(0)} \right)_{i<j} \in \xSkew(n)$.

$\bullet$ For $p \in \{0, \ldots, N-1\}$,

\begin{align*}
b_{ij}^{(p)} & = \frac{1}{f(\mu_{i \to j}^{(p)})^2} \left( (U^{(p)})^T U^{(0)} C^{(0)} (U^{(0)})^T U^{(p)} \right)_{ij} \quad \text{for } 1 \leq i < j \leq n\\
\mu^{(p+1)} & = \mu^{(p)} + h(\dt{\mu})^{(p)} \\
(\dt{\mu}_k)^{(p+1)} & = (\dt{\mu}_k)^{(p)} +  \frac{h}{n}\sum_{\substack{l = 1 \\ l \neq k}}^n \sum_{1 \leq i \leq l < j \leq n} f(\mu_{i \to j}^{(p)}) \: \partial_l f (\mu_{i \to j}^{(p)}) (b_{ij}^{(p)})^2  \\
& \qquad \, \qquad - \frac{h}{n} (n-1) \sum_{1 \leq i \leq k < j \leq n} f(\mu_{i \to j}^{(p)}) \: \partial_k f (\mu_{i \to j}^{(p)}) (b_{ij}^{(p)})^2 \quad \text{for } k \in \{1, \ldots, n\} \\
U^{(p+1)} & = U^{(p)} \exp \left(h B^{(p)} \right)
\end{align*}

\subsubsection*{Numerical examples}
We test the implementation for $n=3$ and $f(\nu) = \frac{1}{4}|\nu|$, with a time-step $h = 0.001$ and for different initial data.
The numerical results are collected in Figures~\ref{figNumGeod1} to \ref{figNumGeod-1}, we now explain the common features displayed in the subfigures.
For each numerical geodesic we first represent the evolution of the weights with respect to time, for $i \in \{1,2,3\}$, $t \mapsto \mu_i(t)$ (subfigure $(a)$), $t \mapsto \lambda_i(t)$ (subfigure $(b)$), $t \mapsto \dt{\mu}_i(t)$ (subfigure $(c)$). We also check that $\mu_1 + \mu_2 + \mu_3$ is constant equal to $1$ while $\dt{\mu}_1 + \dt{\mu}_2 + \dt{\mu}_3$ is constant equal to $0$. We then represent the image of the curve $t \mapsto \mu(t)$ in the simplex $\Delta(3)$ (subfigure $(d)$): as $\Delta(3)$ is included in the plane $x + y + z = 1$, we are able to obtain a $2D$ view thanks to the transformation $(\mu_1, \mu_2, \mu_3) \mapsto \mu_1 (-1,0) + \mu_2 (1,0) + \mu_3 (0, \sqrt{3})$.

%Recalling that the resulting geodesic is arclength parametrized, we represent the evolution of the numerical arclength 
%\[
%h \sqrt{| { \mu^{(p)}}^\prime |^2 + \sum_{1\leq i<j \leq n} f \left(\mu^{(p)}_{i \to j} \right)^2 (b_{ij}^{(p)})^2}
%\]
%with respect to time $t_p$ (subfigure $(d)$).
Representing the evolution of the flag $t \mapsto U(t)$ is more intricate,
for the sake of visualisation, we associate with $(\mu, U) \in \Delta(3) \times \xO(3)$ an ellipsoid $\cE ll = \cE ll(\mu,U)$ as follows: we compute $\lambda_1 = \mu_1 + \frac{\mu_2}{2} + \frac{\mu_3}{3}$, $\lambda_2 = \frac{\mu_2}{2} + \frac{\mu_3}{3}$, $\lambda_3 = \frac{\mu_3}{3}$ and we consider the ellipsoid of axes directed by the columns of $U$: $u_1$, $u_2$, $u_3$ with respective semi-axis lengths $\sqrt{3 \lambda_1}$, $\sqrt{3 \lambda_2}$, $\sqrt{3 \lambda_3}$. For instance, if $(\mu,U) = \left( (1,0,0), U \right)$ corresponds to a line, $\cE ll(\mu,U)$ is a segment directed by $u_1$, if now $(\mu,U) = \left( (0,1,0), U \right)$ corresponds to a plane, $\cE ll(\mu,U)$ is a disk of radius $\sqrt{3/2}$ in the plane spanned by $(u_1,u_2)$, and finally, if $(\mu,U) = \left( (0,0,1), U \right)$ corresponds to $\R^3$, $\cE ll(\mu,U)$ is a sphere of radius $1$. In subfigures $(e),(f),(g)$ we show projections of the ellipsoids $\cE ll \left(\mu^{(p)}, U^{(p)} \right)$ at times $t_p$ (evenly selected in $[0,L]$) in the $3$ planes $x = 0$, $y = 0$ and $z=0$. To complete the visualisation of the evolution of $U$, we compute at each time $t_p$ the principal angles $\theta_1 = \arccos \left(u_1^{(0)} \cdot u_1^{(p)} \right)$ between $\xspan (u_1^{(0)})$ and $\xspan(u_1^{(p)})$, and $\theta_2$ between $\xspan ( u_1^{(0)}, u_2^{(0)} )$ and $\xspan ( u_1^{(p)}, u_2^{(p)} )$. We compute $\theta_2$ through the singular value decomposition $\sigma_1 \geq \sigma_2 = 0$ of $(u_1^{(0)} , u_2^{(0)})^T (u_1^{(p)} , u_2^{(p)})$: $\theta_2 = \arccos(\sigma_1)$
(see \cite{BjorckGolub}).
The time evolution of $\theta_1$ and $\theta_2$ is displayed in subfigure $(h)$.

\begin{center}
\setcounter{subfigure}{0}
\begin{figure}[!htp]
\subcaptionbox{$t \mapsto \mu(t)$}{\includegraphics[width=0.24\textwidth]{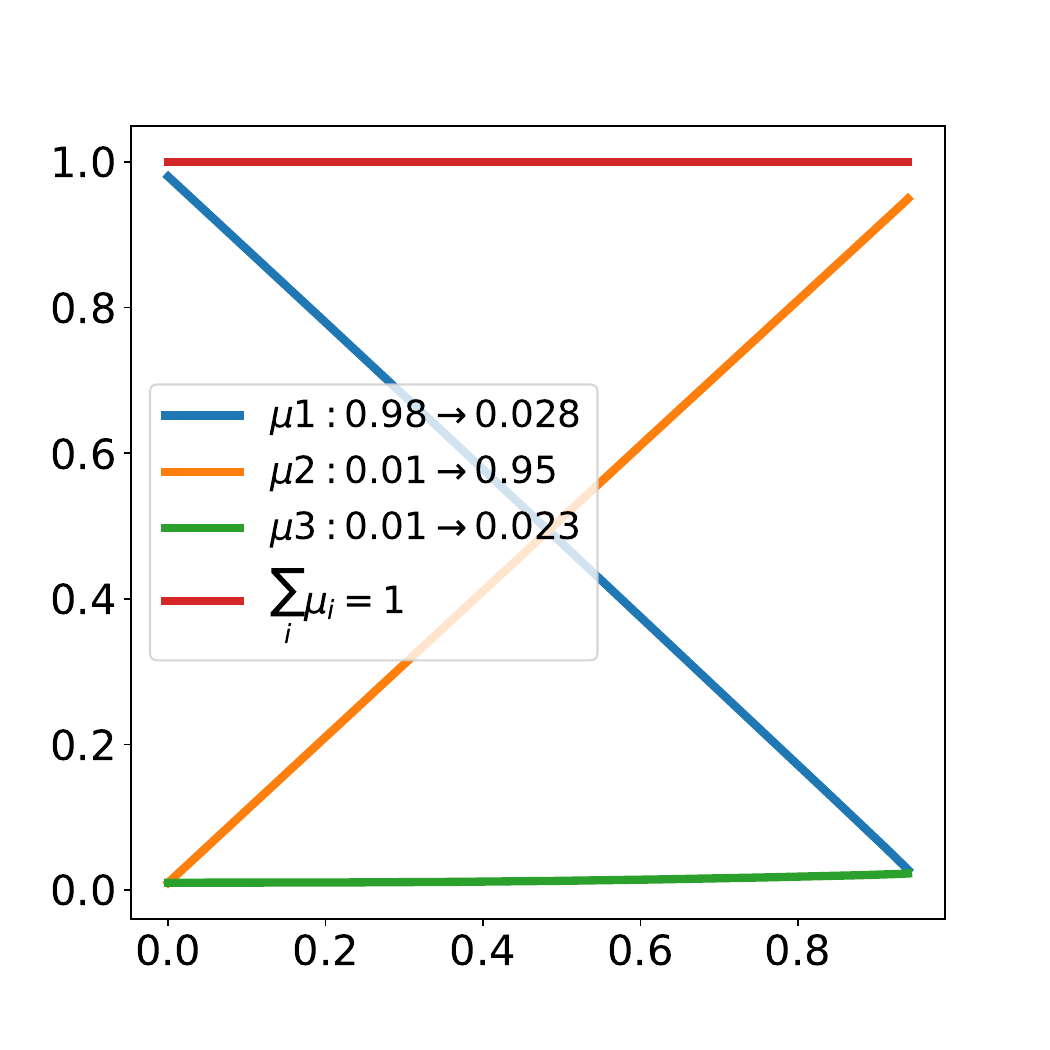}}
\subcaptionbox{$t \mapsto \lambda(t)$}{\includegraphics[width=0.24\textwidth]{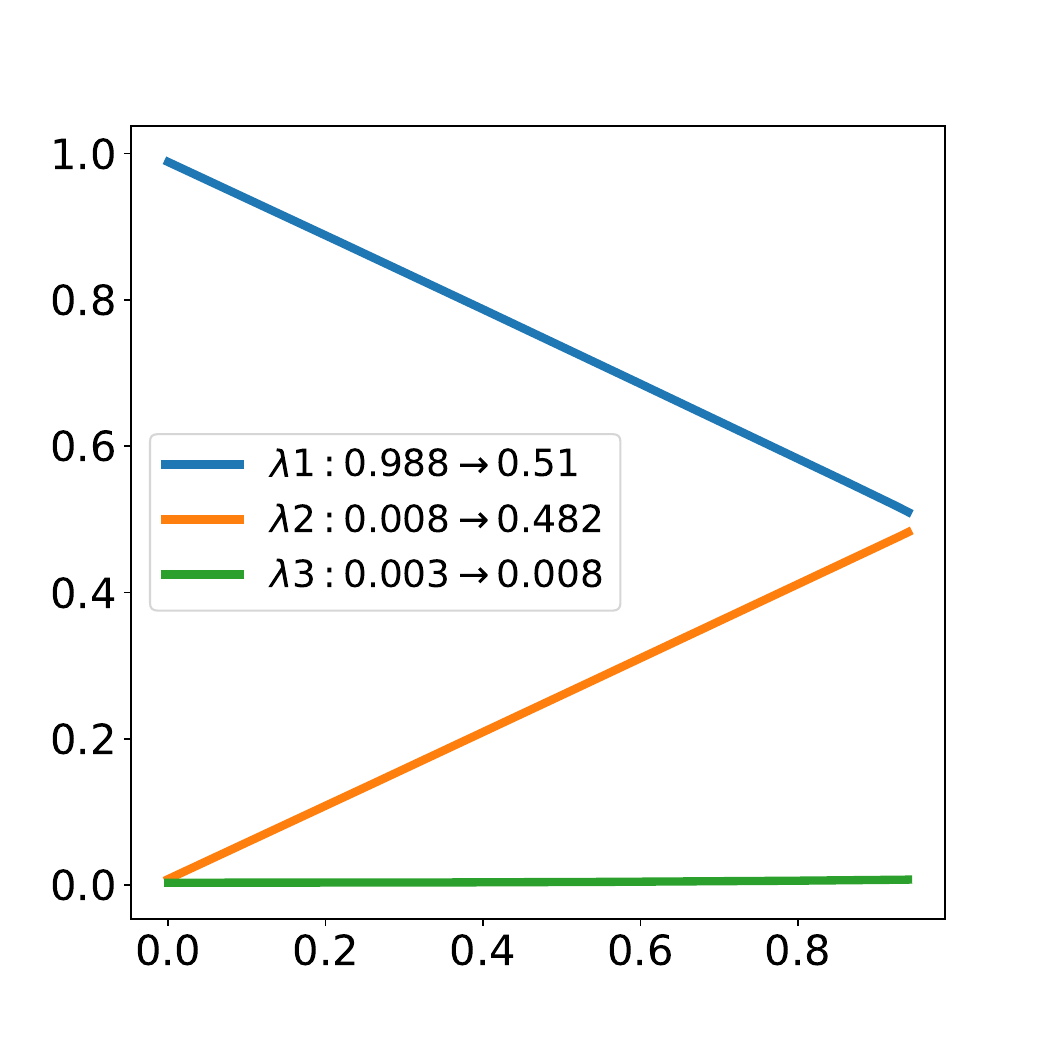}}
\subcaptionbox{$t \mapsto \dt{\mu}(t)$}{\includegraphics[width=0.24\textwidth]{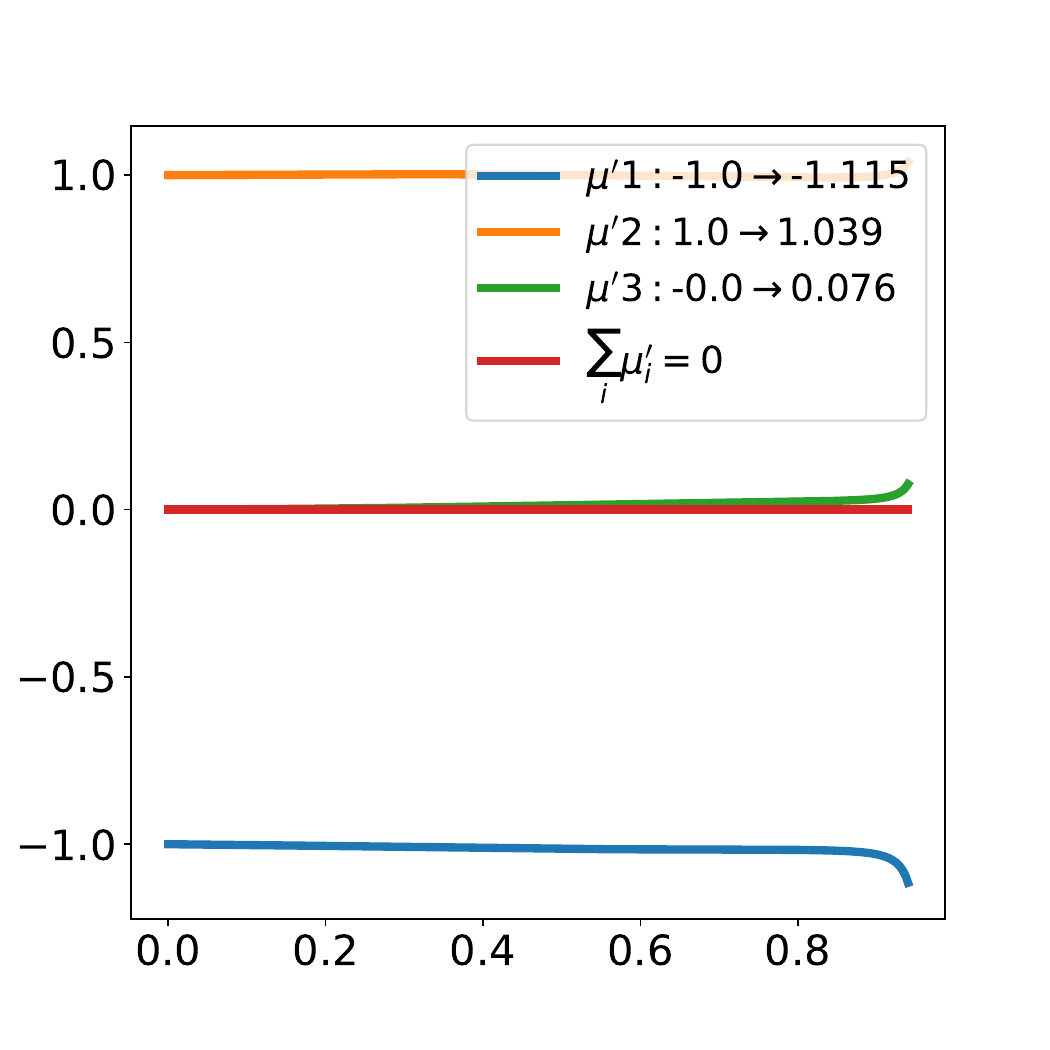}}
\subcaptionbox{Image of $\mu$ in $\Delta(3)$.}{\includegraphics[width=0.24\textwidth]{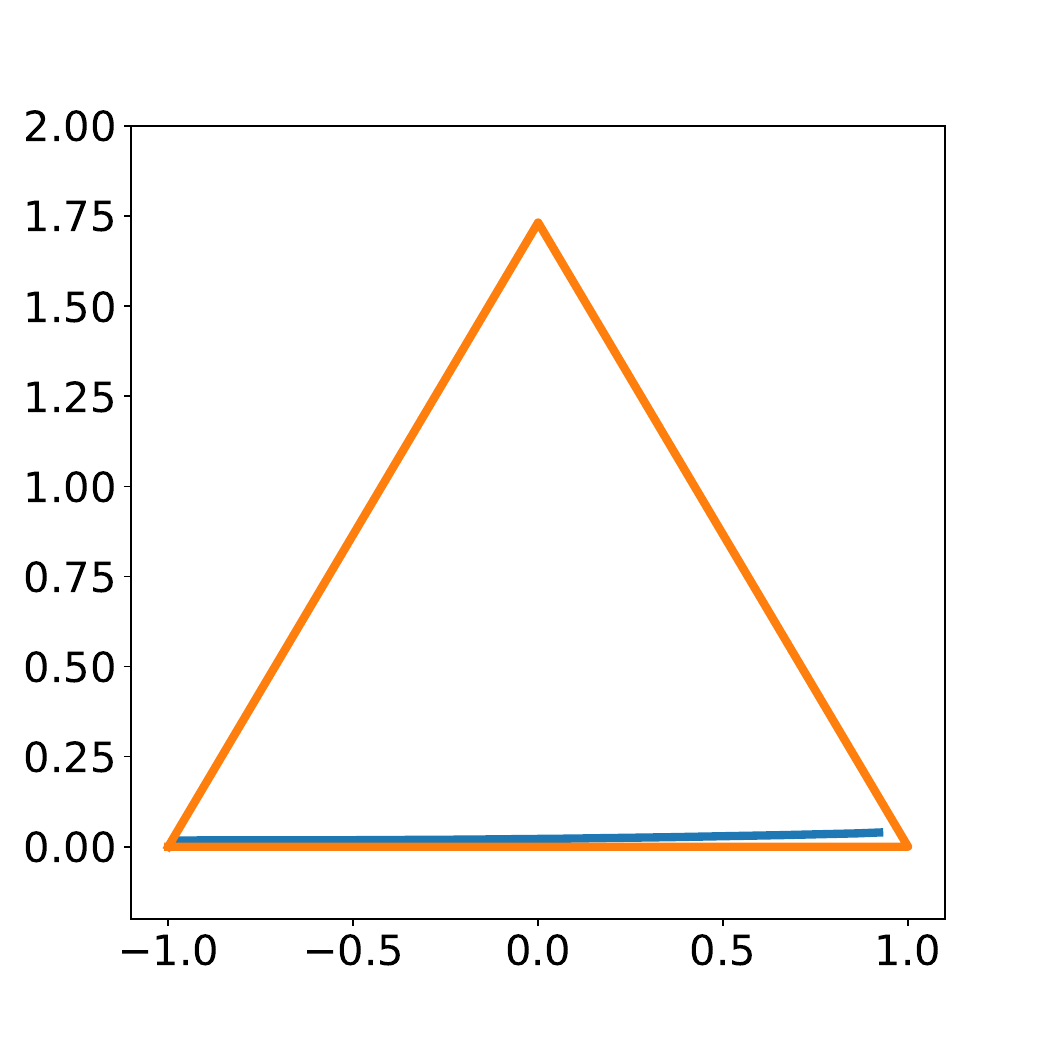}}
%\subcaptionbox{$t \mapsto \sqrt{g_{\gamma(t)} \left( \gamma^\prime(t), \gamma^\prime(t)\right)}$}{\includegraphics[width=0.24\textwidth]{arclengthLineToPlane.pdf}}
%
\subcaptionbox{Plane $x = 0$}{\includegraphics[width=0.24\textwidth]{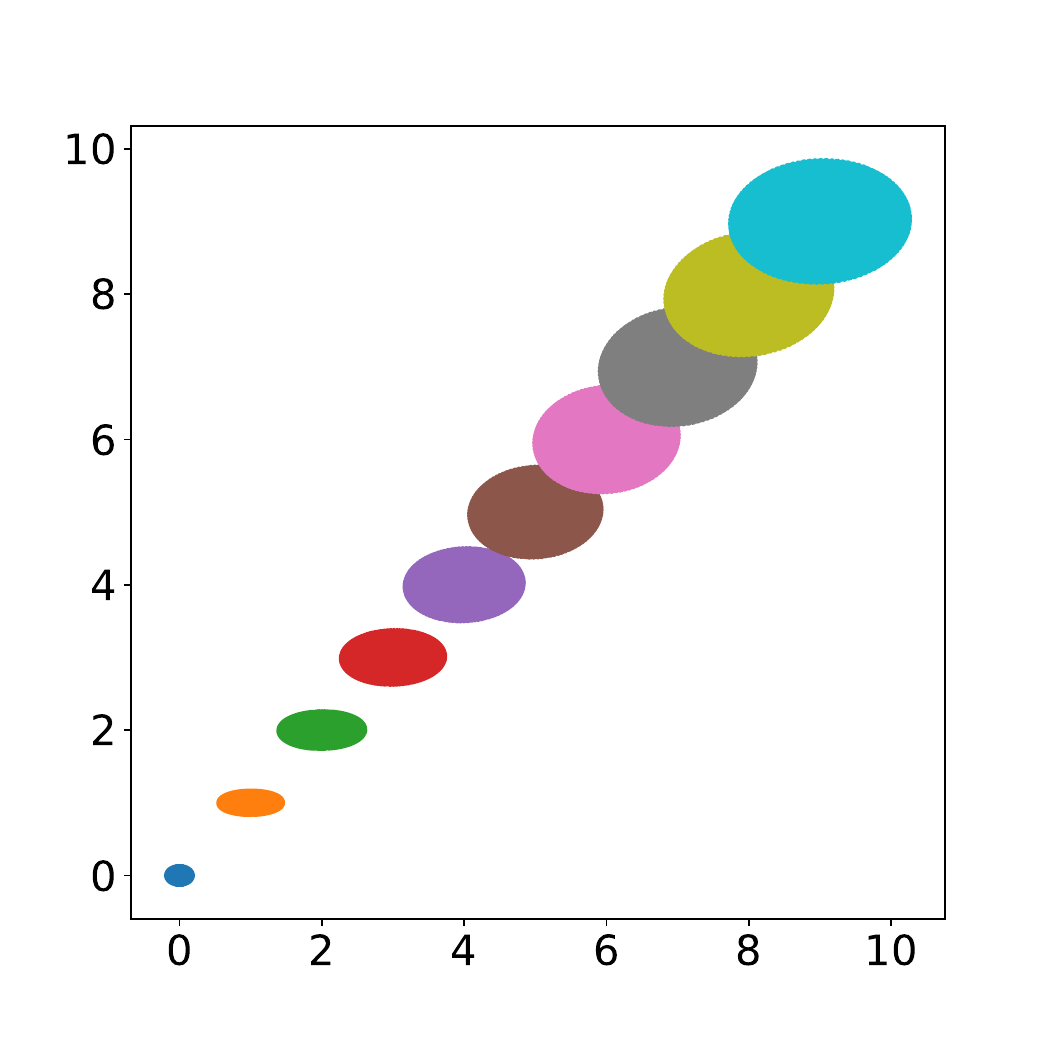}}
\subcaptionbox{Plane $y=0$}{\includegraphics[width=0.24\textwidth]{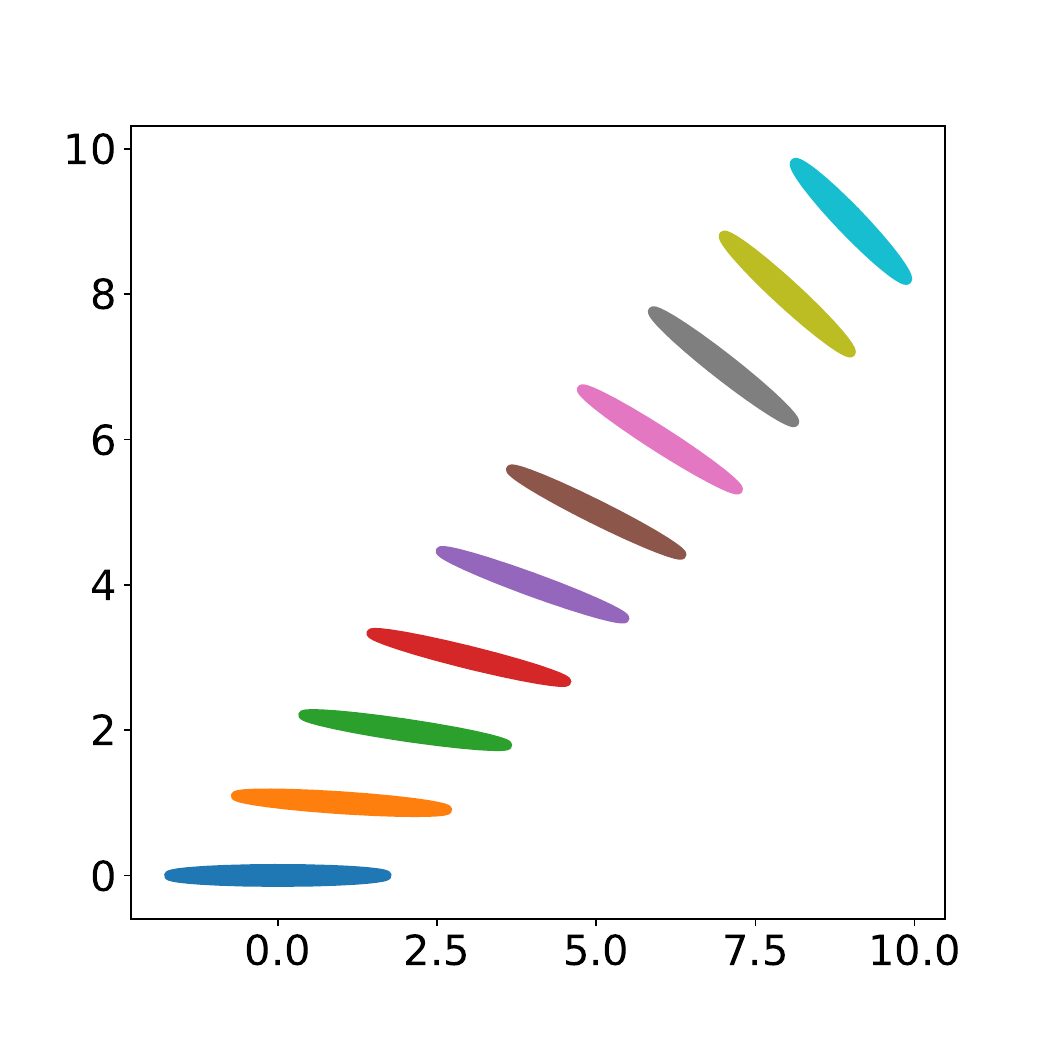}}
\subcaptionbox{Plane $z = 0$}{\includegraphics[width=0.24\textwidth]{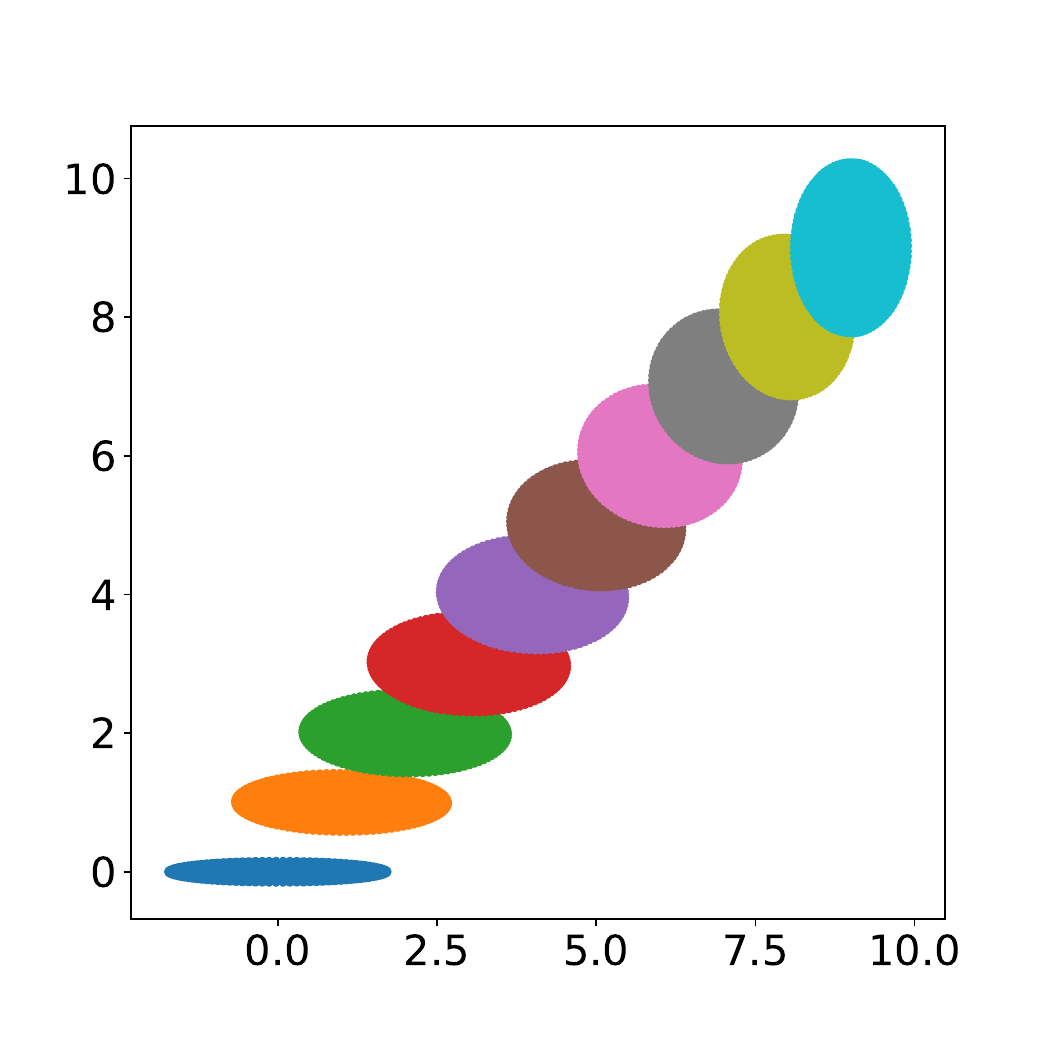}}
\subcaptionbox{Principal angles.}{\includegraphics[width=0.24\textwidth]{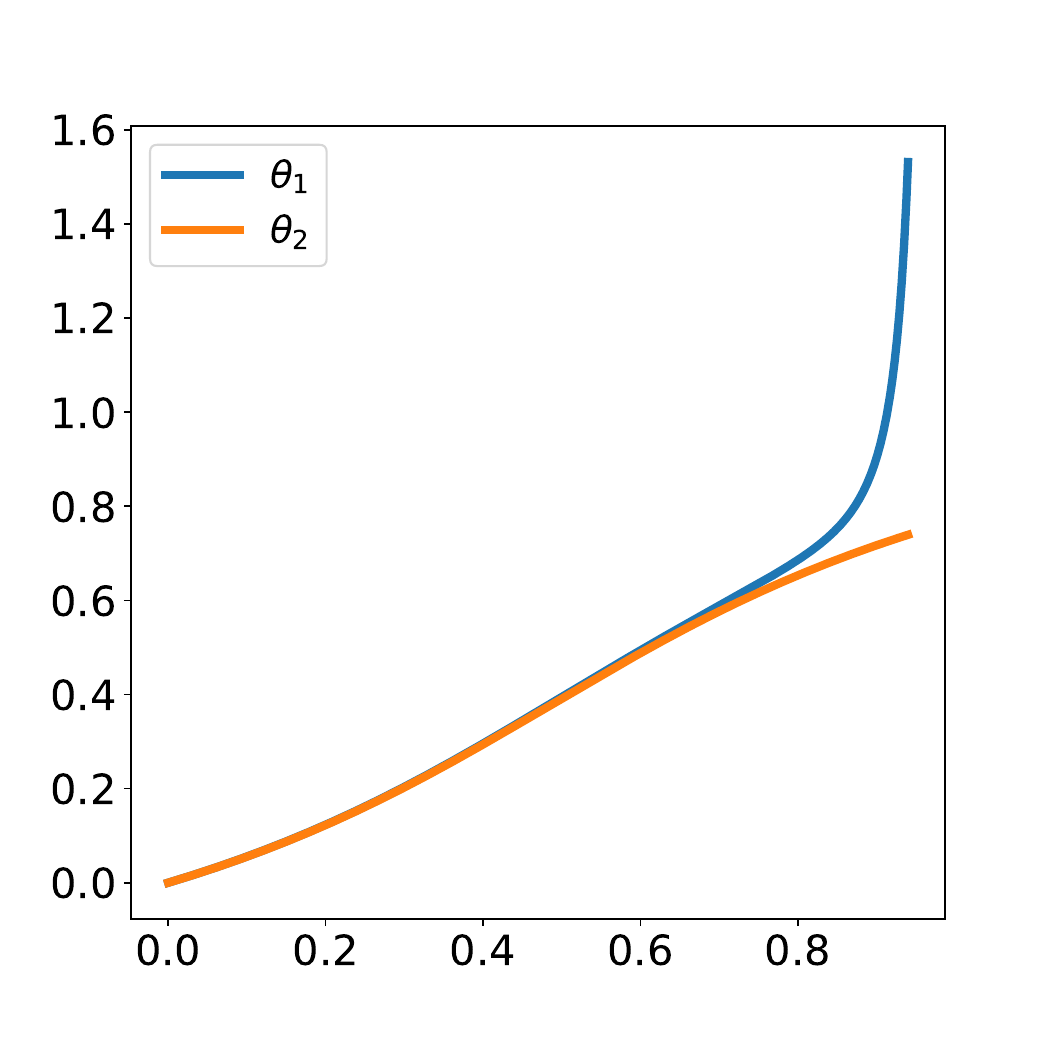}}
\caption{We observe an evolution between a weighted flag close to a line at initial time and close to a plane at final time. Initial data: $\mu^{(0)} = (0.98, 0.01, 0.01) $, $U^{(0)} = I_3$, ${\dot{\mu}}^{(0)} = (-1, 1, 0)$, $(b_{12}, b_{23}, b_{13})^{(0)} = (0.05, 0, 0.5)$. Final point: $\mu^{(N)} = (0.028, 0.95,  0.023)$, $U^{(N)}${\tiny $~=~\left(
\begin{array}{ccc}
0.039 & 0.738 & 0.674\\
-0.997 & 0.072 & -0.021\\
-0.064 & -0.671 & 0.739
\end{array} \right) $}. \label{figNumGeod1}}
\end{figure}
%%%%%%%%%%%%%%%%%%%%%%%%%%%%%%%%%%%%%%%%%%%%%%%%%%%%%%
%
%%%%%%%%%%%%%%%%%%%%%%%%%%%%%%%%%%%%%%%%%%%%%%%%%%%%%
%
\setcounter{subfigure}{0}
\begin{figure}[!htp] 
\subcaptionbox{$t \mapsto \mu(t)$}{\includegraphics[width=0.24\textwidth]{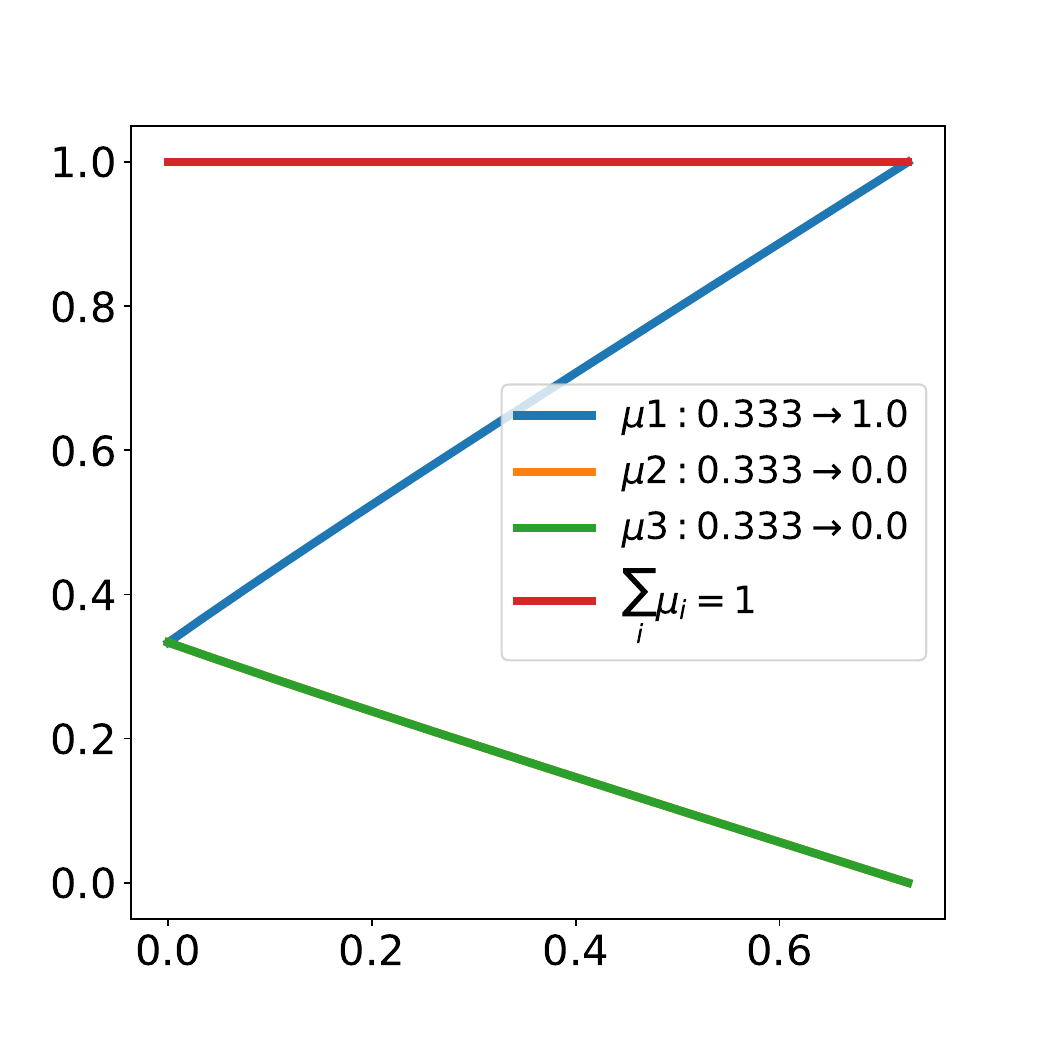}}
\subcaptionbox{$t \mapsto \lambda(t)$}{\includegraphics[width=0.24\textwidth]{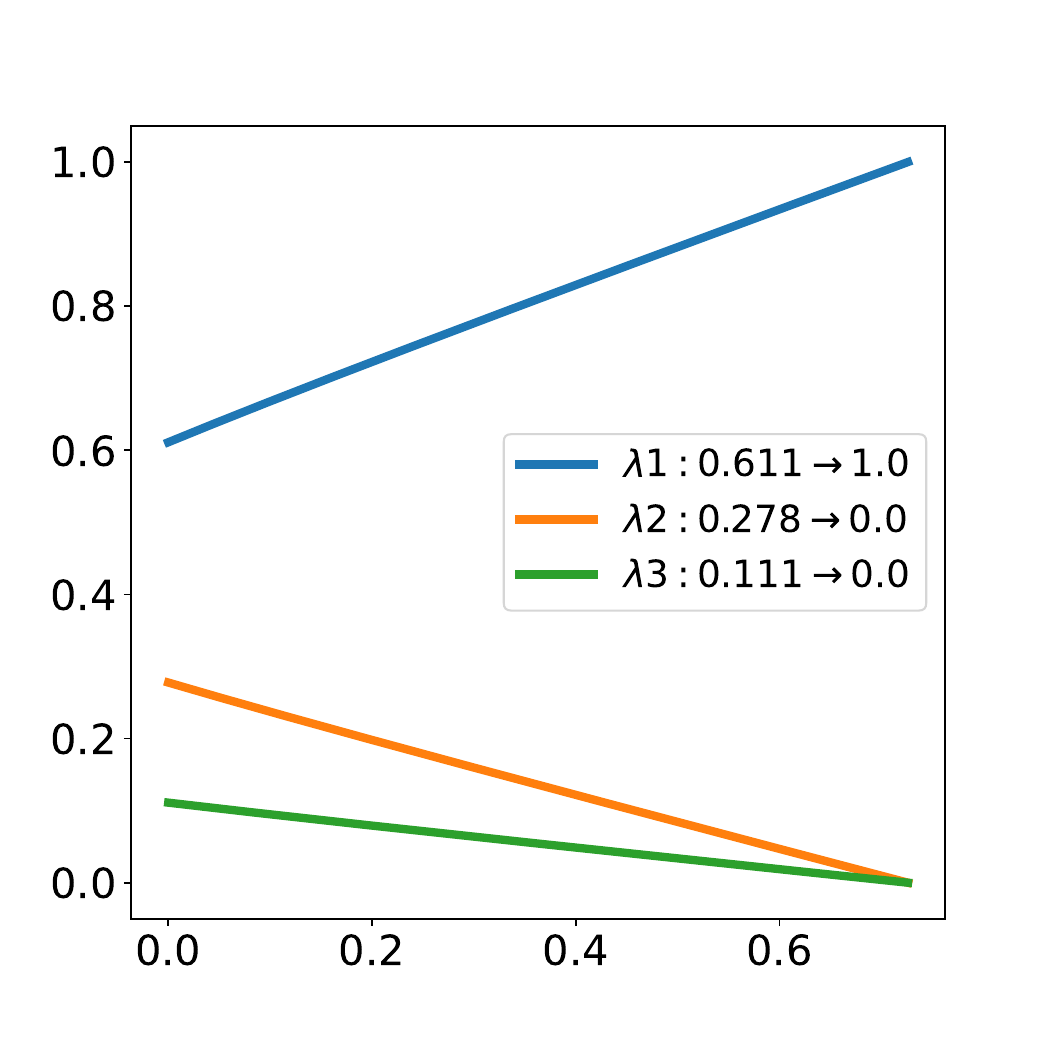}}
\subcaptionbox{$t \mapsto \dt{\mu}(t)$}{\includegraphics[width=0.24\textwidth]{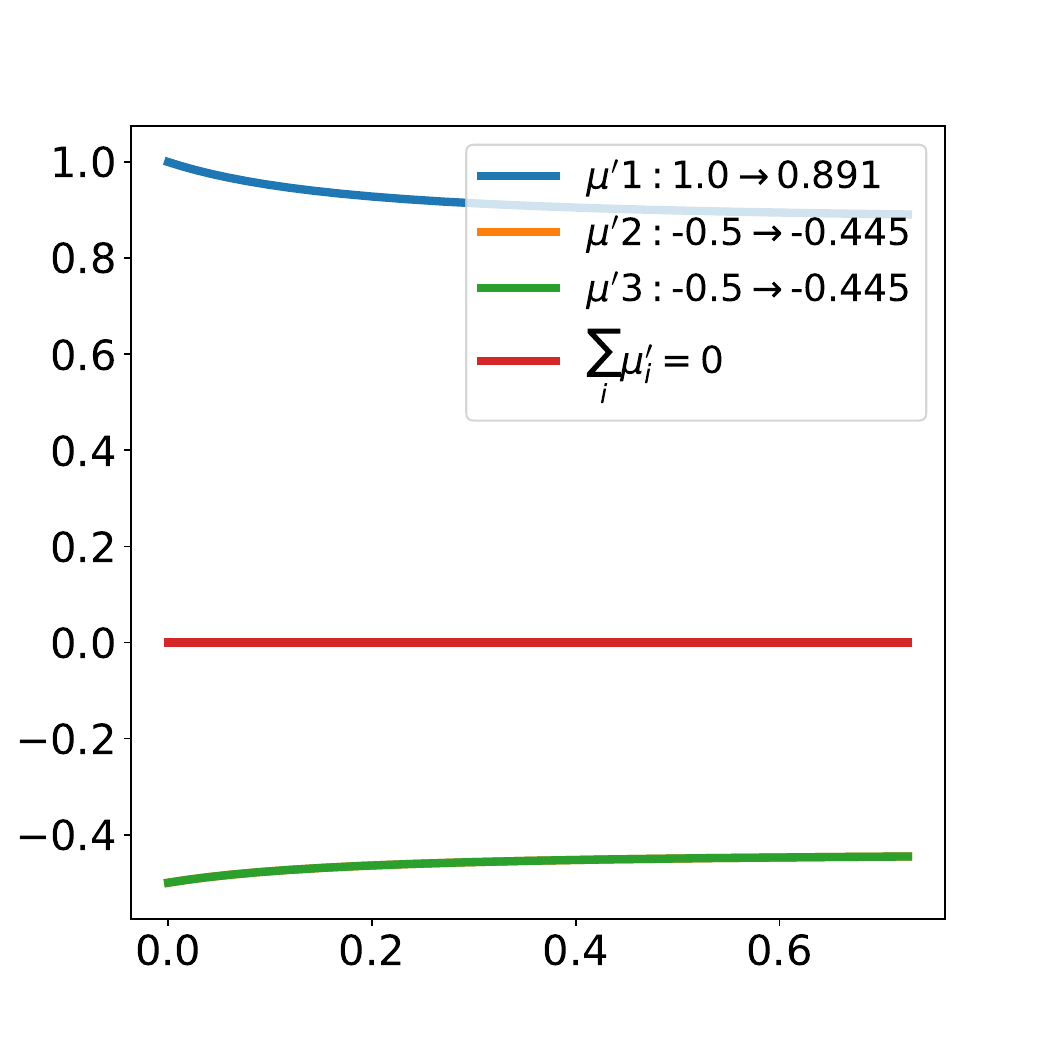}}
\subcaptionbox{Image of $\mu$ in $\Delta(3)$.}{\includegraphics[width=0.24\textwidth]{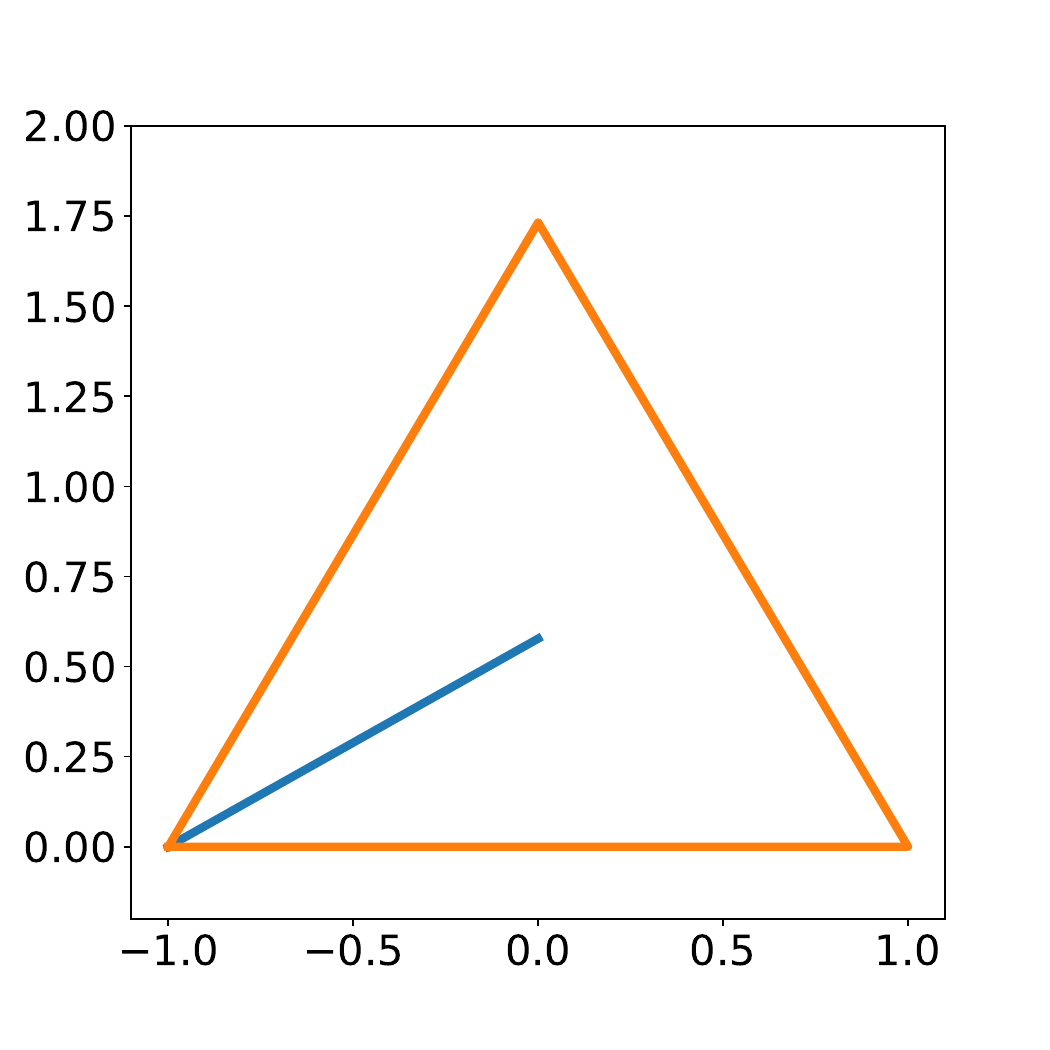}}
%\subcaptionbox{$t \mapsto \sqrt{g_{\gamma(t)} \left( \gamma^\dt(t), \gamma^\dt(t)\right)}$}{\includegraphics[width=0.24\textwidth]{arclengthPlaneLineR3ToLine.pdf}}
%
\subcaptionbox{Plane $x = 0$}{\includegraphics[width=0.24\textwidth]{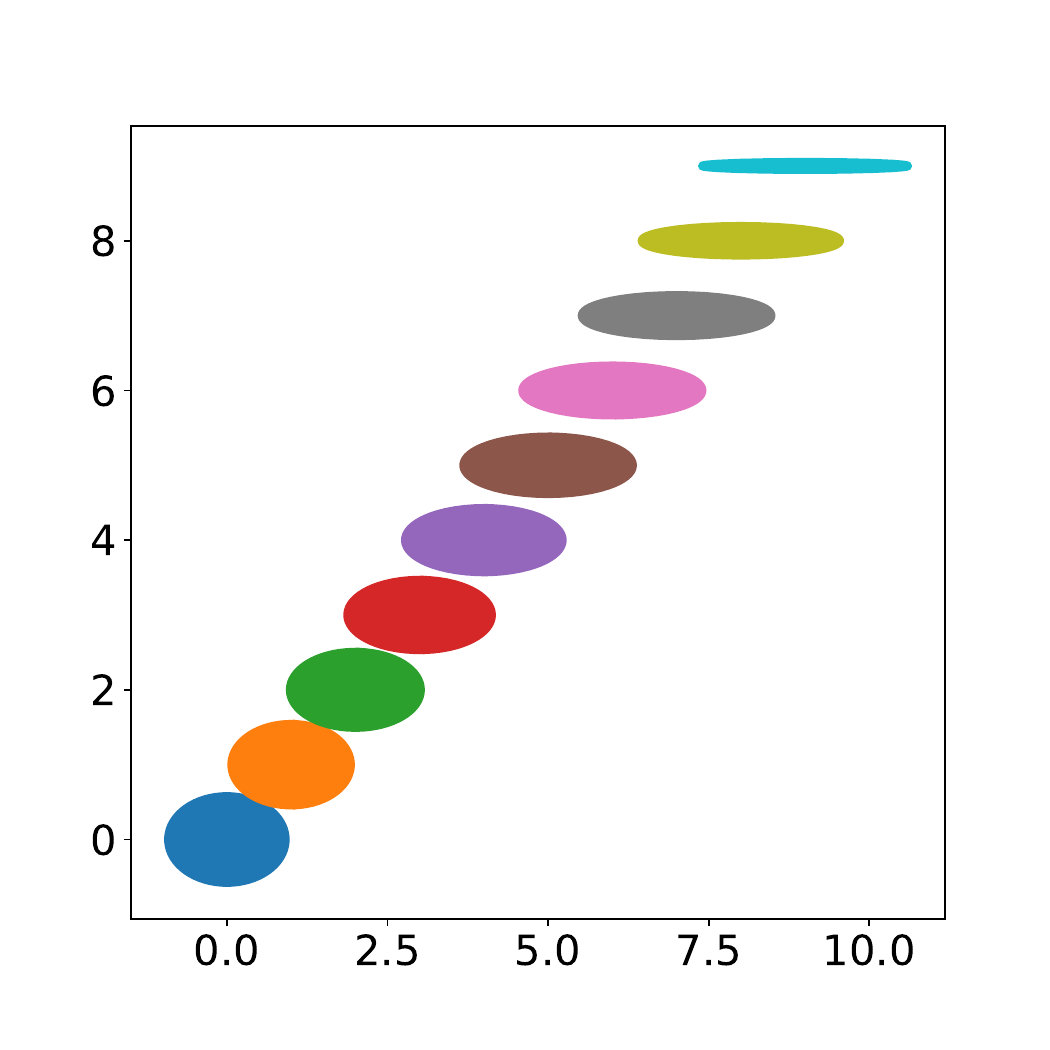}}
\subcaptionbox{Plane $y=0$}{\includegraphics[width=0.24\textwidth]{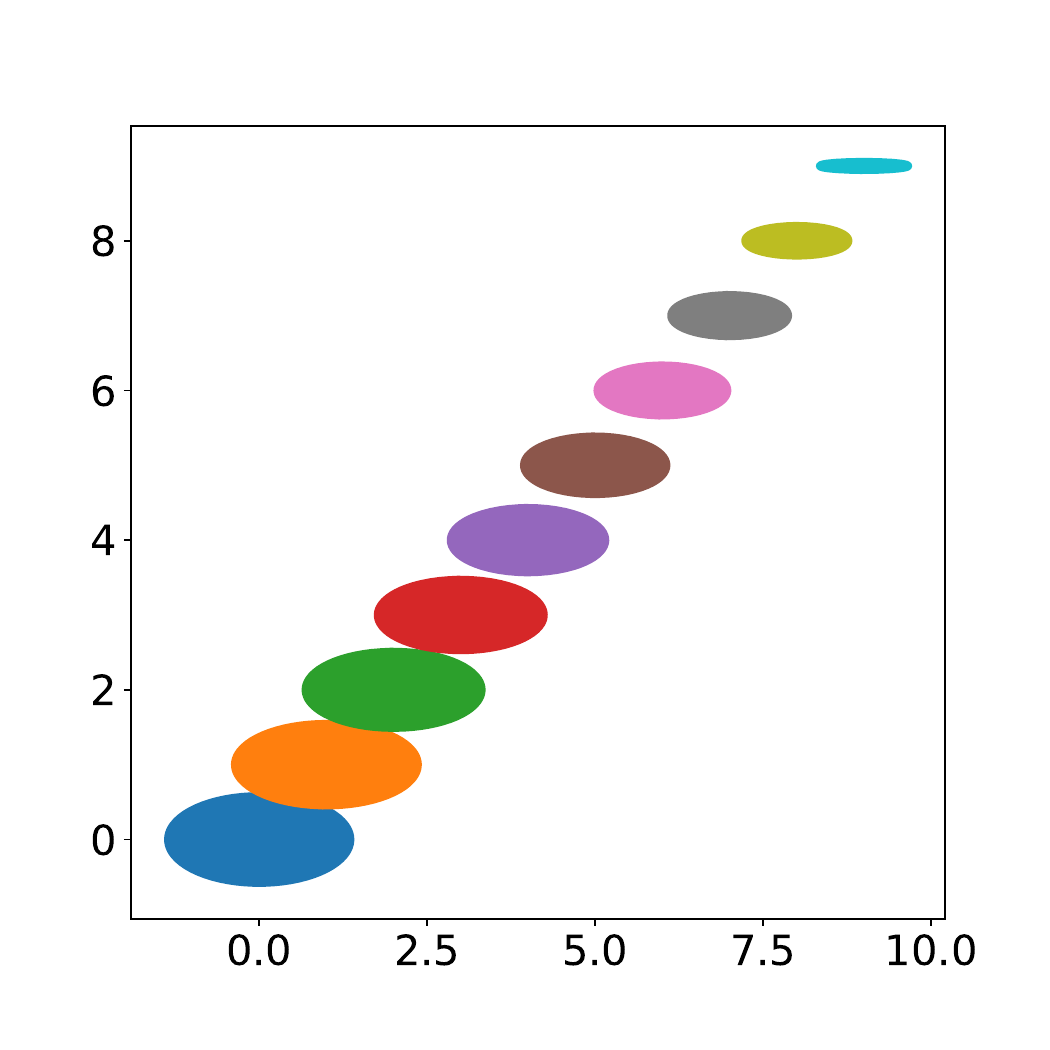}}
\subcaptionbox{Plane $z = 0$}{\includegraphics[width=0.24\textwidth]{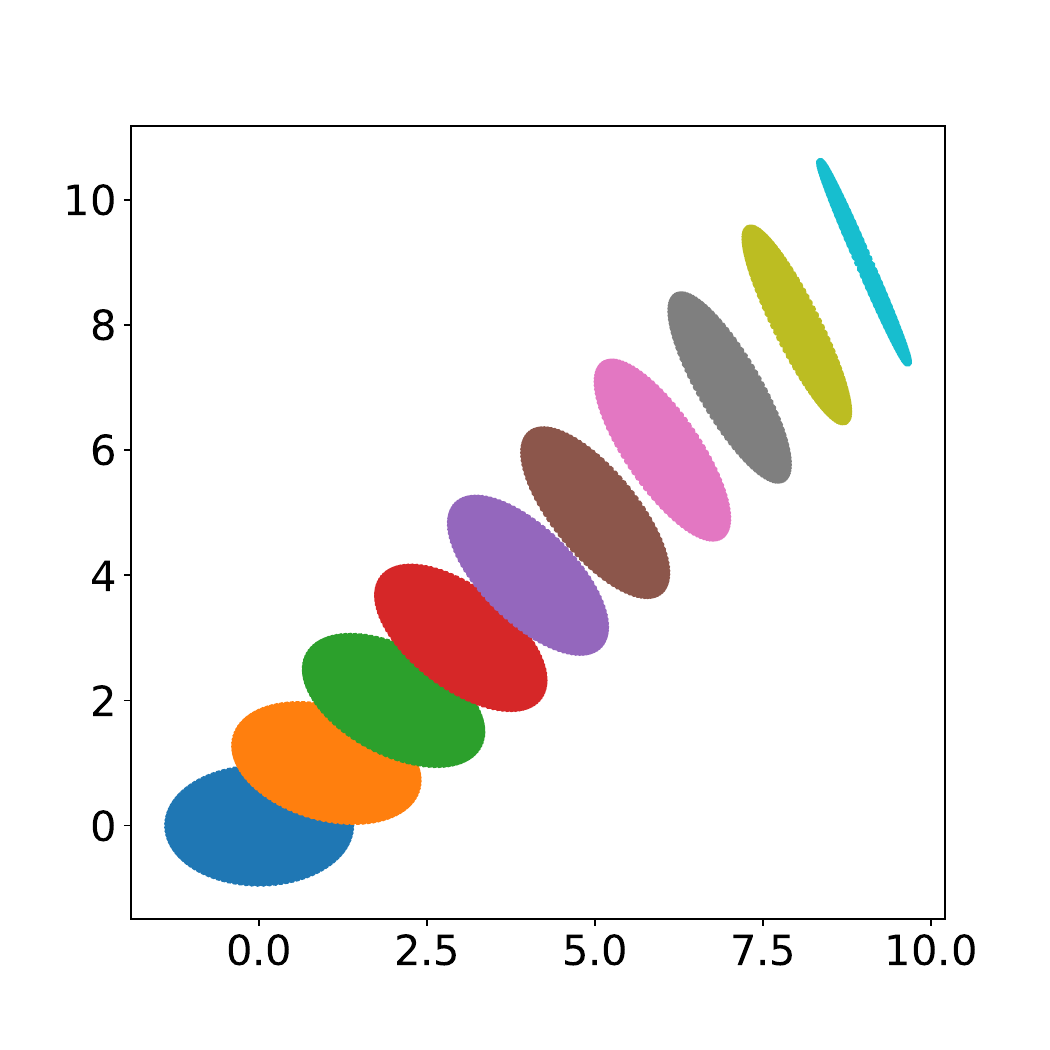}}
\subcaptionbox{Principal angles.}{\includegraphics[width=0.24\textwidth]{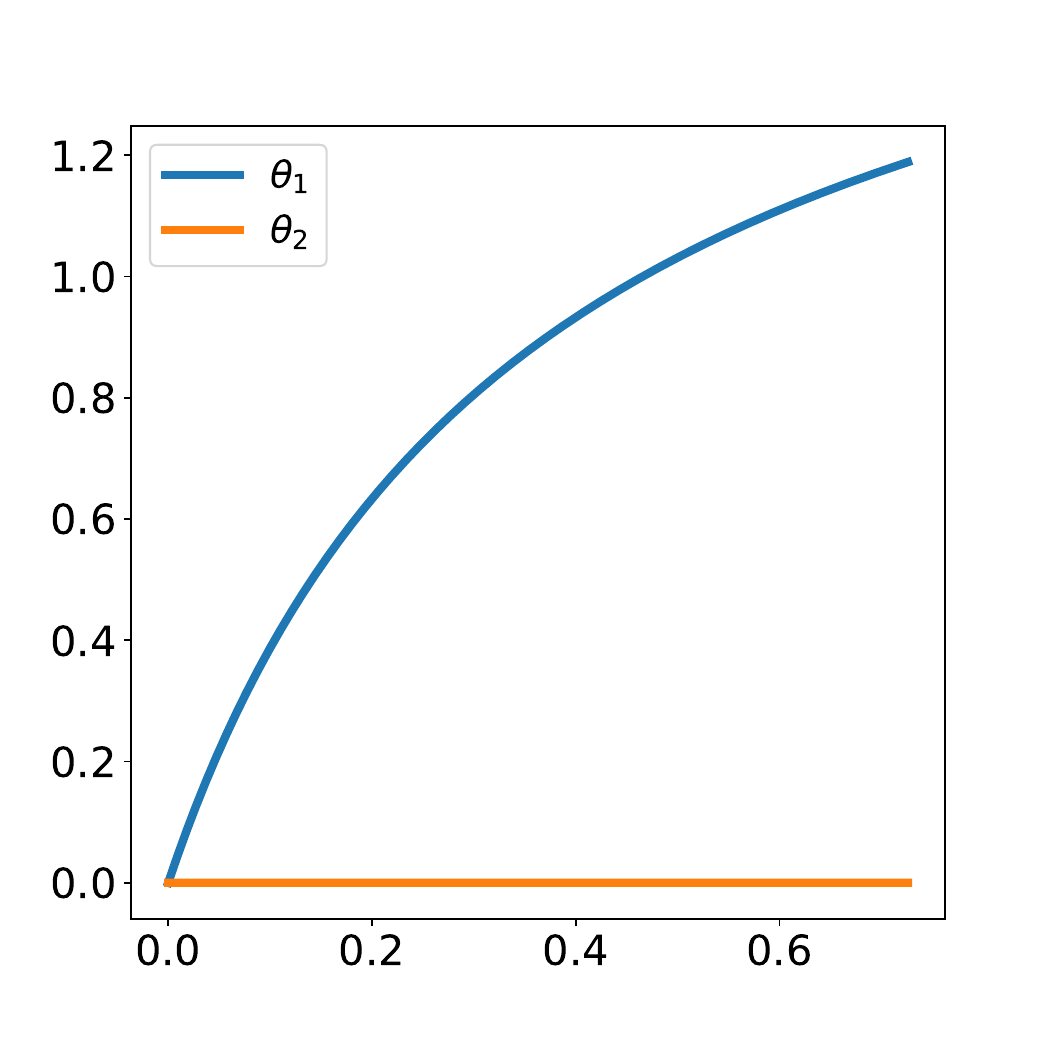}}
\caption{We observe an evolution between a complete weighted flag at initial time and close to a line at final time. Initial data: $\mu^{(0)} = (1/3, 1/3, 1/3) $, $U^{(0)} = I_3$, ${\dot{\mu}}^{(0)} = ( 1,  -0.5, -0.5)$, $(b_{12}, b_{23}, b_{13})^{(0)} = (5, 0, 0)$. Final point: $\mu^{(N)} = (1, 0,  0)$, $U^{(N)}${\tiny $~=~\left(
\begin{array}{ccc}
0.373 & 0.928 & 0\\
-0.928 & 0.373 & 0 \\
0  &  0  & 1
\end{array} \right) $}.}
\end{figure}
\end{center}
%%%%%%%%%%%%%%%%%%%%%%%%%%%%%%%%%%%%%%%%%%%%%%%%%%%%%%
%
%%%%%%%%%%%%%%%%%%%%%%%%%%%%%%%%%%%%%%%%%%%%%%%%%%%%%
%
\begin{center}
\setcounter{subfigure}{0}
\begin{figure}[!htp] 
\subcaptionbox{$t \mapsto \mu(t)$}{\includegraphics[width=0.24\textwidth]{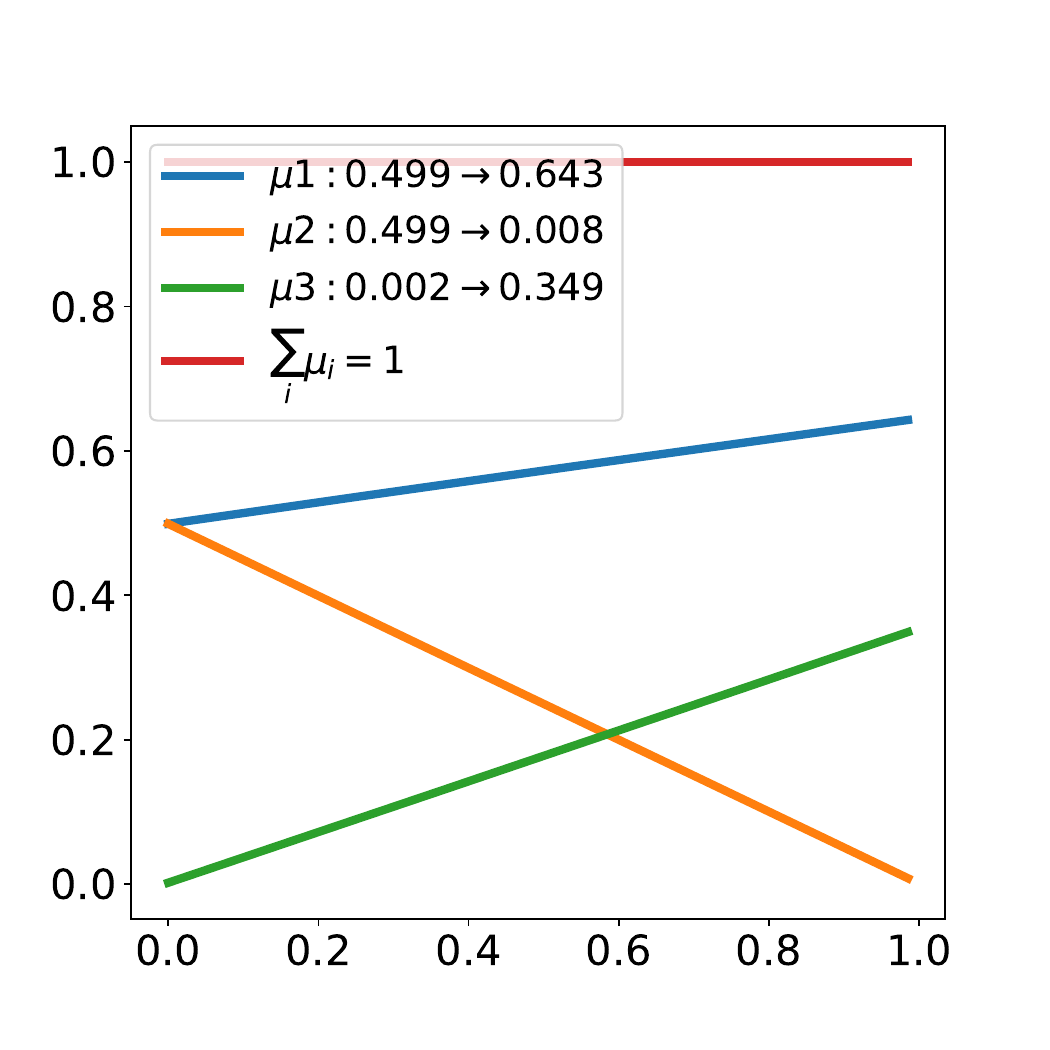}}
\subcaptionbox{$t \mapsto \lambda(t)$}{\includegraphics[width=0.24\textwidth]{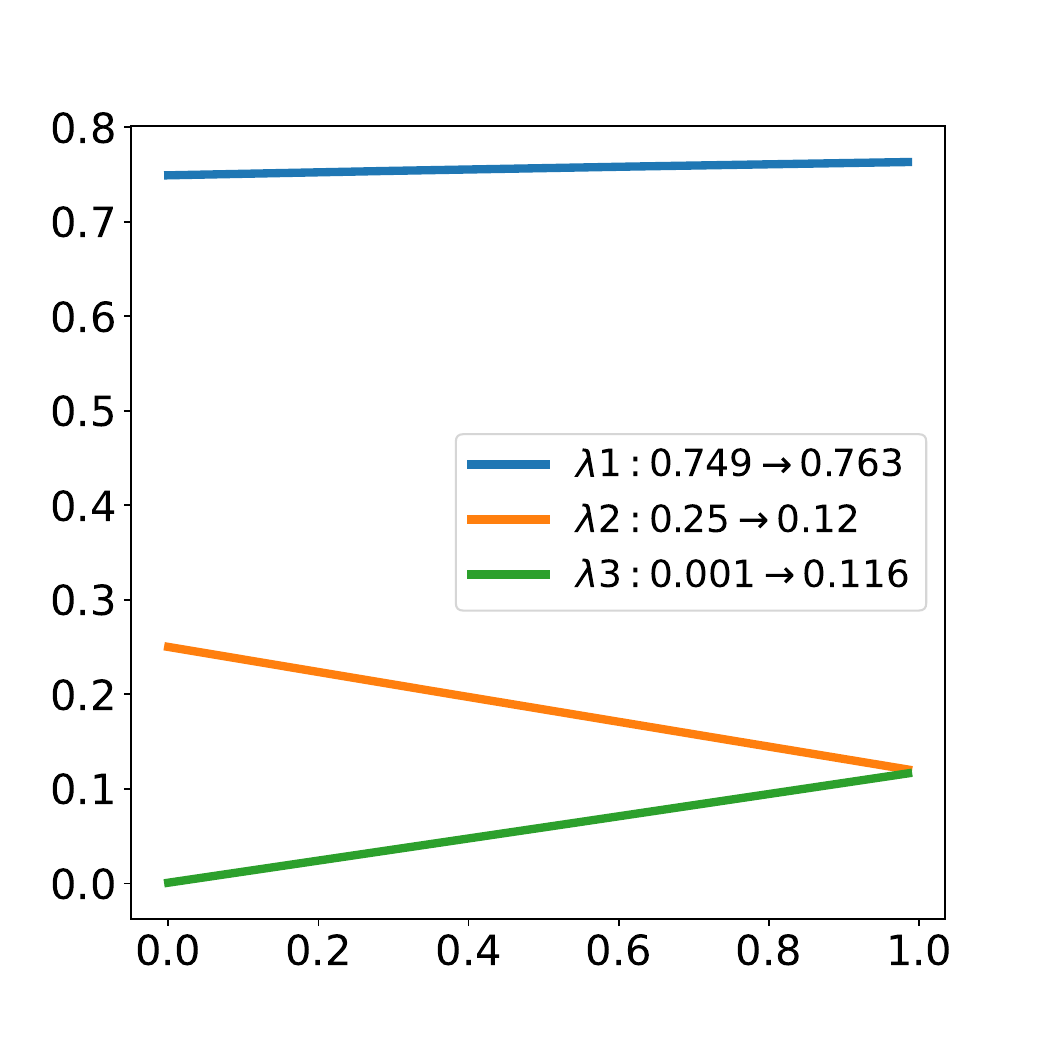}}
\subcaptionbox{$t \mapsto \dt{\mu}(t)$}{\includegraphics[width=0.24\textwidth]{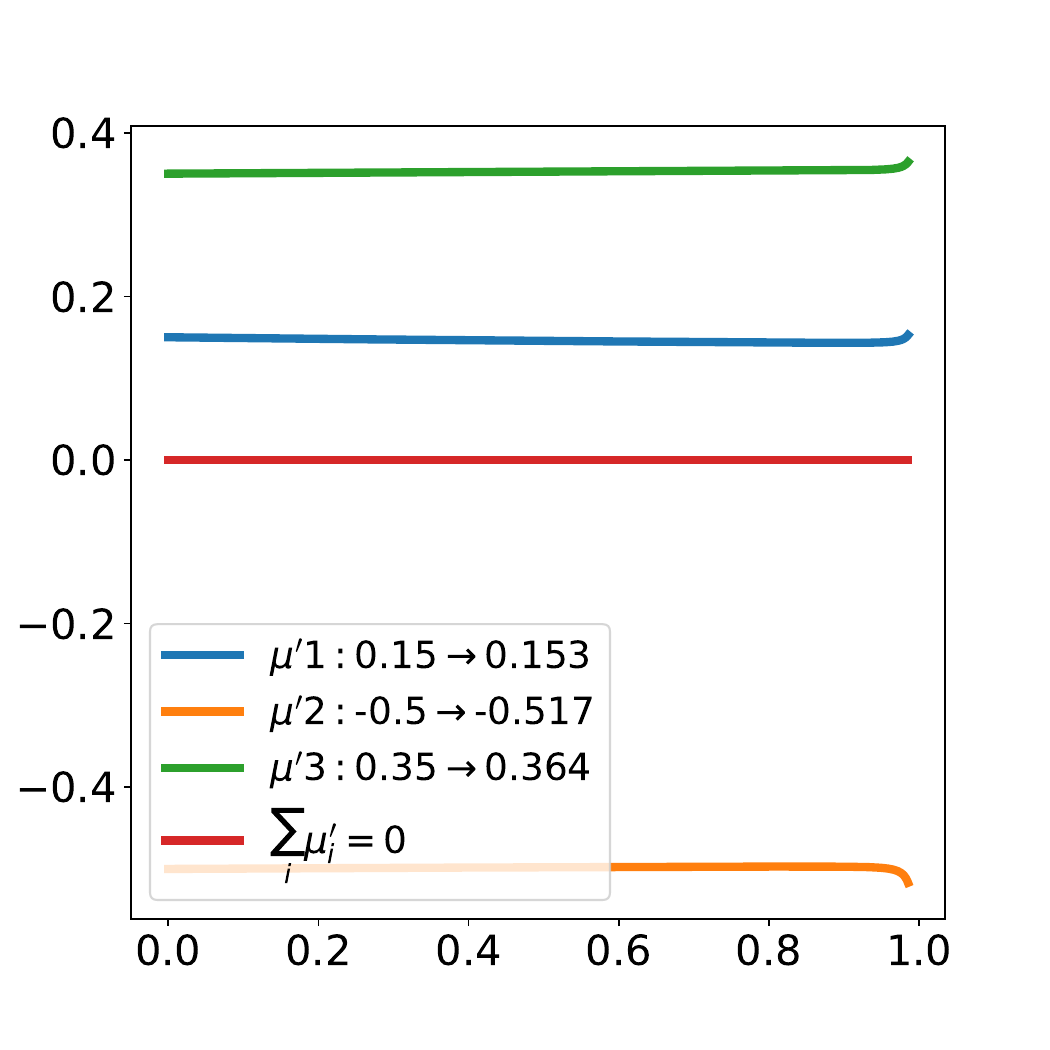}}
\subcaptionbox{Image of $\mu$ in $\Delta(3)$.}{\includegraphics[width=0.24\textwidth]{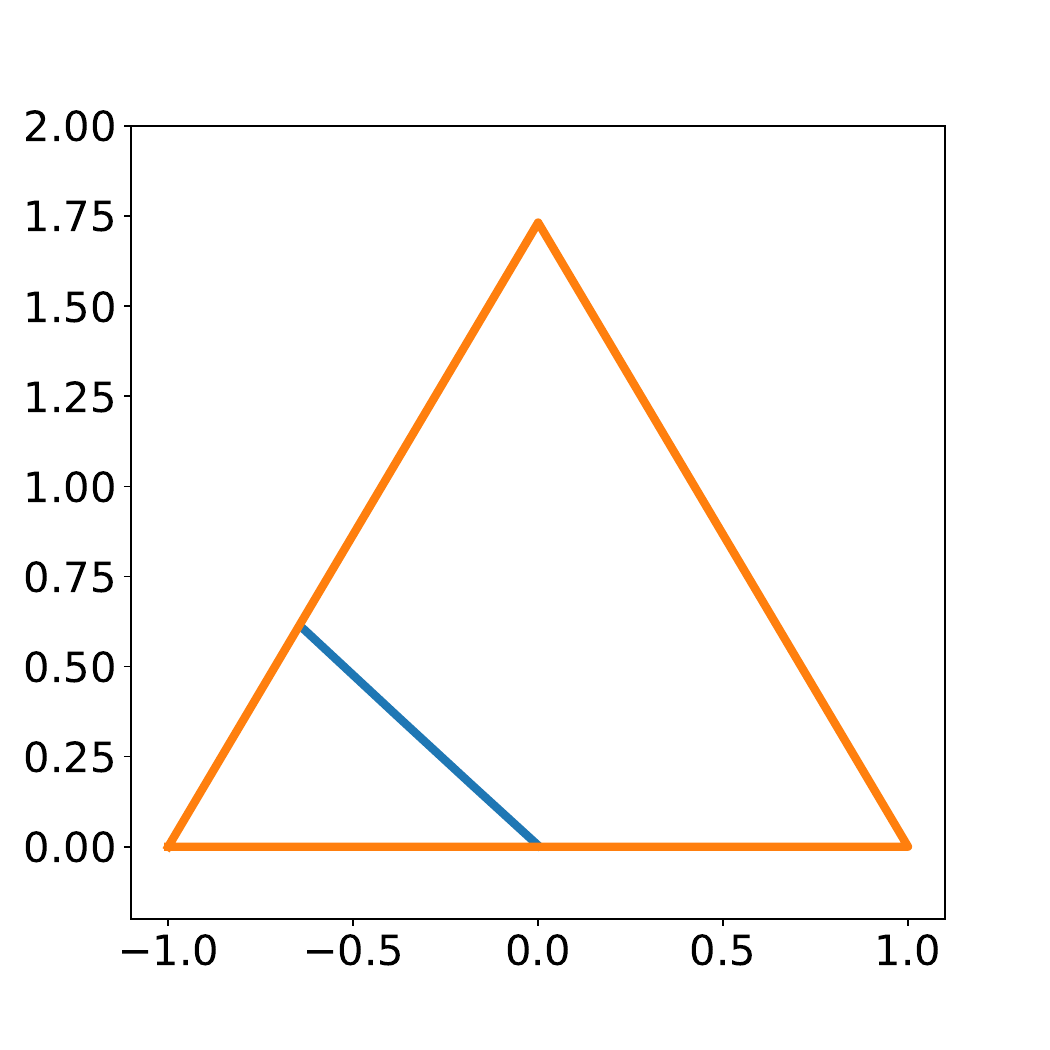}}
%\subcaptionbox{$t \mapsto \sqrt{g_{\gamma(t)} \left( \gamma^\dt(t), \gamma^\dt(t)\right)}$}{\includegraphics[width=0.24\textwidth]{arclengthLinePlaneToLineR3.pdf}}
%
\subcaptionbox{Plane $x = 0$}{\includegraphics[width=0.24\textwidth]{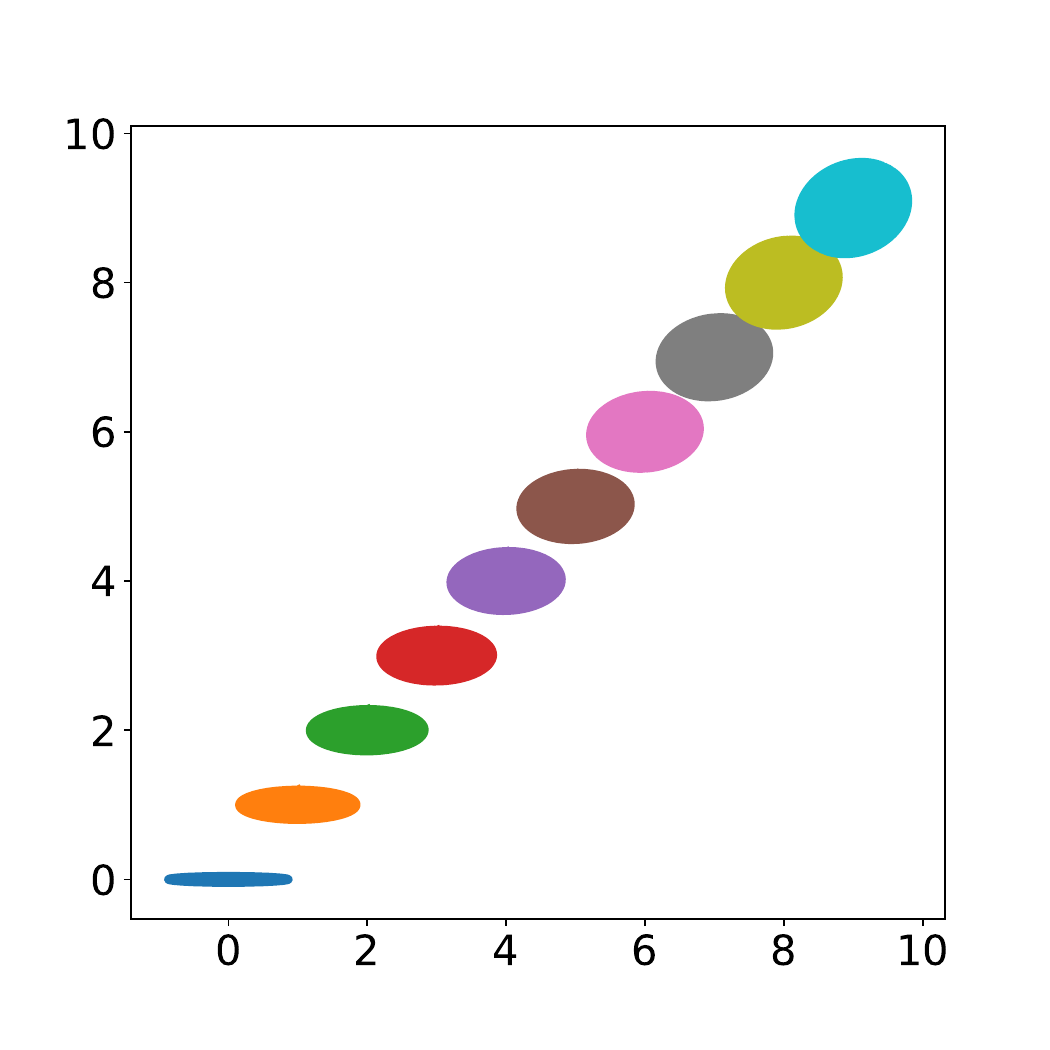}}
\subcaptionbox{Plane $y=0$}{\includegraphics[width=0.24\textwidth]{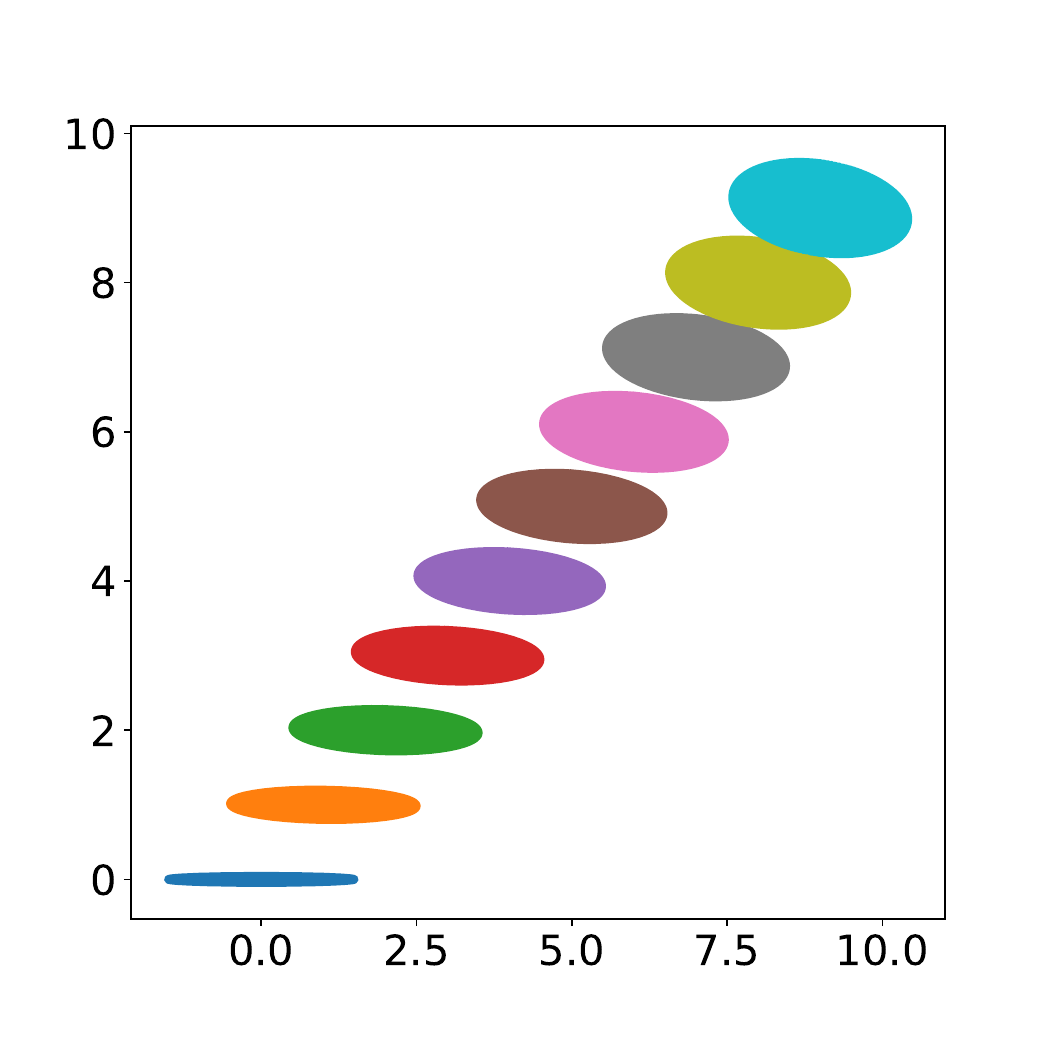}}
\subcaptionbox{Plane $z = 0$}{\includegraphics[width=0.24\textwidth]{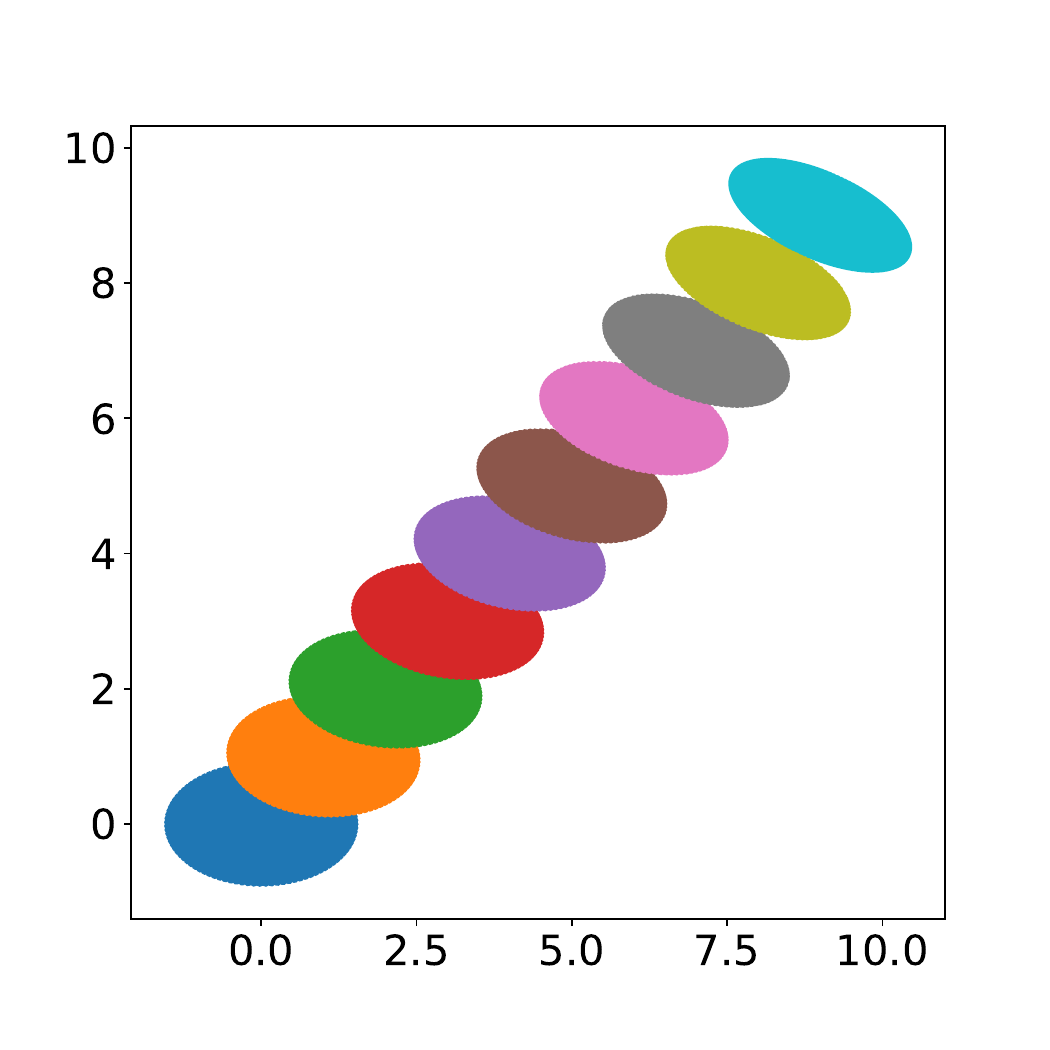}}
\subcaptionbox{Principal angles.}{\includegraphics[width=0.24\textwidth]{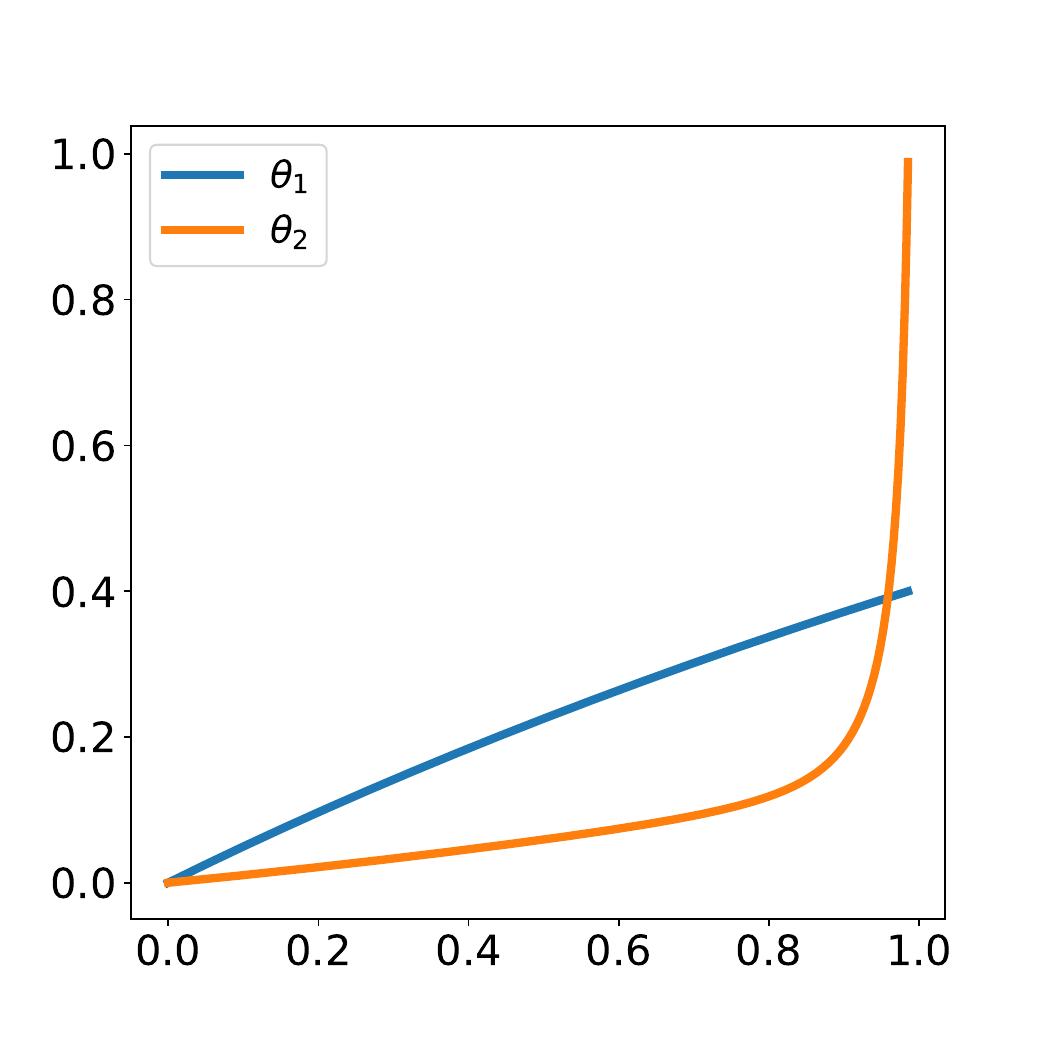}}
\caption{We observe an evolution between a weighted flag close to a line-plane at initial time and close to a line-$\R^3$ at final time. Initial data: $\mu^{(0)} = (0.499, 0.499, 0.002) $, ${\dot{\mu}}^{(0)} = (0.15, -0.5, 0.35)$, $U^{(0)} = I_3$, $(b_{12}, b_{23}, b_{13})^{(0)} = (0.5, 0, 0.1)$. Final point: $\mu^{(N)} = (0.643, 0.008, 0.349)$, $U^{(N)}${\tiny $~=~\left(   
\begin{array}{ccc}
0.921 & 0.118 & 0.371\\
-0.372 & 0.549 & 0.749 \\
-0.115 & -0.827 & 0.55
\end{array} \right) $}.}
\end{figure}
\end{center}
%%%%%%%%%%%%%%%%%%%%%%%%%%%%%%%%%%%%%%%%%%%%%%%%%%%%%%
%
%%%%%%%%%%%%%%%%%%%%%%%%%%%%%%%%%%%%%%%%%%%%%%%%%%%%%
%
\begin{center}
\setcounter{subfigure}{0}
\begin{figure}[!htp] 
\subcaptionbox{$t \mapsto \mu(t)$}{\includegraphics[width=0.24\textwidth]{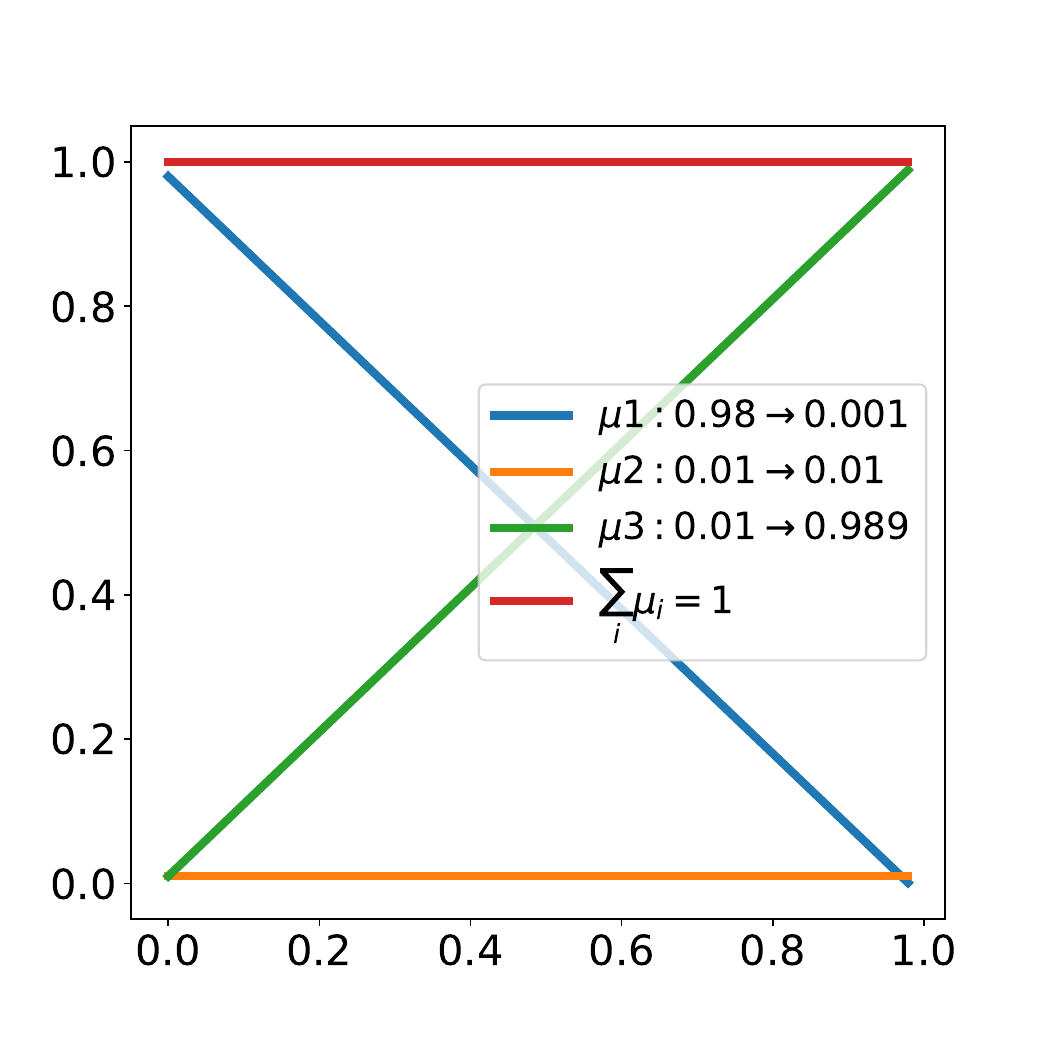}}
\subcaptionbox{$t \mapsto \lambda(t)$}{\includegraphics[width=0.24\textwidth]{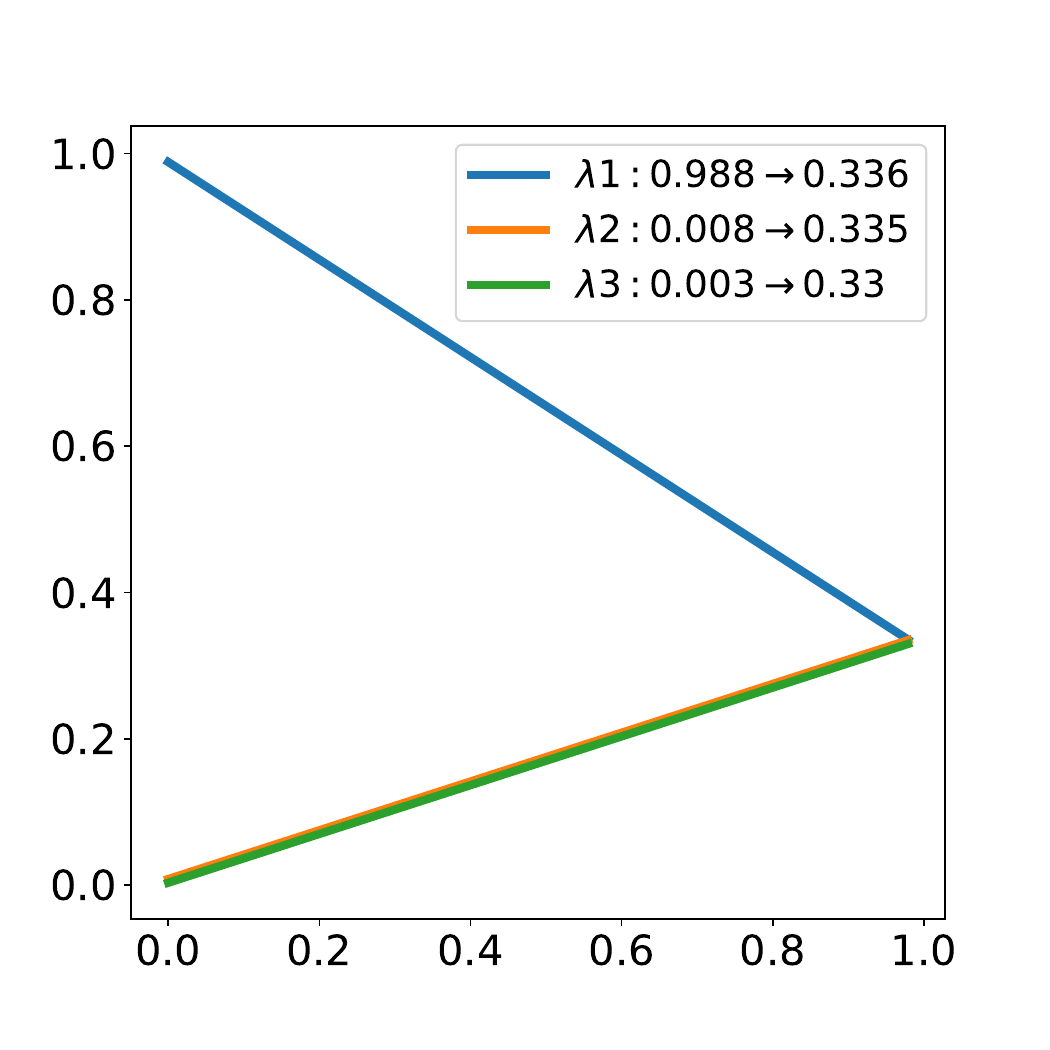}}
\subcaptionbox{$t \mapsto \dt{\mu}(t)$}{\includegraphics[width=0.24\textwidth]{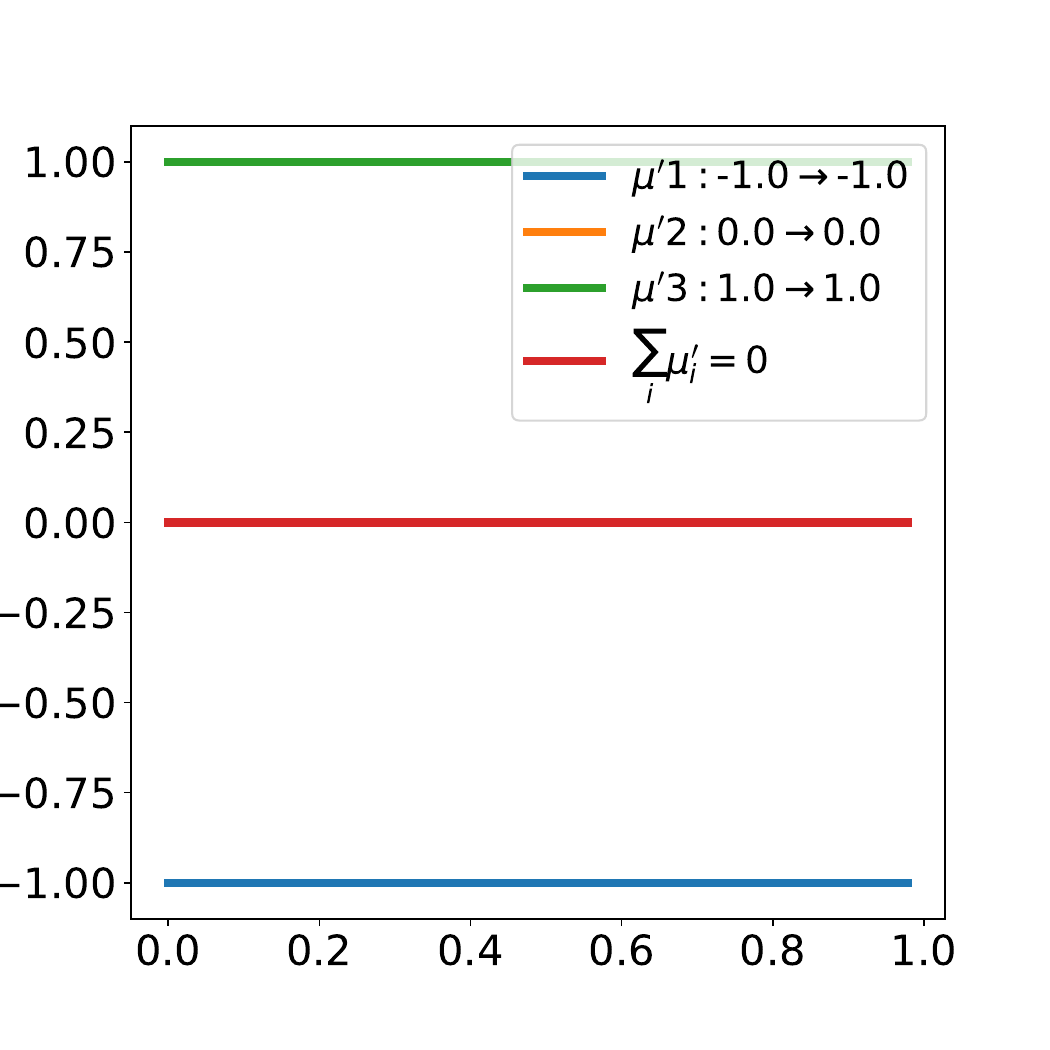}}
\subcaptionbox{Image of $\mu$ in $\Delta(3)$.}{\includegraphics[width=0.24\textwidth]{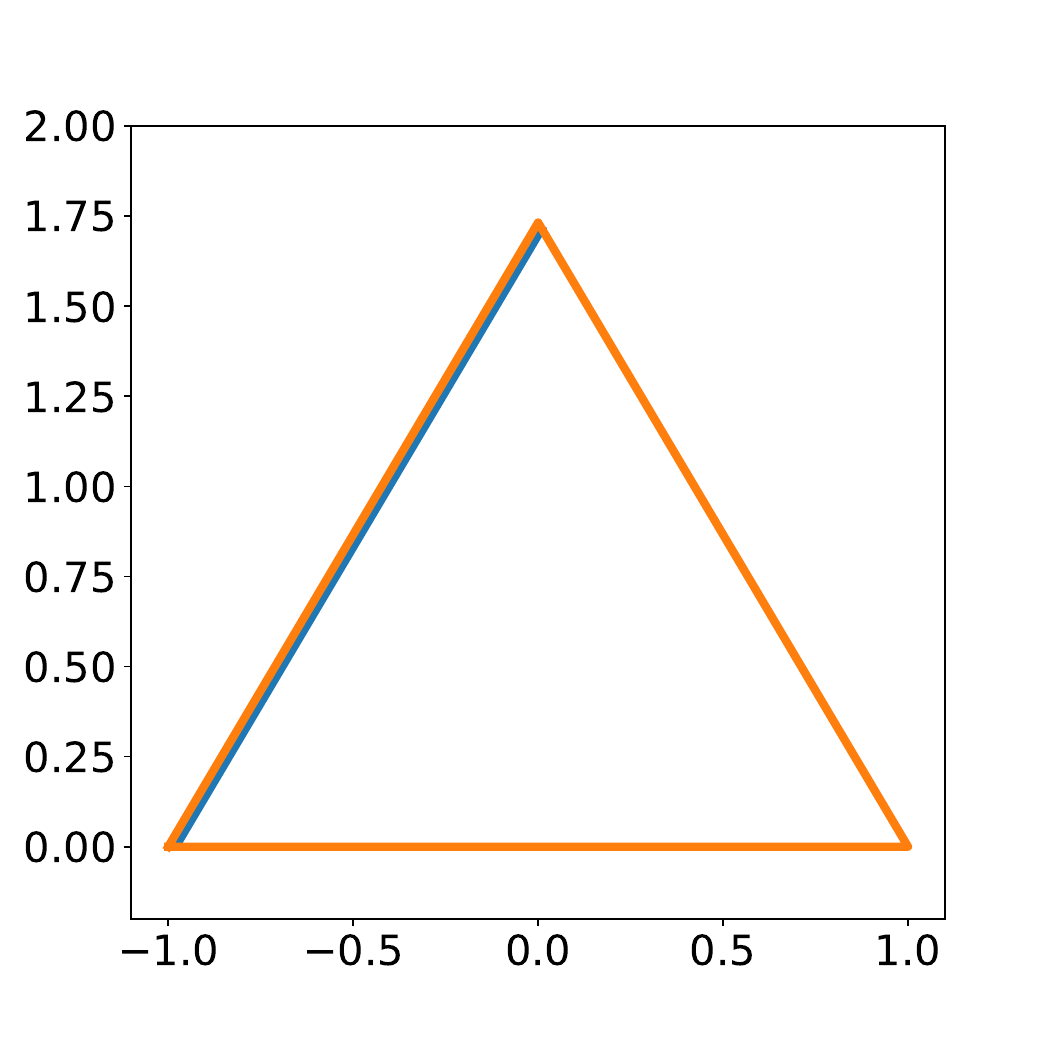}}
%\subcaptionbox{$t \mapsto \sqrt{g_{\gamma(t)} \left( \gamma^\dt(t), \gamma^\dt(t)\right)}$}{\includegraphics[width=0.24\textwidth]{arclengthLineToR3.pdf}}
%
\subcaptionbox{Plane $x = 0$}{\includegraphics[width=0.24\textwidth]{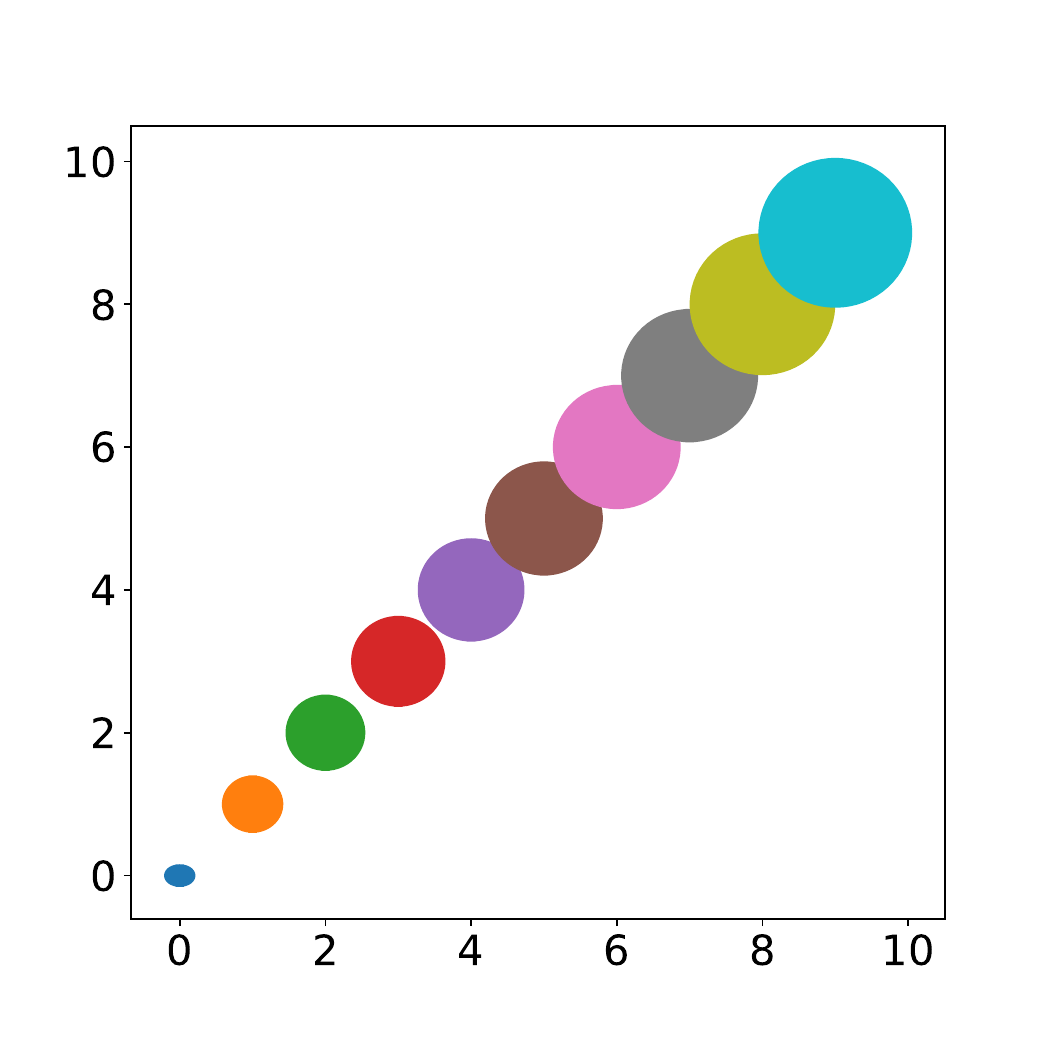}}
\subcaptionbox{Plane $y=0$}{\includegraphics[width=0.24\textwidth]{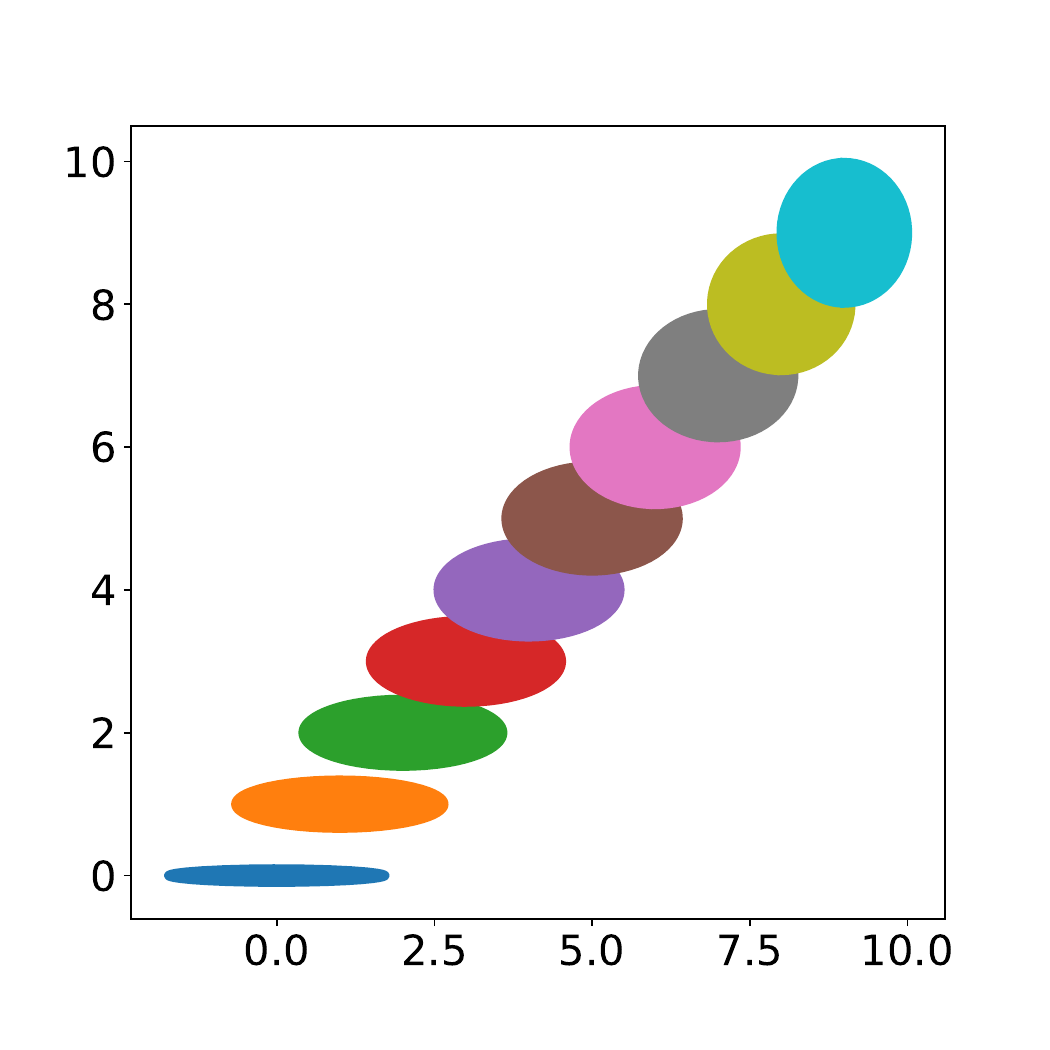}}
\subcaptionbox{Plane $z = 0$}{\includegraphics[width=0.24\textwidth]{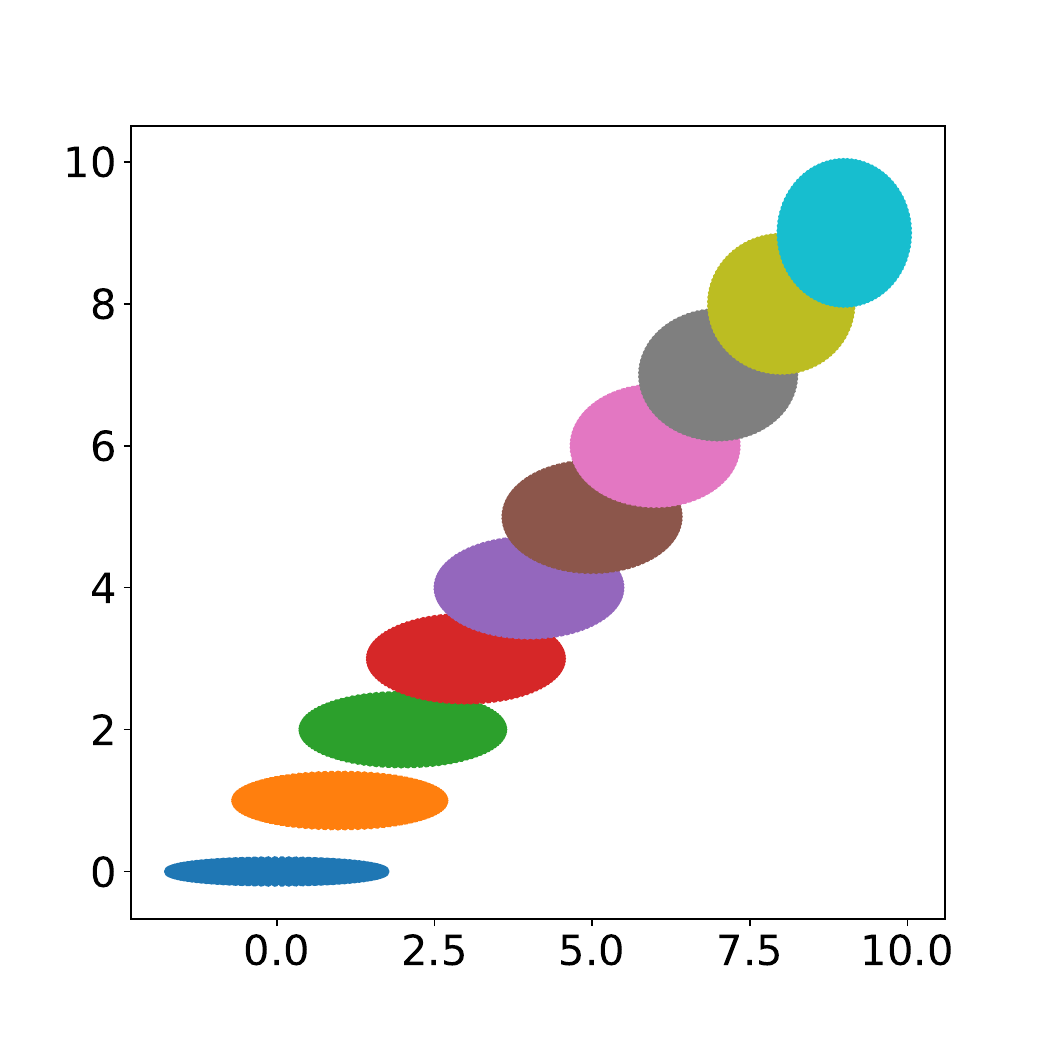}}
\subcaptionbox{Principal angles.}{\includegraphics[width=0.24\textwidth]{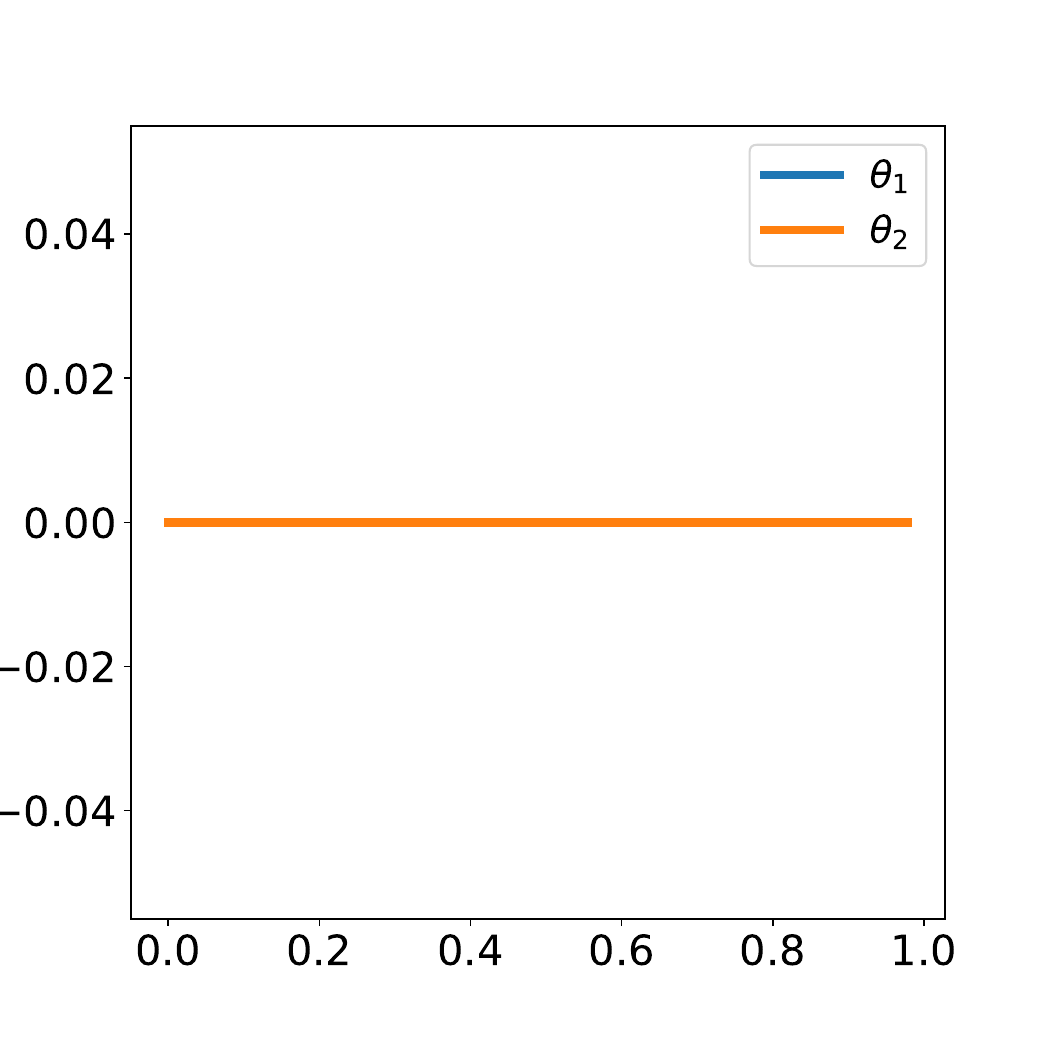}}
\caption{We observe an evolution between a weighted flag close to a line at initial time and close to $\R^3$ at final time. Initial data: $\mu^{(0)} = (0.98, 0.01, 0.01) $, $U^{(0)} = I_3$, ${\dot{\mu}}^{(0)} = (-1, 0, 1)$, $(b_{12}, b_{23}, b_{13})^{(0)} = (0, 0, 0)$. Final point: $\mu^{(N)} = (0.001, 0.01,  0.989)$, $ U^{(N)} = I_3$. \label{figNumGeod-1}}
\end{figure}
\end{center}

Coming back to Example~\ref{ex:euclDist}, we finally draw a comparison with euclidean shortest paths in $\xSym_1^+(n)$ endowed with the usual scalar product $(A,B) \mapsto \tr (A B^T)$. With the initial and final points $(\mu^{(0)}, U^{(0)})$ and $(\mu^{(N)}, U^{(N)})$ obtained in Figure~\ref{figNumGeod1}, we associate $A^{(0)} = U^{(0)} \xdiag(\lambda^{(0)}) {U^{(0)}}^T$, $A^{(N)} = U^{(N)} \xdiag(\lambda^{(N)}) {U^{(N)}}^T \in \xSym_1^+(n)$ and for $p \in {0, 1, \ldots, N}$, we define $A^{(p)} = \left(1 - \frac{t_p}{t_N} \right) A^{(0)} + \frac{t_p}{t_N} A^{(N)}$ and then $(\mu^{(p)}, U^{(p)}) \in \Delta(n) \times O(n)$ (of course not uniquely defined while the projection on $\cW\cF(n)$ is) such that $A^{(p)} = U^{(p)} \xdiag(\lambda^{(p)}) {U^{(p)}}^T$. In Figure~\ref{figEuclGeod1}, we compare the numerical geodesic obtained in Figure~\ref{figNumGeod1} with the euclidean geodesic having the same endpoints. We observe that the euclidean trajectory from a near line to a near plane does not stay at most $2$--dimensional: $\lambda_3$ ($\mu_3$) is not close to $0$ for all times, as it was in Figure~\ref{figNumGeod1}.
\begin{center}
\setcounter{subfigure}{0}
\begin{figure}[!htp] 
\subcaptionbox{$t \mapsto \mu(t)$}{\includegraphics[width=0.26\textwidth]{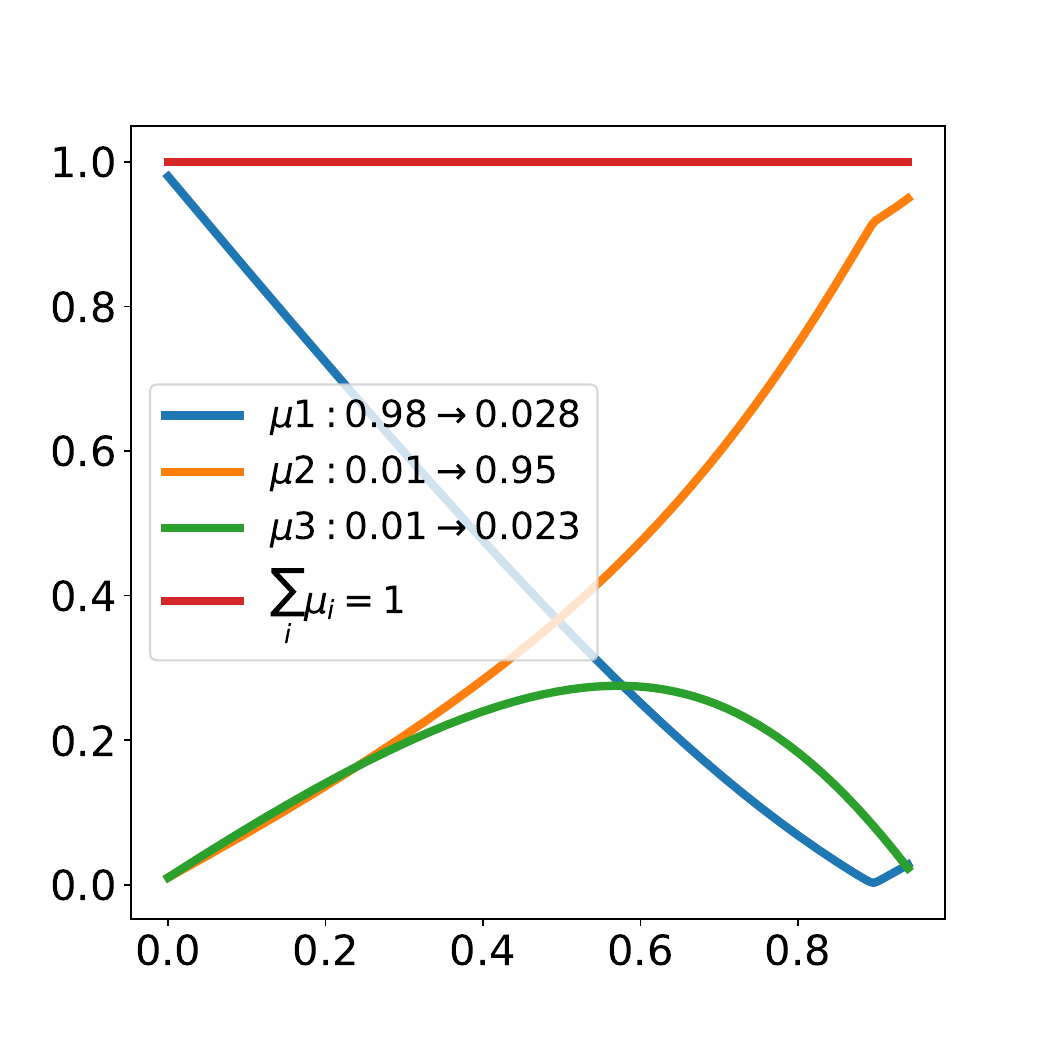}}
\subcaptionbox{$t \mapsto \lambda(t)$}{\includegraphics[width=0.26\textwidth]{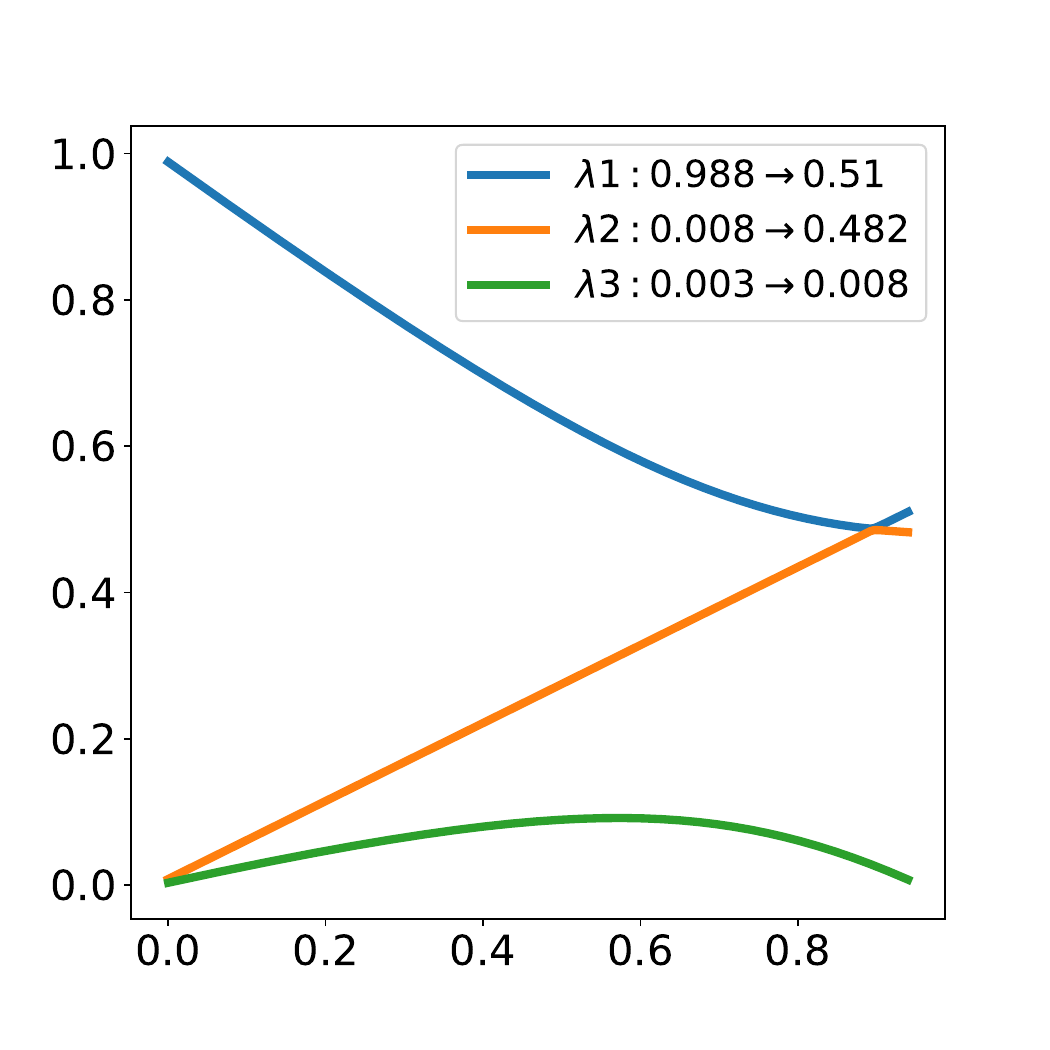}}
%\subcaptionbox{$t \mapsto \mu^\prime(t)$}{\includegraphics[width=0.24\textwidth]{nuprimeLineToPlane.pdf}}
\subcaptionbox{Image of $\mu$ in $\Delta(3)$.}{\includegraphics[width=0.26\textwidth]{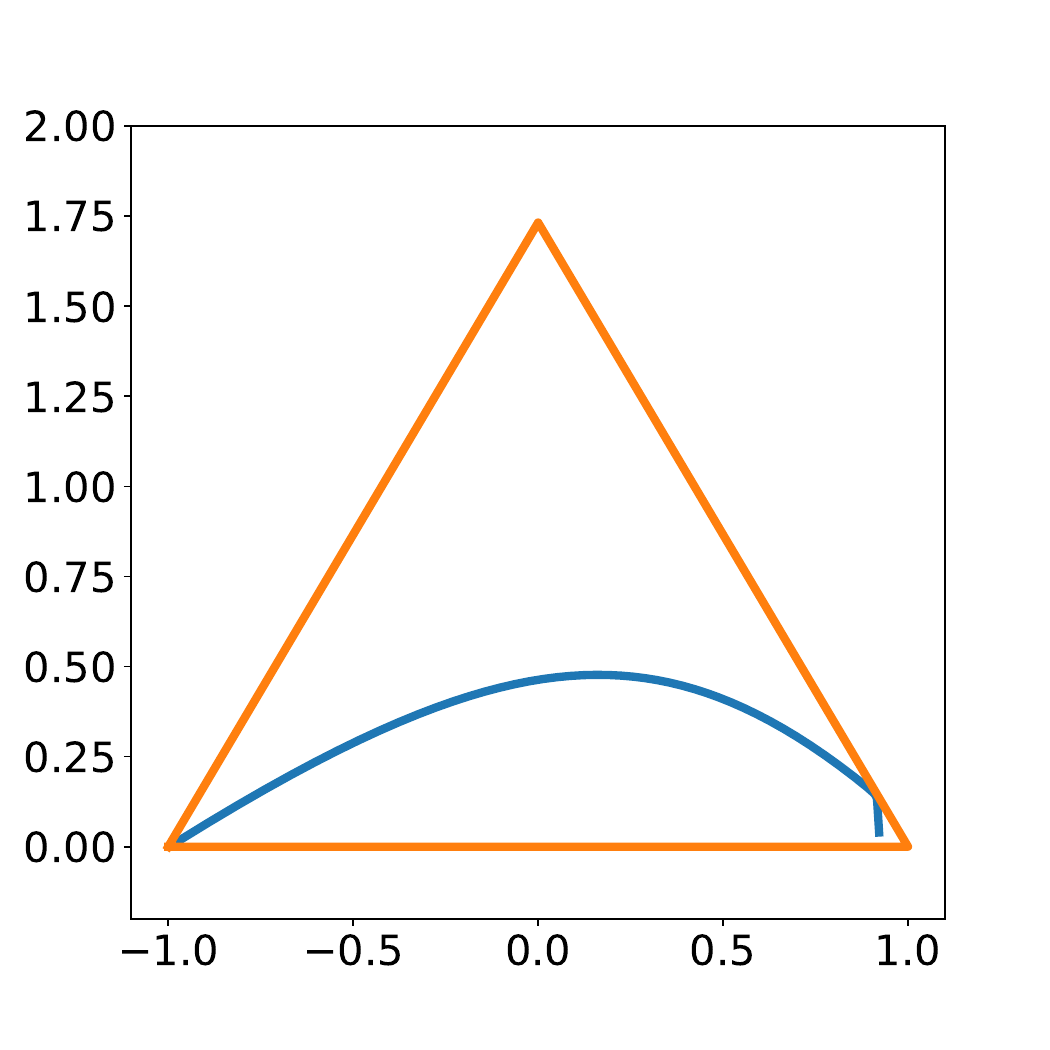}}
%\subcaptionbox{$t \mapsto \sqrt{g_{\gamma(t)} \left( \gamma^\dt(t), \gamma^\dt(t)\right)}$}{\includegraphics[width=0.24\textwidth]{arclengthLineToPlaneEucl.pdf}}
%
\subcaptionbox{Plane $x = 0$}{\includegraphics[width=0.24\textwidth]{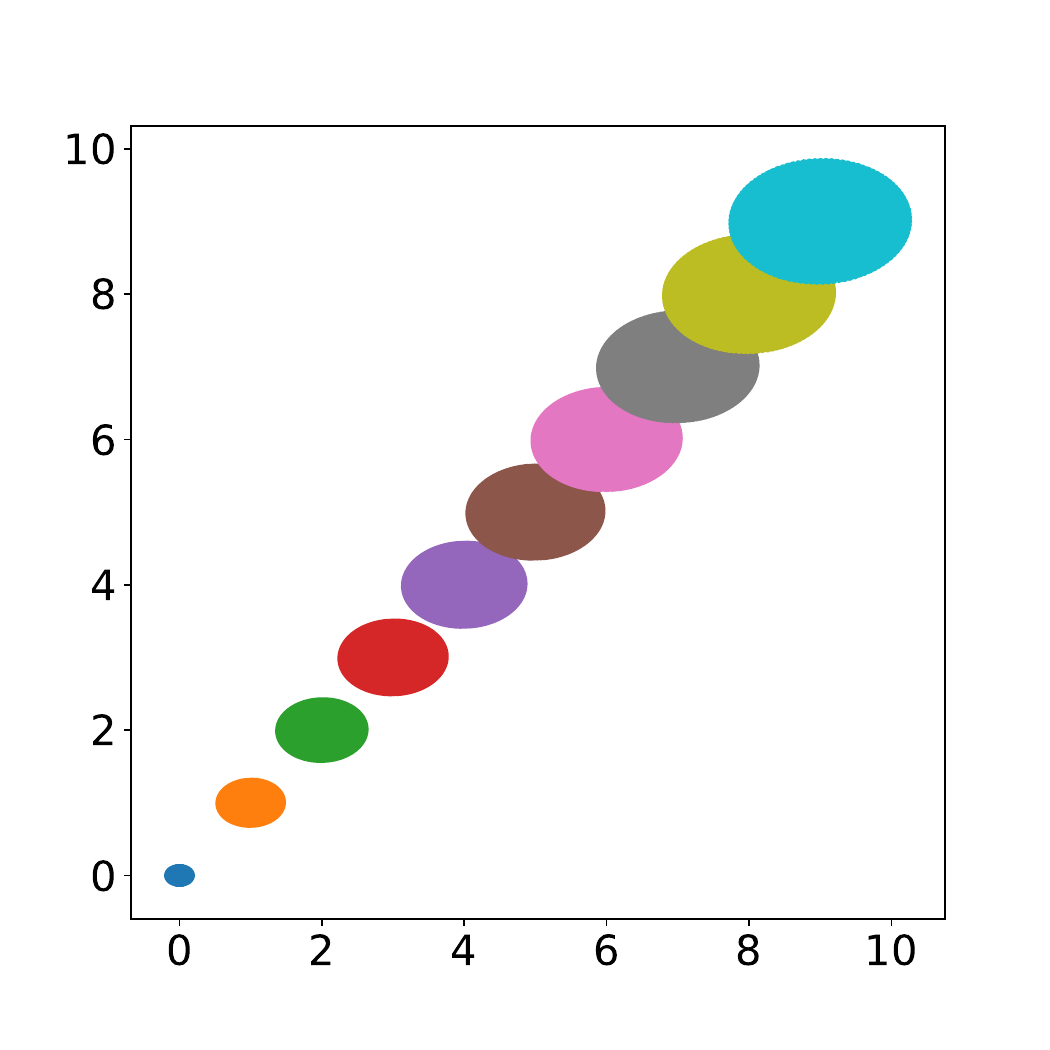}}
\subcaptionbox{Plane $y=0$}{\includegraphics[width=0.24\textwidth]{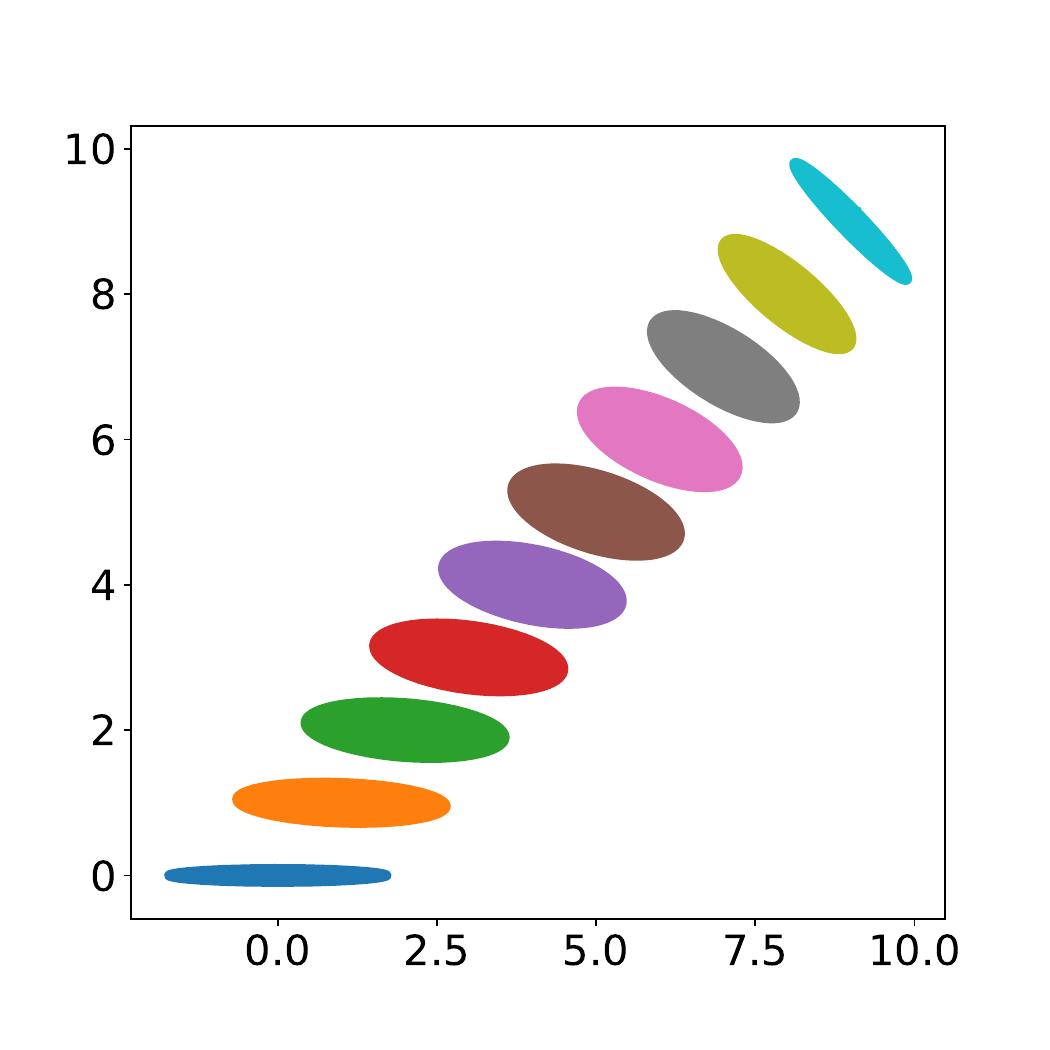}}
\subcaptionbox{Plane $z = 0$}{\includegraphics[width=0.24\textwidth]{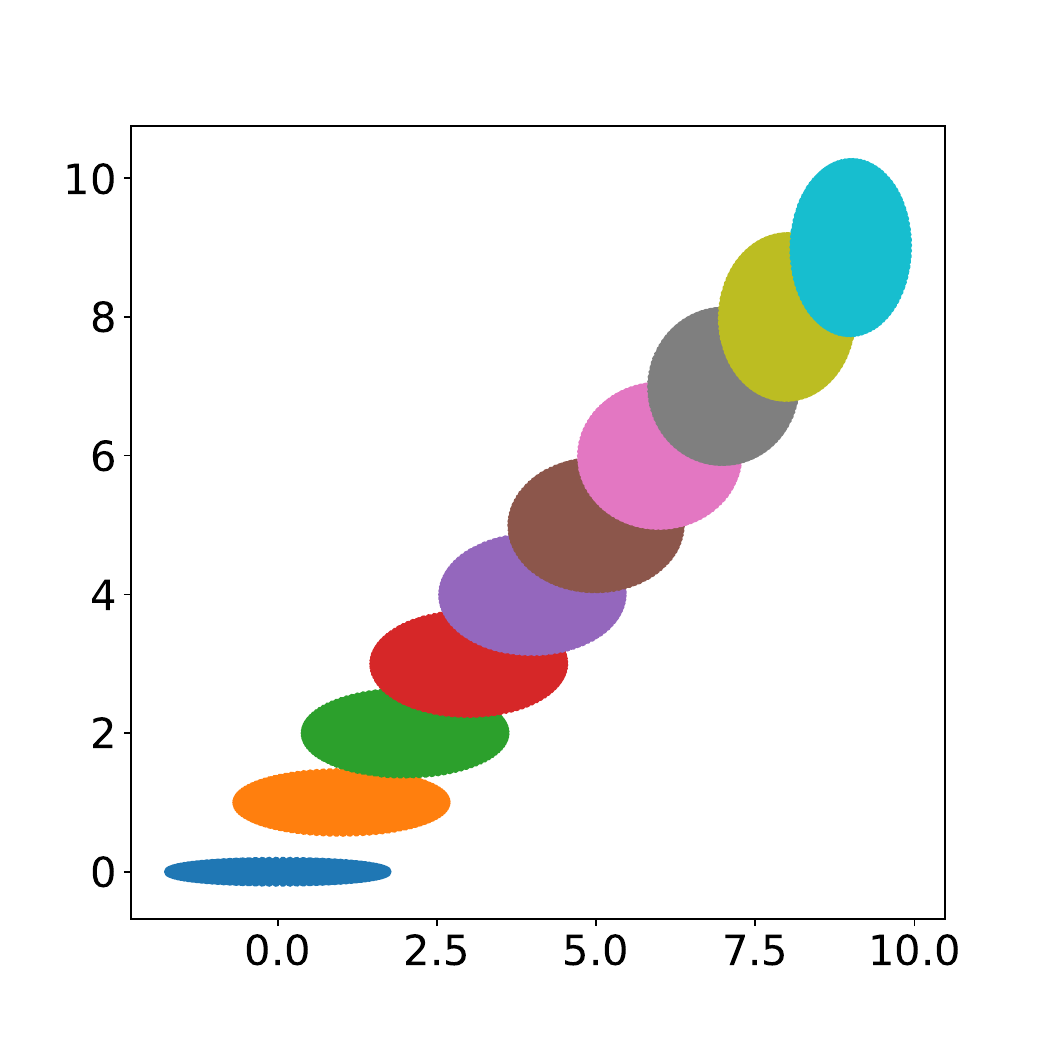}}
\subcaptionbox{Principal angles.}{\includegraphics[width=0.24\textwidth]{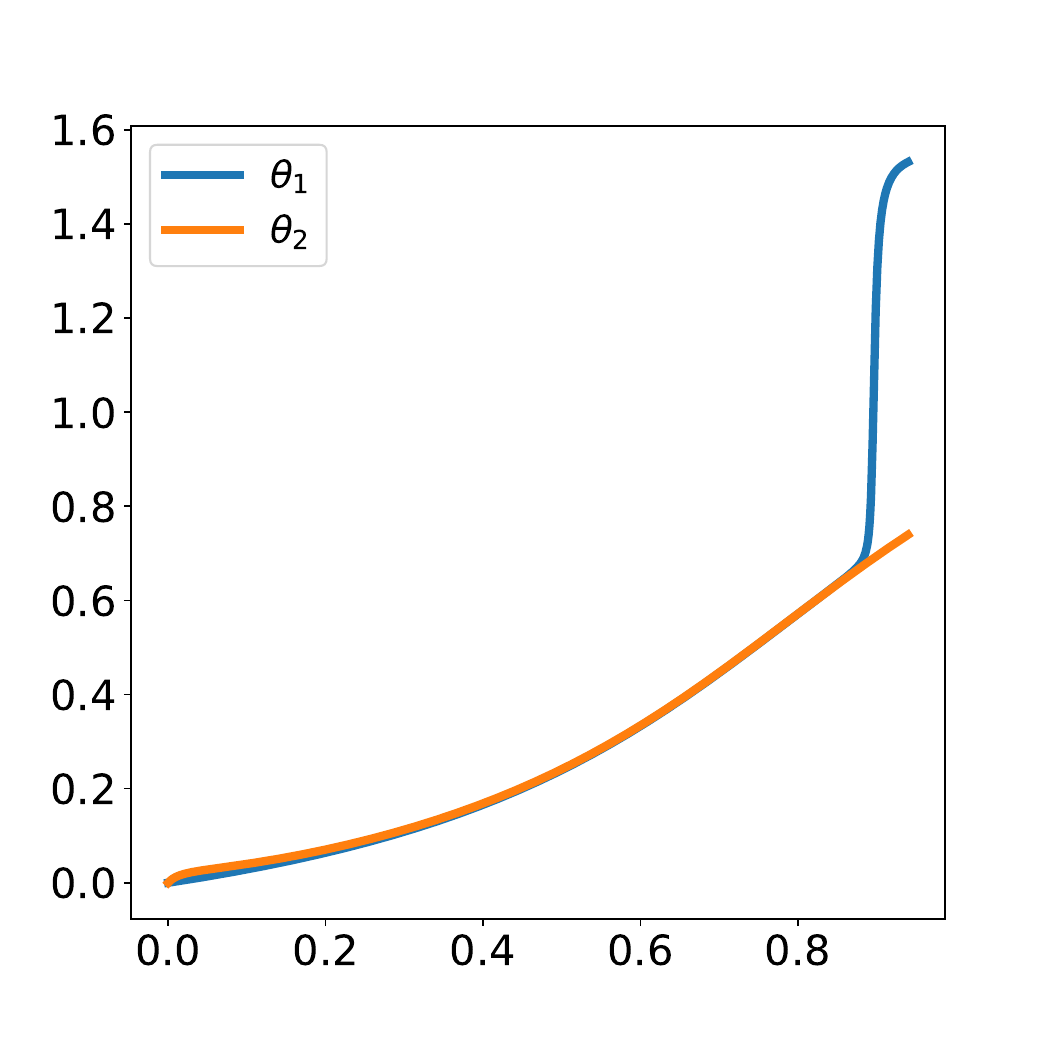}}
\caption{We represent the euclidean geodesic between the same initial and final points as in Figure~\ref{figNumGeod1} that is a weighted flag close to a line at initial time and close to a plane at final time. Initial point: $\mu^{(0)} = (0.98, 0.01, 0.01) $, $U^{(0)} = I_3$ and final point: $\mu^{(N)} = (0.028, 0.95,  0.023)$,
$U^{(N)}${\tiny $~=~\left(
\begin{array}{ccc}
0.039 & 0.738 & 0.674\\
-0.997 & 0.072 & -0.021\\
-0.064 & -0.671 & 0.739                                                                                                                                                                                                                                                        \end{array} \right) $}
. \label{figEuclGeod1}}
\end{figure}
\end{center}

\begin{remk}[Convexity of $\mu_n$.]
For $k = n$, Proposition~\ref{prop:shortestDeltanOn}$(i)$ rewrites
\[
\ddt{\mu}_n = \frac{1}{n}\sum_{l=1}^{n-1} \sum_{1 \leq i \leq l < j \leq n} f(\mu_{i \to j}) \: \partial_l f (\mu_{i \to j}) b_{ij}^2  \: ,
\]
and we can infer that $\mu_n$ is convex if $f$ and its partial derivatives are nonnegative in $[0,1]^n$, which is true for $f(\mu) = |\mu|$ for instance. In particular, $\mu_n$ is maximal at the endpoints of the geodesic, contrary to what is observed for the euclidean geodesic in Figure~\ref{figEuclGeod1}. In particular, if the endpoints are nearly $(n-1)$--dimensional, the whole geodesic remains nearly $(n-1)$--dimensional.
\end{remk}

\newpage
\section{Flagfolds}
\label{section:flagfoldsDef}

In this section, we introduce the notion of \emph{flagfolds} in a very similar way as varifolds are defined. As we explain in Section~\ref{subsec:flagfoldsMotiv}, they model generalized surfaces whose dimension is not a priori fixed while varifolds require an a priori dimension in the definition itself: $d$--varifolds are measures in $\R^n \times \G$ for a given integer $1 \leq d \leq n$. Flagfolds allow to model multi-dimensional shapes and they also give a framework to consider diffuse approximations of thinner shapes as illustrated in Example~\ref{xmpl:diffusedFlagfold}. The definition of flagfolds (Definition~\ref{dfn:flagfolds}) is given in Section~\ref{subsec:flagfoldsDfn} together with some basic notions. Section~\ref{subsec:flagfoldsEmbeddingVarifolds} then focuses on embedding classical $d$--varifolds into flagfolds (see \eqref{eq:varifoldToFlagfold}) as well as projecting flagfolds onto $d$--varifolds (see \eqref{eq:flagfoldsToVarifolds}). Such correspondences are naturally built from the embedding \eqref{eq:Gembedding} of $\G$ into $\cW\cF(n)$. We exhibit in Definition~\ref{dfn:flagfoldsRectifiable} particular flagfolds whose support is contained in the union of $d$--Grassmannians (for $d = 1, \ldots, n$) and such that the restriction to a given $\G$ induces a rectifiable $d$--varifold.

\subsection{Motivations and basic facts about varifolds}
\label{subsec:flagfoldsMotiv}

\emph{Flagfolds} aim at generalizing varifolds to handle objects whose dimension is not fixed, meaning that it can change in space, or at different scales or is not a priori known, which is particularly interesting to model concrete data. 
First of all, let us be more precise about varifolds, we will limit ourselves to a few basic definitions and we refer to \cite{simon} for a complete introduction to varifolds.

\begin{definition}[$d$--Varifold]
\label{dfn:varifold}
A $d$--varifold in an open set $\Omega \subset \R^n$ is a Radon measure in $\Omega \times \G$.
\end{definition}
While such a measure $V$ couples spatial (in $\R^n$) and directional (in $\G$) information, it is possible to push-forward $V$ via the projection ${\rm pr_1} : \Omega \times \G \mapsto \Omega$, $(x,P) \mapsto x$ and only select the spatial part of $V$. 
The resulting Radon measure $\V = {\rm pr_1}_\# V$ is called the \emph{mass} of $V$ and satisfies by definition
\begin{equation}
\label{eq:varifoldMass}
\V (A) = V ( {\rm pr_1}^{-1}(A)) =
V(A \times \G) \: .
\end{equation}
for all Borel sets $A \subset \R^n$.

As one can imagine from Definition~\ref{dfn:varifold}, varifolds can be wild objects (from a geometric perspective at least). We introduce in Definition~\ref{dfn:varifoldRectif} the class of rectifiable $d$--varifolds, based on the notion of $d$--rectifiable sets that is a notion of regularity naturally arising in geometric measure theory. The $d$--dimensional Hausdorff measure in $\R^n$ is denoted by $\cH^d$, and the space of continuous compactly supported function between two metric spaces by $\xC_c (X, Y)$ and $\xC_c(X)$ if $Y = \R$.  
\begin{definition}[Rectifiable sets]
A set $M \subset \R^n$ is said to be countably $d$--rectifiable if there exist a $\cH^d$--negligible set $M_0 \subset \R^n$ and a family of Lipschitz functions $F_j : \R^d \rightarrow \R^n$, $j \in \N$, satisfying
\[
 M \subset M_0 \cup \bigcup_{j \in \N} F_j(\R^d) \: .
\]
If moreover $M$ is $\cH^d$--measurable and for all compact sets $K \subset \R^n$, $\cH^d(M \cap K) < \infty$ (i.e. $\cH^d_{| M}$ is a Radon measure) we say that $M$ is $d$--rectifiable.
\end{definition}

Rectifiability is the regularity notion allowing to extend the notion of tangent in the sense of geometric measure theory:

\begin{proposition}[Approximate tangent space]
\label{prop:approxTgtPlane}
Let $M \subset \R^n$ be a $d$--rectifiable set, then for $\cH^d$--a.e. $x \in M$ there exists a unique $d$--dimensional plane $P_x \in \G$ such that for all $\phi \in \xC_c(\R^n)$,
\begin{equation*}
\frac{1}{\eta^d} \int_{M} \phi \left( \frac{y-x}{\eta} \right) \: d \cH^d(y) \xrightarrow[\eta \to 0_+]{} \int_{P_x} \phi(y) \: d \cH^d(y) \: .
\end{equation*}
The $d$--dimensional plane $P_x$ is called the approximate tangent space to $M$ at $x$, it is denoted $T_x M$ and furthermore satisfies, for all Borel function $\theta : M \rightarrow \R_+$ locally integrable with respect to $\cH^d_{| M}$,
\begin{equation*}
\frac{1}{\eta^d} \int_{M} \phi \left( \frac{y-x}{\eta} \right) \theta(y) \: d \cH^d(y) \xrightarrow[\eta \to 0_+]{} \theta(x) \int_{P_x} \phi(y) \: d \cH^d(y) \: .
\end{equation*}
\end{proposition}
\begin{definition}[Rectifiable $d$--varifold]
\label{dfn:varifoldRectif}
	Given an open set $\Omega \subset \R^n$, let $M \subset \Omega$ be a $d$--rectifiable set and $\theta$ be a non negative function with $\theta > 0$ $\cH^d$--almost everywhere in $M$. A rectifiable $d$--varifold $V= v(M,\theta)$ in $\Omega$ is a non-negative Radon measure on $\Omega \times G_{d,n}$ of the form $V= \theta \mathcal{H}^d_{| M} \otimes \delta_{T_x M}$ i.e.
	\[
	\int_{\Omega \times G_{d,n}} \varphi (x,T) \, dV(x,T) = \int_M \varphi (x, T_x M) \, \theta(x) \, d \mathcal{H}^d (x) \quad \forall \varphi \in \xC_c (\Omega \times G_{d,n})
	\] where $T_x M$ is the approximate tangent space at $x$ which exists $\mathcal{H}^d$--almost everywhere in $M$. The function $\theta$ is called the \emph{multiplicity} of the rectifiable varifold. If additionally $\theta(x)\in \N$ for $\mathcal{H}^d$--almost every $x\in M$, we say that $V$ is an \emph{integral} varifold.
\end{definition}

Let us derive a simple consequence from Definition~\ref{dfn:varifoldRectif}.
Given a Radon measure $\mu$ in $\R^n$ and a nonnegative smooth and even kernel $\omega : \R \rightarrow \R_+$ with support in $(-1,1)$, the following local normalized covariance matrix, for $\eta > 0$:
\begin{equation}
\label{eq:covMeasure}
\Sigma_\eta (\mu , x) = \frac{\displaystyle \int_{\R^n} \omega \left( \frac{|y-x|}{\eta} \right) \frac{y-x}{\eta} \otimes \frac{y-x}{\eta} \: d \mu(y) }{\displaystyle \displaystyle \int_{\R^n} \omega \left( \frac{|y-x|}{\eta} \right) \left| \frac{y-x}{\eta} \right|^2 \: d \mu(y)} \in \xSym_+^1(n) 
\end{equation}
naturally appears when performing local PCA.
Note that for $y \in \R^n$, $y \otimes y$ is the rank $1$ matrix of $(i,j)$--coefficient $y_i y_j$ while for $\mu_1$, $\mu_2$ measures, $\mu_1 \otimes \mu_2$ hereafter refers to the product of measures. In the case of a rectifiable $d$--varifold $V = v(M,\theta)$, we infer from Proposition~\ref{prop:approxTgtPlane}, applied with $\phi(y) = \omega(|y|) y \otimes y$ and $\phi(y) = \omega(|y|) |y|^2$, that for $\V$--a.e. $x$,
\begin{equation*}
\Sigma_\eta (\V , x) \xrightarrow[\eta \to 0]{} \frac{\displaystyle \int_{T_x M} \omega \left( |y| \right) y \otimes y \: d \cH^d (y) }{\displaystyle \displaystyle \int_{T_x M} \omega \left( |y| \right) | y |^2 \: d \cH^d (y)} = \frac{1}{d} \Pi_{T_x M} \: .
\end{equation*}
%
%{\bf ToDo:} add reference or computation details.

As varifolds are measures, they naturally encompass both continuous and discrete shapes. Considering for instance a finite set of points $\{x_i\}_{i = 1 \ldots N}$ in $\R^n$ provided with positive weights $(m_i)_{i = 1 \ldots N}$ and directions $(P_i)_{i = 1 \ldots N}$ in $\G$, one can associate the \emph{point cloud varifold} (\cite{Blanche-rectifiability,BuetLeonardiMasnou})
\begin{equation}
\label{eq:pointCloudVarifold}
 V = \sum_{i = 1}^N m_i \delta_{(x_i,P_i)} \: .
\end{equation}
Nevertheless, such a definition requires to know directions in $\G$ and in particular their dimension $d$ is fixed (a priori known and common to all points). From a practical point of view, data rarely come with a natural notion of direction (e.g. diffusion MRI data), but one often estimates a direction $P_i$ by performing a local PCA in a neighbourhood of $x_i$: the natural directional information first arise as a (renormalized to trace $1$) covariance matrix $S_i \in \xSym_+^1(n)$ and is subsequently truncated to an a priori dimension $d$ so as to define $P_i \in \G$. 
Directly defining the measure
\[
W = \sum_{i = 1}^N m_i \delta_{(x_i, S_i)}
\]
sounds better adapted to practical data and avoids prescribing a fixed dimension. Such a measure $W$ in $\R^n \times \xSym_+^1(n)$ is an example of what we call \emph{flagfolds}.

\emph{Flagfolds} furthermore allow to consider diffuse ($n$--dimensional) approximations of lower $d$--dimensional objects (that converge in the sense of flagfolds). Note that using varifolds to perform such diffuse approximations is possible: the mass of the varifold can be absolutely continuous with respect to the Lebesgue measure in $\R^n$, but the information conveyed by the Grassmannian $\G$ part is $d$--dimensional (by definition) and thus lack coherence with the dimension of the mass measure. 
Let us detail a concrete example.

\begin{xmpl}[Diffused approximations]
\label{xmpl:diffusedFlagfold}
We consider a smooth object and construct successive thinner diffused approximations. For the sake of simplicity, let $D = \xspan(e_3)$ be a straight line in $\R^3$ and $D_\epsilon = \{ x = (x_1, x_2, x_3) \in \R^3 \: : \: x_1^2 + x_2^2 \leq \epsilon^2 \}$ be the cylinder of axis $D$ and radius $\epsilon > 0$. Then $D_\epsilon$ converge to $D$ in the sense of Hausdorff distance or, from a measure point of view, $\mu_\epsilon = ( \pi \epsilon^2)^{-1} \cL^3_{| D_\epsilon}$ converges to $\cH^1_{| D}$.
If we want to incorporate a directional information, let us review some possibilities. We assume that we infer the directional information from the diffused object, for instance computing the covariance matrix $\Sigma_\epsilon(x) = \Sigma_{\eta(\epsilon)}(\mu_\epsilon, x)$ (see \eqref{eq:covMeasure}) at a larger scale $\eta(\epsilon) \gg \epsilon$. A first option is then to truncate the information thanks to some criterion and replace the covariance matrix with an orthogonal projector. In our case, we should obtain at every point the projector $\Pi_D$ onto $D$ and the $1$--varifold $V_\epsilon = \mu_\epsilon \otimes \delta_D$ converges to $V = \cH^1_{| D} \otimes \delta_D$.
%It is also possible to replace $\delta_D$ with some probability measure $\nu_\epsilon$ on $\G$ that takes into account the truncated information
%
With the notion of flagfold, what we propose is to directly consider measures in $\R^n \times \xSym_+^1(n) \simeq \R^n \times \cW\cF(n)$ and define
\[
W_\epsilon = \mu_\epsilon \otimes \delta_{\Sigma_\epsilon} \xrightarrow[\epsilon \to 0]{} \cH^1_{| D} \otimes \delta_{\Pi_D} \: , \quad \text{as measures in } \R^n \times \cW\cF(n) \: ,
\]
the above convergence holding as soon as $\Sigma_\epsilon(x) \xrightarrow[\epsilon \to 0]{\cW\cF(n)}  \Pi_D$ uniformly enough w.r.t. $x \in D_\epsilon$. 
%A natural candidate $\Sigma_\epsilon(x)$ being a local covariance matrix computed at a scale larger than $\epsilon > 0$. 
For instance, taking $\omega = \one_{(-1,1)}$ for the sake of simplicity, easy computations give in a ball of radius $\eta = \eta(\epsilon) \gg \epsilon$ around $0$:
\begin{align*}
\Sigma_\epsilon (0) = \Sigma_\eta (\mu_\epsilon, 0) & = \frac{\displaystyle \int_{B(0,\eta)} y \otimes y \: d\mu_\epsilon(y) }{ \displaystyle \int_{B(0,\eta)} |y|^2 \: d\mu_\epsilon(y)} = \frac{
\left[ \begin{array}{ccc}
\frac{1}{2}\left(\frac{\epsilon}{\eta}\right)^2 & 0 & 0 \\ 
0 & \frac{1}{2}\left(\frac{\epsilon}{\eta}\right)^2 & 0 \\ 
0 & 0 & \frac{2}{3}
\end{array} \right]
+ O \left( \left(\frac{\epsilon}{\eta}\right)^4 \right)
}{\frac{2}{3} + \left(\frac{\epsilon}{\eta}\right)^2 + O \left( \left(\frac{\epsilon}{\eta}\right)^4 \right)} \\
 & = \Pi_D + O \left( \left(\frac{\epsilon}{\eta}\right)^2 \right)
\end{align*}
\end{xmpl}

\subsection{Definition, mass and disintegration}
\label{subsec:flagfoldsDfn}

The definition of flagfolds is very similar to the one of varifolds: we essentially replace $\G$ with $\cW\cF(n)$, and both are compact spaces.

\begin{remk}
\label{remk:WForSym1pos}
We point out that only the topology of $\G$ is involved in the definition and in the basic properties of varifolds. The same is true for flagfolds.
Consequently, one can evenly work with $\xSym_1^+(n)$ or $\cW\cF(n)$ thanks to the homeomorphism of Proposition~\ref{prop:homeo} when dealing with flagfolds. As it will be alternately more convenient to work with $\cW\cF(n)$ or $\xSym_1^+(n)$, we will use both hereafter, with $\cW\cF(n)$ as default choice when stating definitions and properties.
\end{remk}

\begin{definition}[Flagfolds]
\label{dfn:flagfolds}
Let $\Omega \subset \R^n$ be an open set. A \emph{flagfold} in $\Omega$ is a Radon measure in $\Omega \times \cW\cF(n)$.
\end{definition}

\noindent We note that from the definition of $\cW\cF(n)$, flagfolds are not suitable to model $0$--dimensional objects.

As for varifolds (see \eqref{eq:varifoldMass}), it is possible to define the \emph{mass} of a flagfold $V$ as the push-forward of $V$ onto $\Omega$ via the projection ${\rm pr}_1 : \Omega \times \cW\cF(n) \rightarrow \Omega$, $(x,S) \mapsto x$.

\begin{definition}[Mass]
Let $V$ be a flagfold in $\Omega$. The \emph{mass} of $V$ is the Radon measure $\V = {{\rm pr}_1} _\# V$ that is, for any Borel set $B \subset \Omega$,
\[
\V(B) = V(B \times \cW\cF(n)) \: .
\]
\end{definition}
Simple examples of flagfold consist in Dirac masses $\delta_{(x,S)}$ for $(x,S) \in \Omega \times \cW\cF(n)$ and more generally, weighted sums of the form
\begin{equation}
\label{eq:flagfoldPC}
W  = \sum_{i=1}^N m_i \delta_{(x_i,S_i)} \: ,
\end{equation}
with $\{x_i \}_{i=1 \ldots N}$ in $\R^n$, $\left( S_i \right)_{i=1 \ldots N}$ in $\cW\cF(n)$ and $\left( m_i \right)_{i=1 \ldots N}$ in $\R^\ast_+$. We refer to such measures as \emph{point cloud flagfolds}. 
As previously mentioned, if the data set does not come with the directional information, it is possible to infer $(S_i)_i$ thanks to a local PCA around $x_i$, setting for instance 
\[
S_i = \Sigma_\eta(\|W\|,x_i) = \frac{\displaystyle \sum_{i=j}^N m_j \omega \left( \frac{|x_j-x_i|}{\eta} \right) \frac{x_j-x_i}{\eta} \otimes \frac{x_j-x_i}{\eta}}{\displaystyle \sum_{j=1}^N m_j \omega \left( \frac{|x_j-x_i|}{\eta} \right) \left| \frac{x_j-x_i}{\eta} \right|^2} \: ,
\]
as in \eqref{eq:covMeasure} for well-chosen $\eta > 0$. 

Note that for a point cloud flagfold $W$ of the form \eqref{eq:flagfoldPC}, the mass of $W$ satisfies $\displaystyle \| W \| = \sum_{i=1}^N m_i \delta_{x_i}$ and it is possible to split the spatial and directional parts of $W$, at least pointwisely since $m_i  \delta_{(x_i,S_i)} = m_i \delta_{x_i} \otimes \delta_{S_i}$.
Applying a standard disintegration result (see \cite{ambrosio} 2.28-29), it is always possible to decompose a flagfold as a generalized product of the form $\| W \| \otimes w_x$ as described in the following proposition. See \cite[4.1]{MenneSharrer2018} for the case of varifolds.
\begin{proposition}[Young-measure representation]\label{prop:disintegration} Given a $d$--flagfold $W$ on $\Omega$, there exists a family of probability measures $\{w_{x}\}_{x}$ on $\G$ defined for $\|W\|$-almost all $x\in \Omega$,  such that $W = \|W\| \otimes \{w_{x}\}_{x}$, that is,
	\[
	W(\phi) = \int_{x\in\Omega}\int_{S\in\G}\phi(x,S)\, d w_{x}(S)\, d\|W\|(x)
	\]
	for all $\phi\in \xC_{c}(\Omega\times\G)$.
\end{proposition}
\noindent We give the counterpart of rectifiability for flagfolds in Definition~\ref{dfn:flagfoldsRectifiable}.

\subsection{Dimension of a weighted flag and embedding of varifolds into flagfolds}
\label{subsec:flagfoldsEmbeddingVarifolds}

Notice that with the approach we propose, no a priori on the dimension of the shape is been used. However, even after renormalizing the covariance matrix to trace $1$, one can wonder whether there is some information left about the dimension of the object. This is sound since we manage to embed all $d$--Grassmannians $\G$, $d = 1, \ldots, n$ into $\cW\cF(n)$. 

\subsubsection*{Dimension of a weighted flag.} Let us define such a dimension $\overline{d} : \cW\cF(n) \rightarrow [1 , + \infty[$ function on weighted flags, additionally requiring that it is invariant with respect to the flag variable, that is $\overline{d}(\mu,W) = \overline{d}(\mu)$ for all $(\mu, W) \in \cW\cF(n)$.
By construction, there is a natural embedding of $\G$ into $\cW\cF(n)$ through the identification 
\begin{equation} \label{eq:Gembedding}
\G \simeq \{ (\mu, W) \in \cW\cF(n) \: : \: \mu_d= 1 \}
\end{equation}
In other words the $d$--Grassmannian $\G$ is sent on weighted flags whose weight $\mu$ has exactly one non-zero coefficient $\mu_d = 1$, corresponding to $d$ largest eigenvalues that are equal to $\frac{1}{d}$ and the $n-d$ other are zero (see bijection \eqref{eqMuk}), $\G$ hence lies at the corners of the simplex $\Delta(n)$ and the dimension is well-defined as $d$ if $\mu_d = 1$: $\overline{d}((0, \ldots, 0, \underbrace{1}_{d^{th}}, 0, \ldots, 0)) = d$. It is then possible to write any $\mu \in \Delta(n)$ as the simple convex sum
\[
\mu = (\mu_1, \ldots, \mu_n) = \sum_{k=1}^n \mu_k \, (0, \ldots, 0, \underbrace{1}_{k^{th}}, 0, \ldots, 0)
\]
and extend the dimension into a linear function in $\Delta(n)$ by
\begin{equation} \label{eq:dimensionMu}
\overline{d} (\mu,W) = \overline{d} (\mu) := \sum_{k=1}^n \mu_k \overline{d} (0, \ldots, 0, 1, 0, \ldots, 0) = \sum_{k=1}^n k \mu_k \: .
\end{equation}
%As we are averaging all the dimensions of the subspaces of the flag, this weighted mean $\overline{d}(\mu)$ can be associated with the deviation $\displaystyle \sum_{k=1}^n \mu_k (k - \overline{d}(\mu) )^2$.

\subsubsection*{Embedding of Grassmannians into weighted flags.}
We identify hereafter $P \in \G$ with the orthogonal projector $\Pi_P \in \xSym_+(n)$.
Similarly, it is more convenient to use the representation $\xSym_+^1(n)$ for weighted flags (see Proposition~\ref{prop:homeo} giving the homeomorphism $\cW\cF(n) \simeq \xSym_+^1 (n)$). We hence adopt the following notations: given $S \in \xSym_+(n)$,
\begin{equation}
\label{eq:notationEigen}
\left\lbrace
\begin{array}{l}
\lambda_1(S) \geq \ldots \geq \lambda_n(S) \geq  0 \text{ are the eigenvalues of } S \,(\lambda_{n+1}(S) = 0), \\
\mu_1(S), \ldots, \mu_n(S) \text{ are then defined through \eqref{eq:lambdaToMu},}\\
F_k(S) \text{ is the eigenspace corresponding to } \lambda_k(S) , \\
E_k(S) = F_1(S) + \ldots + F_k(S) \: .
 \end{array} 
 \right.
\end{equation}
With such notations, we have $\displaystyle S = \sum_{k=1}^n (\lambda_k(S) - \lambda_{k+1}(S)) \Pi_{E_k(S)} = \sum_{k=1}^n  \frac{\mu_k(S)}{k} \Pi_{E_k(S)}$.

\noindent It is then very natural to embed $\G$ into weighted flags:
\begin{equation}
\label{eq:GtoWF}
\begin{array}{cccl}
i : & \G  & \rightarrow & \xSym_+^1 (n) \\
& P  & \mapsto &  \frac{1}{d}\Pi_{P} \
\end{array}
\end{equation}
Conversely,
it is very natural to associate a $d$--plane with a weighted flag $S$ satisfying $\mu_d(S) = 1$:
\[
\begin{array}{ccl}
\{ S \in \xSym_+^1(n) \: : \: \mu_d(S) = 1\}  & \rightarrow & \G  \\
S  & \mapsto &  E_d(S)
\end{array}
\]
and as we did for the dimension $\overline{d}$, it is possible to extend the previous mapping, linearly with respect to $\mu$, into the following left inverse to the embedding $i$:
\begin{equation}
\label{eq:WFtoG}
\begin{array}{lcl}
 \xSym_+^1(n)  & \rightarrow & \xSym_+(n) \\
 \displaystyle S = \sum_{k=1}^n \frac{\mu_k(S)}{k} \Pi_{E_k(S)} & \mapsto & \overline{S} := \displaystyle \sum_{k=1}^n \mu_k(S) \Pi_{E_k(S)}
\end{array}
\end{equation}
For all $P \in \G$, $\overline{i(P)} = \Pi_P$.
Additionally note that for all $S \in \xSym_1^+(n)$, $\overline{d} (S) = \tr \overline{S} \in [1,n]$ and $\overline{S}$ has the same eigenspaces as $S$ with eigenvalues $\lambda_i(\overline{S}) = \sum_{k=i}^n \mu_k(S)$ while we recall that $\lambda_i(S) = \sum_{k=i}^n \frac{1}{k} \mu_k(S)$.

\subsubsection*{Embedding of varifolds into flagfolds}
Embedding \eqref{eq:GtoWF} and retraction \eqref{eq:WFtoG} give rise to a similar correspondence  between $d$--varifolds and flagfolds as we describe hereafter.

\noindent Given a $d$--varifold $V_d$, it is possible to define the associated flagfold $\widehat{V_d} = ({\rm Id}, i)_\# V_d$: for $\psi \in \xC_c (\Omega \times \cW\cF(n))$,
\begin{equation} \label{eq:varifoldToFlagfold}
\int_{\Omega \times \cW\cF(n)} \psi (x,S) \: d \widehat{V_d} (x,S) = \int_{\Omega \times \G} \psi \left( x, \frac{1}{d} \Pi_P \right) \: dV_d (x,P) \: .
\end{equation}
Notice that by definition $\| \widehat{V_d} \| = \| V_d\|$.
Conversely, given a flagfold $W$, it is possible to define a $d$--varifold $V_d$, for $d \in \{1, \ldots, n\}$ as follows, for $\psi \in \xC_c(\Omega \times \G)$,
\begin{equation} \label{eq:flagfoldsToVarifolds}
\int_{\Omega \times \G} \psi(x,P) \: dV_d(x,P) = \int_{\Omega \times \cW\cF(n)} \mu_d(S) \psi(x, E_d(S)) \: d W(x,S) \: .
\end{equation} 
The previous mappings \eqref{eq:varifoldToFlagfold} and \eqref{eq:flagfoldsToVarifolds} define a one-to-one correspondence between $d$--varifolds and flagfolds with support in $\G \subset \cW\cF(n)$. More generally, the mapping
\[
\begin{array}{ccl}
\left\lbrace \left( V_d \right)_{d=1 \ldots n} \: : \: V_d \text{ is a } d\text{--varifold} \right\rbrace & \rightarrow & \{ \text{flagfolds in } \Omega \} \\
(V_1 , \ldots, V_n) & \mapsto & \sum_{d=1}^n \widehat{V_d} \: .
\end{array}
\]
is one-to-one, and onto $
\left\lbrace W \text{ flagfold} \: : \: {\rm spt} \: W \subset \Omega \times \cup_{d=1}^n \G \right\rbrace$. Note that for such flagfolds supported in $\Omega \times \cup_{d=1}^n \G$, $\widehat{V_d} = W_{| \Omega \times \G}$ and $W = \sum_{d=1}^n \widehat{V_d}$. We point out that $W = \sum_{d=1}^n \widehat{V_d}$ does not hold in general, see Example~\ref{xmpl:flagfold1}, nevertheless the following identity is true for any flagfold
\begin{equation*}
     \| W \| = \sum_{d=1}^n \| \widehat{V_d} \| = \sum_{d=1}^n \|V_d\| \: .
\end{equation*}
Indeed, as $\sum_{d=1}^n \mu_d(S) = 1$, let $\phi \in \xC_c (\Omega)$,
\begin{align*}
\| W \| (\phi) & = \int_{\Omega \times \cW\cF(n)} \phi (x) \: dW(x,S) = \sum_{k=1}^n \int_{\Omega \times \cW\cF(n)} \phi(x) \mu_k(S)  \: d W(x,S) \\
& = \sum_{k=1}^n \int_{\Omega \times {\rm G}_{k,n}} \phi(x) \: dV_k(x,S) = \sum_{k=1}^n  \| V_k \|(\phi) \: .
\end{align*}

\subsubsection*{Rectifiable flagfolds}
There is a special case of flagfold whose support is included in $\Omega \times \cup_{d=1}^n \G$ and we additionally require that each $d$--varifold $V_d$ is a rectifiable $d$--varifold. We call such flagfolds \emph{rectifiable}. They exactly model unions of regular (in the sense of rectifiability) sets of different dimensions.

\begin{definition}[Rectifiable flagfolds]
\label{dfn:flagfoldsRectifiable}
We say that a flagfold $W$ in $\Omega$ is a \emph{rectifiable flagfold} if there exists a family of rectifiable $d$--varifolds $V_d = (M_d, \theta_d)$ for $d = 1, \ldots, n$ (see Definition~\ref{dfn:varifoldRectif}), such that
\[
W = \sum_{d=1}^n \widehat{V_d} \quad %\text{with} \quad V_d = v(M_d,\theta_d) = \theta_d \cH^d_{| M_d} \otimes \delta_{T_x M_d} 
\: .
\]
\end{definition}

\noindent Equivalently, $W$ is a rectifiable flagfold if and only if $W$ is supported in $\Omega \times \cup_{d=1}^n \G$ and for $d=1, \ldots,n$, $V_d$ defined as in \eqref{eq:flagfoldsToVarifolds} is a rectifiable $d$--varifold.

\section{First variation of a flagfold}
\label{section:flagfoldFirstVariation}

This section is dedicated to extend the so called \emph{first variation} from varifolds to flagfolds and investigate analogous regularity properties based on Allard's rectifiability theorem. More precisely, Section~\ref{subsec:pushforward} first focuses on the notion of push-forward of a flagfold and provides Definition~\ref{dfn:pushforwardFlagfold} that is consistent with the push-forward of a varifold. The first variation $\delta W(X)$ of a flagfold $W$ is subsequently defined in Definition~\ref{dfn:flagfoldFirstVariation} as the first variation of the total mass of the flagfold under the diffeomorphism flow induced by $X \in \xC_c^1(\Omega,\R^n)$. By such a construction, it naturally extends the notion of first variation of a varifold introduced by Allard in \cite{Allard72}.
In this same paper \cite{Allard72}, it is proved that for a $d$--varifold $V$ with positive lower $d$--dimensional density, assuming that the first variation $\delta V$ is a Radon measure (while $\delta V$ is generally speaking a distribution of order $1$) implies that $V$ is a rectifiable $d$--varifold.
Section~\ref{subsec:monotonicity} explores the implications of a similar control of the first variation $\delta W$ of a flagfold $W$. Following the same path as for varifolds, Proposition~\ref{prop:monotonicity} adapts a usual monotonicity formula (see for instance \cite[§3,chap.8]{simon}) though now taking into account the varying dimension $\overline{d}$: it implies the monotonicity of $r \mapsto r^{-d} \| W \| (B(x,r))$ for all $0 < d < \overline{d}(x)$ a.e. assuming that $\delta W$ is a Radon measure. Theorem~\ref{thm:rectifiability} eventually states that a flagfold $W$ with positive lower densities and locally bounded first variation has integer dimension a.e. and consequently identifies with a sum of classical varifolds $W = \widehat{V_1} + \ldots + \widehat{V_n}$. The rectifiability of the flagfold is proved though only under the additional assumption that the first variations of $V_1, \ldots, V_n$ are locally bounded,  leaving the issue open under the more general assumption that $W$ has locally bounded first variation.

\subsection{Push-forward of a flagfold and variation of the mass}
\label{subsec:pushforward}
%Let $X \in \xC_c^1(\Omega)$ and $\phi_t(x) = x + t X(x)$ a diffeomorphism of $\Omega$ for $t$ small enough. 
First of all, let us recall the definition of push-forward $\phi^\# V_d$ of a $d$--varifold $V$ by a diffeomorphism $\Phi$ of $\Omega$:
\[
\int_{\Omega \times \G} \psi  \: d ( \Phi^\# V )  = \int_{\Omega \times \G} \psi \left( \Phi(x), D\Phi(x) \cdot P \right) J\Phi(x,P) \: dV(x,P), \quad \text{for } \psi \in \xC_c(\Omega \times \G)
\]
where $D\Phi(x) \in \xM_n(\R)$ is the Jacobian of $\Phi$ at $x$ and, if $(\tau_1, \ldots, \tau_d)$ is an orthonormal basis of $P$, we have $D\Phi(x) \cdot P = \xspan (D\Phi(x) \tau_1, \ldots , D\Phi(x) \tau_d ) \in \G$ and
\begin{equation}
\label{eq:jacobian}
J \Phi(x,P) := \sqrt{\det \left( D\Phi(x)_{| P}^T D\Phi(x)_{| P} \right)} = \sqrt{\det M_x} \quad \text{with} \quad (M_x)_{ij} = D\Phi(x)\tau_i \cdot D\Phi(x)\tau_j
\end{equation}
is the $d$--Jacobian of $D\Phi(x)_{| P}$ (see \cite{simon}). We point out that this is not the standard notion of push-forward for the Radon measure $V$ since we have
\[
\Phi^\# V = J \Phi(x,P) \, \left[ (x,P) \mapsto (\Phi(x), D \Phi(x) \cdot P) \right]_\# V \: .
\]
The push-forward of varifolds is defined to be coherent with the area formula holding for $d$--rectifiable sets (see \cite[11.2]{Maggi}).
We now propose to extend the notion of push-forward to flagfolds, consistently with the embedding of varifolds into flagfolds that is, if $V$ is a $d$--varifold and $\widehat{V}$ is the flagfold associated through \eqref{eq:varifoldToFlagfold}, we require that both push-forward coincide: $\Phi^\# \widehat{V} = \widehat{\Phi^\# V}$. As a result, if $W = \sum_{d=1}^n \widehat{V_d}$ is a rectifiable flagfold in the sense of Definition~\ref{dfn:flagfoldsRectifiable} then the push-forward of $W$ is obtained by pushing-forward each rectifiable $d$--varifold $V_d = v(M_d, \theta_d)$ and summing the contributions over $d$. Besides such natural consistency, we perform the extension by linearity with respect to the $d$--dimensional slices of the flagfold as stated in \eqref{eq:pushforwardSum}. 

\begin{definition}[Push-forward of flagfolds]
\label{dfn:pushforwardFlagfold}
Let $W$ be a flagfold in $\Omega$ and let $\Phi$ be a diffeomorphism of $\Omega$.
Using the notations introduced in \eqref{eq:notationEigen} for the eigen decomposition of $S \in \xSym_+^1(n)$,
we define the flagfold $\Phi^\# W$ by
\begin{multline*}
\int_{\Omega \times \cW\cF(n)} \psi  \: d ( \Phi^\# W )  = \sum_{k=1}^n \int_{\Omega \times \cW\cF(n)} \mu_k(S) \psi \left( \Phi(x), i(D\Phi(x) \cdot E_k(S))\right) J\Phi(x, E_k(S)) \: dW(x,S), \\
\text{ for } \psi \in \xC_c(\Omega \times \cW\cF(n)) \: ,
\end{multline*}
where we recall that for $P \in \G$, $i(P) = \frac{1}{d} \Pi_P$.
\end{definition}

Let $V$ be a $d$--varifold in $\Omega$ and $\widehat{V} = ({\rm id}, i)_\# V$ be the flagfold associated through \eqref{eq:varifoldToFlagfold}. Let $P \in \G$ then $i(P) = \frac{1}{d} \Pi_P$ and $\mu_d(i(P)) = 1$, $E_d(i(P)) = P$ and for $k \neq d$, $\mu_k(i(P)) = 0$. Consequently,
for $\psi \in \xC_c(\Omega \times \cW\cF(n))$, we have by definitions,
\begin{align*}
\int_{\Omega \times \cW\cF(n)} \psi  \: d ( \Phi^\# \widehat{V} )  & = \sum_{k=1}^n \int_{\Omega \times \cW\cF(n)} \mu_k(S) \psi \left( \Phi(x), i \left( D\Phi(x) \cdot E_k(S)\right) \right) J\Phi(x, E_k(S)) \: d \widehat{V} (x,S)  \\
& = \sum_{k=1}^n \int_{\Omega \times \G} \hspace*{-7pt} \mu_k(i(P)) \psi \left( \Phi(x), i \left( D\Phi(x) \cdot E_k(i(P)) \right) \right) J\Phi(x, E_k(i(P))) \: d V (x,P) \\
& = \int_{\Omega \times \G}  \psi \left( \Phi(x), i \left( D\Phi(x) \cdot P \right) \right) J\Phi(x, P) \: d V (x,P)  \\
& = \int_{\Omega \times \G} \psi(x, i(P)) \: d (\Phi^\# V  )(x,P)  = \int_{\Omega \times \cW\cF(n)}  \psi  \:  d (\widehat{ \Phi^\# V } ) \: ,
\end{align*}
which shows the required consistency $\Phi^\# \widehat{V} = \widehat{\Phi^\# V}$. Let moreover $W$ be a flagfold in $\Omega$ and for $d = 1, \ldots, n$, let $V_d$ be the $d$--varifold associated with $W$ through \eqref{eq:flagfoldsToVarifolds}. Again by definitions, we have for $d = 1, \ldots, n$ and $\psi \in \xC_c(\Omega \times \cW\cF(n))$,
\begin{align*}
\int_{\Omega \times \cW\cF(n)} \psi  \: d (\widehat{\Phi^\# V_d})  & =  \int_{\Omega \times \G} \psi \left( \Phi(x), i \left(  D\Phi(x) \cdot P \right) \right) J\Phi(x, P) \: d V_d (x,P)  \\
& = \int_{\Omega \times \cW\cF(n)} \mu_d(S) \psi \left( \Phi(x), i \left( D\Phi(x) \cdot E_d(S) \right) \right) J\Phi(x, E_d(S)) \: d W (x, S) \: ,
\end{align*}
and summing from $d = 1$ to $n$, we obtain
\begin{equation}
\label{eq:pushforwardSum}
\Phi^\# W = \sum_{d=1}^n \widehat{ \Phi^\# V_d } \: .
\end{equation}
Note that the previous equality holds even in the case $W \neq \sum_{d=1}^n \widehat{  V_d }$.

We can now extend the first variation from varifolds to flagfolds. In the case of a $d$--varifold $V$, one computes the variation of the total mass $\V (\Omega)$ under a diffeomorphism $\Phi_t$ as follows (see \cite[Chap. 2]{simon}). Let $X \in \xC_c^1(\Omega)$ and $\Phi_t(x) = x + t X(x)$ be a diffeomorphism of $\Omega$ for $t$ small enough, then
\begin{equation*}
\left. \frac{d}{dt} \left\| \Phi_t^\# V \right\|(\Omega) \right|_{t=0} = \int_{\Omega \times \G} \xdiv_P X(x) \: dV(x,P)  \quad \text{with} \quad \xdiv_P X(x) := \tr \left( \Pi_P DX(x) \right) \: .
\end{equation*}
The first variation of $V$ is then defined by (see \cite{Allard72,simon})
\[
\begin{array}{lcccl}
\delta V & : & \xC_c^1 (\Omega, \R^n) & \rightarrow & \R \\
         &   &             X          & \mapsto     & \displaystyle \int_{\Omega \times \G} \xdiv_P X(x) \: dV(x,P)  \: .
\end{array}
\]
As it is done when defining the first variation $\delta V$ of a varifold $V$, we now compute the variation of the total mass of a flagfold $W$:
\begin{equation}
\left. \frac{d}{dt} \left\| \Phi_t^\# W \right\|(\Omega) \right|_{t=0} = \int_{\Omega \times \cW\cF(n)} \tr \left( \overline{S} DX(x) \right) \: dW(x,S)   \: .
\end{equation}
\noindent Indeed, for $(x,P) \in \Omega \times \G$, we recall that $\displaystyle \left.  \frac{d}{dt} J\Phi_t(x,P) \right|_{t=0} = \tr \left(\Pi_P DX(x) \right) = \xdiv_{P} X(x)$, therefore,
\begin{align*}
\left. \frac{d}{dt} \left\| \Phi_t^\# W \right\|(\Omega) \right|_{t=0} & = \left. \frac{d}{dt} \sum_{d=1}^n \int_{\Omega \times \cW\cF(n)}  J\Phi_t(x, E_d(S)) \mu_d( S)  \: dW(x,S) \right|_{t=0} \\
& = \sum_{d=1}^n \int_{\Omega \times \cW\cF(n)} \left.  \frac{d}{dt} J\Phi_t(x, E_d(S)) \right|_{t=0} \mu_d( S)  \: dW(x,S) \\
& = \int_{\Omega \times \cW\cF(n)} \sum_{d=1}^n  \mu_d(S) \tr \left( \Pi_{E_d(S)} DX(x) \right) \: d W(x,S) \\
& = \int_{\Omega \times \cW\cF(n)} \tr \left( \overline{S} DX(x) \right) \: dW(x,S)  \: .
\end{align*}
where we have used $\overline{S} = \sum_{d=1}^n \mu_d(S) \Pi_{E_d(S)}$ and linearity of the trace.

\begin{definition}[First Variation of a flagfold]
\label{dfn:flagfoldFirstVariation}
Let $W$ be a flagfold in $\Omega$. We define the \emph{first variation} of $W$ as the continuous linear form
\[
\begin{array}{lcccl}
\delta W & : & \xC_c^1 (\Omega, \R^n) & \rightarrow & \R \\
         &   &             X          & \mapsto     & \displaystyle \int_{\Omega \times \cW\cF(n)} \tr \left( \overline{S} DX(x) \right) \: dW(x,S) \\
         &   &                        &             & \displaystyle = \sum_{d=1}^n \int_{\Omega \times \cW\cF(n)} \mu_d(S) \xdiv_{E_d(S)} X(x) \: d W(x,S) \: .
\end{array}
\]
\end{definition}
\noindent We will use the notation $\xdiv_{\overline{S}} X(x) := \tr \left( \overline{S} DX(x) \right)$. We say that $W$ has \emph{locally bounded first variation} if $\delta W$ identifies with a Radon measure, that is: for all $K \subset \Omega$ compact set, there exists $c_K > 0$ such that for all $X \in \xC_c^1(\Omega, \R^n)$ with support in $K$,
\[
| \delta W(X) | \leq c_K \| X \|_\infty \: .
\]

\begin{remk}
As a consequence of both the consistency of the push-forward of flagfold with respect to the embedding of varifolds and \eqref{eq:pushforwardSum}, we note the following correspondences concerning the first variation of flagfolds and varifolds:
\begin{enumerate}[$\bullet$]
    \item Let $V$ be a $d$--varifold, then $\delta \widehat{V} = \delta V$.
    \item Let $W$ be a flagfold and let $(V_d)_{d = 1 \ldots n}$ be the $d$--varifolds associated through \eqref{eq:flagfoldsToVarifolds}, then
    \begin{equation}
        \delta W = \sum_{d=1}^n \delta V_d 
    \end{equation}
    even though $W = \widehat{V_1} + \ldots + \widehat{V_n}$ does not hold in general.
\end{enumerate}
\end{remk}

We conclude this section with two simple examples of flagfolds in $\R^3$ that are not (associated with) varifolds. The canonical basis of $\R^3$ is $(e_1, e_2, e_3)$ and we fix the weighted flag 
\begin{equation*}
S = \left( \begin{array}{ccc}
    \frac{2}{3} & 0 & 0 \\
    0 & \frac{1}{3} & 0 \\
    0 &  0  &  0
\end{array} \right) \in \xSym_+^1(3) \quad \text{i.e.} \quad \left\lbrace 
\begin{array}{l}
    (\lambda_1, \lambda_2, \lambda_3) = \left( \lambda_1(S), \lambda_2(S), \lambda_3(S) \right) = \left( \frac{2}{3}, \frac{1}{3}, 0 \right) \\
    (\mu_1, \mu_2, \mu_3) = \left( \mu_1(S), \mu_2(S), \mu_3(S) \right) = \left( \frac{1}{3}, \frac{2}{3}, 0 \right) \\
     E_1 = E_1(S) = \xspan(e_1) , \,   E_2 = E_2(S) = \xspan(e_1,e_2)  
\end{array}
\right. \: .
\end{equation*}

\begin{xmpl}
\label{xmpl:flagfold1}
    We first consider $W = \cL^2_{|E_2} \otimes \delta_S$. Applying definitions, we have $\| W \| = \cL^2_{|E_2}$ and
\begin{equation*}
 V_k = \mu_k \cL^2_{|E_2} \otimes \delta_{E_k} \text{ for } k = 1,2 \quad  \text{while } V_3 = 0 \quad \Rightarrow \quad \widehat{V_k} = \mu_k \cL^2_{|E_2} \otimes \delta_{i(E_k)} \text{ for } k = 1,2 \quad \text{and } \widehat{V_3} = 0 \: . 
\end{equation*}
Considering the Borel set $C = B(0,1) \times \{ i(E_1) \} = B(0,1) \times \{ \Pi_{E_1} \} \subset \R^3 \times \cW\cF(3)$, we observe that $\widehat{V_1}(C) = \mu_1 \cL^2 (B(0,1) \cap E_2) = \frac{\pi}{3} $ and $\widehat{V_2}(C) = \widehat{V_3}(C) $ as well as $W(C) = 0$ and in such a simple case, $W = \widehat{V_1} + \widehat{V_2} + \widehat{V_3}$ does not hold. We also note that for $k = 1,2$,
\begin{equation*}
\| V_k \| = \mu_k \cL^2_{|E_2} \quad \Rightarrow \quad \Theta_\ast^k (x , \| V_k \| ) = \liminf_{\rho \to 0_+} \frac{\|V_k\|(B(x,\rho))}{\rho^k} = \left\lbrace
\begin{array}{ccc}
   \lim_{\rho \to 0_+} \mu_1\pi \rho = 0  & \text{if} & k = 1  \\
    \mu_2 \pi > 0 & \text{if} & k = 2
\end{array}
\right. \: .
\end{equation*}
Let us compute the first variation $\delta W = \delta V_1 + \delta V_2$. 
Let $X \in \xC_c^1 (\R^3, \R^3)$, applying Fubini theorem we obtain
\begin{align*}
    \delta V_1 (X) = \mu_1 \int_{x \in E_2} \xdiv_{E_1} X(x) \: d \cL^2(x) = \mu_1 \int_{x_2 \in \R} \underbrace{ \int_{x_1 \in \R} \partial_1 X_1 (x_1,x_2,x_3) \: d x_1 }_{= 0} d x_2 = 0 \: .
\end{align*}
Note that $V_1$ is not a rectifiable $1$--varifold ($E_2$ is a $2$--rectifiable set but it is not a $1$--rectifiable set) though $\delta V_1 = 0$, we point out that Allard's rectifiability theorem requires $\Theta_\ast^1 (x , \| V_1 \| ) > 0$ that fails to be true here.
As $V_2$ is a $2$--varifold associated with a plane whose mean curvature is zero (and the multiplicity $\mu_2$ is constant), we have that $\delta V_2 = 0$ (this can also be checked by similar computations) and eventually $\delta W = 0$. 
\end{xmpl}

\begin{xmpl}
    We now consider $W = \cL^1_{|E_1} \otimes \delta_S$. Applying definitions, we have $\| W \| = \cL^1_{|E_1}$ and
\begin{equation*}
 V_k = \mu_k \cL^1_{|E_1} \otimes \delta_{E_k} \text{ for } k = 1,2 \quad  \text{while } V_3 = 0 \quad \Rightarrow \quad \widehat{V_k} = \mu_k \cL^1_{|E_1} \otimes \delta_{i(E_k)} \text{ for } k = 1,2 \quad \text{and } \widehat{V_3} = 0 \: . 
\end{equation*}
As in the previous example $W = \widehat{V_1} + \widehat{V_2} + \widehat{V_3}$ does not hold and we now have
\begin{equation*}
\| V_k \| = \mu_k \cL^1_{|E_1} \quad \Rightarrow \quad \Theta_\ast^k (x , \| V_k \| ) = \left\lbrace
\begin{array}{ccc}
   2 \mu_1 > 0  & \text{if} & k = 1  \\
    \lim_{\rho \to 0_+}  \frac{2 \mu_2}{\rho} = +\infty & \text{if} & k = 2
\end{array}
\right. \: .
\end{equation*}
Let us compute the first variation $\delta W = \delta V_1 + \delta V_2$. This time $V_1$ is a $1$--varifold associated with a line whose mean curvature is zero and thus $\delta V_1 = 0$.
Let $X \in \xC_c^1 (\R^3, \R^3)$,
\begin{align*}
    \delta V_2 (X) & = \mu_2 \int_{x \in E_1} \xdiv_{E_2} X(x) \: d \cL^1(x) = \mu_2 \int_{x_1 \in \R} \partial_1 X_1 (x_1,x_2,x_3) + \partial_2 X_2 (x_1,x_2,x_3) \: d x_1 \\
    & = \mu_2 \int_{x_1 \in \R}  \partial_2 X_2 (x_1,x_2,x_3) \: d x_1  \: ,
\end{align*}
and we can infer that $\delta W = \delta V_2$ is not a Radon measure.
\end{xmpl}

\subsection{Monotonicity formula and consequences}
\label{subsec:monotonicity}

Let $W$ be a flagfold in $\Omega$ and let $(V_1, \ldots, V_n)$ be the associated varifolds defined in \eqref{eq:flagfoldsToVarifolds}.
We recall that $\displaystyle \|W \| = \sum_{k=1}^n \|V_k\|$ and $\delta W = \sum_{k=1}^n \delta V_k$ even though $W = \sum_{k=1}^n V_k$ does not hold in general, as we illustrated in the simple case of Example~\ref{xmpl:flagfold1}. In this section, we give sufficient conditions, involving lower $d$--dimensional density for $d = 1, \ldots, n$ and the first variation of $W$, that allow to recover $W$ as the sum of the aforementioned $d$--varifolds: loosely speaking, such conditions ensure that the weighted flags of mixed dimensions have zero measure.

We recall that we defined the continuous function $\displaystyle\overline{d} : \cW\cF(n) \rightarrow [1,n]$, $\overline{d}(S) = \tr \overline{S} = \sum_{k=1}^n k \mu_k(S)$ in \eqref{eq:dimensionMu}. Given a flagfold $W = \| W \| \otimes w^y$ (see Proposition~\ref{prop:disintegration}) in $\Omega \subset \R^n$, we define space $\mu : \Omega \subset \Delta(n)$ and dimension $\overline{d} : \Omega \rightarrow [1,n]$ by
\begin{equation} \label{eq:mux}
\mu (y) := \int_{S \in \cW\cF(n)} \mu(S) \: d w^y(S)
\quad \text{and} \quad
\overline{d}(y) := \int_{S \in \cW\cF(n)} \overline{d}(S) \: d w^y(S) = \sum_{k=1}^n k \mu_k(y)
\end{equation}
that are well-defined for $\| W \|$--a.e. $y \in \Omega$ and moreover in $\xL^1_{loc}(\|W\|)$.
%Before that, we use Proposition~\ref{prop:disintegration} to define for $\| W\|$--a.e. $x \in \Omega$,
%\begin{equation} 
%\mu_d(x) = \int_{S \in \cW\cF(n)} \mu_d(S) \: d w^x(S), \quad d = 1, \ldots ,n \: .
%\end{equation}
With such a notation, we have $\displaystyle \| V_d \| = \| \widehat{V}_d \| = \mu_d(y) \| W \|$.

\begin{proposition}[A monotonicity formula] \label{prop:monotonicity}
Let $W$ be a flagfold in $\Omega$ with locally bounded first variation: we identify $\delta W$ with a vector Radon measure in $\Omega$ and denote by $| \delta W|$ its total variation.
Let $x \in \Omega$ and assume that $\Lambda \geq 0$, $0 < d^\ast \leq n$ and $\rho^\ast > 0$ satisfy for all $0 < \rho \leq \rho^\ast$,
\begin{equation} \label{eq:monotonicityHyp}
\int_{B_\rho(x)} (\overline{d}(y) - d^\ast) \, d \|W\|(y) \geq 0 \quad \text{and} \quad | \delta W |(B_\rho(x)) \leq \Lambda \|W\|(B_\rho(x)) \: .
\end{equation}
Then $\rho \mapsto e^{\Lambda \rho} \rho^{-d^\ast} \| W \|(B_\rho(x))$ is non-decreasing. Moreover, in the particular case where $\delta W = 0$, then for all $0 < \sigma < \xi < {\rm dist} (x, \R^n \setminus \Omega)$,
\begin{multline} \label{eq:monotonicity}
\frac{\|W\|(B_\xi(x))}{\xi^{d^\ast}} -  \frac{\|W\|(B_\sigma(x))}{\sigma^{d^\ast}} = 
 \int_{B_\xi(x) \setminus B_\sigma(x) \times \cW\cF(n)} \frac{\overline{S}^\perp(y-x) \cdot (y-x)}{|y-x|^{d^\ast + 2}} \: dW(y,S) \\
+ \int_\sigma^\xi \rho^{- d^\ast-1} \int_{B_\rho(x)} \left( \overline{d}(y) - d^\ast \right)  \: d \|W \|(y) \: d \rho
\end{multline}
\end{proposition}

\begin{remk}
Proposition~\ref{prop:monotonicity} differs from the classical monotonicity formula (see \cite[§3,chap.8]{simon}, \cite[Thm 2.1]{DeLellis18}) because of the additional term involving $\overline{d}(S) - d^\ast$, (and then $\overline{d}(y) - d^\ast$) while for $d$-varifolds, for all $S \in \G$, $\overline{d}(S) = \tr \Pi_S = d$ and it is possible to chose $d^\ast = d$ to fulfill assumption \eqref{eq:monotonicityHyp}. For the sake of clarity, we give the whole proof of Proposition~\ref{prop:monotonicity} though the reader can focus on the term involving $\overline{d}(S) - d^\ast$.
\end{remk}

\begin{proof}
{\bf Step 1}: As in the classical monotonicity formula for $d$--varifolds (see \cite{simon}), we test the first variation $\delta W$ with a radial test function $X = X_x \in \xC_c^1 (\Omega, \R^n)$ defined as
\begin{equation} \label{eq:monotonicity8}
X(y) = \phi \left( \frac{|y-x|}{\rho} \right) (y-x) \quad \text{with} \quad \left\lbrace
\begin{array}{l}
\phi : \R \rightarrow \R_+ \text{ of class } \xC^1\\
\phi \text{ is symmetric and supported in } [-1,1]\\
\dt{\phi} \leq 0 \text{ in } [0,+\infty[
\end{array} \right. \: .
\end{equation}
Given $S \in \cW\cF(n)$, we first compute
\begin{align}
\xdiv_{\overline{S}} X(y) & = \phi \left( \frac{|y-x|}{\rho} \right) \tr \left( \overline{S} \right) + \frac{1}{\rho} \dt{\phi} \left( \frac{|y-x|}{\rho} \right) \frac{\overline{S} (y-x) \cdot (y-x)}{|y-x|} \nonumber \\
& = \overline{d}(S) \phi \left( \frac{|y-x|}{\rho} \right) + \frac{|y-x|}{\rho} \dt{\phi} \left( \frac{|y-x|}{\rho} \right) \left( 1 - \frac{(I_n - \overline{S})(y-x)}{|y-x|} \cdot \frac{(y-x)}{|y-x|} \right) \: . \label{eq:monotonicity1}
\end{align}
As $\displaystyle \frac{|y-x|}{\rho} \dt{\phi} \left( \frac{|y-x|}{\rho} \right) = - \rho \frac{d}{d \rho} \left\lbrace \phi \left( \frac{|y-x|}{\rho} \right) \right\rbrace$ we can rewrite
\[
\int \frac{|y-x|}{\rho} \dt{\phi} \left( \frac{|y-x|}{\rho} \right) \: dW(y,S) = - \rho \frac{d}{d\rho} \mathcal{I}(\rho,\phi) \quad \text{with} \quad \mathcal{I}(\rho,\phi):= \int \phi \left( \frac{|y-x|}{\rho} \right) \: dW(y,S) \: .
\]
and denoting by $\overline{S}^\perp := I_n - \overline{S}\in \xSym_+(n)$, 
\begin{multline} \label{eq:monotonicity7}
\int \frac{|y-x|}{\rho} \dt{\phi} \left( \frac{|y-x|}{\rho} \right)\frac{\overline{S}^\perp(y-x)}{|y-x|} \cdot \frac{(y-x)}{|y-x|} \: dW(y,S) = - \rho \frac{d}{d\rho} \mathcal{J}(\rho,\phi) \leq 0 \\ 
\text{with} \quad \mathcal{J}(\rho,\phi):= \int \phi \left( \frac{|y-x|}{\rho} \right)\frac{\overline{S}^\perp(y-x)}{|y-x|} \cdot \frac{(y-x)}{|y-x|} \: dW(y,S) \: .
\end{multline}
%We denote by $\overline{S}^\perp := I_n - \overline{S} \in \xSym_+(n)$ so that
%\begin{equation} \label{eq:monotonicity3}
%\mathcal{J} (\rho, \phi) := -\int \frac{|y-x|}{\rho} \phi^\prime \left( \frac{|y-x|}{\rho} \right) \frac{\overline{S}^\perp(y-x)}{|y-x|} \cdot \frac{(y-x)}{|y-x|}  \: dW(y,S) \geq 0 \: .
%\end{equation}
We can integrate \eqref{eq:monotonicity1} with respect to $W$, leading to
\begin{equation} \label{eq:monotonicity2}
\delta W (X) - \rho \frac{d}{d\rho} \mathcal{J}(\rho,\phi) - \int \left( \overline{d}(S) - d^\ast \right) \phi \left( \frac{|y-x|}{\rho} \right) \: dW(y,S) = d^\ast \mathcal{I}(\rho,\phi) - \rho \frac{d}{d\rho} \mathcal{I}(\rho,\phi) \: .
\end{equation}

{\bf Step 2}: We assume $\delta W = 0$ and prove \eqref{eq:monotonicity}.  Notice that
\[
\frac{d}{d \rho} \left\lbrace \rho^{-d^\ast} \mathcal{I}(\rho,\phi) \right\rbrace = - d^\ast \rho^{-d^\ast - 1} \mathcal{I}(\rho,\phi) + \rho^{-d^\ast} \frac{d}{d\rho} \mathcal{I}(\rho,\phi) = - \rho^{-d^\ast - 1} \left( d^\ast \mathcal{I}(\rho,\phi) - \rho \frac{d}{d\rho} \mathcal{I}(\rho,\phi)\right) \: .
\]
%We start with rewriting
%\[
%\mathcal{J} (\rho, \phi) = \rho \frac{d}{d\rho} \int  \phi \left( \frac{|y-x|}{\rho} \right) \frac{\overline{S}^\perp(y-x)}{|y-x|} \cdot \frac{(y-x)}{|y-x|}  \: dW(y,S) 
%\]
%which we insert in \eqref{eq:monotonicity2}, and after multiplication by $\rho^{-d^\ast-1}$, we obtain
%\begin{multline} \label{eq:monotonicity6}
%-\rho^{-d^\ast} \frac{d}{d\rho} \left\lbrace \int  \phi \left( \frac{|y-x|}{\rho} \right) \frac{\overline{S}^\perp(y-x)}{|y-x|} \cdot \frac{(y-x)}{|y-x|}  \: dW(y,S) \right\rbrace
%\\ - \rho^{- d^\ast-1} \int \left( \overline{d}(S) - d^\ast \right) \phi \left( \frac{|y-x|}{\rho} \right) \: dW(y,S) 
%= \frac{d}{d\rho} \left\lbrace \rho^{-d^\ast} \mathcal{I}(\rho,\phi) \right\rbrace \: .
%\end{multline}
Therefore, we first multiply \eqref{eq:monotonicity2} by $- \rho^{-d^\ast-1}$ and integrate between $0 < \sigma < \xi < {\rm dist}(x, \R^n \setminus \Omega)$ and integrate by parts,
\begin{multline} \label{eq:monotonicity9}
  \xi^{-d^\ast} \mathcal{J}(\xi,\phi) -\sigma^{-d^\ast} \mathcal{J}(\sigma,\phi)  + \int_{\rho= \sigma}^\xi  d^\ast \rho^{-d^\ast-1} \mathcal{J} (\rho, \phi) \: d\rho
\\ 
 + \int_{\Omega \times \cW\cF(n)} \left( \overline{d}(S) - d^\ast \right) \int_{\rho=\sigma}^\xi \rho^{- d^\ast-1} \phi \left( \frac{|y-x|}{\rho} \right) \: d\rho \: dW(y,S) 
=   \xi^{-d^\ast} \mathcal{I}(\xi,\phi)  - \sigma^{-d^\ast} \mathcal{I}(\sigma,\phi) \: .
\end{multline}
We choose a sequence $(\phi_h)_{h \in \N}$ of functions satisfying \eqref{eq:monotonicity8} and pointwise converging to $\one_{]-1,1[}$ and let $n \to +\infty$ in \eqref{eq:monotonicity9}.
%we infer \eqref{eq:monotonicity} by dominated convergence. Indeed,
By dominated convergence, for all $\sigma \leq \rho \leq \xi$,
\begin{enumerate}[(i)]
\item $\displaystyle \rho^{-d^\ast}  \cI (\rho, \phi_h) \xrightarrow[h \to +\infty]{} \rho^{-d^\ast} \|W\|(B_\rho(x))
$,
\item $\displaystyle \cJ (\rho, \phi_h) \xrightarrow[h \to +\infty]{} \cJ (\rho) := \int_{B_\rho(x) \times \cW\cF(n)} \frac{\overline{S}^\perp(y-x)}{|y-x|} \cdot \frac{(y-x)}{|y-x|} \: dW(y,S)$ and 
\begin{align*}
\int_{\rho= \sigma}^\xi & d^\ast  \rho^{-d^\ast-1}  \mathcal{J} (\rho, \phi_h) \: d\rho \xrightarrow[h \to +\infty]{} 
\int_{\rho= \sigma}^\xi  d^\ast  \rho^{-d^\ast-1}  \mathcal{J} (\rho) \: d\rho \\
& = d^\ast \int_{B_\xi(x) \times \cW\cF(n)} \frac{\overline{S}^\perp(y-x)}{|y-x|} \cdot \frac{(y-x)}{|y-x|} \int_{\rho= \max \{ \sigma, |y-x| \}}^\xi \rho^{-d^\ast-1} \: d\rho \: dW(y,S) \\
& = \int_{B_\xi(x) \times \cW\cF(n)}  \frac{\overline{S}^\perp(y-x)}{|y-x|} \cdot \frac{(y-x)}{|y-x|} \left(  \max \{ \sigma, |y-x| \}^{- d^\ast} - {\xi}^{-d^\ast} \right) \: dW(y,S) \\
& = - {\xi}^{-d^\ast} \cJ(\xi) + {\sigma}^{-d^\ast} \cJ(\sigma) + \int_{B_\xi(x) \setminus B_\sigma(x) \times \cW\cF(n)} \frac{\overline{S}^\perp(y-x)}{|y-x|} \cdot \frac{(y-x)}{|y-x|}   |y-x|^{- d^\ast} \: dW(y,S)
\end{align*}
\item $\displaystyle \int_{\Omega \times \cW\cF(n)} \left( \overline{d}(S) - d^\ast \right) \int_{\rho=\sigma}^\xi \rho^{- d^\ast-1} \phi_h \left( \frac{|y-x|}{\rho} \right) \: d\rho \: dW(y,S)$ similarly converges to 
\begin{multline*} 
\int_{\Omega \times \cW\cF(n)} \left( \overline{d}(S) - d^\ast \right) \int_{\rho=\sigma}^\xi \rho^{- d^\ast-1} \one_{B_\rho(x)} \: d\rho \: dW(y,S) \\
 = - \frac{\xi^{-d^\ast}}{d^\ast}\int_{B_\xi(x) \times \cW\cF(n)} \left( \overline{d}(S) - d^\ast \right) \: dW(y,S) + \frac{\sigma^{-d^\ast}}{d^\ast}\int_{B_\sigma(x) \times \cW\cF(n)} \left( \overline{d}(S) - d^\ast \right) \: dW(y,S) \\
+ \frac{1}{d^\ast} \int_{B_\xi(x) \setminus B_\sigma(x) \times \cW\cF(n)} \left( \overline{d}(S) - d^\ast \right) |y-x|^{-d^\ast} \: dW(y,S)
\end{multline*}
Formula \eqref{eq:monotonicity} follows from $(i)$--$(iii)$ and \eqref{eq:monotonicity9}.
\end{enumerate}

{\bf Step 3}: We are left with the monotonicity of $\rho \mapsto e^{\Lambda \rho} \rho^{-d^\ast} \| W \|(B_\rho(x))$.
As in Step $2$, we multiply \eqref{eq:monotonicity2} by $- \rho^{-d^\ast-1}$, then using \eqref{eq:monotonicity7} and
\[
- \delta W(X) \geq - \int \underbrace{\phi \left( \frac{|y-x|}{\rho} \right)}_{\leq \one_{B_\rho(x)}} \underbrace{|y-x|}_{\leq \rho} \: d | \delta W|(y) \geq - \rho |\delta W| (B_\rho(x)) \: ,
\]
we obtain
\begin{equation} \label{eq:monotonicity4}
- \rho^{- d^\ast} |\delta W| (B_\rho(x)) + \rho^{- d^\ast-1} \int \left( \overline{d}(S) - d^\ast \right) \phi \left( \frac{|y-x|}{\rho} \right) \: dW(y,S) \leq \frac{d}{d\rho} \left\lbrace \rho^{-d^\ast} \mathcal{I}(\rho,\phi) \right\rbrace \: .
\end{equation}
We integrate the previous inequality \eqref{eq:monotonicity4} between $0 < \sigma < \xi \leq \rho^\ast$ and let $\phi \rightarrow \one_{]-1,1[}$, we get by dominated convergence
\begin{multline*}
%- \int_{\sigma}^\xi \rho^{- d^\ast} |\delta W| (B_\rho(x)) \: d\rho & \leq
- \int_{\sigma}^\xi \rho^{- d^\ast} |\delta W| (B_\rho(x)) \: d\rho + \int_\sigma^\xi \rho^{- d^\ast-1} \int_{B_\rho(x)} \left( \overline{d}(y) - d^\ast \right)  \: d \|W \|(y) \: d \rho \\
 \leq  \xi^{-d^\ast} \| W \| (B_\xi(x)) - \sigma^{-d^\ast} \| W \| (B_\sigma(x)) \: .
\end{multline*}
Thanks to assumptions \eqref{eq:monotonicityHyp} we eventually obtain the following integral relation
\begin{equation} \label{eq:monotonicity5}
\xi^{-d^\ast} \| W \| (B_\xi(x)) - \sigma^{-d^\ast} \| W \| (B_\sigma(x)) \geq - \Lambda \int_{\rho = \sigma}^\xi \rho^{- d^\ast} \| W \|(B_\rho(x)) \: d\rho  \: .
\end{equation}
The proof can be concluded as for the classical monotonicity formula (we give the details for the sake of completeness).\\
Let $f : ]0,\rho_\ast[ \rightarrow \R_+$, $\rho \mapsto \rho^{-d^\ast} \|W\|(B_\rho(x))$, $f$ has bounded variation (since $\rho \mapsto \| W \|(B_\rho(x))$ is non-decreasing) and we denote by $Df$ its distributional derivative. The integral inequation \eqref{eq:monotonicity5} then implies that
$
Df + \Lambda f 
$ is a positive Radon measure.

\noindent Indeed, let $\psi \in \xC_c^1 (]0,\rho_\ast[, \R_+)$, $h$ small enough so that $\psi (\cdot - h) \in \xC_c^1 (]0,\rho_\ast[, \R_+)$ as well, and $\xi = \sigma + h$ in \eqref{eq:monotonicity5}:
\[
- \Lambda \int_{\sigma = 0}^{\rho^\ast} \psi(\sigma) \frac{1}{h}\int_{\rho = \sigma}^{\sigma + h} f(\rho) \: d\rho  \: d \sigma \leq \int_{\sigma = 0}^{\rho^\ast} \frac{f(\sigma + h) - f(\sigma)}{h} \psi(\sigma) \: d \sigma = \int_{\sigma = 0}^{\rho^\ast} f(\sigma) \frac{\psi(\sigma-h) - \psi(\sigma)}{h} \: d \sigma
\]
Letting $h \to 0$, we conclude thanks to dominated convergence that
\[
\left\langle - \Lambda f , \psi \right\rangle = - \Lambda \int_{\sigma = 0}^{\rho^\ast} \psi(\sigma) f(\sigma)  \: d \sigma \leq - \int_{\sigma = 0}^{\rho^\ast} f(\sigma) \dt{\psi}(\sigma) \: d \sigma = \left\langle Df , \psi \right\rangle
\]
Eventually, $g (\rho) := e^{\Lambda \rho} f(\rho)$ has bounded variation as well and 
\[
Dg = e^{\Lambda \rho}\left( Df + \Lambda f \right)
\]
that is a {\it positive} Radon measure and consequently $g$ is non-decreasing.

\end{proof}

\begin{remk}
Notice that the assumption \eqref{eq:monotonicityHyp} for all $0 < \rho < \rho^\ast$, $\displaystyle \int_{B_\rho(x)} (d(y) - d^\ast) \, d \|W\|(y) \geq 0$ could be replaced with $\displaystyle \int_{B_\rho(x)} (d(y) - d^\ast) \, d \|W\|(y) \geq - \Lambda^\prime \rho \| W \|(B_\rho(x))$, leading to the monotonicity of $\rho \mapsto e^{(\Lambda + \Lambda^\prime)\rho} \rho^{-d^\ast} \|W\|(B_\rho(x))$.
\end{remk}

\begin{proposition} \label{prop:zeroDensity}
Let $W$ be a flagfold in $\Omega$ with locally bounded first variation, then for $\|W\|$--a.e. $x \in \Omega$,
\[
\Theta^d (x, \|W\|) := \lim_{r \to 0_+} \frac{\|W\|(B_r(x))}{\omega_d r^d} = 0, \quad \forall 1 \leq d < \overline{d}(x) \: ,
\]
and in the particular case where $\overline{d}(x) = 1$, then $\displaystyle \Theta^1 (x, \|W\|)$ exists and is finite.
\end{proposition}

\begin{proof}
It is an easy consequence of Proposition~\ref{prop:monotonicity}. First notice that thanks to differentiation properties of Radon measures and $\overline{d} \in \xL^1(\|W\|)$, for $\| W \|$--a.e. $x \in \Omega$,
\[
\lim_{\rho \to 0_+} \frac{1}{\|W\|(B_\rho(x))} \int_{B_\rho(x)} \overline{d}(y) \: d \|W\|(y) = \overline{d}(x) \: .
\]
Therefore, given $1 \leq d^\ast < \overline{d}(x)$ (we start with the case where $\overline{d}(x) > 1$), there exists $\rho^\ast = \rho^\ast(x) > 0$ such that for all $0 < \rho < \rho^\ast$,
\begin{equation} \label{eq:zeroDensity1}
\int_{B_\rho(x)} \overline{d}(y) \: d \|W\|(y) \geq d^\ast \|W\|(B_\rho(x)) \quad \text{i.e.} \quad \int_{B_\rho(x)} (\overline{d}(y) - d^\ast) \: d \|W\|(y) \geq 0 \: .
\end{equation}
In addition, for $\| W \|$--a.e. $x \in \Omega$,
\[
\lim_{\rho \to 0_+} \frac{|\delta W|(B_\rho(x))}{\|W\|(B_\rho(x))} \text{ exists and is finite,}
\]
and for such an $x \in \Omega$, there exists consequently $\Lambda(x) \geq 0$ such that for all $0 < \rho < \rho^\ast(x)$,
\begin{equation} \label{eq:zeroDensity2}
|\delta W| (B_\rho(x)) \leq \Lambda(x) \|W\|(B_\rho(x)) \: .
\end{equation}
Thanks to \eqref{eq:zeroDensity1} and \eqref{eq:zeroDensity2}, assumptions \eqref{eq:monotonicityHyp} are satisfied $\|W\|$--a.e. in $\Omega$, whence Proposition~\ref{prop:monotonicity} ensures that for $\| W \|$--a.e. $x \in \Omega$, the application $\rho \mapsto e^{\Lambda(x) \rho} \rho^{-d^\ast} \|W\|(B_\rho(x))$ is non-decreasing in a neighbourhood $]0,\rho^\ast(x)[$ of $0_+$ and has a finite limit when $\rho \to 0_+$. We infer
\begin{equation} \label{eq:zeroDensity3}
\Theta^{d^\ast} (x, \|W\|) \text{ exists and is finite.}
\end{equation}
Let $1 \leq d < d^\ast < \overline{d}(x)$, from \eqref{eq:zeroDensity3} we know that
\[
\frac{\|W\|(B_\rho(x))}{\rho^{d}} = \underbrace{ \frac{\|W\|(B_\rho(x))}{\rho^{d^\ast}} }_{\xrightarrow[\rho \to 0_+]{} \Theta^{d^\ast}(x,\|W\|) < \infty} \underbrace{\rho^{d^\ast - d}}_{\xrightarrow[\rho \to 0_+]{} 0} \xrightarrow[\rho \to 0_+]{} 0 \: .
\]
We are left with the case where $\overline{d}(x) = 1$, note that for any $S \in \cW\cF(n)$, $\overline{d}(S) \geq 1$ hence for $\| W \|$--a.e. $y \in \Omega$,
$\displaystyle
\overline{d}(y) = \int_{S \in \cW\cF(n)} \overline{d}(S) \: d w^y(S) \geq w^y (\cW\cF(n)) = 1   = \overline{d}(x)
$
and thus Proposition~\ref{prop:monotonicity} applies with $d^\ast = \overline{d}(x) = 1$ since
\[
\int_{B_\rho(x)} (\overline{d}(y) - 1) \: d \|W\|(y) \geq 0 \: .
\]
Eventually $\Theta^1 (x, \|W\|)$ exists and is finite due to the monotonicity of $\rho \mapsto e^{\Lambda(x)\rho} \rho^{-1} \|W\|(B_\rho(x))$.
\end{proof}

We can finally state a structure theorem for a flagfold $W$ with locally bounded first variation, furthermore requiring that the $d$--lower density of each $\| V_d \|$, $d =1, \ldots, n$ defined in \eqref{eq:flagfoldsToVarifolds} do not vanish.

\begin{theorem}
\label{thm:rectifiability}
Let $W$ be a flagfold in $\Omega$ with locally bounded first variation and assume that for all $d = 1, \ldots, n$ and for $\| V_d \|$--a.e. $x \in \Omega$,
\[
\Theta_\ast^d (x, \| V_d \|) := \liminf_{\rho \to 0_+} \frac{\| V_d \|(B_\rho(x))}{\rho^d}  > 0 \: ,
\]
where $V_d$ has been defined in \eqref{eq:flagfoldsToVarifolds}.
Then, for $\| W \|$--a.e. $x \in \Omega$,
\[
\overline{d}(x) \in \N \cap [1,n] \quad \text{and} \quad \mu_{\overline{d}(x)} = 1 \: ,
\]
and $W = \widehat{V_1} + \ldots + \widehat{V_n}$ identifies with a sum of classical $d$--varifolds.

\noindent In particular, if we assume that $\delta V_1, \ldots, \delta V_n$ are Radon measures, then $W$ is a rectifiable flagfold.
\end{theorem}

\begin{proof}
We introduce for $d = 1, \ldots, n$,
\[
N_d = \{ x \in \Omega \: : \: \Theta_\ast^d (x, \| V_d \|) = 0\}  \quad \text{and} \quad M_d = \{ x \in \Omega \: : \: \mu_d(x) > 0 \} \: ,
\]
where $M_d$ is defined up to a set of $\|W\|$--null measure.

{\bf Step $1$}: We already know that $\| V_d \|(N_d) = 0$. The key point is to show that for all $i \neq j$,
\begin{equation} \label{eq:rectifibility0}
\| W \| (M_i \cap M_j) = 0 \: .
\end{equation}
We prove it by induction on $i \in \{1, \ldots, n\}$.

\noindent Let $x \in M_1 \cap M_j$ for $j \in \{ 2, \ldots, n\}$. First, notice that $\overline{d}(x) > 1$. Indeed, recalling \eqref{eq:dimensionMu} and \eqref{eq:mux}
\begin{align*}
\overline{d} (x) & = \int_{S \in \cW\cF(n)} \sum_{k=1}^n k \mu_k(S) \: d w^x(S) = \sum_{k=1}^n k \int_{S \in \cW\cF(n)} \mu_k(S) \: d w^x(S) = \sum_{k=1}^n k \mu_k(x) \\
& \geq  \underbrace{j \mu_j(x)}_{> \mu_j(x)} + \sum_{\substack{k = 1 \\ k \neq j}}^n \mu_k(x) > \sum_{k=1}^n \mu_k(x) = 1
\end{align*}
Therefore, thanks to Proposition~\ref{prop:zeroDensity} we know that $\Theta^1(x, \|W\|) = 0$. As $\| V_1 \| = \mu_1(x) \| W \| \leq \|W\|$, this directly implies $\Theta^1(x, \|V_1\|) = 0$ and $M_1 \cap M_j \subset N_1$ whence $\| V_1 \| (M_1 \cap M_j) = 0$. Moreover, when restricted to $M_1$, $\| W \|_{| M_1} = \frac{1}{\mu_1(x)} \| V_1 \|_{| M_1} << \| V_1 \|_{| M_1}$ and therefore
$
\| W \| (M_1 \cap M_j) = 0$. 

\noindent Let $i_0 \in \{ 2, \ldots, n\}$, assume that for all $i \leq i_0 -1$, for all $j \neq i$,  $\| W \| (M_{i} \cap M_j) = 0$ and let us prove that $\| W \| (M_{i_0} \cap M_j) = 0$ for all $j \neq i_0$.

\noindent Let $j \neq i_0$. If $j < i_0$ then $\| W \| (M_j \cap M_{i_0}) = 0$ by induction assumption.\\
If $j > i_0$, on one hand, let $x \in 
\displaystyle
M_{i_0} \cap M_j \cap \bigcap_{i=1}^{i_0 -1 } (\R^n \setminus M_{i}) 
$, and notice as previously that $\overline{d}(x) > i_0$. Indeed, $\mu_1(x) = \ldots = \mu_{i_0 -1}(x) = 0$ so that
 \begin{align*}
\overline{d} (x)  =   \sum_{k=1}^n k \mu_k(x) =  \sum_{k=i_0}^n k \mu_k(x)
 \geq  \underbrace{j \mu_j(x)}_{> i_0 \mu_j(x)} + \sum_{\substack{k = i_0 \\ k \neq j}}^n i_0 \mu_k(x) > i_0 \sum_{k=i_0}^n \mu_k(x) = i_0 \: .
\end{align*}
Consequently $\Theta^{i_0} (x, \|W\|) = \Theta^{i_0} (x, \| V_{i_0} \|) = 0$ and since moreover $\| V_{i_0} \| (N_{i_0}) = 0$ and $\| W \|_{| M_{i_0}} = \frac{1}{\mu_{i_0}(x)} \| V_{i_0} \|_{| M_{i_0}} \ll \| V_{i_0} \|_{| M_{i_0}}$ then
\begin{equation} \label{eq:rectifibility1}
M_{i_0} \cap M_j \cap \bigcap_{i=1}^{i_0 -1 } (\Omega \setminus M_{i}) \subset N_{i_0} \quad \Rightarrow \quad \| W \| \left( M_{i_0} \cap M_j \cap \bigcap_{i=1}^{i_0 -1 } (\Omega \setminus M_{i}) \right) = 0 \: .
\end{equation}
On the other hand 
\begin{equation} \label{eq:rectifibility2}
\| W \| \left( M_{i_0} \cap M_j \cap \bigcup_{i=1}^{i_0 -1 }  M_{i} \right) \leq \| W \| \left( \bigcup_{i=1}^{i_0-1} M_i \cap M_{i_0} \right) = 0 \text{ by induction assumption.}
\end{equation}
From \eqref{eq:rectifibility1} and \eqref{eq:rectifibility2} we infer that $\| W \|(M_{i_0} \cap M_j) = 0$ and conclude the induction.

\noindent {\bf Step $2$}: We draw consequences of \eqref{eq:rectifibility0}. First of all
\[
M_d = \{ x \in \Omega \: : \: \mu_d(x) > 0 \} \quad \text{and} \quad M_d^\prime = \{ x \in \Omega \: : \: \mu_d(x) = 1 \}
\]
are equal up to a set of $\| W \|$--measure zero. Moreover, the points $\mu \in \Delta(n)$ such that $\mu_d = 1$ for some $d$ and $\mu_i = 0$ for $i \neq d$ are exactly the vertices, that is the extremal points of $\Delta(n)$, hence they cannot arise as non trivial convex combinations of elements of $\Delta(n)$. In other words, if $x \in M_d^\prime$,
\[
1 = \mu_d(x) = \int_{S \in \cW\cF(n)} \underbrace{ \mu_d(S) }_{\leq 1} \: d w^x(S) \leq w^x(S) = 1
\]
with equality if and only if for $w^x$--a.e. $S \in \cW\cF(n)$, $\mu_d(S) = 1$ i.e. $S \in \G$ via \eqref{eq:Gembedding} $\G \simeq \{ S \in \cW\cF(n) \: : \: \mu_d(S) = 1 \}$ i.e. $S = \frac{1}{d} \Pi_{E_d(S)} = i(E_d(S))$. It follows that for $\psi \in \xC_c(\Omega \times \cW\cF(n))$,
\[
\int_{S \in \cW\cF(n)} \mu_d(S) \psi (x, i(E_d(S))) \: d w^x(S) = \int_{S \in \cW\cF(n)}  \psi (x, S) \: d w^x(S)
\]
and consequently, integrating with respect to $x \in M_d$ ($=M_d^\prime$ up to a set of $\|W\|$-null measure) and applying definitions of $V_d$ and $\widehat{V_d}$,
\begin{align*}
 \int_{M_d \times \cW\cF(n)}   \psi (x, S) d W(x, S)  & =  \int_{M_d \times \cW\cF(n)}   \mu_d(S) \psi (x, i(E_d(S))) d W(x, S) \\
 & = \int_{M_d \times \G} \psi(x, i(P)) \: dV_d(x,P) 
 = \int_{M_d \times \cW\cF(n)} \psi(x,S) \widehat{V_d} (x,S) \: .
\end{align*}
Therefore $W_{|M_d \times \cW\cF(n)} = \widehat{V_d}_{|M_d \times \cW\cF(n)}$. Moreover $\| \widehat{V_d} \| (\Omega \setminus M_d) = \int_{\Omega \setminus M_d} \mu_d(x) \: d \|W\|(x) = 0 $ and $\| W \| \left( \Omega \setminus \cup_{d=1}^n M_d \right) = 0$ and we can conclude that
\[
W = \sum_{d=1}^n W_{|M_d \times \cW\cF(n)} = \sum_{d=1}^n \widehat{V_d} \: .
\]

\noindent {\bf Step $3$}: If $\delta V_1, \ldots, \delta V_n$ are Radon measures, so is their sum $\delta W$. It remains to prove that $V_d$ is a rectifiable $d$--varifold provided that $\delta V_d$ is Radon, which is exactly given by Allard's rectifiability theorem.
\end{proof}

%There is still a question: can we infer the final conclusion only with the assumption $\delta W$ Radon. Is it possible to prove that actually $\delta V_d$ are automatically Radon measures and apply again Allard's regularity theorem. Another strategy would be to follow \cite{DePhilippis}, while the blow up in the first variation should provide the fact that tangent measures of $V_d$ are invariant along $d$ directions and the probability measure on $\G$ is a Dirac mass, the fact that $\Theta^{d, \ast} (x , \| V_d \|) < + \infty$ for $\| V_d \|$--a.e. $x \in \Omega$ is not clear.

%\subsection{Rectifiability theorem}
%
%dimension is integer (Marstrand theorem, tangent measure are $d$--uniform + property w.r.t. restriction)
%
%? characterization through blow up
%
%\section{Approximation properties}

\bibliographystyle{alpha} 
\bibliography{flagifoldsBib}
%\nocite{*}
     
\end{document}